\documentclass[11pt]{book}

\usepackage{amsmath}
\usepackage{amssymb}
\usepackage[mathcal]{euscript}
\usepackage{epsfig}
\usepackage{epic,eepic}

\newcommand{\itbf}{\itshape\bfseries}

\newcommand{\eps}{\epsilon}

\newcommand{\ee}{{\bf e}}
\newcommand{\bn}{{\boldsymbol n}}
\newcommand{\be}{{\boldsymbol e}}
\newcommand{\ba}{{\boldsymbol a}}
\newcommand{\bb}{{\boldsymbol b}}

\newcommand{\bbZ}{{\mathbb Z}}
\newcommand{\bbR}{{\mathbb R}}
\newcommand{\bbN}{{\mathbb N}}
\newcommand{\bbC}{{\mathbb C}}
\newcommand{\bbS}{{\mathbb S}}
\newcommand{\bbP}{{\mathbb P}}
\newcommand{\bbQ}{{\mathbb Q}}

\newcommand{\cA}{{\mathcal A}}

\newcommand{\cC}{{\mathcal C}}
\newcommand{\cD}{{\mathcal D}}
\newcommand{\cE}{{\mathcal E}}
\newcommand{\cF}{{\mathcal F}}
\newcommand{\cG}{{\mathcal G}}
\newcommand{\cH}{{\mathcal H}}
\newcommand{\cK}{{\mathcal K}}

\newcommand{\cP}{{\mathcal P}}
\newcommand{\cS}{{\mathcal S}}
\newcommand{\cT}{{\mathcal T}}

\newcommand{\cX}{{\mathcal X}}

\newcommand{\PP}{{\mathcal B}}
\newcommand{\OH}{{\mathcal O}}
\newcommand{\cn}{{\mathbb L}^{N+1}}

\newcommand{\gh}{{\mathfrak h}}

\newcommand{\veps}{{\mathbf \eps}}

\newcommand{\et}{{\textstyle\frac{\epsilon}{2}}}

\newcommand{\ra}{\to}
\newcommand{\rest}{\!\!\upharpoonright}

\def\To#1{\overset{#1}\longmapsto}
\def\I{{\rm I}}
\def\II{{\rm II}}
\def\III{{\rm III}}
\def\IV{{\rm IV}}
\def\V{{\rm V}}

\newtheorem{thm}{Theorem}[chapter]

\newtheorem{lem}[thm]{Lemma}

\newtheorem{dfn}[thm]{Definition}

\begin{document}
\thispagestyle{empty} \quad \vspace{1cm}
\begin{center}{\huge Discrete Differential Geometry. \\
Consistency as Integrability}\\ \vspace{2cm}
{\Large Alexander
I.\,
Bobenko, Yuri B.\, Suris}
\end{center}
\vspace{1.5cm}
\begin{center}
\includegraphics[width=0.6\textwidth]{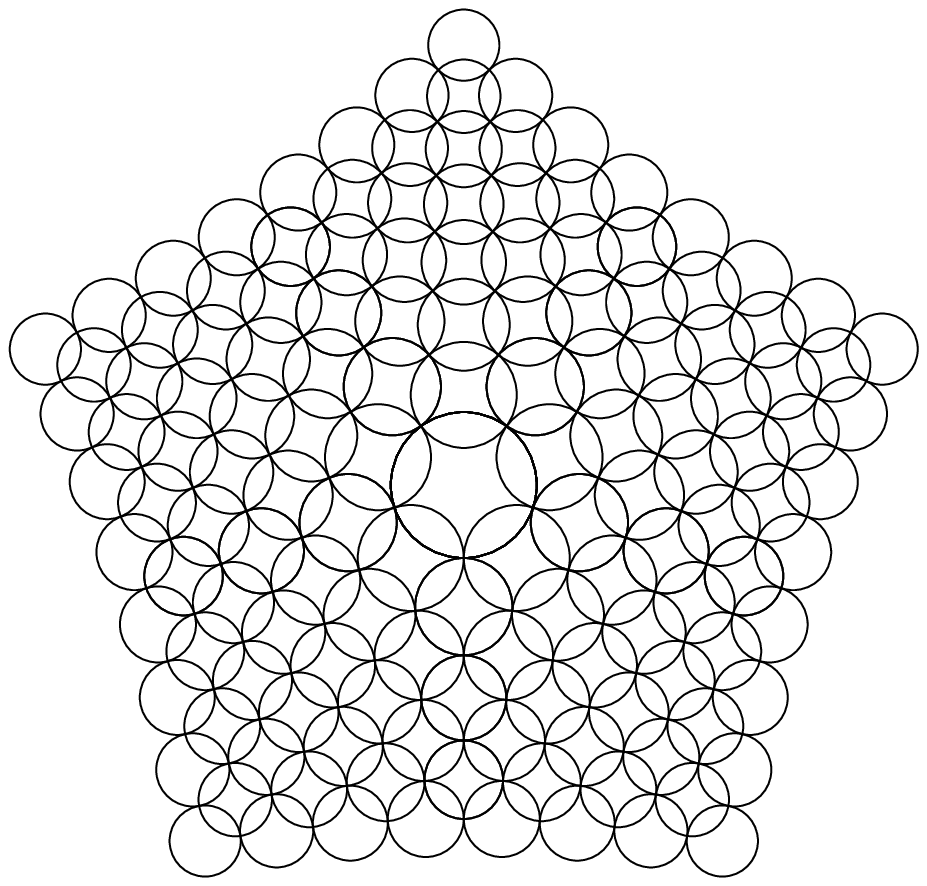}
\end{center}

 \vfill
 Supported by the DFG Forschergruppe
``Polyhedral Surfaces'' and the DFG Research Center
\textsc{Matheon} ``Mathematics for key technologies'' in Berlin.

\tableofcontents
\setcounter{chapter}{-1}
\chapter{Introduction} \label{Sect: intro}

\begin{figure}
\begin{center}
    \begin{tabular}{l c c c c}
\parbox[c]{0.3\textwidth}{\centerline{Differential  Geometry}}
& &
\parbox[c]{0.3\textwidth}{\center{Discrete  Differential \\ Geometry}}
& &
\parbox[c]{0.2\textwidth}{\centerline{Integrability}}
\vspace*{6mm}
\\
\parbox[c]{0.3\textwidth}{
    \makebox[0pt][l]{\hspace{10pt}\includegraphics[width=0.20\textwidth]
    {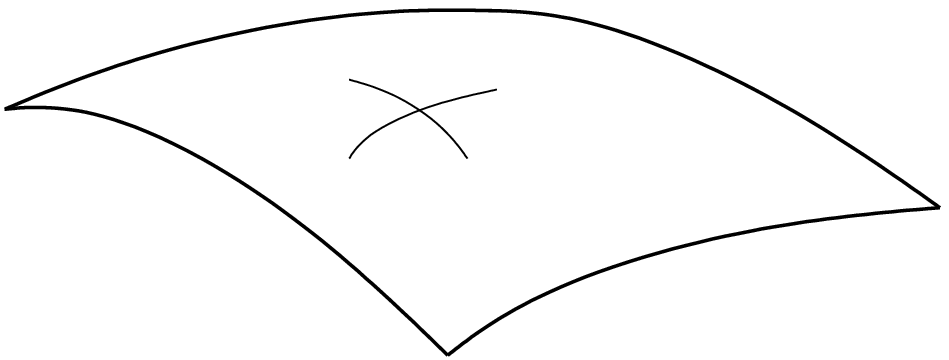}}
    \centerline{surfaces}}
    &
    \parbox[c]{0.05\textwidth}{\boldmath{$\Longleftarrow$}}
    &
\parbox[c]{0.2\textwidth}{
    \makebox[0pt][l]{\hspace{5pt}\setlength{\unitlength}{0.03em}
\begin{picture}(160,160)(0,0)
 \put(0,0){\circle*{15}}    \put(150,0){\circle*{15}}
 \put(0,150){\circle*{15}}  \put(150,150){\circle{15}}
 \path(0,0)(150,0)       \path(0,0)(0,150)
 \path(150,0)(150,142.5)   \path(0,150)(142.5,150)
\end{picture}}

    \vspace{10pt}
    \centerline{discrete nets}}
    &
    \parbox[c]{0.05\textwidth}{\boldmath{$\Longrightarrow$}}
    &\parbox[c]{0.2\textwidth}{integrable equations}
\vspace*{2mm}
\\
         & &\boldmath{$\Uparrow$}  & &
\vspace*{2mm}
\\
\parbox[c]{0.3\textwidth}{
    \makebox[0pt][l]{\hspace{0.5pt}\includegraphics[width=0.25\textwidth]
    {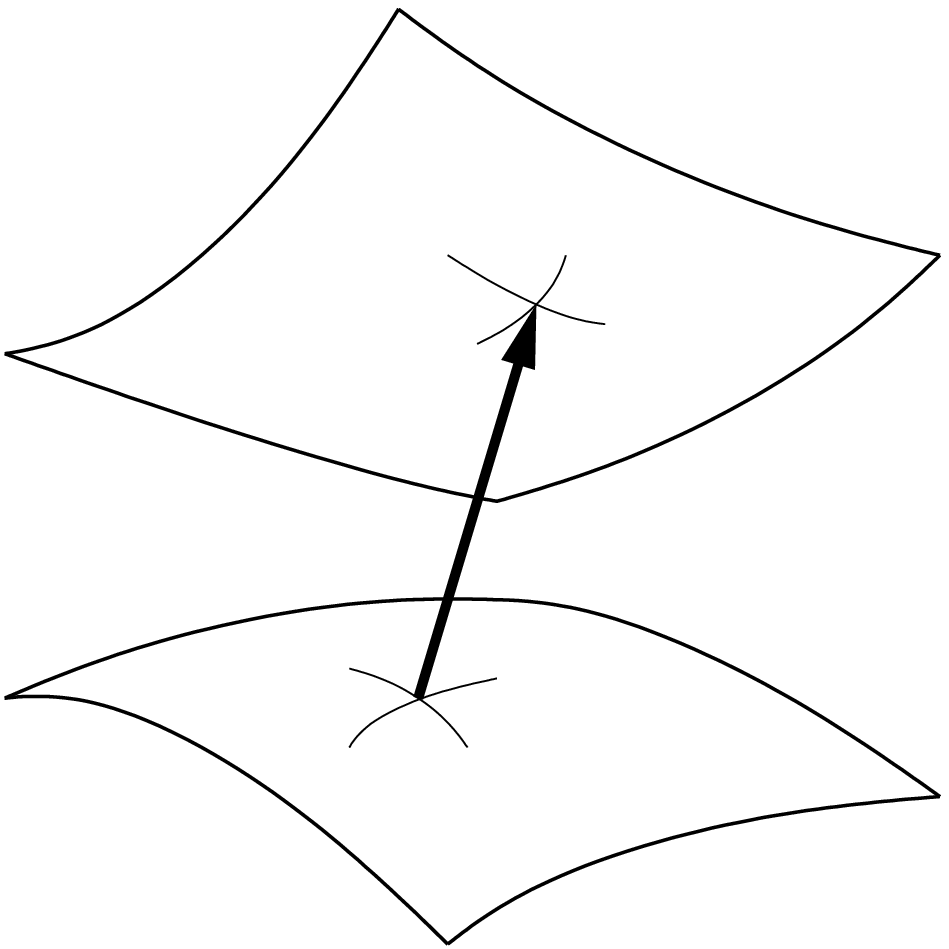}}
    \centerline{surfaces and transformations}}
    &
    \parbox[c]{0.05\textwidth}{\boldmath{$\Longleftarrow$}}
    &
\parbox[c]{0.3\textwidth}{
    \makebox[0pt][l]{\hspace{10pt}\setlength{\unitlength}{0.035em}
\begin{picture}(200,210)(0,0)
 \put(0,0){\circle*{15}}    \put(150,0){\circle*{15}}
 \put(0,150){\circle*{15}}  \put(150,150){\circle{15}}
 \put(50,50){\circle*{15}} \put(50,200){\circle{15}}
 \put(200,50){\circle{15}}
 \put(193,193){\tiny$\Box$}
 \path(0,0)(150,0)       \path(0,0)(0,150)
 \path(150,0)(150,142.5)   \path(0,150)(142.5,150)
 \path(0,150)(44.7,194.7)    \path(155.3,155.3)(193,193)
 \path(57.5,200)(193,200)
 \path(200,193)(200,57.5) \path(150,0)(194.7,44.7)
 \dashline[+30]{10}(0,0)(50,50)
 \dashline[+30]{10}(50,50)(50,192.5)
 \dashline[+30]{10}(50,50)(192.5,50)
 \end{picture}}

    \vspace{10pt}
    \centerline{CONSISTENCY}}
    &
    \parbox[c]{0.05\textwidth}{\boldmath{$\Longrightarrow$}}
    &\parbox[c]{0.2\textwidth}{zero curvature representation,
                                B\"acklund-Dar\-boux transformations}
\vspace*{2mm}
\\
         &  &\boldmath{$\Downarrow$}  & &
\vspace*{2mm}
\\
    \parbox[l]{0.3\textwidth}{
    \makebox[0pt][l]{\hspace{-5pt}\includegraphics[width=0.35\textwidth]
    {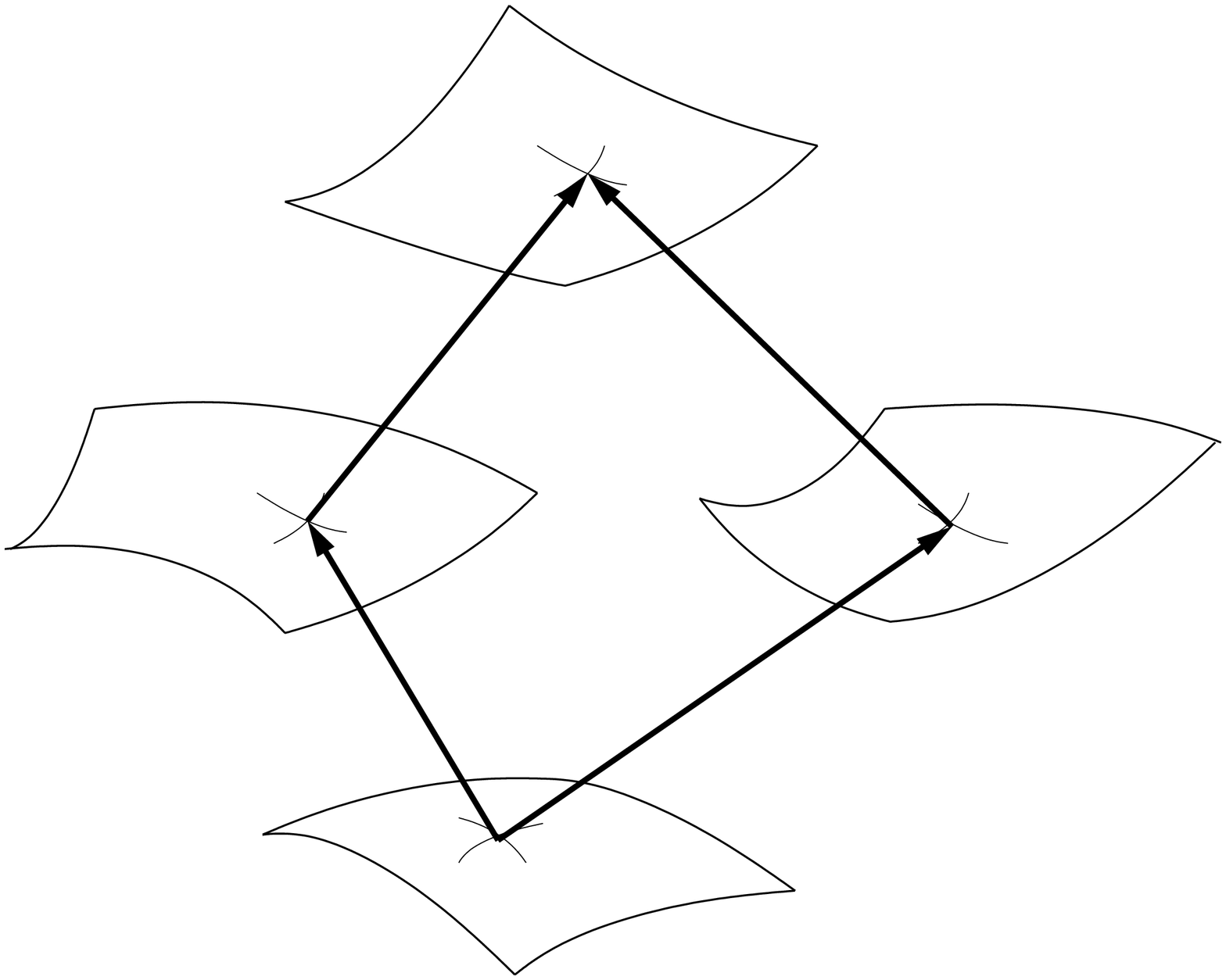}}
    \centerline{Bianchi permutability}}
    &
    \parbox[c]{0.05\textwidth}{\boldmath{$\Longleftarrow$}}
    &\parbox[c]{0.3\textwidth}{
    \makebox[0pt][l]{\hspace{10pt}
\setlength{\unitlength}{0.03em}
\begin{picture}(200,220)(-100,-90)

 \drawline(15,-20)(50,0)(50,47)
 \drawline(47,50)(0,50)(-35,30)(-35,-20)(15,-20)(15,30)(47,48.5)
 \drawline(15,30)(-35,30)
 \drawline(-35,-20)(0,0)(0,50)
 \drawline(0,0)(50,0)
 \drawline(30,-90)(131,-32)
 \drawline(135,-26)(135,116)
 \drawline(131,120)(-11,120)
 \drawline(-19,118)(-120,60)(-120,-90)(30,-90)(30,56)
 \drawline(34,62)(135,120)
 \drawline(26,60)(-120,60)
 \drawline(-120,-90)(-15,-30)(-15,116)
 \drawline(-15,-30)(131,-30)
  \drawline(0,0)(-15,-30)
  \drawline(-35,-20)(-120,-90)
  \drawline(50,0)(132,-29)
  \drawline(0,50)(-14,116)
  \drawline(15,-20)(30,-90)
  \drawline(-35,30)(-120,60)
  \drawline(53,51.5)(135,120)
  \drawline(15,30)(26,56)

  \put(-35,-20){\circle*{8}}       
  \put(15,-20){\circle*{8}}        
  \put(0,0){\circle*{8}}           
  \put(-35,30){\circle*{8}}        
  \put(50,0){\circle*{8}}          
  \put(0,50){\circle*{8}}          
  \put(15,30){\circle*{8}}         
  \put(50,50){\circle*{8}}          

  \put(-120,-90){\circle*{8}}     
  \put(-15,-30){\circle*{8}}       
  \put(30,-90){\circle*{8}}       
  \put(-120,60){\circle*{8}}      
  \put(135,-30){\circle*{8}}       
  \put(-15,120){\circle*{8}}       
  \put(30,60){\circle*{8}}         
  \put(135,120){\circle*{8}}       
\end{picture}
}

    \vspace{10pt}
    \centerline{multi-consistency}}
    &
    \parbox[c]{0.05\textwidth}{\boldmath{$\Longrightarrow$}}
    &\parbox[c]{0.2\textwidth}{hierarchies of commuting flows}
\end{tabular}
\end{center}
\caption{Fundamental consistency principle of the discrete
differential geometry as a conceptual basis of the differential
geometry of special surfaces and of the integrability.}
\label{fig:consistency}
\end{figure}

 A new field of {\em discrete differential geometry}
is presently emerging on the border between differential and
discrete geometry. Whereas classical differential geometry
investigates smooth geometric shapes (such as surfaces), and
discrete geometry studies geometric shapes with finite number of
elements (such as polyhedra), the discrete differential geometry
aims at the development of discrete equivalents of notions and
methods of smooth surface theory. Current interest in this field
derives not only from its importance in pure mathematics but also
from its relevance for other fields like computer graphics. Recent
progress in discrete differential geometry has lead, somewhat
unexpectedly, to a better understanding of some fundamental
structures lying in the basis of the classical differential
geometry and of the theory of integrable systems (this is
schematically presented on Fig.~\ref{fig:consistency}). The goal
of this book is to give a systematic presentation of current
achievements in this field.

The classical period of development of surface theory resulted in
the beginning of the 20-th century in an enormous amount of
knowledge on numerous special classes of surfaces, coordinate
systems and their transformations, which is summarized in
extensive volumes by Darboux \cite{Da1, Da2}, Bianchi \cite{Bi}
etc. One can say that the local differential geometry of special
classes of surfaces and coordinate systems has been completed
during this heroic period. Mathematicians of that era have found
most (if not all) geometries of interest and knew nearly
everything about their properties. It was observed that such
special geometries as minimal surfaces, surfaces with constant
curvature, isothermic surfaces, orthogonal and conjugate
coordinate systems, Ribaucour sphere congruences, Weingarten line
congruences etc. have many similar features. Among others we
mention here B\"acklund and Darboux type transformations, with
remarkable permutability properties investigated mainly by
Bianchi, and existence of special deformations within the class
(associated family). Geometers realized that there should exist a
unifying fundamental structure behind all these common properties
of quite different geometries. And they were definitely looking
for this structure \cite{Jonas, E2}.

Much later, after advent of the solitons theory in the the last
quarter of the 20-th century, these common similar features were
recognized to be associated to {\em integrability} of the
underlying differential equations. However, the current status of
the notion of integrability remains unsatisfactory from the point
of view of a mathematician. There exists no commonly accepted
definition of the integrability (as the title of the volume ``What
is integrability?'' \cite{Z} clearly demonstrates). Different
scientists suggest different properties as defining ones. Usually,
one just refers to some common features as Darboux transformations
etc., exactly like the differential geometers of the classical
period did.

A progress in understanding of the unifying fundamental structure
the classical differential geometers were looking for, and
simultaneously in understanding of the very nature of
integrability, came from the efforts to {\em discretize} all these
theories. It turns out that many sophisticated properties of
differential-geometric objects find their simple explanation
within the discrete differential geometry. The early period of its
development is documented in the works of Sauer and
Wunderlich~\cite{Sauer, W}. The modern period began with the work
by Bobenko and Pinkall~\cite{BP1, BP2} and by Doliwa and
Santini~\cite{DS1, CDS}. A closely related development of the
spectral theory of difference operators on graphs was initiated by
Novikov with collaborators \cite{ND,No1,No2}, see also \cite{DN}
for a further development of a discrete complex analysis on
triangulated manifolds. Discrete differential geometry deals with
multidimensional discrete nets (i.e., maps from the regular cubic
lattice $\bbZ^m$ into $\bbR^N$ or some other suitable space)
specified by certain geometric properties. In this setting,
discrete surfaces appear as two-dimensional layers of
multidimensional discrete nets, and their transformations
correspond to shifts in the transversal lattice directions. A
characteristic feature of the theory is that all lattice
directions are on equal footing with respect to the defining
geometric properties. Due to this symmetry, discrete surfaces and
their transformations become indistinguishable. We associate such
a situation with the {\em multidimensional consistency} (of
geometric properties, resp. of equations which serve for their
analytic description). The multidimensional consistency, and
therefore the existence and construction of multidimensional
discrete nets, relies just on certain incidence theorems of
elementary geometry.

\begin{figure}[htbp]
\begin{center}
\includegraphics[width=0.40\textwidth]{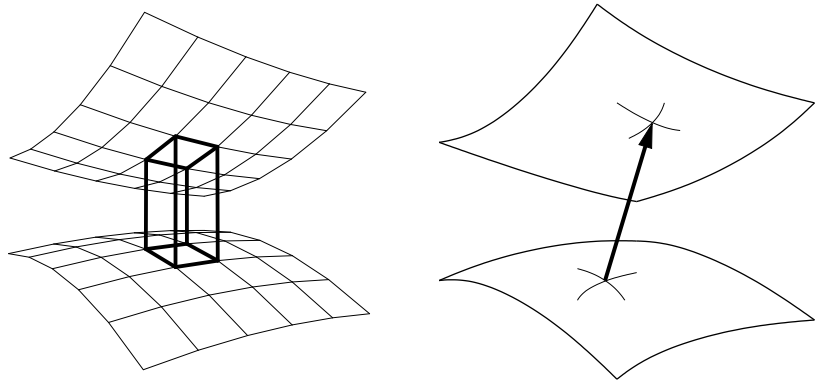}\qquad
\includegraphics[width=0.40\textwidth]{GlatteTrafo.eps}
\end{center}
\caption{From the discrete master theory to the classical theory:
surfaces and their transformations appear by refining two of three
net directions.} \label{fig:Backlund}
\end{figure}

Conceptually, one can think of passing to a continuum limit by
refining mesh size in some of the lattice directions. In these
directions the net converges to smooth surfaces whereas those
directions that remain discrete correspond to transformations of
the surfaces (see Fig.~\ref{fig:Backlund}). Differential geometric
properties of special classes of surfaces and their
transformations follow in this way from (and find their simple
explanation in) the elementary geometric properties of the
original multidimensional discrete nets. In particular, difficult
classical theorems about permutability of the B\"acklund-Darboux
type transformations (Bianchi permutability) for various
geometries follow directly from the symmetry of the underlying
discrete nets, and are therefore built in the very core of the
theory. Thus the pass from differential geometry to elementary
geometry via discretization (or, in an opposite direction, the
derivation of the differential geometry from the discrete
differential geometry) leads to enormous conceptual
simplifications, and the true roots of the classical theory of
special classes of surfaces are found in various incidence
theorems of elementary geometry. However, these elementary roots
become deeply hidden in the classical differential geometry, since
the continuum limit from the discrete master theory to the
classical one is inevitably accompanied by a break of the symmetry
among the lattice directions, which always yields essential
structural complications.

To give a rigorous justification of this philosophy, one needs to
prove convergence of the procedure just described. Theorems of
this kind were lacking in the literature until recently. The first
results of this kind have been proven by the authors in a common
work with Matthes \cite{BMS1, BMS2}. This makes the general
philosophy of the discrete differential geometry to a firmly
established mathematical truth for several important classes of
surfaces and coordinate systems, like conjugate nets, orthogonal
nets, including general surfaces parametrized along curvature
lines, surfaces with constant negative Gaussian curvature, and
general surfaces parametrized along asymptotic lines. For some
other classes, like isothermic surfaces, the convergence results
still wait to be rigorously established.

But finding simple discrete explanations for complicated
differential geometric theories is not the only outcome of this
development. It is well known that differential equations which
analytically describe interesting special classes of surfaces are
integrable (in the sense of the theory of integrable systems),
and, conversely, many of interesting integrable systems admit a
differential-geometric interpretation. Having identified the roots
of the integrable differential geometry in the multidimensional
consistency of discrete nets, we are led to a new (geometric)
understanding of the integrability itself. First, we adhere to the
viewpoint that the central role in this theory belongs to {\em
discrete integrable systems}. In particular, all the great variety
of integrable differential equations can be derived from several
fundamental discrete systems by performing different continuous
limits. Further, and more important, we come to a constructive and
almost algorithmic definition of integrability of discrete
equations as their multidimensional consistency, introduced by the
authors in \cite{BS1} (and independently in \cite{Nij}). It turns
out that this definition captures enough structure to yield such
traditional attributes of integrable equations as zero curvature
representations and B\"acklund-Darboux transformations (which, in
turn, serve as the basis for applying analytic methods like
inverse scattering, finite gap integration, Riemann-Hilbert
problems, etc.). A continuous counterpart (and consequence) of the
multidimensional consistency is the well-known fact that
integrable systems never appear alone but are organized into
hierarchies of commuting flows.

The geometric way of thinking about the discrete integrability has
also led to introducing novel concepts into the latter. One of the
reasons to consider discrete integrable systems on the regular
square lattice $\bbZ^2$ is the desire to have a proper model for
{\em parametrized surfaces}. However, an immanent and important
feature of various parametrizations of surfaces is the existence
of distinguished points, where the combinatorics of coordinate
lines change (like umbilic points, where the combinatorics of the
curvature lines is special). This compels us to introduce {\em
quad-graphs}, which are cell decompositions of topological
two-manifolds with quadrilateral faces. Their elementary building
blocks are still quadrilaterals but are attached to one another in
a manner which can be more complicated than in $\bbZ^2$. This
notion has been introduced into the context of discrete
differential geometry in \cite{BP3}, and a systematic development
of the theory of integrable systems on quad-graphs has been
undertaken by the authors in \cite{BS1}.

The structure of the book follows the logic of this introduction.
We start in Chapter~\ref{Sect: smooth diffe geo} with an overview
of some classical results of the surface theory, focusing on
transformations of surfaces. The geometries considered here
include general conjugate and orthogonal nets in spaces of
arbitrary dimension, asymptotic nets on general surfaces, as well
as special classes of surfaces, like isothermic ones and surfaces
with constant negative Gaussian curvature. There are no proofs in
this chapter: on one hand, these tedious analytic proofs can be
found in the original literature, and on the other hand, the
discrete approach which we develop in the subsequent chapters will
lead to conceptually transparent and technically much simpler
proofs.

In Chapter~\ref{Chap: discr geom} we define and investigate
discrete analogs of the classical geometries discussed in the
previous chapter, focusing on the idea of multi-dimensional
consistency of discrete nets. It turns out that all these discrete
analogs are reductions of discrete conjugate nets, which are
multidimensional nets with the combinatorics of $\bbZ^m$,
consisting of planar quadrilaterals. Imposing additional
constraints on the geometry of elementary quadrilaterals, one
comes to discrete orthogonal nets, discrete asymptotic nets,
discrete isothermic surfaces, discrete surfaces with constant
negative Gaussian curvature etc.

Then in Chapter~\ref{Sect: approx} we develop an approximation
theory for hyperbolic difference systems, which is applied to
derive the classical theory of smooth surfaces as a continuum
limit of the discrete theory. We prove that discrete nets of
Chapter~\ref{Chap: discr geom} approximate the corresponding
smooth geometries of Chapter~\ref{Sect: smooth diffe geo}
simultaneously with their transformations. Bianchi's and
Eisenhart's permutability theorems for transformations appear in
this approach as simple corollaries.

In Chapter~\ref{Sect: consist} we formulate the concept of
multi-dimensional consistency as a defining principle of
integrability. We derive basic features of integrable systems,
such as the zero curvature representation, B\"acklund-Darboux
transformations, from the consistency principle. At this point,
the theory makes an interesting conceptual turn. First, we
generalize the underlying combinatorial structure of the
two-dimensional theory from the regular square lattice to
arbitrary quad-graphs, i.e., cell decomposition of surfaces with
all quadrilateral faces. Introducing these generalized
combinatorics has been partly motivated by the desire to have
proper discrete models for geometrically significant special
points of parametrized surfaces, like umbilic points in the
curvature line parametrization. But then the multidimensional
consistency allows us to regard integrable systems on quad-graphs
as systems on regular square lattices $\bbZ^m$ restricted to
quadrilateral surfaces. At this point, interesting interrelations
with the discrete geometry and the combinatorial analysis come
onto the scene.

Finally, in Chapters~\ref{Chap:_discr_analysis-linear},
\ref{Chap:_discr_analysis-nonlinear} these ideas are applied to
discrete complex analysis. We study Laplace operators on graphs,
discrete harmonic and holomorphic functions. The linear discrete
complex analysis appears here as a linearization of the theory of
circle patterns. The consistency principle allows us to single out
distinguished cases where we obtain more detailed analytic results
(like Green's function and isomonodromic special functions).

Essential parts of this book are based on results obtained jointly
with Vsevolod Adler, Tim Hoffmann, Daniel Matthes, Christian
Mercat and Ulrich Pinkall. We warmly thank them for inspiring
collaboration.

\chapter{Classical differential geometry}
\label{Sect: smooth diffe geo}

In this chapter we discuss some classical results of the
differential geometry of nets (parametrized surfaces and
coordinate systems) in $\bbR^N$, mainly concentrated around the
topics of transformations of nets and of their permutability
properties. This classical area was very popular in the
differential geometry of the 19th and of the first quarter of the
20th century, and is well documented in the fundamental treatises
\cite{Bi, Da1, Da2, E1, E2, Tzi} and others. Our presentation
mainly follows these classical treatments, of course with
modifications which reflect our present points of view. We do not
trace back the exact origin of the concrete classical results:
often enough this turns out to be a complicated task in the
history of mathematics which still waits for its competent
investigation.

For the classes of nets described by essentially two-dimensional
systems (special classes of surfaces such as surfaces with a
constant negative Gaussian curvature or isothermic surfaces), the
permutability theorems, mainly due to Bianchi, are dealing with a
{\em quadruple} of surfaces (depicted as vertices of a so-called
{\em Bianchi quadrilateral}). Given three surfaces of such a
quadruple, the fourth one is uniquely defined, see Theorems
\ref{thm:ks permut} and \ref{thm:iso permut}.

For the classes of nets described by essentially three-dimensional
systems (conjugate nets; orthogonal nets, including general
surfaces parametrized by curvature lines; Moutard nets; general
surfaces parametrized by asymptotic lines), the situation is
somewhat different. The corresponding permutability theorems
(Theorems \ref{thm:conj permut}, \ref{thm:ortho permut},
\ref{thm:mp permut}, and \ref{thm:as permut}) consist of two
parts. The first parts present the traditional view and are
dealing with Bianchi quadrilaterals. In our opinion, this is not
the proper setting in the three-dimensional context, and the
non-uniqueness of the fourth net in these theorems reflects this.
The natural setting for permutability is given in the second
parts, where the permutability is associated with an {\em octuple}
of nets, depicted as vertices of a combinatorial cube, so that the
eighth net is uniquely determined by other seven ones. We found
the first instance of this kind of statement in \cite{E2}, \S 24
(``extended theorem of permutability'' for conjugate nets), and
propose therefore to term such an octuple of nets as an {\em
Eisenhart cube}; in the modern literature on integrable systems,
this kind of permutability theorems in classical differential
geometry was (re-)discovered in \cite{GTs} (for orthogonal nets).
Our discrete philosophy makes the origin of such permutability
theorems quite transparent.

A few remarks on notations: we denote independent variables of a
net $f:\bbR^m\to\bbR^N$ by $u=(u_1,\ldots,u_m)\in\bbR^m$, and we
set $\partial_i=\partial/\partial u_i$. We write $\PP_{i_1\ldots
i_s}=\big\{u\in\bbR^m: u_i=0,\quad i\neq i_1,\ldots,i_s\big\}$ for
$s$-dimensional coordinate planes (coordinate axes, if $s=1$). We
always suppose that the dimension of the ambient space $N\ge 3$.

\section[Conjugate nets]{Conjugate nets and their transformations}
 \label{Sect: conj}

{\bf Conjugate nets.} This classical notion can be defined as follows.
\begin{dfn}\label{dfn:cn}
  A map $f:\bbR^m\to\bbR^N$ is called an $m$-dimensional
  {\em conjugate net} in $\bbR^N$, if
  $\partial_i\partial_j f\in{\rm span}(\partial_i f,\partial_j f)$
  at any $u\in \bbR^m$ and for all pairs $1\leq i\neq j\leq m$.
\end{dfn}
A parametrized surface in the three-space ($m=2$, $N=3$) is a
conjugate net, if its second fundamental form is diagonal in this
parametrization. Such a parametrization exists for a general
surface in the three-space. It is important to note that
Definition \ref{dfn:cn}, as well as Definition \ref{dfn:jp} below,
are dealing with projectively invariant notions only, and thus
belong to the projective differential geometry. In this setting
the space $\bbR^N$ where conjugate net lives should be interpreted
as an affine part of $\bbP(\bbR^{N+1})$.

>From Definition \ref{dfn:cn} there follows that conjugate nets are
described by the following (linear) differential equations:
\begin{equation} \label{eq:cn property}
    \partial_i\partial_j f=c_{ji}\partial_i f+c_{ij}\partial_j f,
    \quad i\neq j,
\end{equation}
with some functions $c_{ij}:\bbR^m\to\bbR$. Compatibility of
these equations is expressed by the following system of
(nonlinear) differential equations:
\begin{equation}\label{eq:cn c}
   \partial_i c_{jk}=c_{ij}c_{jk}+c_{ji}c_{ik}-c_{jk}c_{ik},
   \quad i\neq j\neq k\neq i,
\end{equation}
which thus split off from eqs. (\ref{eq:cn property}) for $f$.
System  (\ref{eq:cn property}), (\ref{eq:cn c}) is hyperbolic (see
Sect. \ref{Sect: consist}); the following data define a well-posed
Goursat problem for this system and determine a conjugate net $f$
uniquely:
\newcounter{Q}
\begin{list}{(Q$_\arabic{Q}$)}
{\usecounter{Q}\setlength{\leftmargin}{1cm}}
\item values of $f$ on the coordinate axes $\PP_i$ for $1\le i\le m$, i.e.,
 $m$ smooth curves $f\rest_{\PP_i}$ with a common intersection point $f(0)$;
\item values of $c_{ij}$, $c_{ji}$ on the coordinate planes $\PP_{ij}$ for
 all $1\leq i<j\leq m$, i.e., $m(m-1)$ smooth real-valued functions
 $c_{ij}\rest_{\PP_{ij}}\,$ of two variables.
\end{list}
{\bf Alternative analytic description of conjugate nets.}
Given the functions $c_{ij}$, define functions $h_i:\bbR^m\to\bbR$
as solutions of the system of differential equations
\begin{equation}\label{eq:cn h}
\partial_ih_j=c_{ij}h_j\,, \qquad i\neq j.
\end{equation}
Compatibility of this system is a consequence of (\ref{eq:cn c}).
Define vectors $v_i=h_i^{-1}\partial_i f$. It follows from
(\ref{eq:cn property}) and (\ref{eq:cn h}) that these vectors
satisfy the following differential equations:
\begin{equation}\label{eq:cn v}
\partial_iv_j=\frac{h_i}{h_j}c_{ji}v_i\,, \qquad i\neq j.
\end{equation}
Thus, defining the {\em rotation coefficients} as
\begin{equation}\label{eq: cn rot}
\beta_{ji}=\frac{h_i}{h_j}c_{ji}\,,
\end{equation}
we end up with the following system:
\begin{eqnarray}
    \partial_if  & = & h_i v_i\,,                          \label{eq:DX}\\
    \partial_i v_j & = & \beta_{ji} v_i\,,\quad i\neq j,   \label{eq:DV}\\
    \partial_i h_j & = & h_i\beta_{ij}\,,\quad i\neq j.    \label{eq:DH}
\end{eqnarray}
Rotation coefficients satisfy a
closed system of differential equations, which follow from eqs.
(\ref{eq:cn c}) upon substitution (\ref{eq: cn rot}):
\begin{equation}\label{eq:DB}
\partial_i\beta_{kj} = \beta_{ki}\beta_{ij}\,,
    \quad i\neq j\neq k\neq i.
\end{equation}
Eqs. (\ref{eq:DB}), known as the {\em Darboux system}, can be
regarded as compatibility conditions of the linear differential
equations (\ref{eq:DV}).

Observe an important difference between two descriptions of
conjugate nets: while the functions $c_{ij}$ describe the {\em
local geometry} of a net, this is not the case for the rotation
coefficients $\beta_{ij}$. Indeed, to define the latter, one needs
to find $h_i$ as solutions of differential  equations (\ref{eq:cn h}).
\smallskip

{\bf Transformations of conjugate nets.} The most general class of
transformations of conjugate nets was introduced by Jonas and Eisenhart.
\begin{dfn}\label{dfn:jp}
  A pair of $m$-dimensional conjugate nets
  $f,f^+:\bbR^m\to\bbR^N$ is called a {\em Jonas pair}, if
  three vectors $\partial_i f$, $\partial_i f^+$ and $\delta f= f^+-f$
  are co-planar at any point $u\in\bbR^m$ of the definition domain and
  for any $1\leq i\leq m$. The net $f^+$ is called a {\em Jonas
  transform} of the net $f$.
\end{dfn}
This definition yields that Jonas transformations are described by the
following (linear) differential equations:
\begin{equation}\label{eq:jp property}
\partial_i f^+=a_i\partial_i f+b_i(f^+-f).
\end{equation}
Of course, functions $a_i,b_i:\bbR^m\to\bbR$ have to satisfy
(nonlinear) differential equations, which express the compatibility
of eqs. (\ref{eq:jp property}) with (\ref{eq:cn property}):
\begin{eqnarray}
\partial_ia_j & = & (a_i-a_j)c_{ij}+b_i(a_j-1),       \label{eq:jp dA}\\
\partial_ib_j & = & c_{ij}^+b_j+c_{ji}^+b_i-b_jb_i,   \label{eq:jp dB}\\
a_jc_{ij}^+ & = & a_ic_{ij}+b_i(a_j-1).               \label{eq:jp C+}
\end{eqnarray}
Following data determine a Jonas transform $f^+$ of a given
conjugate net $f$ uniquely:
\newcounter{J}
\begin{list}{(J$_\arabic{J}$)}
{\usecounter{J}\setlength{\leftmargin}{1cm}}
\item a point $f^+(0)$;
\item values of $a_i$, $b_i$ on the coordinate axes
 $\PP_i$ for $1\le i\le m$,
 i.e., $2m$ smooth real-valued functions $a_i\rest_{\PP_i}$,
 $b_i\rest_{\PP_i}$ of one variable.
\end{list}

Observe a remarkable conceptual similarity between Definitions
\ref{dfn:cn} and \ref{dfn:jp}. Indeed, one can interpret the
condition of Definition \ref{dfn:cn} as planarity of infinitesimal
quadrilaterals $(f(u),f(u+\eps_i e_i),f(u+\eps_i e_i+\eps_j e_j),
f(u+\eps_j e_j))$, while the condition of Definition \ref{dfn:jp}
can be interpreted as planarity of infinitesimally narrow
quadrilaterals $(f(u),f(u+\eps_ie_i),f^+(u+\eps_ie_i),f^+(u))$.
\smallskip

{\bf Classical formulation of the Jonas transformation.} Our
formulation of Jonas transformations is rather different from the
classical one, which can be found, e.g., in \cite{E2}. The latter
is based on the formula
\begin{equation}\label{eq:jp Eis}
f^+=f-\frac{\phi}{\psi}\,g,
\end{equation}
whose data are: an additional solution $\phi:\bbR^m\to\bbR$ of eq.
(\ref{eq:cn property}), a parallel to $f$ net $g:\bbR^m\to\bbR^N$,
and the function $\psi:\bbR^m\to\bbR$, associated to $\phi$ in the
same way as $g$ is related to $f$. We now demonstrate how to
identify these ingredients within our approach and how they are
specified by the initial data (J$_{1,2}$).

There follows from eqs. (\ref{eq:jp dA})--(\ref{eq:jp C+}):
\begin{equation}\label{eq:jp dB/A}
\partial_i\Big(\frac{b_j}{a_j}\Big)=c_{ij}\frac{b_j}{a_j}+c_{ji}
\frac{b_i}{a_i}-\frac{b_jb_i}{a_ja_i}.
\end{equation}
The symmetry of the right-hand sides of eqs. (\ref{eq:jp dB/A}),
(\ref{eq:jp dB}) yields the existence of the functions
$\phi,\phi^+:\bbR^m\to\bbR$ such that
\begin{equation}\label{eq:jp phis}
\frac{\partial_i\phi}{\phi}=\frac{b_i}{a_i}\,,\qquad
\frac{\partial_i\phi^+}{\phi^+}=b_i\,,\qquad 1\le i\le m.
\end{equation}
These equations define $\phi,\phi^+$ uniquely up to respective constant
factors, which can be fixed by requiring $\phi(0)=\phi^+(0)=1$.
An easy computation based on eqs. (\ref{eq:jp dB/A}), (\ref{eq:jp dB})
shows that the functions $\phi,\phi^+$ satisfy the following equations:
\begin{eqnarray}
\partial_i\partial_j\phi & = & c_{ij}\partial_j\phi+c_{ji}\partial_i\phi,
\label{eq:jp phi eq}  \\
\partial_i\partial_j\phi^+ & = & c_{ij}^+\partial_j\phi^+
+c_{ji}^+\partial_i\phi^+,   \label{eq:jp phi+ eq}
\end{eqnarray}
for all $1\le i\neq j\le m$. Thus, a Jonas transformation yields
some additional scalar solutions $\phi$ and $\phi^+$  of the
equations describing the nets $f$ and $f^+$, respectively. It is
clear that the solution $\phi$ is directly specified by the
initial data (J$_2$). Indeed, these data yield the values $\phi$
along the coordinate axes, through integrating the first equations
in (\ref{eq:jp phis}); these values determine the solution of eq.
(\ref{eq:jp phi eq}) uniquely.

Further, introduce the quantities
\begin{equation}\label{eq:jp aux0}
g=\frac{f^+-f}{\phi^+}\,,\qquad \psi=-\frac{\phi}{\phi^+}\,.
\end{equation}
Then a direct computation based on eqs. (\ref{eq:jp property}),
(\ref{eq:jp dA})--(\ref{eq:jp C+}), and (\ref{eq:jp phis}) shows that
the following equations hold:
\begin{eqnarray}
\partial_i g & = & \alpha_i\partial_i f,             \label{eq:jp aux1}\\
\partial_i\psi & = & \alpha_i\partial_i\phi,         \label{eq:jp aux2}
\end{eqnarray}
where
\begin{equation}\label{eq:jp aux3}
\alpha_i=\frac{a_i-1}{\phi^+}\,.
\end{equation}
Thus, $g$ is a parallel net to $f$, and $\psi$ is an associated to $\phi$
function, in Eisenhart's terminology. Another computation leads to the
relation
\begin{equation}\label{eq:jp aux4}
\partial_i\alpha_j=c_{ij}(\alpha_i-\alpha_j).
\end{equation}
The same argument as above shows that the data (J$_2$) yield the
values of $\phi^+$, and thus the values of $\alpha_i$, on the
coordinate axes $\PP_i$. This uniquely specifies the solutions
$\alpha_i$ of the compatible linear system (\ref{eq:jp aux4}).
This, in turn, allows for a unique determination of the solutions
$g$, $\psi$ of eqs. (\ref{eq:jp aux1}), (\ref{eq:jp aux2}) with
the initial data $g(0)=f^+(0)-f(0)$ and $\psi(0)=1$ (here the data
(J$_1$) enter into the construction). Thus, the classical formula
(\ref{eq:jp Eis}) is recovered.
\smallskip

One can iterate Jonas transformations and
obtain a sequence $f$, $f^+$, $(f^+)^+$, ... , of conjugate nets.
We will see that this can be interpreted as generating a conjugate net
of dimension $M=m+1$, with $m$ continuous directions and one discrete
direction. The most remarkable property of Jonas transformations is
the following permutability theorem.

\begin{thm}[Permutability of Jonas transformations]
\label{thm:conj permut}\quad

 1)
  Let  $f$ be an $m$-dimensional conjugate net, and let
  $f^{(1)}$ and $f^{(2)}$ be its two Jonas transforms.
  Then there exists a two-parameter family of conjugate nets
  $f^{(12)}$ that are Jonas transforms of both $f^{(1)}$ and $f^{(2)}$.
  Corresponding points of the four conjugate nets $f$, $f^{(1)}$, $f^{(2)}$
  and $f^{(12)}$ are co-planar.

2)
  Let  $f$ be an $m$-dimensional conjugate net. Let
  $f^{(1)}$, $f^{(2)}$ and $f^{(3)}$ be its three Jonas transforms,
  and let three further conjugate nets $f^{(12)}$, $f^{(23)}$ and $f^{(13)}$
  be given such that $f^{(ij)}$ is a simultaneous Jonas transform
  of $f^{(i)}$ and $f^{(j)}$. Then there exists
  generically a unique conjugate net $f^{(123)}$
  that is a Jonas transform of $f^{(12)}$, $f^{(23)}$ and $f^{(13)}$.
  The net $f^{(123)}$ is uniquely defined by the condition that any its
  point is co-planar with the corresponding points of $f^{(i)}$, $f^{(ij)}$
  and $f^{(ik)}$ for any permutation $(ijk)$ of $(123)$.
\end{thm}

The situations described in this theorem can be interpreted as conjugate
nets of dimension $M=m+2$, resp. $M=m+3$, with $m$ continuous and two
(resp. three) discrete directions.

The theory of discrete conjugate nets allows one to put all
directions on an equal footing and to unify the theories of smooth
nets and of their transformations. Moreover, we will see that both
these theories may be seen as a {\em continuum limit} (in some
precise sense) of the fully discrete theory, if the mesh sizes of
all or some of the directions becomes infinitely small (see
Fig.~\ref{fig:Backlund}). This way of thinking is the guiding idea
and the philosophy of the {\em discrete differential geometry}.

\section[Orthogonal nets]{Orthogonal nets and their transformations}
 \label{Sect: subsect os}

{\bf Orthogonal nets.} An important subclass of conjugate nets
is fixed in the following definition.

\begin{dfn} \label{dfn:os}
  A conjugate net $f:\bbR^m\to\bbR^N$
  is called an $m$-dimensional {\em O-net (orthogonal net)} in $\bbR^N$,
  if there holds $\partial_i f \perp \partial_j f\,$ at any $u\in \bbR^m$
  and for all $1\leq i\neq j\leq m$.
  Such a net is called an {\em orthogonal coordinate system} if $m=N$.
\end{dfn}
Two-dimensional $(m=2)$ orthogonal nets are called {\it O-surfaces}.
An O-surface in $\bbR^3$ is nothing but a surface parametrized
along its curvature lines, or, otherwise said, parametrized so
that both the first and the second fundamental forms are diagonal. Such
a parametrization exists for a general surface in $\bbR^3$ in the
neighborhood of a non-umbilic point. Note that Definition \ref{dfn:os} is
dealing only with notions which are invariant under M\"obius
transformations. Thus orthogonal nets (as well as their Ribaucour
transformations from Definition \ref{dfn:rp} below) belong to the
M\"obius differential geometry. It will be important to preserve
this symmetry group under discretization.

For an analytic description of an orthogonal net $f:\bbR^m\to\bbR^N$,
introduce metric coefficients $h_i=|\partial_i f|$ and (pairwise orthogonal)
unit vectors $v_i=h_i^{-1}\partial_i f$. Then there hold
eqs. (\ref{eq:DX})--(\ref{eq:DB}), supplemented by the {\em orthogonality
constraint}
\begin{equation}\label{eq:DL}
    \partial_i\beta_{ij}+\partial_j\beta_{ji}
    =-\langle \partial_i v_i,\partial_j v_j\rangle\,,\quad i\neq j.
\end{equation}
Indeed, eq. (\ref{eq:DV}) holds since $f$ is a conjugate net and  $v_j$
are orthonormal, and serves as a definition of {\em rotation coefficients}
$\beta_{ji}$. Eq. (\ref{eq:DH}) is a direct consequence of (\ref{eq:DX}),
(\ref{eq:DV}). To derive eq. (\ref{eq:DL}), one considers the identity
$\partial_i\partial_j\langle v_i,v_j\rangle=0$.
So, a distinctive feature of orthogonal nets among general
conjugate ones is that the rotation coefficients $\beta_{ji}$
reflect the {\em local geometry}. In the same spirit, a solution
to the system (\ref{eq:cn h}) is given by the locally defined
metric coefficients $h_i$.

Eq. (\ref{eq:DL}) is an {\em admissible constraint} for the system
(\ref{eq:DX})--(\ref{eq:DB}). This has the following meaning:
eq. (\ref{eq:DL}) involves two independent variables $i$, $j$ only,
and it is therefore sensible to ask for it to be fulfilled
on the coordinate plane $\PP_{ij}$. One can easily check that if a
solution to the system (\ref{eq:DX})--(\ref{eq:DB}) fulfills eq.
(\ref{eq:DL}) on all coordinate planes $\PP_{ij}$ for $1\le i<j\le m$,
then it is fulfilled everywhere on $\bbR^m$.

System (\ref{eq:DX})--(\ref{eq:DB}), (\ref{eq:DL}) is not
hyperbolic, therefore it is less clear what data form a well-posed
problem for it. It can be shown (see Sect. \ref{subsect: os
approx}) that the following data can be used to determine an
orthogonal net $f$ uniquely:
\newcounter{O}
\begin{list}{(O$_\arabic{O}$)}
{\usecounter{O}\setlength{\leftmargin}{1cm}}
\item  values of $f$ on the coordinate axes $\PP_i$ for $1\le i\le m$,
 i.e., $m$ smooth curves $f\rest_{\PP_i}$, intersecting pairwise
 orthogonally at $f(0)$;
\item  $m(m-1)/2$ smooth functions $\gamma_{ij}:\PP_{ij}\to\bbR$
 for all $1\le i<j\le m$, which have the meaning of $\gamma_{ij}=\frac{1}{2}
 (\partial_i\beta_{ij}-\partial_j\beta_{ji})\rest_{\PP_{ij}}$.
\end{list}
\medskip

{\bf Ribaucour transformations of orthogonal nets.}
\begin{dfn} \label{dfn:rp}
  A pair of $m$-dimensional orthogonal nets $f,f^+:\bbR^m\to\bbR^N$
  is called a {\em Ribaucour pair}, if the corresponding coordinate curves
  of $f$ and $f^+$ envelope one-parameter families of circles, i.e.
  if at any $u\in\bbR^m$ and for any $1\le i\le m$ the straight lines
  spanned by the vectors $\partial_i f$, $\partial_i f^+$ at the corresponding
  points $f$, $f^+$ are interchanged by the reflection in the affine
  hyperplane orthogonal to $\delta f=f^+-f$, which interchanges $f$ and
  $f^+$. The net $f^+$ is called a {\em Ribaucour transform} of $f$.
\end{dfn}

To describe a Ribaucour transformation analytically, we write:
  \begin{equation}  \label{eq:Riba}
    \partial_i f^+= r_i\bigg(\partial_i f-2\,
    \frac{\langle\partial_i f,\delta f\rangle}
    {\langle \delta f,\delta f\rangle}\delta f\bigg),
  \end{equation}
with some functions $r_i:\bbR^m\to\bbR$ which obviously coincide (up to
a sign) with the quotients of the corresponding metric coefficients,
$r_i^2=(h_i^+/h_i)^2$. Further, denote $\ell=|\delta f|$ and
introduce the unit vector $y=\ell^{-1}\delta f$, so that $f^+=f+\ell y$.
Then, in the case $r_i>0$, we find:
\begin{equation}\label{eq:rp v}
v_i^+=v_i-2\langle v_i,y\rangle y,  \qquad
\partial_i y={\textstyle\frac{1}{2}}\theta_i(v_i^++v_i),
\end{equation}
with the functions $\theta_i:\bbR^m\to\bbR$ defined as $\theta_i
=(r_i-1)h_i/\ell$.
Eqs. (\ref{eq:rp v}) imply equations for the metric coefficients:
\begin{equation}\label{eq:rp h}
h_i^+=h_i+\theta_i\ell, \qquad
\partial_i\ell=-\langle v_i,y\rangle(h_i^++h_i).
\end{equation}
(In the case $r_i<0$ one has to change the sign of the quantities $v_i^+$,
$h_i^+$ in eqs. (\ref{eq:rp v}), (\ref{eq:rp h}).)
Compatibility of the system (\ref{eq:rp v}) yields that
$\theta_i$ have to satisfy certain differential equations:
\begin{equation}\label{eq:rp beta}
\beta_{ij}^+=\beta_{ij}-2\langle v_i,y\rangle\theta_j,  \qquad
\partial_i\theta_j={\textstyle\frac{1}{2}}\theta_i(\beta_{ij}^++\beta_{ij}).
\end{equation}
Following data determine a Ribaucour transform $f^+$ of a given
orthogonal net $f$ uniquely:
\newcounter{R}
\begin{list}{(R$_\arabic{R}$)}
{\usecounter{R}\setlength{\leftmargin}{1cm}}
\item a point $f^+(0)$;
\item values of $\theta_i$ on the coordinate axes $\PP_i$ for $1\le i\le m$,
 i.e., $m$ smooth functions $\theta_i\rest_{\PP_i}$ of one variable.
\end{list}

According to the general philosophy, iterating Ribaucour
transformations can be interpreted as adding an additional
(discrete) dimension to an orthogonal net. The situation arising
by adding two or three discrete dimensions is described in the
following fundamental theorem.

\begin{thm}[Permutability of Ribaucour transformations]
\label{thm:ortho permut}\quad

 1)
  Let  $f$ be an $m$-dimensional orthogonal net, and let
  $f^{(1)}$ and $f^{(2)}$ be its two Ribaucour transforms.
  Then there exists a one-parameter family of orthogonal nets
  $f^{(12)}$ that are Ribaucour transforms of both $f^{(1)}$ and $f^{(2)}$.
  Corresponding points of the four orthogonal nets $f$, $f^{(1)}$, $f^{(2)}$
  and $f^{(12)}$ are concircular.

2)
  Let  $f$ be an $m$-dimensional orthogonal net. Let
  $f^{(1)}$, $f^{(2)}$ and $f^{(3)}$ be its three Ribaucour
  transforms, and let three further orthogonal nets $f^{(12)}$,
  $f^{(23)}$ and $f^{(13)}$
  be given such that $f^{(ij)}$ is a simultaneous Ribaucour transform
  of $f^{(i)}$ and $f^{(j)}$. Then there exists
  generically a unique orthogonal net $f^{(123)}$
  that is a Ribaucour transform of $f^{(12)}$, $f^{(23)}$ and $f^{(13)}$.
  The net $f^{(123)}$ is uniquely defined by the condition that
  the corresponding points of $f^{(i)}$, $f^{(ij)}$, $f^{(ik)}$ and
  $f^{(123)}$ are concircular for any permutation $(ijk)$ of $(123)$.
\end{thm}

The theory of discrete orthogonal nets will unify the theories
of smooth orthogonal nets and of their transformations.
\medskip

{\bf M\"obius-geometric description of orthogonal nets.} Since
orthogonal nets belong to the M\"obius differential geometry, it
is useful to describe them with the help of the corresponding
apparatus (a sketch of which is given in Appendix to Chapter
\ref{Sect: smooth diffe geo}).
This has major conceptual and technical advantages.
First, this description linearizes the invariance group of
orthogonal nets, i.e., the M\"obius group of the sphere $\bbS^N$
(which can be considered as a compactification of $\bbR^N$ by a
point at infinity). Further, using the Clifford algebra model of
the M\"obius differential geometry enables us to give a frame
description of orthogonal nets, which turns out to be a key
technical device.

In this formalismus, the ambient space for points and hyperspheres
of the conformal $N$-sphere is the projectivized Minkowski space

$\bbP(\bbR^{N+1,1})$. The standard basis of the Minkowski space
$\bbR^{N+1,1}$ is denoted by $\{\ee_1, \ldots, \ee_{N+2}\}$. We
denote also $\ee_0=\frac{1}{2}(\ee_{N+2}+\ee_{N+1})$  and
$\ee_\infty=\frac{1}{2}(\ee_{N+2}-\ee_{N+1})$. The points of the
conformal $N$-sphere are elements of the projectivized light cone
$\bbP(\cn)$, i.e., straight line generators of $\cn$. The
Euclidean space $\bbR^N$ is identified, via
\begin{equation}\label{eq: pi0}
\pi_0:\; \bbR^N\ni f  \;\mapsto \;\hat{f}=f+\ee_0+|f|^2\ee_\infty\in
\bbQ_{\,0}^N\,,
\end{equation}
with the section $\bbQ_{\,0}^N$ of the cone $\cn$ by the affine
hyperplane $\{\xi_0=1\}$, where $\xi_0$ is the $\ee_0$-component
of $\xi\in\bbR^{N+1,1}$ in the basis
$\{\ee_1,\ldots,\ee_N,\ee_0,\ee_\infty\}$. Orientation preserving
Euclidean motions of $\bbR^N$ are represented as conjugations by
elements of $\cH_\infty$, the isotropy subgroup of $\ee_\infty$ in
${\rm Spin}^+(N+1,1)$.

It is not difficult to derive the following nice characterization of
orthogonal nets, due to Darboux (its second half follows directly from
eq. (\ref{eq: pi0})).
\begin{thm}\label{prop: os in K}
A conjugate net $f:\bbR^m\to\bbR^N$ is orthogonal, if and only if $|f|^2$
satisfies the same equation (\ref{eq:cn property}) as $f$ does, in other
words, if the corresponding $\hat{f}=\pi_0\circ f:\bbR^m\to\bbQ_{\,0}^N$
is a conjugate net in $\bbR^{N+1,1}$.
\end{thm}
As easily seen, metric coefficients $h_i=|\partial_i f|$ satisfy also
$h_i=|\partial_i\hat{f}|$. Hence, vectors
$\hat{v}_i=h_i^{-1}\partial_i\hat{f}=v_i+2\langle f,v_i\rangle\ee_\infty$
have the (Lorentz) length 1. Since
$\langle\hat{f},\hat{f}\rangle=0$,
one readily finds that $\langle\hat{f},\hat{v}_i\rangle=0$ and
$h_i=-\langle\partial_i\hat{v}_i,\hat{f}\rangle$.

\begin{thm}[Spinor frame of an O-net]\label{thm: os frame}
For an orthogonal net $f:\bbR^m\to\bbR^N$, i.e., for the corresponding
conjugate net $\hat{f}:\bbR^m\to\bbQ_{\,0}^N$, there exists a function
$\psi:\bbR^m\to\cH_\infty$ (called {\em a frame of} $\hat{f}$), such that
\begin{eqnarray}
   \hat{f} & = & \psi^{-1}\ee_0\psi,               \label{eq:fpoint}\\
   \hat{v}_i & = & \psi^{-1}\ee_i\psi,\quad 1\leq i\leq m,
                                                   \label{eq:fdirection}
\end{eqnarray}
and satisfying the system of differential equations:
\begin{equation}\label{eq: frame os}
   \partial_i\psi=-\ee_i\psi\hat{s}_i,\quad
   \hat{s}_i={\textstyle\frac{1}{2}}\partial_i\hat{v}_i,
   \quad 1\le i\leq m.
\end{equation}
\end{thm}
Note that for an orthogonal coordinate system $(m=N)$ the frame
$\psi$ is uniquely determined at any point by the requirements
(\ref{eq:fpoint}) and (\ref{eq:fdirection}).
\medskip

It is readily seen that the unit tangent vectors $\hat{v}_i$
satisfy eq. (\ref{eq:DV}) with the same rotation coefficients
$\beta_{ji}=\langle\partial_i\hat{v}_j,\hat{v}_i\rangle=
-\langle\partial_i\hat{v}_i,\hat{v}_j\rangle$. With the help of the
frame $\psi$ we extend the set of vectors $\{\hat{v}_i: 1\le i\le m\}$
to an orthonormal basis $\{\hat{v}_k: 1\le k\le N\}$ of
$T_{\hat{f}}\bbQ_{\,0}^N$:
\begin{equation}\label{eq:fbasis}
   \hat{v}_k=\psi^{-1}\ee_k\psi,\quad 1\leq k\leq N.
\end{equation}
Correspondingly, we extend the set of rotation coefficients according
to the formula
\[
 \beta_{ki}=\langle\partial_i\hat{v}_k,\hat{v}_i\rangle=
 -\langle\partial_i\hat{v}_i,\hat{v}_k\rangle=
 -\langle\partial_i\hat{v}_i,\psi^{-1}\ee_k\psi\rangle,
 \quad 1\le i\le m,\quad 1\le k\le N.
\]
Recall that we also have:
\[
 h_i=-\langle\partial_i\hat{v}_i,\hat{f}\rangle=
 -\langle\partial_i\hat{v}_i,\psi^{-1}\ee_0\psi\rangle,\quad 1\le i\le m.
\]
Thus, introducing vectors $S_i=\psi\hat{s}_i\psi^{-1}$, we have
the following expansion with respect to the vectors $\ee_k$:
\begin{equation}\label{eq: os S}
 S_i=\psi\hat{s}_i\psi^{-1}=
 {\textstyle\frac{1}{2}}\psi(\partial_i\hat{v}_i)\psi^{-1}=
 -{\textstyle\frac{1}{2}}\,\sum_{k\neq i}\beta_{ki}\ee_k+h_i\ee_\infty.
\end{equation}
It is easy to see that eq. (\ref{eq:DB}) still holds, if the range
of the indices is extended to all pairwise distinct with $1\le
i,j\le m$ and $1\le k\le N$, and that the orthogonality constraint
(\ref{eq:DL}) can be now put as
\begin{equation}\label{eq:ddL}
 \partial_i \beta_{ij}+\partial_j\beta_{ji}=
  -\sum_{k\neq i,j}\beta_{ki}\beta_{kj}\,.
\end{equation}
The system consisting of (\ref{eq:DB}), (\ref{eq:ddL}) carries the
name of the {\em Lam\'e system}.

\section[Moutard nets]{Moutard nets and their transformations}
 \label{Subsect: mn}

We introduce now Moutard nets \cite{Mou} without a geometric
motivation, but they will play an extremely important role in the
subsequent geometric considerations.
\begin{dfn}\label{def:M-net}
A map $f:\bbR^2\to\bbR^N$ is called an {\em M-net (Moutard net)},
if it satisfies the {\em Moutard differential equation}
\begin{equation}\label{eq:Mou}
\partial_1\partial_2 f=q_{12} f
\end{equation}
with some $q_{12}:\bbR^2\to\bbR$.
\end{dfn}
On the first sight, the notion of M-net is not related to that of
a conjugate net. In particular, there do not exist $M$-dimensional
M-nets with $M\ge 3$. However, the relation is easily established:
if $\nu:\bbR^2\to\bbR$ is any solution of the same Moutard
equation (\ref{eq:Mou}) (for instance, any component of the vector
$f$), then $y=\nu^{-1}f:\bbR^2\to\bbR^N$ is a special conjugate
net in $\bbR^N$:
\[
\partial_1\partial_2 y=-(\partial_2\log\nu)\partial_1 y
                       -(\partial_1\log\nu)\partial_2 y.
\]
Such nets $y$ are called {\em conjugate nets with equal
invariants}, and they were intensively studied in the classical
projective differential geometry \cite{Tzi}. Thus, in a projective
space the class of M-nets coincides with the class of conjugate
nets with equal invariants.

Following data determine an M-net $f$ uniquely:
\newcounter{M}
\begin{list}{(M$_\arabic{M}$)}
{\usecounter{M}\setlength{\leftmargin}{1cm}}
\item values of $f$ on the coordinate axes $\PP_1$, $\PP_2$, i.e., two
 smooth curves $f\rest_{\PP_i}$ with a common intersection point $f(0)$;
\item  a smooth function $q_{12}:\bbR^2\to\bbR$, having the
 meaning of the coefficient of the Moutard equation.
\end{list}

\begin{dfn}\label{dfn:mp}
Two M-nets $f,f^+:\bbR^2\to\bbR^N$ are called {\em Moutard transforms}
of one another, if they satisfy (linear) differential equations
\begin{eqnarray}
\partial_1 f^++\partial_1 f & = & p_1(f^+-f),  \label{eq:mp 1}\\
\partial_2 f^+-\partial_2 f & = & p_2(f^++f),  \label{eq:mp 2}
\end{eqnarray}
with some functions $p_1,p_2:\bbR^2\to\bbR$
(or similar equations with all plus and minus signs interchanged).
\end{dfn}
The functions $p_1$, $p_2$, specifying the Moutard transform,
have to satisfy (nonlinear) differential equations, which
express compatibility of eqs. (\ref{eq:mp 1}), (\ref{eq:mp 2})
with eq. (\ref{eq:Mou}):
\begin{eqnarray}
\partial_1 p_2=\partial_2 p_1 & = & -q_{12}+p_1p_2,  \label{eq:Mou 1}\\
q_{12}^+ & = & -q_{12}+2p_1p_2.                      \label{eq:Mou 2}
\end{eqnarray}
Following data determine a Moutard transform $f^+$ of a given M-net $f$:
\newcounter{T}
\begin{list}{(MT$_\arabic{T}$)}
{\usecounter{T}\setlength{\leftmargin}{1.1cm}}
\item a point $f^+(0)\in\bbR^N$;
\item values of the functions $p_i$ on the coordinate axes $\PP_i$ for
 $i=1,2$, i.e., two smooth functions $p_i\rest_{\PP_i}$ of one variable.
\end{list}
{\bf Classical formulation of the Moutard transformation.}
Due to the first equation in (\ref{eq:Mou 1}), for any Moutard transformation
there exists a function $\theta:\bbR^2\to\bbR$, unique up to a constant
factor, such that
\begin{equation}\label{eq:Mou a}
p_1=-\frac{\partial_1\theta}{\theta},\qquad
p_2=-\frac{\partial_2\theta}{\theta}.
\end{equation}
The last equation in (\ref{eq:Mou 1}) implies that $\theta$
satisfies eq. (\ref{eq:Mou}). This scalar solution of eq. (\ref{eq:Mou})
can be specified by its values on the coordinate axes $\PP_i\,$ ($i=1,2$),
which are readily obtained from the data (MT$_2$) by integrating the
corresponding eqs. (\ref{eq:Mou a}).
This establishes a bridge to the classical formulation of the Moutard
transformation (see, e.g., \cite{Mou, Tzi}).
They used to specify a Moutard transform $f^+$ of the solution $f$ of the
Moutard equation (\ref{eq:Mou}) by an additional scalar solution $\theta$
of this equation, via eqs. (\ref{eq:mp 1}), (\ref{eq:mp 2}) with
(\ref{eq:Mou a}). From these equations one can conclude that $f^+$ solves
the Moutard equation (\ref{eq:Mou}) with the transformed potential
\begin{equation}\label{eq:Mou q}
q_{12}^+=q-2\partial_1\partial_2\log\theta=
\frac{\partial_1\partial_2 \theta^+}{\theta^+},\qquad
\theta^+=\frac{1}{\theta}.
\end{equation}
In our formulation, the origin of the function $\theta$ becomes clear:
it comes from $p_1$, $p_2$ by integrating the system (\ref{eq:Mou a}).
Eq. (\ref{eq:Mou q}) is then
nothing but an equivalent form of eq. (\ref{eq:Mou 2}).

\begin{thm}[Permutability of Moutard transformations]
\label{thm:mp permut}\quad

 1)
  Let  $f$ be an M-net, and let $f^{(1)}$ and $f^{(2)}$ be its two
  Moutard transforms. Then there exists a one-parameter family of
  M-nets $f^{(12)}$ that are Moutard transforms of both $f^{(1)}$
  and $f^{(2)}$.

2)
  Let  $f$ be an M-net. Let $f^{(1)}$, $f^{(2)}$ and $f^{(3)}$ be its
  three Moutard transforms, and let three further M-nets $f^{(12)}$,
  $f^{(23)}$ and $f^{(13)}$ be given such that $f^{(ij)}$ is a
  simultaneous Moutard transform of $f^{(i)}$ and $f^{(j)}$. Then there
  exists generically a unique M-net $f^{(123)}$ that is a Moutard
  transform of $f^{(12)}$, $f^{(23)}$ and $f^{(13)}$.
\end{thm}

\section[Asymptotic nets]{Asymptotic nets and their transformations}
 \label{Subsect: an}

\begin{dfn} \label{dfn:an}
  A map $f:\bbR^2\to\bbR^3$ is called an {\em A-surface (a surface
  parametrized along asymptotic lines)}, if at any point the vectors
  $\partial_1^2 f$, $\partial_2^2 f$  lie in the tangent plane to the
  surface $f$, spanned by $\partial_1 f$, $\partial_2 f$.
\end{dfn}
Thus, the second fundamental form of an A-surface in $\bbR^3$ is
off-diagonal. Such a paramerization exists for a general surface
with a negative Gaussian curvature.
Definition \ref{dfn:an}, like the definition of conjugate nets, contains
projectively invariant notions only. Therefore A-surfaces belong
actually to the geometry of the three-dimensional projective space.
In our presentation, however, we will use for convenience additional
structures of $\bbR^3$ (Euclidean structure and the cross-product).
For the projective interpretation of these constructions, see
\cite{KP}. A convenient description of A-surfaces is provided
by the {\em Lelieuvre representation} which states:
there exists a unique (up to sign) normal field
$n:\bbR^2\to\bbR^3$ to the surface $f$ such that
\begin{equation}\label{eq:Lel}
\partial_1 f=\partial_1 n\times n,\qquad \partial_2 f=n\times \partial_2 n.
\end{equation}
Cross-differentiation of eq. (\ref{eq:Lel}) reveals that
$\partial_1\partial_2 n\times n=0$, that is, the Lelieuvre normal
field satisfies the {\em Moutard equation}
\begin{equation}\label{eq: as Mou}
\partial_1\partial_2 n=q_{12} n
\end{equation}
with some $q_{12}:\bbR^2\to\bbR$. This reasoning can be reversed: integration
of eqs. (\ref{eq:Lel}) with any solution $n:\bbR^2\to\bbR^3$ of the Moutard
equation generates an A-surface $f:\bbR^2\to\bbR^3$.
\begin{thm}
A-surfaces in $\bbR^3$ are in a one-to-one correspondence, via the Lelieuvre
representation (\ref{eq:Lel}), with M-nets in $\bbR^3$.
\end{thm}
An A-surface $f$ is reconstructed uniquely (up to a translation) from its
Lelieuvre normal field $n$. In its turn, an M-net $n$ is uniquely determined
by the initial data (M$_{1,2}$), which we denote in this context by (A$_{1,2}$):
\newcounter{A}
\begin{list}{(A$_\arabic{A}$)}
{\usecounter{A}\setlength{\leftmargin}{1cm}}
\item values of the Lelieuvre normal field on the coordinate axes
 $\PP_1$, $\PP_2$, i.e., two smooth curves $n\rest_{\PP_i}$ with
 a common intersection point $n(0)$;
\item  a smooth function $q_{12}:\bbR^2\to\bbR$, having the
 meaning of the coefficient of the Moutard equation for $n$.
\end{list}

\begin{dfn}\label{dfn:wp}
A pair of A-surfaces $f,f^+:\bbR^2\to\bbR^3$ is called a {\em
Weingarten pair}, if, for any $u\in\bbR^2$, the line
$[f(u),f^+(u)]$ is tangent to both surfaces $f$ and $f^+$
at the corresponding points. The surface
$f^+$ is called a {\em Weingarten transform} of the surface $f$.
\end{dfn}
It can be demonstrated that the Lelieuvre normal fields of a Weingarten
pair $f$, $f^+$ of A-surfaces satisfy (with the suitable choice of their
signs)
the following relation:
\begin{equation}\label{eq:Wei}
f^+-f=n^+\times n.
\end{equation}
Differentiating the last equation and using the
Lelieuvre formulas (\ref{eq:Lel}) for $f$ and for $f^+$, one easily
sees that the normal fields of a Weingarten pair are related by
(linear) differential equations:
\begin{eqnarray}
\partial_1 n^++\partial_1 n & = & p_1(n^+-n),  \label{eq:wp 1}\\
\partial_2 n^+-\partial_2 n & = & p_2(n^++n),  \label{eq:wp 2}
\end{eqnarray}
with some functions $p_1,p_2:\bbR^2\to\bbR$. Thus:
\begin{thm}
The Lelieuvre normal fields $n$, $n^+$ of a Weingarten pair $f$, $f^+$
of A-surfaces are Moutard transforms of one another.
\end{thm}
A Weingarten transform $f^+$ of a given A-surface $f$ is reconstructed
from a Moutard transform $n^+$ of the Lelieuvre normal field $n$. The data
necessary for this are the data (MT$_{1,2}$) for $n$:
\newcounter{W}
\begin{list}{(W$_\arabic{W}$)}
{\usecounter{W}\setlength{\leftmargin}{1cm}}
\item a point $n^+(0)\in\bbR^3$;
\item values of the functions $p_i$ on the coordinate axes $\PP_i$ for
 $i=1,2$, i.e., two smooth functions $p_i\rest_{\PP_i}$ of one variable.
\end{list}
Following statement is a direct consequence of Theorem \ref{thm:mp permut}.

\begin{thm}[Permutability of Weingarten transformations]
\label{thm:as permut}\quad

 1)
  Let  $f$ be an A-surface, and let $f^{(1)}$ and $f^{(2)}$ be its two
  Weingarten transforms. Then there exists a one-parameter family of
  A-surfaces $f^{(12)}$ that are Weingarten transforms of both $f^{(1)}$
  and $f^{(2)}$.

2)
  Let  $f$ be an A-surface. Let $f^{(1)}$, $f^{(2)}$ and $f^{(3)}$ be its
  three Weingarten transforms, and let three further A-surfaces $f^{(12)}$,
  $f^{(23)}$ and $f^{(13)}$ be given such that $f^{(ij)}$ is a simultaneous
  Weingarten transform of $f^{(i)}$ and $f^{(j)}$. Then there exists
  generically a unique A-surface $f^{(123)}$
  that is a Weingarten transform of $f^{(12)}$, $f^{(23)}$ and $f^{(13)}$.
  The net $f^{(123)}$ is uniquely defined by the condition that any its
  point lies in the tangent planes to $f^{(12)}$, $f^{(23)}$
  and $f^{(13)}$ at the corresponding points.
\end{thm}

\section[Surfaces with constant Gaussian curvature]
{Surfaces with constant negative Gaussian \\ curvature and
their transformations}
 \label{Subsect: Ksurf}

Up to now, we discussed special classes of coordinate systems in
space, or special parametrizations of a general surface. Now, we
turn to the discussion of several special classes of surfaces. The
distinctive feature of these classes is the existence of
transformations with certain permutability properties.

\begin{dfn} \label{dfn:ks}
  An A-surface $f:\bbR^2\to\bbR^3$ is called a {\em K-surface} (or a
  {\em pseudospheric surface}), if its Gaussian curvature $K$
  is constant, i.e., does not depend on $u\in\bbR^2$.
\end{dfn}

K-surfaces constitute one of the most prominent examples of
integrability in the differential geometry. One of the approaches
to their analytical study is based on the investigation of the
angle between asymptotic lines which is governed by the famous
sine-Gordon equation. This approach was transferred to the
discrete setting in \cite{BP1}, see also a presentation in
\cite{BMS1} based on the notion of consistency. In the present
paper, we take an alternative route, based on the study of the
Gauss map of K-surfaces. Following are the classical
characterization results.

\begin{thm}\label{Th: ks}\quad

 1) An A-surface $f:\bbR^2\to\bbR^3$ is a K-surface, if and only if the
 functions $|\partial_i f|=\alpha_i$ depend on $u_i$ only ($i=1,2$).

 2) The Lelieuvre normal field $n:\bbR^2\to\bbR^3$ of a K-surface with
 $K=-1$ takes values in the sphere $\bbS^2\subset\bbR^3$, thus coinciding
 with the Gauss map. Conversely, any M-net in the unit sphere $\bbS^2$
 is the Gauss map and the Lelieuvre normal field of a K-surface with $K=-1$.
 There holds: $|\partial_i n|=\alpha_i$ for $i=1,2$, with the
 same functions $\alpha_i=\alpha_i(u_i)$ as in 1).
\end{thm}

Thus, K-surfaces are in a one-to-one correspondence with M-nets in $\bbS^2$
(otherwise said, with Lorentz-harmonic $\bbS^2$-valued functions).
It is important to observe that the coefficient $q_{12}$  of the Moutard
equation satisfied by a a Lorentz-harmonic $\bbS^2$-valued function
$n:\bbR^2\to\bbR^3$ is completely determined by $n$, more precisely,
by its first order derivatives:
\begin{equation}\label{eq: Mou n q}
 q_{12}=\langle \partial_1\partial_2 n,n\rangle
 =-\langle \partial_1 n,\partial_2 n\rangle.
\end{equation}
Therefore, following data determine the Gauss map $n$ of a K-surface $f$:
\begin{list}{(K)}
{\setlength{\leftmargin}{0.8cm}}
\item values of the Gauss map on the coordinate axes
 $\PP_1$, $\PP_2$, i.e., two smooth curves $n\rest_{\PP_i}$ in $\bbS^2$
 intersecting at a point $n(0)$.
\end{list}
The K-surface $f$ is reconstructed from $n$ uniquely, up to a
translation, via formulas (\ref{eq:Lel}).

Historically the first class of surface transformations with
remarkable permutability properties was introduced by B\"acklund.

\begin{dfn}[B\"acklund transformation] \label{dfn:bt}
  A Weingarten pair of K-surfaces $f,f^+:\bbR^2\to\bbR^N$ forms a
  {\em B\"acklund pair}, if the distance $|f^+-f|$ is constant,
  i.e., does not depend on $u\in\bbR^2$.
\end{dfn}
The Gauss maps $n,n^+$ of a B\"acklund pair of K-surfaces $f,f^+$
are related by the Moutard transformation (\ref{eq:wp 1}),
(\ref{eq:wp 2}). From these equations there follows easily that
for a B\"acklund pair the quantity $\langle n,n^+\rangle$ is
constant; thus, the intersection angle of the tangent planes at
the corresponding points of a B\"acklund pair is constant.
Moreover, eq. (\ref{eq:Wei}) yields that this constant angle is
related to the constant distance between $f$ and $f^+$ via
\[
|f^+-f|^2=1-\langle n,n^+\rangle^2.
\]
The fact that $n,n^+\in\bbS^2$ allows one to express the
coefficients $p_1$, $p_2$ in eqs. (\ref{eq:wp 1}), (\ref{eq:wp 2}) in terms
of the solutions themselves:
\begin{eqnarray}
p_1 & =& \frac{\langle \partial_1 n,n^+\rangle-
\langle n,\partial_1 n^+\rangle}{2-2\langle n,n^+\rangle}=
\frac{\langle \partial_1 n,n^+\rangle}{1-\langle n,n^+\rangle},
\label{eq: ks Mou p1}\\
p_2 & = & \frac{\langle n,\partial_2 n^+\rangle
-\langle \partial_2 n,n^+\rangle}{2+2\langle n,n^+\rangle}=
\frac{-\langle \partial_2 n,n^+\rangle}{1+\langle n,n^+\rangle}.
\label{eq: ks Mou p2}
\end{eqnarray}
With these expressions, eqs. (\ref{eq:wp 1}), (\ref{eq:wp 2}) become a
compatible system of first order differential equations for $n^+$, therefore
the following data determine a B\"acklund transform $f^+$ of the given
K-surface $f$ uniquely:
\begin{list}{(B)}
{\setlength{\leftmargin}{0.8cm}}
\item a point $n^+(0)\in\bbS^2$.
\end{list}

\noindent
Permutability of B\"acklund transformations is due to Bianchi:
\begin{thm}[Permutability of B\"acklund transformations]
\label{thm:ks permut}\quad

  Let  $f$ be a K-surface, and let $f^{(1)}$ and $f^{(2)}$ be its two
  B\"acklund transforms. Then there exists a unique K-surface $f^{(12)}$
  which is simultaneously a B\"acklund transform of $f^{(1)}$ and of
  $f^{(2)}$. The fourth surface $f^{(12)}$ is uniquely defined by the
  properties $|f^{(12)}-f^{(1)}|=|f^{(2)}-f|$ and
  $|f^{(12)}-f^{(2)}|=|f^{(1)}-f|$, or, in terms of the Gauss maps,
  $\langle n^{(1)},n^{(12)}\rangle=\langle n,n^{(2)}\rangle$ and
  $\langle n^{(2)},n^{(12)}\rangle=\langle n,n^{(1)}\rangle$. Equivalently,
  there holds one of the relations
  \[
  n^{(12)}-n \;||\; n^{(1)}-n^{(2)}\qquad or \qquad
  n^{(12)}+n \;||\; n^{(1)}+n^{(2)}.
  \]
\end{thm}

We will see how the theory of discrete K-surfaces unifies the
theories of smooth K-surfaces and of their B\"acklund transformations.

\section[Isothermic surfaces]
{Isothermic surfaces and their transformations}
 \label{Subsect: isotherm}

Classically, theory of isothermic surfaces and their transformations was
considered as one of the highest achievements of the local differential
geometry.
\begin{dfn} \label{dfn:is}
  An O-surface $f:\bbR^2\to\bbR^N$ is called an {\em I-surface (isothermic
  surface)}, if
  its first fundamental form is conformal, possibly upon a
  re-parametrization
  $u_i\to\varphi_i(u_i)\,$  $(i=1,2)$ of the dependent variables,
  i.e., if at any point $u\in\bbR^2$ of the definition domain there holds
  $|\partial_1 f|^2/|\partial_2 f|^2=\alpha_1(u_1)/\alpha_2(u_2)$.
\end{dfn}
In other words, isothermic surfaces are characterized by the relations
$\partial_1\partial_2 f\!\in{\rm span}(\partial_1 f,\partial_2 f)$ and
\begin{equation}\label{eq:is prop}
\langle\partial_1 f,\partial_2 f\rangle=0,\quad
|\partial_1 f|^2=\alpha_1 s^2,\quad
|\partial_2 f|^2=\alpha_2 s^2,
\end{equation}
with some $s:\bbR^2\to\bbR_+$ and with the functions $\alpha_i$ depending
on $u_i$ only $(i=1,2$). These conditions may be equivalently represented as
\begin{equation}\label{eq:is prop1}
\partial_1\partial_2 f=(\partial_2\log s)\partial_1 f+
(\partial_1 \log s)\partial_2 f,\qquad
\langle\partial_1 f,\partial_2 f\rangle=0.
\end{equation}
The following property of isothermic surfaces can actually serve as their
another characterization.
\begin{thm}[Dual I-surface]
Let $f:\bbR^2\to\bbR^N$ be an isothermic surface. Then the $\bbR^N$-valued
one-form $df^*$ defined by
\begin{equation}\label{eq: is dual}
\partial_1 f^*=\alpha_1\frac{\partial_1 f}{|\partial_1 f|^2}=
\frac{\partial_1 f}{s^2},\qquad
\partial_2 f^*=-\alpha_2\frac{\partial_2 f}{|\partial_2 f|^2}=
-\frac{\partial_2 f}{s^2},
\end{equation}
is closed. The surface $f^*:\bbR^2\to\bbR^N$, defined (up to a translation)
by the integration of this one-form, is isothermic, with
\begin{equation}\label{eq:is prop dual}
\langle\partial_1 f^*,\partial_2 f^*\rangle=0,\quad
|\partial_1 f^*|^2=\alpha_1s^{-2},\quad
|\partial_2 f^*|^2=\alpha_2s^{-2}.
\end{equation}
The surface $f^*$ is called {\em dual} to the surface $f$, or the {\em
Christoffel transform} of the surface $f$.
\end{thm}
Another important class of transformations of isothermic surfaces build the
Darboux transformations.
\begin{dfn}[Darboux transformation] \label{dfn:dt}
  A Ribaucour transform \linebreak
  $f^+:\bbR^2\to\bbR^N$ of a given isothermic surface
  $f:\bbR^2\to\bbR^N$ is called a {\em Darboux transform}, if its first
  fundamental form is likewise conformal, possibly upon a reparametrization
  of the dependent variables, i.e., if at any point $u\in\bbR^2$ of the
  definition domain there holds
  $|\partial_1 f^+|^2/|\partial_2 f^+|^2=\alpha_1(u_1)/\alpha_2(u_2)$.
\end{dfn}

Introduce the corresponding function $s^+:\bbR^2\to\bbR_+$ for the
surface $f^+$, and denote $r=s^+/s:\bbR^2\to\bbR_+$. Thus,
$|\partial_i f^+|^2/|\partial_i f|^2=(s^+/s)^2=r^2$ for $i=1,2$. Comparing
this with the definition (\ref{eq:Riba}) of Ribaucour transformations we see
that one of the two possibilities holds:
\[
{\rm (i)}\;\; r_1=r_2=r,\quad {\rm or\quad(ii)}\;\; r_1=-r_2=-r.
\]
It can be demonstrated that in the case
(i) the surface $f^+$ is with necessity a M\"obius transformation of $f$;
we will not consider this trivial case further. In the case (ii) one gets
proper Darboux transformations. An important property of the
Darboux transformations is the following: the quantity
\begin{equation}\label{eq: Darboux int}
c=\frac{\ell^2}{ss^+}=\frac{|f^+-f|^2}{ss^+}
\end{equation}
is constant, i.e., does not depend on $u\in\bbR^2$. It is called a
{\em parameter} of the Darboux transformation. Following data determine
a Darboux transform $f^+$ of a given isothermic surface $f$ uniquely:
\newcounter{D}
\begin{list}{(D$_\arabic{D}$)}
{\usecounter{D}\setlength{\leftmargin}{1cm}}
\item a point $f^+(0)$;
\item a real number $c$, designated to be the constant (\ref{eq: Darboux int}).
\end{list}

As usual, we regard iterating a Darboux transformation as adding a third
(discrete) dimension to a two-dimensional isothermic  net.
The main classical result on Darboux transformations is the following
theorem, which assures that one can add several discrete
dimensions in a consistent way.

\begin{thm}[Permutability of Darboux transformations]
\label{thm:iso permut}\quad

  Let  $f$ be an isothermic surface, and let
  $f^{(1)}$ and $f^{(2)}$ be its two Darboux transforms, with parameters
  $c_1$ and $c_2$, respectively. Then there exists a unique isothermic
  surface $f^{(12)}$ which is simultaneously a Darboux transform of $f^{(1)}$
  with the parameter $c_2$ and a Darboux transform of $f^{(2)}$ with the
  parameter $c_1$. The surface $f^{(12)}$ is uniquely defined by the
  condition that the corresponding points of the four isothermic surfaces
  are concircular, and have a constant cross-ratio
  \[
  q\big(f,f^{(1)},f^{(12)},f^{(2)}\big)=\frac{c_1}{c_2}\,.
  \]
\end{thm}

{\bf Remark.} The real cross-ratio of four concircular points
$a,b,c,d\in\bbR^N$ may be defined as
\begin{equation}\label{eq:q}
q(a,b,c,d)=(a-b)(b-c)^{-1}(c-d)(d-a)^{-1},
\end{equation}
where the points are interpreted as elements of the Clifford algebra
${\cal C}\ell(\bbR^N)$. In more down-to-earth terms, since the four points
are co-planar, we may identify the plane where they lie with the complex
plane, and interpret in the above formula the symbols $a,b,c,d$ as complex
numbers.
\medskip

The theory of discrete isothermic surfaces unifies the theories
of smooth isothermic surfaces and of their Darboux transformations.
\smallskip

{\bf M\"obius-geometric characterization of isothermic surfaces
and their Darboux transformations.} It is easily checked that
conditions (\ref{eq:is prop}) are invariant with respect to affine
transformations of $\bbR^N$, as well as with respect to the
inversion $f\to f/\langle f,f\rangle$. In other words, the notion
of isothermic surfaces belongs to the M\"obius differential
geometry. The same holds for their Darboux transformations.
Therefore, it is useful to characterize these notions within the
M\"obius-geometric formalism. (However, the notion of the dual
surface, or Christoffel transformation, is essentially based on
the Euclidean structure of the ambient space $\bbR^N$.)

To find such a characterization, note first of all that eqs.
(\ref{eq:is prop1}) are equivalent to
\[
\partial_1\partial_2\hat{f}=(\partial_2\log s)\partial_1\hat{f}+
(\partial_1 \log s)\partial_2\hat{f}
\]
for the image $\hat{f}:\bbR^2\to\bbQ_{\,0}^N$ of $f$ in the
quadric $\bbQ_{\,0}^N\subset\cn$.

\begin{thm}\label{Th: iso Mob}
The lift $\hat{s}=s^{-1}\hat{f}:\bbR^2\to\cn$ of an isothermic surface
$f:\bbR^2\to\bbR^N$ to the light cone of
$\bbR^{N+1,1}$ satisfies the {\em Moutard equation}
\begin{equation}\label{eq: is Mou}
\partial_1\partial_2 \hat{s}=q_{12} \hat{s},
\end{equation}
with $q_{12}=s\partial_1\partial_2(s^{-1})$.

Conversely, given an M-net $\hat{s}:\bbR^2\to\cn$ in the light cone,
define $s:\bbR^2\to\bbR_+$ and $f:\bbR^2\to\bbR^N$ by
\[
\hat{s}=s^{-1}(f+\ee_0+|f|^2\ee_\infty)
\]
(so that $s^{-1}$ is the $\ee_0$-component, and $s^{-1}f$ is the
$\bbR^N$-part of $\hat{s}$ in the basis
$\ee_1,\ldots,\ee_N,\ee_0,\ee_\infty$). Then $f$ is an isothermic
surface.
\end{thm}
Note that the functions $\alpha_i=\langle\partial_i\hat{s},
\partial_i\hat{s}\rangle\,$ $(i=1,2$) depend on $u_i$ only and coincide
with the namesake functions from the definition (\ref{eq:is prop}).

Thus, we see that {\em I-surfaces are in a one-to-one
correspondence with M-nets in $\cn$, i.e., with Lorentz-harmonic
$\cn$-valued functions}.

Let us address the problem of minimal data which determine an
isothermic surface (i.e., an M-net in $\cn$) uniquely. Guided by
an analogy with the case of K-surfaces, one is tempted to think
that two arbitrary curves $\hat{s}\rest_{\PP_i}$ in $\cn$ would be
such data. However, as a consequence of the fact that now we are
dealing with the light cone $\cn=\{\langle
\hat{s},\hat{s}\rangle=0\}$ rather than with the sphere
$\bbS^2=\{\langle n,n\rangle=1\}$ as a quadric where M-nets live,
we cannot find an expression for $q_{12}$ in terms of the first
derivatives of $\hat{s}$ anymore; instead, one has:
\[
q_{12}=\frac{\langle \partial_1^3\hat{s},\partial_2\hat{s}\rangle}
{\langle \partial_1\hat{s},\partial_1\hat{s}\rangle}
=\frac{\langle \partial_2^3\hat{s},\partial_1\hat{s}\rangle}
{\langle \partial_2\hat{s},\partial_2\hat{s}\rangle}.
\]
This shows that the coordinate curves $\hat{s}\rest_{\PP_i}$ are
not arbitrary but rather subject to certain further conditions. We
leave the question on correct initial data for an isothermic
surface open.

Darboux pairs of isothermic surfaces are characterized in terms of their
lifts as follows.
\begin{thm}
The lifts $\hat{s},\hat{s}^+:\bbR^2\to\cn$ of a Darboux pair of
isothermic surfaces $f,f^+:\bbR^2\to\bbR^N$ are related by a {\em
Moutard transformation}, i.e., there exist two functions
$p_1,p_2:\bbR^2\to\bbR$ such that
\begin{equation}\label{eq: iso Mou}
\partial_1\hat{s}^++\partial_1\hat{s}=p_1(\hat{s}^+-\hat{s}),\qquad
\partial_2\hat{s}^+-\partial_2\hat{s}=p_2(\hat{s}^++\hat{s}).
\end{equation}
Conversely, for an M-net $\hat{s}$ in $\cn$, any Moutard transform
$\hat{s}^+$ with values in $\cn$ is a lift of a Darboux transform $f^+$
of the isothermic surface $f$.
\end{thm}
Note that the quantity $\langle\hat{s},\hat{s}^+\rangle$ is constant
(does not depend on $u\in\bbR^2$),
and is related to the parameter $c$ of the Darboux transformation:
$\langle\hat{s},\hat{s}^+\rangle=-c/2$. The formulas
\begin{equation}\label{eq: iso Mou p}
p_i=\frac{\langle\hat{s},\partial_i\hat{s}^+\rangle-
\langle \partial_i\hat{s},\hat{s}^+\rangle}{2\langle\hat{s},\hat{s}^+\rangle}=
-\frac{\langle \partial_i\hat{s},\hat{s}^+\rangle}
{\langle\hat{s},\hat{s}^+\rangle},\qquad i=1,2,
\end{equation}
make it apparent that a Moutard transform $\hat{s}^+$ is completely
determined by prescribing its value $\hat{s}^+(0)$ at one point. Indeed,
eqs. (\ref{eq: iso Mou}) with coefficients (\ref{eq: iso Mou p}) form
a compatible system of first order differential equations for
$\hat{s}^+:\bbR^2\to\cn$. Of course, data (D$_{1,2}$) are encoded in
$\hat{s}^+(0)$ in a straightforward manner.
\medskip

We summarize the considerations of these chapter in the following table:
\begin{eqnarray*}
 \mbox{\rm O-net}\;\; f\;\;{\rm in}\;\;\bbR^N & \leftrightarrow &
  {\rm conjugate\;\;net}\;\;\hat{f}\;\;{\rm in}\;\;
  \bbQ_{\,0}^N\simeq \bbP(\cn),\\
 \mbox{\rm A-net}\;\; f\;\;{\rm in}\;\;\bbR^3 & \leftrightarrow &
  \mbox{\rm M-net}\;\;n\;\;{\rm in}\;\;\bbR^3,\\
 \mbox{\rm K-net}\;\; f\;\;{\rm in}\;\;\bbR^3 & \leftrightarrow &
  \mbox{\rm M-net}\;\;n\;\;{\rm in}\;\;\bbS^2,\\
 \mbox{\rm I-net}\;\; f\;\;{\rm in}\;\;\bbR^N & \leftrightarrow &
  \mbox{\rm M-net}\;\;\hat{s}\;\;{\rm in}\;\;\cn.
\end{eqnarray*}
The next chapter will be devoted to discretizing all these relations.

\section*{Appendix: M\"obius-geometric formalism}
\label{Sect: Moebius}
\addcontentsline{toc}{section}{Appendix: M\"obius-geometric formalism}

The classical model of the $N$-dimensional M\"obius geometry,
which allows for a linear representation of M\"obius
transformations, lives in the Min\-kowski space $\bbR^{N+1,1}$,
i.e., in an $(N+2)$-dimensional space with the basis $\{\ee_1,
\ldots, \ee_{N+2}\}$, equipped with the Lorentz scalar product in
which $\ee_i$ are pairwise orthogonal and
\begin{equation*}
  \langle \ee_i,\ee_i \rangle =\left\{\begin{array}{rl}
  1, & 1\le i\le N+1,\\
  -1, & i=N+2.\end{array}\right.
\end{equation*}
(Although we use the same symbol for the Lorentz scalar product in
$\bbR^{N+1,1}$ and for the Euclidean scalar products in $\bbR^N$ and
in $\bbR^{N+1}$, its
concrete meaning should be always clear from the context.)

{\bf Points.} The space of points in the M\"obius geometry is
$\bbP(\cn)$ -- the space of straight line generators of the {\it
light cone}
\begin{equation}\label{eq:Lightcone}
  \cn=\left\{\xi\in\bbR^{N+1,1}:\;\langle \xi,\xi \rangle=0\right\}.
\end{equation}
The sphere $\bbS^N\subset\bbR^{N+1}$ (where we regard $\bbR^{N+1}$ as
spanned by $\ee_1,\ldots,\ee_{N+1}$) is identified with a section of
$\cn$ by the affine hyperplane $\langle \xi,\ee_{N+2}\rangle=-1$:
\[
\bbS^N \simeq \bbQ_{\,1}^N=\{\xi\in\cn:\; \xi_{N+2}=1\},
\]
which is a copy of $\bbS^N$ shifted by $\ee_{N+2}$:
\begin{equation}\label{eq:pi1}
\pi_1: \; \bbS^N\ni y \;\; \mapsto \;\;\hat{y}
=\hat{y}_{\,\rm Sph}= y+\ee_{N+2}\in\bbQ_{\,1}^N\,.
\end{equation}
Similarly, the Euclidean space $\bbR^N$ may be identified
with the section of $\cn$ by the affine hyperplane
$\langle \xi,\ee_{N+2}-\ee_{N+1}\rangle=-1$,
\[
 \bbR^N  \simeq  \bbQ_{\,0}^N=\{\xi\in\cn: \;\xi_{N+1}+\xi_{N+2}=1\},
\]
(Euclidean metric $d\xi_1^2+\ldots+d\xi_N^2$ being
induced from the ambient $\bbR^{N+1,1}$):
\begin{eqnarray}\label{eq:pi0}
\pi_0:\; \bbR^N\ni x & \mapsto & \hat{x}=\hat{x}_{\,\rm Euc}=
x+\tfrac{1}{2}(1-|x|^2)\,\ee_{N+1}+\tfrac{1}{2}(1+|x|^2)\,\ee_{N+2}
\qquad\nonumber\\
& & \quad =x+\ee_0+|x|^2\ee_\infty\in\bbQ_{\,0}^N\,.
\end{eqnarray}
Here the following notations are introduced:
\begin{equation*}
  \ee_0=\tfrac{1}{2}(\ee_{N+2}+\ee_{N+1})\quad {\rm and} \quad
  \ee_\infty=\tfrac{1}{2}(\ee_{N+2}-\ee_{N+1}).
\end{equation*}
Thus, the space $\bbR^N$ is modelled as a paraboloid in an
$(N+1)$-dimensional affine subspace through $\ee_0$ spanned by
$\ee_1,\ldots,\ee_N,\ee_\infty$. An important property of the Euclidean
identification (\ref{eq:pi0}) is:
\begin{equation}
\langle \hat{x}_1,\hat{x}_2\rangle=-\tfrac{1}{2}|x_1-x_2|^2,\qquad
\forall x_1,x_2\in\bbR^N.
\end{equation}
Note that the correspondence between $\bbQ_{\,1}^N$ and $\bbQ_{\,0}^N$
along the straight line generators of $\cn$ induces the stereographic
projection $\sigma:\bbS^N\to\bbR^N$,
\begin{equation*}
 y=\sigma^{-1}(x)=\frac{2}{1+|x|^2}\,x+\frac{1-|x|^2}{1+|x|^2}\,\ee_{N+1}.
\end{equation*}
In particular, the generators of $\cn$ through the points $\ee_0$
and $\ee_\infty$ correspond to the north pole $y_0=\ee_{N+1}$ and
the south pole $y_\infty=-\ee_{N+1}$ on $\bbS^N$, and to the zero
and the point at infinity in $\bbR^N$, respectively.
\medskip

{\bf Spheres.} A hypersphere $S$ in the conformal $N$-sphere is
the (non-empty) intersection of $\bbP(\cn)$ with a projectivized
hyperplane. Thus, $S$ can be put into a correspondence with the
point $S\in\bbP(\bbR^{N+1,1})$ polar to the above mentioned
hyperplane with respect to the light cone. This point is
space-like, i.e., any its representative $\hat{s}\in\bbR^{N+1,1}$
has $\langle \hat{s},\hat{s}\rangle>0$. There are various choices
of this representative which have nice geometric interpretations.

Fixing the $\ee_{N+2}$-component of $\hat{s}$ leads to the choice
\begin{equation}\label{eq:sph1}
\hat{s}=\hat{s}_{\,\rm Sph}=s+\ee_{N+2},\quad {\rm with}\quad
 s\in\bbR^{N+1}, \;\;
\langle s,s\rangle >1.
\end{equation}
This is related to a description of the hypersphere
$S\subset\bbS^N$ as the intersection of $\bbS^N$ with the
hyperplane $\langle s,y\rangle=1$ in $\bbR^{N+1}$. Indeed, the
latter equation is equivalent to $\langle
\hat{s},\hat{y}\rangle=0$ for $\hat{y}$ from eq. (\ref{eq:pi1}).
The point $s$ lies outside of $\bbS^N$, and $S\subset\bbS^N$ is
the contact set of $\bbS^N$ with the tangent cone to $\bbS^N$ with
the apex $s$. Also, $S\subset\bbS^N$ is the intersection of
$\bbS^N$ with the orthogonal $N$-sphere in $\bbR^{N+1}$ with the
center $s$ and the radius $\rho=(\langle s,s\rangle-1)^{1/2}$. See
Fig.~\ref{two circles}.

Fixing the $\ee_0$-component of $\hat{s}$ leads to the choice
\begin{equation}\label{eq:sph0}
\hat{s}=\hat{s}_{\,\rm Euc}=c+\ee_0+(|c|^2-r^2)\ee_\infty,\quad
 {\rm where}\quad c\in\bbR^N,
\end{equation}
related to the Euclidean description of the hypersphere
$S\subset\bbR^N$. Indeed, $\langle\hat{s},\hat{x}\rangle=0$ for
$\hat{x}$ from eq. (\ref{eq:pi0}) is equivalent to $|x-c|^2=r^2$.
This is the equation for points $x$ of the sphere $S$ with the
Euclidean center $c$ and the Euclidean radius $r$.

Still another choice is to fix the Lorentz norm of $\hat{s}$:
\begin{equation}\label{eq:sph can}
\hat{s}=\hat{s}_{\,\rm
M\ddot{o}b}=\pm\frac{\kappa}{\rho}(s+\ee_{N+2})=
\pm\frac{\kappa}{r}\big(c+\ee_0+(|c|^2-r^2)\ee_\infty\big)
\in\cn_{\,\kappa},
\end{equation}
where for any $\kappa>0$ the following quadric is introduced:
\begin{equation}\label{eq:Lkappa}
\cn_{\,\kappa}=\{\xi\in\bbR^{N+1,1}:\;\langle \xi,\xi\rangle=\kappa^2\}.
\end{equation}
Actually, this represents {\em oriented hyperspheres}, each choice
of the sign $\pm$ corresponding to one of the two possible
orientations of a given hypersphere. For any two (oriented)
hyperspheres $S_1$, $S_2$ the scalar product of their
representatives $\hat{s}_{\,\rm M\ddot{o}b}$ is a M\"obius
invariant: if $\kappa=1$, then
\[
\langle \hat{s}_1,\hat{s}_2\rangle=
\frac{1}{\rho_1\rho_2}\big(\langle s_1,s_2\rangle-1\big)=
\frac{1}{2r_1r_2}\big(r_1^2+r_2^2-|c_1-c_2|^2\big),
\]
is the cosine of the intersection angle of $S_1$, $S_2$, if they
intersect, and the inversive distance between $S_1$, $S_2$, otherwise.

\begin{figure}[htbp]
\begin{center}
\setlength{\unitlength}{0.05em}
\begin{picture}(500,320)(-250,-160)
\dashline[+30]{10}(0,0)(150,-62.132)
\dashline[+30]{10}(0,0)(150,0)
\dashline[+30]{10}(150,0)(150,-62.132)
\dashline[+30]{10}(0,0)(106.07,-106.07)
\dashline[+30]{10}(150,-62.132)(106.07,-106.07)
\path(88.767,-150)(211.233,150)
\thicklines
\put(0,0){\circle{300}}
\put(150,0){\circle*{5}}
\put(106.07,-106.07){\circle*{5}}
\put(150,-62.132){\circle{124.264}}
\put(150,-62.132){\circle*{5}}
\put(150,-80){$s$}
\put(163,-40){$\rho$}
\put(165,9){$y\in S$}
\put(118,118){$\bbS^N$}
\put(220,120){$\langle s,y\rangle=1$}
\end{picture}
\caption{Hypersphere $S\subset\bbS^N$, with an orthogonal $N$-sphere
through $S$}\label{Moeb circles}
\end{center}
\end{figure}
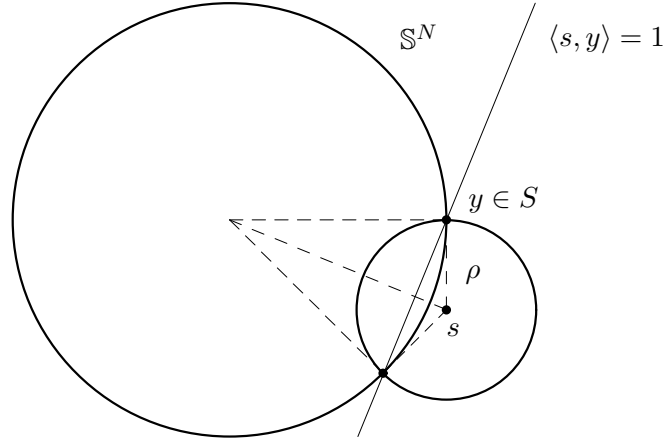
\medskip

{\bf Transformations.} Elements of the group ${\rm O}^+(N+1,1)$,
i.e., Lorentz transformations preserving the time-like direction,
induce M\"obius transformations on $\bbQ_{\,1}^N\simeq\bbS^N$.
Therefore we identify ${\rm O}^+(N+1,1)$ with the group ${\rm
M}(N)$ of M\"obius transformations of $\bbS^N$. Similarly, Lorentz
transformations which fix $\ee_\infty$, induce Euclidean motions
on $\bbQ_{\,0}^N\simeq\bbR^N$, and therefore we identify the
corresponding isotropy subgroup ${\rm O}^+_\infty(N+1,1)$ with
${\rm E}(N)$, the group of Euclidean motions of $\bbR^N$.

It is convenient to work with spinor representations of these
groups. Recall that the Clifford algebra ${\cal C}\ell(N+1,1)$ is an
algebra over $\bbR$ with  generators
$\ee_1,\ldots,\ee_{N+2}\in\bbR^{N+1,1}$ subject to the relation
\begin{equation*}  \label{eq:cliffrule}
  \xi\eta+\eta\xi=-2\langle \xi,\eta\rangle{\mathbf 1}=
  -2\langle \xi,\eta\rangle, \qquad
  \forall \xi,\eta\in\bbR^{N+1,1}.
\end{equation*}
This implies that $\xi^2=-\langle \xi,\xi \rangle$, therefore any
vector $\xi\in\bbR^{N+1,1}\setminus\cn$ has an inverse
$\xi^{-1}=-\xi/\langle \xi,\xi \rangle$. The multiplicative group
generated by invertible vectors is called the {\it Clifford
group}. We need its subgroup generated by unit space-like vectors:
\begin{equation*}
  \cG={\rm Pin}^+(N+1,1)=\{\psi=\xi_1\cdots \xi_n:\; \xi_i^2=-1\},
\end{equation*}
and its subgroup generated by vectors orthogonal to $\ee_\infty$:
\begin{equation*}
  \cG_\infty={\rm Pin}_\infty^+(N+1,1)=\{\psi=\xi_1\cdots \xi_n:\;
  \xi_i^2=-1,\;\langle \xi_i,\ee_\infty\rangle=0\}.
\end{equation*}
These groups act on $\bbR^{N+1,1}$ by twisted conjugations:
$A_\psi(\eta)=(-1)^n\psi^{-1}\eta\psi$.
In particular, for a vector $\xi$ with $\xi^2=-1$ one has:
\begin{equation*}
  A_\xi(\eta)=-\xi^{-1}\eta\xi=\xi\eta\xi=
  \eta-2\langle\xi,\eta\rangle \xi,
\end{equation*}
which is the reflection in the hyperplane orthogonal to $\xi$.
Thus, $\cG$ is generated by reflections, while $\cG_\infty$ is
generated by reflections which fix $\ee_\infty$, and therefore
leave $\bbQ_{\,0}^N$ invariant. Actually, $\cG$ is a double cover
of ${\rm O}^+(N+1,1)\simeq{\rm M}(N)$, while $\cG_\infty$ is a
double cover of ${\rm O}^+(N+1,1)_\infty\simeq{\rm E}(N)$.
Orientation preserving transformations from $\cG$, $\cG_\infty$
form the subgroups
\[
  \cH={\rm Spin}^+(N+1,1),\qquad
  \cH_\infty={\rm Spin}_\infty^+(N+1,1),
\]
which are singled out by the condition that the number $n$ of
vectors $\xi_i$ in the multiplicative representation of their
elements $\psi=\xi_1\cdots\xi_n$ is even. The Lie algebras of the
Lie groups $\cH$ and $\cH_\infty$ consist of bivectors:
\begin{eqnarray*}
  \gh & \!= \!& {\rm spin}(N+1,1)\;=\;{\rm span}\Big\{\ee_i\ee_j\,:\,
  i,j\in\{0,1,\ldots,N,\infty\},\;i\neq j\Big\},\\
  \gh_\infty & \!=\! & {\rm spin}_\infty(N+1,1)\;=\;
  {\rm span}\Big\{\ee_i\ee_j\,:\,
  i,j\in\{1,\ldots,N,\infty\},\;i\neq j\Big\}.
\end{eqnarray*}

\chapter{Discrete differential geometry} \label{Chap: discr geom}

For functions on $\bbZ^M$, we define translation and difference
operators in a standard manner:
\[
  (\tau_i f)(u)=f(u+e_i),\quad
  (\delta_i f)(u)=f(u+e_i)-f(u),
\]
where $e_i$ is the $i$-th coordinate vector of $\bbZ^M$. We use the
same notation for (discrete) $s$-dimensional coordinate planes,
\[
\PP_{i_1\ldots i_s}=\big\{u\in\bbZ^M: \;u_i=0\quad{\rm for}\quad i\neq
i_1,\ldots,i_s\big\},
\]
as in the continuous case.

\section{Discrete conjugate nets}
 \label{Sect: subsect discr conj}

The following definition is due to Sauer \cite{Sauer} for $M=2$, and to
Doliwa and Santini \cite{DS1} for general $M$. See \cite{D1, D2, D3, D5, DSM,
MDS} for further relevant developments.
\begin{dfn}\label{dfn:dcn}
  A map $f:\bbZ^M\to\bbR^N$ is called an
  $M$-dimensional {\em Q-net (quadrilateral net, or
  discrete conjugate net)} in $\bbR^N$, if any of its elementary
  quadrilaterals is planar, i.e., if at any $u\in\bbZ^M$
  and for all pairs $1\leq i\neq j\leq M$ the four points
  $f$, $\tau_i f$, $\tau_j f$, and $\tau_i\tau_j f$ are co-planar.
\end{dfn}
Note that this definition actually belongs to the projective geometry,
as it should. To understand what restrictions does this condition impose,
we consider various values of $M$.
\smallskip

{\itbf M=2:} {\bf discrete surface parametrized by conjugate
lines.} Suppose two coordinate lines, $f\rest_{\PP_1}$ and
$f\rest_{\PP_2}\,,$ on a Q-surface $f:\bbZ^2\to\bbR^3$ are given.
To extend the surface into the quadrant $\bbZ_+^2$, say, one
proceeds by induction whose step consists of choosing $f_{12}$ in
the plane spanned by $f$, $f_1$ and $f_2$, provided the latter
three points are known (here we write $f_i$, $f_{ij}$ for $\tau_i
f$, $\tau_i\tau_j f$, etc.). The planarity condition is equivalent
to the relation
\[
\delta_1\delta_2 f=c_{21}\delta_1 f+c_{12}\delta_2 f.
\]
So, one has two free real parameters $c_{21}$, $c_{12}$ on each
such step. It is convenient to think of these parameters as
attached to the elementary square $(u,u+e_1,u+e_1+e_2,u+e_2)$ of
the lattice $\bbZ^2$. Thus, one can define a Q-surface $f$ by
prescribing its two coordinate lines $f\rest_{\PP_1},\,$
$f\rest_{\PP_2}\,$, and two real-valued functions $c_{12}, c_{21}$
defined on all elementary squares of $\bbZ^2$.

Actually, the combinatorics of Q-surfaces may well be more
complicated than that of $\bbZ^2$. Indeed, Definition
\ref{dfn:dcn} can be literally extended to maps
$f:V(\cD)\to\bbR^N$, where $V(\cD)$ is the set of vertices of an
arbitrary {\em quad-graph} $\cD$. A quad-graph is a strongly
regular cell decomposition of a surface with all quadrilateral
faces. We will later need also the notation $E(\cD)$ and $F(\cD)$
for the sets of edges, resp. faces, of a quad-graph $\cD$. As we
will show, the integrable nature (multi-dimensional consistency)
of the Q-nets gives an opportunity to realize $\cD$ as a surface
in some $\bbZ^M$ and to work only on this larger (but simpler)
definition domain.

\begin{figure}[htbp]
\begin{center}
\setlength{\unitlength}{0.05em}
\begin{picture}(200,240)(0,0)
 \put(0,0){\circle*{15}}    \put(150,0){\circle*{15}}
 \put(0,150){\circle*{15}}  \put(150,150){\circle*{15}}
 \put(50,50){\circle*{15}} \put(50,200){\circle*{15}}
 \put(200,50){\circle*{15}}
 \put(200,200){\circle{15}}
 \path(0,0)(150,0)       \path(0,0)(0,150)
 \path(150,0)(150,150)   \path(0,150)(150,150)
 \path(0,150)(50,200)    \path(150,150)(194,194)
 \path(50,200)(192.5,200)
 \path(200,192.5)(200,50) \path(200,50)(150,0)
 \dashline[+30]{10}(0,0)(50,50)
 \dashline[+30]{10}(50,50)(50,200)
 \dashline[+30]{10}(50,50)(200,50)
 \put(-30,-5){$f$}
 \put(-35,145){$f_3$} \put(215,45){$f_{12}$}
 \put(165,-5){$f_1$} \put(160,140){$f_{13}$}
 \put(15,50){$f_2$}  \put(10,205){$f_{23}$}
 \put(215,200){$f_{123}$}
\end{picture}
\caption{3D system on an elementary cube}\label{Fig:cube eq}
\end{center}
\end{figure}
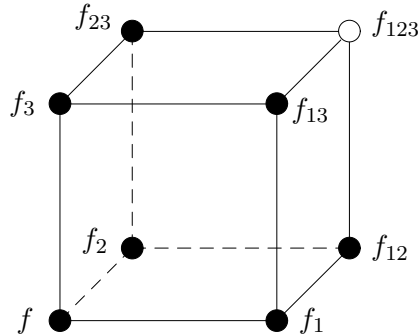

{\itbf M=3:} {\bf basic 3D system.} Suppose that three coordinate
surfaces of a three-dimensional Q-net $f$ are given, that is,
$f\rest_{\PP_{12}}$, $f\rest_{\PP_{23}}$ and $f\rest_{\PP_{13}}$.
Of course, each one of them is a Q-surface. To extend the net into
the octant $\bbZ_+^3$, one proceeds by induction whose step
consists of determining $f_{123}$, provided $f$, $f_i$ and
$f_{ij}$ are known for all $1\leq i\neq j\leq 3$. The point
$f_{123}$ has to lie in three planes $\tau_i\Pi_{jk}$ $(1\leq
i\leq 3)$, where $\tau_i\Pi_{jk}$ is the plane passing through
three points $(f_i,f_{ij},f_{ik})$. This condition determines
$f_{123}$ uniquely. Indeed, all three planes $\tau_1\Pi_{23}$,
$\tau_2\Pi_{13}$ and $\tau_3\Pi_{12}$ belong to the
three-dimensional affine space through the point $f$ spanned by
the vectors $\delta_i f$ $(1\leq i\leq 3)$, and therefore
generically these planes intersect at exactly one point. An
elementary construction step of a three-dimensional Q-net out of
its three coordinate surfaces, i.e., finding the eighth vertex of
an elementary hexahedron out of the known seven vertices, is
symbolically represented on Fig.~\ref{Fig:cube eq}. This is the
picture we have in mind when thinking (and speaking) about
three-dimensional systems. Of course, one can also give an
analytic formulation of this picture. This is done as follows. The
characteristic property of a Q-net is:
\begin{equation} \label{eq:dcn property}
    \delta_i\delta_j f=c_{ji}\delta_i f+c_{ij}\delta_j f, \quad i\neq j.
\end{equation}
Here, as before, functions $c_{ij}$, $c_{ji}$, as well as equation
(\ref{eq:dcn property}) itself, are thought of as defined on
elementary squares of $\bbZ^3$ parallel to the coordinate plane
$\PP_{ij}$. Six such equations, attached to six facets of an
elementary cube of $\bbZ^3$, form the three-dimensional system
encoded on Fig.~\ref{Fig:cube eq}. Here the numbers $\{c_{jk}\}$
on the facets adjacent to $f$ are considered as known, while the
numbers $\{\tau_ic_{jk}\}$ on the facets adjacent to $f_{123}$ are
uniquely defined by the compatibility of eqs. (\ref{eq:dcn
property}) on all six facets. In other words, it is required that
$\delta_i(\delta_j\delta_k f)$ does not depend on the permutation
$(i,j,k)$ of the indices $(1,2,3)$. This compatibility condition
gives:
\begin{equation}\label{eq:dcn c}
  \delta_i c_{jk}=(\tau_k c_{ij})c_{jk}+(\tau_k c_{ji})c_{ik}
  -(\tau_i c_{jk})c_{ik}, \quad  i\neq j\neq k\neq i.
\end{equation}
Note that equations for $c_{jk}$ split off from equations for $f$;
they constitute a system of 6 (linear) equations for 6 unknown variables
$\tau_i c_{jk}$ in terms of the known ones $c_{jk}$.
The resulting map $\{c_{jk}\}\mapsto \{\tau_i c_{jk}\}$ is birational.
Sometimes it is this map that is considered as the fundamental
3D system encoded on Fig.~\ref{Fig:cube eq}.

\begin{figure}[htbp]
\begin{center}
\setlength{\unitlength}{0.07em}
\begin{picture}(200,220)(-100,-90)

 \drawline(15,-20)(50,0)(50,47)
 \drawline(47,50)(0,50)(-35,30)(-35,-20)(15,-20)(15,30)(47,48.5)
 \drawline(15,30)(-35,30)
 \dashline{4}(-35,-20)(0,0)(0,50)\dashline{4}(0,0)(50,0)
 \drawline(30,-90)(131,-32)
 \drawline(135,-26)(135,116)
 \drawline(131,120)(-11,120)
 \drawline(-19,118)(-120,60)(-120,-90)(30,-90)(30,56)
 \drawline(34,62)(131,118)
 \drawline(26,60)(-120,60)
 \dashline{4}(-120,-90)(-15,-30)(-15,116)
 \dashline{4}(-15,-30)(131,-30)
  \dashline{2}(0,0)(-15,-30)
  \dashline{2}(-35,-20)(-120,-90)
  \dashline{2}(50,0)(132,-29)
  \dashline{2}(0,50)(-14,116)
  \dashline{2}(15,-20)(30,-90)
  \dashline{2}(-35,30)(-120,60)
  \dashline{2}(53,51.5)(131,117.5)
  \dashline{2}(15,30)(26,56)

  \put(-35,-20){\circle*{8}}       
  \put(15,-20){\circle*{8}}        
  \put(0,0){\circle*{8}}           
  \put(-35,30){\circle*{8}}        
  \put(50,0){\circle*{8}}          
  \put(0,50){\circle*{8}}          
  \put(15,30){\circle*{8}}         
  \put(50,50){\circle{8}}          

  \put(-120,-90){\circle*{10}}     
  \put(-15,-30){\circle*{9}}       
  \put(30,-90){\circle*{10}}       
  \put(-120,60){\circle*{10}}      
  \put(135,-30){\circle{10}}       
  \put(-15,120){\circle{10}}       
  \put(30,60){\circle{10}}         
  \put(130,115){$\Box$}            

  \put(-45,-15){$f$}
  \put(0,-15){$f_1$}
  \put(-18,6){$f_2$}
  \put(-50,20){$f_3$}
  \put(58,3){$f_{12}$}
  \put(-5,20){$f_{13}$}
  \put(-25,50){$f_{23}$}
  \put(58,44){$f_{123}$}
  \put(-143,-84){$f_4$}
  \put(5,-83){$f_{14}$}
  \put(-12,-41){$f_{24}$}
  \put(-143,49){$f_{34}$}
  \put(145,-28){$f_{124}$}
  \put(10,71){$f_{134}$}
  \put(-9,108){$f_{234}$}
  \put(140,103){$f_{1234}$}
\end{picture}
\caption{4D consistency of 3D systems}\label{hypercube}
\end{center}
\end{figure}
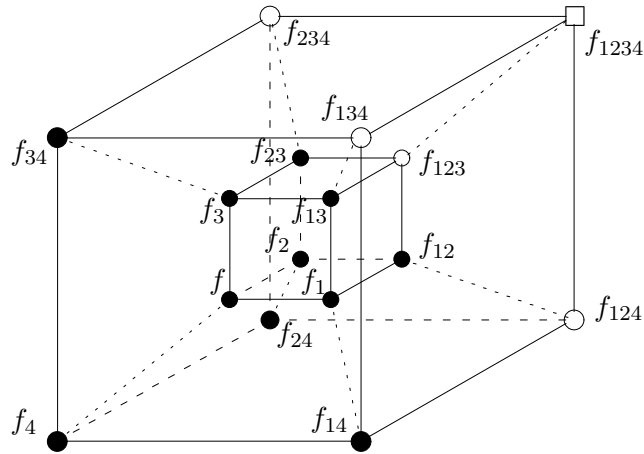

{\itbf M \boldmath$\ge$ 4:} {\bf consistency.} Turning to the case
$M\ge 4$, we see that one can prescribe all two-dimensional
coordinate surfaces of a Q-net, i.e., $f\rest_{\PP_{ij}}$ for all
$1\leq i<j\leq M$. Indeed, these data are clearly independent, and
one can construct the whole net from them. In doing so, one
proceeds by induction, again. The inductive step is essentially
three-dimensional and consists of determining $f_{ijk}$, provided
$f$, $f_i$ and $f_{ij}$ are known. However, this inductive process
works, only if one does not encounter contradictions. To see the
possible source of contradictions, consider in detail the case of
$M=4$; higher dimensions do not add anything new. From $f$, $f_i$
and $f_{ij}$ one determines all $f_{ijk}$ uniquely. After that,
one has, in principle, four different ways to determine
$f_{1234}$, from four 3-dimensional cubes adjacent to this point;
see Fig.~\ref{hypercube}. A remarkable property of Q-nets is that
these four values for $f_{1234}$ automatically coincide. We call
this property the 4D consistency.
\begin{dfn}
A 3D system is called 4D consistent, if it can be imposed on all
three-dimensional facets of an elementary hypercube of $\bbZ^4$.
\end{dfn}

The following fundamental theorem is due to Doliwa and Santini
\cite{DS1}.
\begin{thm}\label{prop:dcn consistency}
The 3D system governing Q-nets is 4D-consistent.
\end{thm}
{\bf Proof.} Actually, the statement of the theorem is about the
properties of the map $\{c_{jk}\}\mapsto\{\tau_ic_{jk}\}$. For
such maps with the fields on 2D plaquettes (here each plaquette
carries two fields) the 4D consistency means that the two values
$\tau_i(\tau_jc_{k\ell})$ and $\tau_j(\tau_ic_{k\ell})$ coincide
for any permutation $(i,j,k,\ell)$ of the indices $(1,2,3,4)$.
However, an algebraic proof of this claim could be hardly
performed without help of a computer system for symbolic
computations. A geometric approach, dealing with this system
augmented by the fields $f\in\bbR^N$ on the vertices, allows us to
perform a conceptual proof, free of computations. Indeed, the map
$\{c_{jk}\}\mapsto\{\tau_ic_{jk}\}$ does not depend on the
dimension $N$ of the space where $f$ lies. Therefore, we are free
to assume that $N\ge 4$. This will enable us to use geometric
``general position'' arguments.

In the construction above, the four values in question are
\[
f_{1234}=\tau_1\tau_2\Pi_{34}\cap\tau_1\tau_3\Pi_{24}\cap
\tau_1\tau_4\Pi_{23}\,,
\]
and three other ones obtained by cyclic shifts of indices. In
general position, the affine space $\Pi_{1234}$ through $f$
spanned by $\delta_i f$ $(1\leq i\leq 4)$, is four-dimensional. We
prove that all six planes $\tau_i\tau_j\Pi_{k\ell}$ in this
four-dimensional space intersect generically at exactly one point
(which will be then $f_{1234}$). It is easy to understand that the
plane $\tau_i\tau_j\Pi_{k\ell}$ is the intersection of two
three-dimensional subspaces $\tau_i\Pi_{jk\ell}$ and
$\tau_j\Pi_{ik\ell}$. Here, of course, the subspace
$\tau_i\Pi_{jk\ell}$ is the one through the four points
$(f_i,f_{ij},f_{ik},f_{i\ell})$, or, equivalently, the one through
$f_i$ spanned by $\delta_j f_i$, $\delta_k f_i$ and $\delta_\ell
f_i$. Now the intersection in question can be alternatively
described as the intersection of the four three-dimensional
subspaces $\tau_1\Pi_{234}$, $\tau_2\Pi_{134}$, $\tau_3\Pi_{124}$
and $\tau_4\Pi_{123}$ of one and the same four-dimensional space
$\Pi_{1234}$. This intersection consists in the generic case of
exactly one point. $\Box$
\medskip

The $M$-dimensional consistency for $M>4$ is defined and proved
analogously. Actually, it follows from the 4-dimensional consistency.
\medskip

On the level of formulas we have for $M\ge 4$ the system
(\ref{eq:dcn property}), (\ref{eq:dcn c}), where now all indices $i,j,k$
vary between 1 and $M$. This system consists of interrelated
three-dimensional building blocks: for any triple
of pairwise different indices $(i,j,k)$ the equations involving
these indices only form a closed subset. The $M$-dimensional consistency
of this system means that all three-dimensional building blocks can be
imposed without contradictions. A set of initial data which determines
a solution of the system (\ref{eq:dcn property}), (\ref{eq:dcn c}),
consists of
\newcounter{dQ}
\begin{list}{(Q$_\arabic{dQ}^\Delta$)}
{\usecounter{dQ}\setlength{\leftmargin}{1cm}}
\item values of $f$ on the coordinate axes $\PP_i$ for
  $1\le i\le M$;
\item values of $c_{ij}$, $c_{ji}$ on all elementary squares of
  the coordinate plane $\PP_{ij}$, for $1\le i<j\le M$.
\end{list}

{\bf Quadratic reduction of Q-nets.} An important observation made
by Doliwa \cite{D2} is that quadrilateral nets can be consistently
restricted to an arbitrary quadric in $\bbR^N$. The importance of
this resides on the fact that many of the geometrically relevant
nets turn out to be reductions of quadrilateral nets to some
quadrics or to intersections of quadrics. The quadratic reduction
is based on the following fundamental claim:
\begin{thm}\label{Th:quadratic}
If seven points $f$, $f_i$ and  $f_{ij}\,$ ($1\le i<j\le 3$) of an
elementary hexahedron of a quadrilateral net belong to a quadric
$\,\bbQ\subset\bbR^N$, then so does the eighth point $f_{123}$.
\end{thm}
{\bf Proof.} A deep reason for this is the fact, well known in the
projective geometry of the 19-th century: for seven points of a
three-dimensional projective space in a general position, one can
find the eighth {\em associated point} which belongs to any
quadric through the original seven points (these quadrics form a
two-parameter linear family, which contains $\bbQ$). Under
conditions of Theorem \ref{Th:quadratic}, we have three
(degenerate) quadrics through the original seven points: the pairs
of planes $\Pi_{jk}\cup\tau_i\Pi_{jk}$ for $i=1,2,3$. Clearly,
their eighth intersection point is
$f_{123}=\tau_1\Pi_{23}\cap\tau_2\Pi_{31}\cap\tau_3\Pi_{12}$, and
this has to be the associated point. $\Box$
\medskip

{\bf Corollary.} {\em If the coordinate surfaces
$f\rest_{\PP_{ij}}$ of a Q-net $f:\bbZ^M\to\bbR^N$ lie in a
quadric $\bbQ$, then so does the whole of $f$.}
\medskip

{\bf Alternative analytic description of Q-nets.} In a complete
analogy with the smooth case, one can give a (non-local)
description of discrete conjugate nets, with somewhat simpler
equations. Given the plaquette functions $c_{ij}$, define
quantities $h_j$, attached to the edges parallel to the $j$-th
coordinate axes, as solutions of the system of difference
equations
\begin{equation}\label{eq:dcn h}
\delta_ih_j=c_{ij}h_j\,, \qquad i\neq j,
\end{equation}
whose compatibility is assured by eq. (\ref{eq:dcn c}). Introduce
vectors $v_j=h_j^{-1}\delta_j f$, attached to the same edges.
There follows from (\ref{eq:dcn property}) and (\ref{eq:dcn h})
that these vectors satisfy the following difference equations:
\begin{equation}\label{eq:dcn v}
\delta_iv_j=\frac{h_i}{\tau_ih_j}c_{ji}v_i\,, \qquad i\neq j.
\end{equation}
Define the {\em discrete rotation coefficients} (attached to plaquettes)
as
\begin{equation}\label{eq: dcn rot}
\beta_{ji}=\frac{h_i}{\tau_ih_j}c_{ji}\,.
\end{equation}
Then we end up with the following system:
\begin{eqnarray}
    \delta_i f  & = & h_i v_i\,,                         \label{eq:dDX}\\
    \delta_i v_j & = & \beta_{ji} v_i\,,\quad i\neq j,   \label{eq:dDV}\\
    \delta_i h_j & = & (\tau_jh_i)\beta_{ij}\,,\quad i\neq j.
    \label{eq:dDH}
\end{eqnarray}
Rotation coefficients satisfy a closed system of difference
equations ({\em discrete Darboux system}), which follows from eqs.
(\ref{eq:dcn c}) upon substitution (\ref{eq: dcn rot}), or
otherwise can be derived as compatibility conditions of the linear
difference equations (\ref{eq:dDV}):
\begin{equation}\label{eq:dDB}
\delta_i\beta_{kj} = (\tau_j\beta_{ki})\beta_{ij}\,,
    \quad i\neq j\neq k\neq i.
\end{equation}
This system, which is considerably simpler than (\ref{eq:dcn c}),
was first derived in \cite{BK} without any relation to geometry. A
geometric interpretation in terms of Q-nets was found in \cite{DS1}.
It should be mentioned that eqs. (\ref{eq:dDH}) and (\ref{eq:dDB})
are implicit, but can be easily solved for the shifted variables,
resulting in
\begin{eqnarray}
    \tau_i h_j & = & \frac{h_j+h_i\beta_{ij}}{1-\beta_{ij}\beta_{ji}}\,,
    \quad i\neq j,  \label{eq:dDH expl}\\
    \tau_i\beta_{kj} & = &  \frac{\beta_{kj}+\beta_{ki}\beta_{ij}}
    {1-\beta_{ij}\beta_{ji}}\,,  \quad i\neq j\neq k\neq i.
    \label{eq:dDB expl}
\end{eqnarray}
The last formula defines an explicit rational 3D map
$\{\beta_{kj}\} \mapsto\{\tau_i\beta_{kj}\}$. Like the map
$\{c_{kj}\} \mapsto\{\tau_ic_{kj}\}$ for the local plaquette
coefficients, the map (\ref{eq:dDB expl}) is 4D-consistent, but
now this can be checked by an easy computation ``by hands''.
\medskip

{\bf Transformations of Q-nets.} A natural generalization of
Definition \ref{dfn:jp} would be the following one.
\begin{dfn}\label{dfn:dcn jp}
  A pair of $m$-dimensional Q-nets $f,f^+:\bbZ^m\to\bbR^N$
  is called a {\em Jonas pair}, if four points $f$, $\tau_i f$,
  $f^+$ and $\tau_i f^+$ are co-planar at any point $u\in\bbZ^m$
  and for any $1\leq i\leq m$. The net $f^+$ is called
  a {\em Jonas transform} of the net $f$.
\end{dfn}
But actually this relation can be re-phrased as follows: set
$F(u,0)=f(u)$ and $F(u,1)=f^+(u)$, then
$F:\bbZ^m\times\{0,1\}\to\bbR^N$ is an $M$-dimensional Q-net,
where $M=m+1$. Thus, in the discrete case there is no difference
between conjugate nets and their Jonas transformations. The
situation of Definition \ref{dfn:dcn jp} is governed by the
equation
\begin{equation}\label{eq:djp property}
\delta_i f^+=a_i\delta_i f+b_i(f^+-f),
\end{equation}
where coefficients $a_i$, $b_i$ are nothing but $a_i=1+c_{Mi}\,,$
$b_i=c_{iM}$. These coefficients are naturally attached to
elementary squares of $\bbZ^M$ parallel to the coordinate plane
$\PP_{iM}$. It is also convenient to think of them as attached to
edges of $\bbZ^m$ parallel to $\PP_i$ (to which the corresponding
``vertical'' squares are adjacent). Equations of the system
(\ref{eq:dcn c}) with one of the indices equal to $M$ give:
\begin{eqnarray}
\delta_ia_j & = & (\tau_jb_i)(a_j-1)+(\tau_ja_i-\tau_ia_j)c_{ij},
  \label{eq: djp a}\\
\delta_ib_j & = & c_{ij}^+b_j+c_{ji}^+b_i-(\tau_ib_j)b_i,
  \label{eq: djp b}\\
a_jc_{ij}^+ & = & (\tau_ja_i)c_{ij}+(\tau_jb_i)(a_j-1).
  \label{eq: djp c}
\end{eqnarray}

Following data are needed to specify a Jonas transform $f^+$ of a
given $m$-dimensional Q-net $f$:
\newcounter{dJ}
\begin{list}{(J$_\arabic{dJ}^\Delta$)}
{\usecounter{dJ}\setlength{\leftmargin}{1cm}}
\item value of $f^+(0)$;
\item values of $a_i$, $b_i$ on all edges of the respective coordinate
 axis $\PP_i$, for $1\le i\le m$.
\end{list}

{\bf Alternative description of discrete Jonas transformations.} We
give here a discrete version of the Eisehart's formulation of the
Jonas transformation. One derives from eqs.
(\ref{eq: djp a})--(\ref{eq: djp c}) the following formulas:
\begin{equation}\label{eq: djp B/A}
\bigg(1+\tau_i\Big(\frac{b_j}{a_j}\Big)\bigg)\bigg(1+\frac{b_i}{a_i}\bigg)=
1+(1+c_{ij})\frac{b_j}{a_j}+(1+c_{ji})\frac{b_i}{a_i},
\end{equation}
\begin{equation}\label{eq: djp B}
\Big(1+\tau_ib_j\Big)\Big(1+b_i\Big)=1+(1+c_{ij}^+)b_j+(1+c_{ji}^+)b_i.
\end{equation}
The symmetry of their right-hand sides implies that there exist
functions $\phi,\phi^+:\bbZ^m\to\bbR$ (associated to points of $\bbZ^m$)
such that
\begin{equation}\label{eq: djp phis}
\frac{\tau_i\phi}{\phi}=1+\frac{b_i}{a_i}\,,\qquad
\frac{\tau_i\phi^+}{\phi^+}=1+b_i\,,\qquad 1\le i\le m.
\end{equation}
These functions are defined uniquely up to constant factors, which
can be fixed by requiring $\phi(0)=\phi^+(0)=1$. Moreover, eqs.
(\ref{eq: djp B/A}), (\ref{eq: djp B}) imply that the functions
$\phi,\phi^+$ satisfy the following equations:
\begin{eqnarray}
\delta_i\delta_j\phi & = & c_{ij}\delta_j\phi+c_{ji}\delta_i\phi,
\label{eq: djp phi eq}  \\
\delta_i\delta_j\phi^+ & = & c_{ij}^+\delta_j\phi^+
+c_{ji}^+\delta_i\phi^+,   \label{eq: djp phi+ eq}
\end{eqnarray}
for all $1\le i\neq j\le m$. Thus, like in the smooth case, a
discrete Jonas transformation yields additional scalar solutions
$\phi$ and $\phi^+$ of the equations describing the nets $f$ and
$f^+$, respectively. The solution $\phi$ is directly specified by
the initial data (J$_2^\Delta$). Introduce the functions
$g:\bbZ^m\to\bbR^N$ and $\psi:\bbZ^m\to\bbR$ by the formulas
(\ref{eq:jp aux0}), so that also the classical representation
(\ref{eq:jp Eis}) remains valid. A direct computation based on
eqs. (\ref{eq:djp property}), (\ref{eq: djp a})--(\ref{eq: djp
c}), and (\ref{eq: djp phis}) shows that the following equations
hold:
\begin{eqnarray}
\delta_i g & = & \alpha_i\delta_i f,             \label{eq: djp aux1}\\
\delta_i\psi & = & \alpha_i\delta_i\phi,         \label{eq: djp aux2}
\end{eqnarray}
where
\begin{equation}\label{eq: djp aux0}
\alpha_i=\frac{a_i-1}{\tau_i\phi^+}\,.
\end{equation}
Next, one checks that the quantities $\alpha_i$ satisfy the equations
\begin{equation}\label{eq: djp aux3}
\delta_i\alpha_j=c_{ij}(\tau_j\alpha_i-\tau_i\alpha_j).
\end{equation}
The same argument as in the smooth case shows that the data (J$_2^\Delta$)
yield the values of $\phi^+$, and thus the values of $\alpha_i$, on the
coordinate axes $\PP_i$. This uniquely specifies the solutions $\alpha_i$
of the compatible linear system (\ref{eq: djp aux3}), which, in turn,
allows for a unique determination of the solutions $g$, $\psi$ of eqs.
(\ref{eq: djp aux1}), (\ref{eq: djp aux2}) with the initial data
$g(0)=f^+(0)-f(0)$ and $\psi(0)=1$.
\smallskip

{\bf Continuous limit.} Observe that eqs. (\ref{eq:dcn property}),
(\ref{eq:dcn c}) are quite similar to eqs. (\ref{eq:cn property}),
(\ref{eq:cn c}) characterizing smooth conjugate nets. We will
demonstrate in Sect.~\ref{Sect: approx} that the status of this
similarity can be raised to that of a mathematical theorem about
approximation of smooth conjugate nets by discrete ones. More
precisely, we will show how to choose initial data for a discrete
system (with $M=m$ and a small mesh size $\eps$) so that it
approximates a given $m$-dimensional smooth conjugate net as
$\eps\to 0$.

Analogously, eqs. (\ref{eq:djp property})--(\ref{eq: djp c}) are
similar to eqs. (\ref{eq:jp property})--(\ref{eq:jp C+}).
Accordingly, initial data of a discrete system with $M=m+1$ can be
choosen so that, keeping one direction discrete, one arrives in
the limit at a given smooth conjugate net and its Jonas transform.

For $M=m+2$ and $M=m+3$, keeping the last two, resp. three
directions discrete, one proves the permutability properties of
Jonas transformations formulated in Theorem \ref{thm:conj permut}.
Thus, permutability of Jonas transformations, which is a
non-trivial theorem of differential geometry, becomes an obvious
consequence of the multidimensional consistency of discrete
conjugate nets, combined with the convergence result mentioned
above.

\section{Discrete orthogonal nets}
 \label{Subsect: discr ortho}

Two-dimensional circular nets ($M=2$) were introduced in \cite{MPS, Nu}
as discrete analogs of the curvature lines parametrized surfaces. A
discretization of triply orthogonal coordinate systems ($M=3$, $N=3$) was
first proposed in \cite{B}. The next crucial step was done in \cite{CDS},
where discrete orthogonal nets were generalized to arbitrary dimensions
by considering them as a reduction of Q-nets. For further developments, see
\cite{DMS, KSch, AKV, DS2, BH}.

\begin{dfn} \label{dfn:dos}
  A map $f:\bbZ^M\to\bbR^N$ is called an $M$-dimensional {\em discrete
  O-net (discrete orthogonal net, or circular net)}
  in $\bbR^N$, if any of its elementary quadrilaterals is circular,
  i.e., if at any $u\in\bbZ^M$ and for all pairs $1\le i\neq j\le M$
  the four points $f$, $\tau_i f$, $\tau_j f$ and $\tau_i\tau_j f$ are
  concircular.
\end{dfn}
It is important to observe the M\"obius invariance of Definition
\ref{dfn:dos}. To understand restrictions imposed on Q-nets by the
circularity condition, we consider again various values of $M$.
\smallskip

{\itbf M=2:} {\bf discrete surfaces parametrized along curvature
lines.} Suppose two coordinate lines $f\rest_{\PP_1}$ and
$f\rest_{\PP_2}$ on a discrete O-surface $f$ are given. An
elementary inductive step for extending the O-surface to the
quadrant $\bbZ_+^2$, consists of choosing $f_{12}$ on the circle
through $f$, $f_1$ and $f_2$. In doing so, one has the freedom of
choosing one real parameter at each such step, for instance, the
cross-ratio of four points $q_{12}=q(f,f_1,f_{12},f_2)$, which is
naturally attached to the elementary square
$(u,u+e_1,u+e_1+e_2,u+e_2)$. Thus, to define a discrete O-surface
$f$, one needs to prescribe  the coordinate lines $f\rest_{\PP_1}$
and $f\rest_{\PP_2}$ and one real-valued function $q_{12}$ defined
on elementary squares.

\begin{figure}[hbtp]
\begin{center}
\setlength{\unitlength}{0.05em}
\begin{picture}(100,220)(-50,-110)
\put(0,0){\circle{200}}
\path(-80,-60)(80,-60)(60,80)(-80,60)(-80,-60)
\path(10,-55)(20,-60)(10,-65)
\path(-85,30)(-80,40)(-75,30)
\path(-7,76)(4,72)(-5,66)
\path(63,27)(66,38)(73,29)
\put(0,-46){$v_{i}$}
\put(30,0){$\tau_i v_{j}$}
\put(-70,0){$v_{j}$}
\put(-10,54){$\tau_j v_{i}$}
\put(-100,-70){$f$}
\put(85,-70){$\tau_i f$}
\put(60,90){$\tau_i\tau_j f$}
\put(-100,70){$\tau_j f$}
\end{picture}
\caption{An elementary quadrilateral of a discrete orthogonal net}
\label{fig:elem circle}
\end{center}
\end{figure}
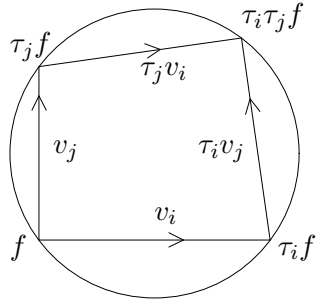

However, it turns out to be technically more convenient to use other
functions on elementary squares of $\bbZ^2$ characterizing the form of
circular quadrilaterals, namely, {\em discrete rotation coefficients}.
Introduce discrete metric coefficients $h_i=|\delta_i f|$
and unit vectors $v_i=h_i^{-1}\delta_i f$, so that
\begin{equation}\label{dos: evolve x}
  \tau_i f=f+h_iv_i.
\end{equation}
Then the rotation coefficients $\beta_{12}$ and $\beta_{21}$ are
defined by the formula
\begin{equation} \label{dos: evolve v}
 \tau_i v_j=\nu_{ji}^{-1}(v_j+\beta_{ji}v_i),
\end{equation}
which holds due to planarity of the elementary quadrilateral, with
$\nu_{ji}$ being appropriate normalization coefficients. From
(\ref{dos: evolve x}), (\ref{dos: evolve v}) there easily follows
a formula for metric coefficients $h_j$:
\begin{equation} \label{dos: evolve h}
  \tau_i h_j=\nu_{ji}^{-1}(h_j+h_i\beta_{ij}).
\end{equation}It is natural to assume that the variables $v_i, h_i$
are attached to the edges of $\bbZ^2$ parallel to the coordinate
axis $\PP_i$, while the rotation coefficients
$\beta_{12},\beta_{21}$ are attached to elementary squares of
$\bbZ^2$. Elementary considerations of the geometry of a circular
quadrilateral on Fig.~\ref{fig:elem circle} show that the discrete
rotation coefficients are related by
\begin{equation} \label{eq: dos beta}
  \beta_{ij}+\beta_{ji}+2\langle v_i,v_j \rangle=0,
\end{equation}
while the normalization coefficients satisfy
\begin{equation} \label{eq: dos n}
  \nu_{ij}^2=\nu_{ji}^2=1-\beta_{ij}\beta_{ji}.
\end{equation}
For an embedded quadrilateral there holds $\nu_{ij}=\nu_{ji}>0$,
while for a non-embedded one there holds $\nu_{ij}=-\nu_{ji}$.
Thus, formula (\ref{eq: dos beta}) expresses the {\em circularity
constraint}, and reflects the fact that the number of independent
plaquette functions necessary to define a discrete O-surface is
equal to 1. For instance, it is enough to prescribe $\beta_{12}$,
or $\beta_{21}$, or $\beta_{12}-\beta_{21}$, or any other function
of $\beta_{12}$, $\beta_{21}$, independent on
$\beta_{12}+\beta_{21}$.

Equations (\ref{dos: evolve v}), (\ref{dos: evolve h}) give the
evolution of the edge variables $v_j$, $h_j$, known the plaquette
variables $\beta_{12}$, $\beta_{21}$. Together with eq. (\ref{dos:
evolve x}), this allows us to reconstruct a discrete O-surface
$f$, provided the following initial data are given: two coordinate
lines $f\rest_{\PP_1}$, $f\rest_{\PP_2}$ (or, equivalently, the
point $f(0)$ and the functions $v_i$, $h_i$ on the edges of the
coordinate axes $\PP_i$, $i=1,2$), and the function $\beta_{21}$
(say) defined on all elementary squares of $\bbZ^2$.

As in the case of Q-surfaces, the combinatorics of discrete
O-surfaces may be more complicated than that of $\bbZ^2$, because
Definition~\ref{dfn:dos} can be literally extended to an arbitrary
quad-graph. This possibility is very important from the geometric
point of view, since vertices of valence different from 4 serve as
a model for umbilic points of a smooth surface parametrized  along
curvature lines.
\smallskip

\begin{figure}[t]
\begin{center}
\epsfig{file=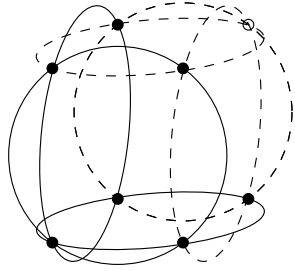,width=55mm} \caption{An elementary
hexahedron of a discrete orthogonal net} \label{fig:Miquel}
\end{center}
\end{figure}

{\itbf M=3:} {\bf basic 3D system.} Suppose that three coordinate
surfaces of a three-dimensional discrete O-net $f$ are given, that
is, $f\rest_{\PP_{12}}$, $f\rest_{\PP_{23}}$ and
$f\rest_{\PP_{13}}$. An inductive extension step consists of
determining $f_{123}$, provided $f$, $f_i$ and $f_{ij}$ are known
for all $1\leq i<j\leq 3$ and satisfy the circular condition on
the three corresponding squares. Such a step is possible due to
the {\em Miquel theorem} from the elementary geometry:
\begin{thm}
Given three circles $C_{ij}$ and seven points $f, f_i, f_{ij}$ in
$\bbR^3$ $(1\le i<j\le 3)$, such that each quadruple
$(f,f_i,f_j,f_{ij})$ lies on $C_{ij}$, define three new circles
$\tau_iC_{jk}$ as those passing through the triples
$(f_i,f_{ij},f_{ik})$, respectively. Then these new circles
intersect at one point:
\[
f_{123}=\tau_1C_{23}\cap\tau_2C_{31}\cap\tau_3C_{12}\,.
\]
\end{thm}
The claim of the Miquel theorem is equivalent to the following
one: if seven points $f, f_i, f_{ij}$ of an elementary hexahedron
of a Q-net $f:\bbZ^3\to\bbR^3$ lie on a two-sphere $S^2$, then so
does the eighth point $f_{123}$. To see the equivalence, note that
the four points $f$, $f_i$ determine the sphere $S^2$ uniquely.
Then $f_{ij}\in S^2$ is equivalent to $f_{ij}\in
C_{ij}=\Pi_{ij}\cap S^2$. Now, the intersection point of the three
circles $\tau_iC_{jk}=\tau_i\Pi_{jk}\cap S^2$ belongs to $S^2$ and
coincides with the unique intersection point of the three planes
\[
f_{123}=\tau_1\Pi_{23}\cap\tau_2\Pi_{31}\cap\tau_3\Pi_{12}.
\]
This construction is illustrated on Fig.~\ref{fig:Miquel}.
Clearly, this is a particular issue of Theorem \ref{Th:quadratic}.
Thus, finding the eighth point of a three-dimensional discrete
O-net out of the seven known ones is a 3D system in the sense of
Fig.~\ref{Fig:cube eq}.

The Miquel theorem  guarantees that circularity constraint
propagates in the construction of a Q-net from its three
coordinate surfaces. This was first observed by Doliwa and Santini
\cite{CDS}:
\begin{thm}
If three coordinate surfaces $f\rest_{\PP_{12}}$,
$f\rest_{\PP_{23}}$, $f\rest_{\PP_{13}}$ of a Q-net
$f:\bbZ^3\to\bbR^3$ are discrete O-surfaces, then the whole of $f$
is a discrete O-net.
\end{thm}
This is a discrete analog of the statement which holds also in the smooth
context, and is a sort of inversion of the classical Dupin theorem.

On the level of formulas, for discrete O-nets the system encoded
on Fig.~\ref{Fig:cube eq} consists of eqs. (\ref{dos: evolve x}),
(\ref{dos: evolve v}), (\ref{dos: evolve h}) for all $1\le i\neq
j\le 3$, with the normalization coefficients given in (\ref{eq:
dos n}), which have to be augmented by evolution equations for the
rotation coefficients. These latter ones represent the solvability
of the system of linear equations (\ref{dos: evolve v}) and are
derived from the requirement
$\tau_i(\tau_jv_k)=\tau_j(\tau_iv_k)$.
%
Like in the case of Q-nets, these equations are decoupled from the
rest ones, and generate a well defined map
$\{\beta_{kj}\}\mapsto\{\tau_i\beta_{kj}\}$:
\begin{equation}\label{dos: evolve beta}
\tau_i\beta_{kj}=(\nu_{ki}\nu_{ij})^{-1}(\beta_{kj}+\beta_{ki}\beta_{ij}).
\end{equation}
However, this evolution of the discrete rotation coefficients
$\beta_{kj}$ can be coupled to the evolution (\ref{dos: evolve v})
of the unit vectors $v_i$ by means of the circularity constraint
(\ref{eq: dos beta}). It follows from the Miquel theorem, and can
be easily checked analytically, that this constraint  propagates
in the coupled evolution.
\smallskip

{\itbf M \boldmath$\ge$ 4:} {\bf consistency.} The 4D consistency
of discrete O-nets is a consequence of the analogous property of
Q-nets, since the O-constraint is compatible with the Q-property.

On the level of formulas we have for $M\ge 4$ the system
(\ref{dos: evolve x}), (\ref{dos: evolve v}), (\ref{dos: evolve
h}), (\ref{dos: evolve beta}), augmented by the circularity
constraint (\ref{eq: dos beta}), where now all indices $i,j,k$
vary between 1 and $M$. For any triple of pairwise different
indices $(i,j,k)$, equations involving these indices solely, form
a closed subset. The $M$-dimensional consistency of this system
means that all three-dimensional building blocks can be imposed
without contradictions. Initial data which allow for a unique
solution of this system consist of:
\newcounter{dO}
\begin{list}{(O$_\arabic{dO}^\Delta$)}
{\usecounter{dO}\setlength{\leftmargin}{1cm}}
\item values of $f$ on the coordinate axes $\PP_i$
(or, what is equivalent, value $f(0)$ and values of $v_i$ and $h_i$ on all
edges of the coordinate axes $\PP_i$) for $1\le i\le M$;
\item values of $M(M-1)/2$ functions $\beta_{ji}$ (say) on all elementary
squares of the coordinate planes $\PP_{ij}$ for $1\le i<j\le M$.
\end{list}

At this point, it will be useful to observe the similarity of eqs.
(\ref{dos: evolve x}), (\ref{dos: evolve v}), (\ref{dos: evolve
h}), (\ref{dos: evolve beta}) and the constraint (\ref{eq: dos
beta}) with the system (\ref{eq:DX})--(\ref{eq:DB}) and the
constraint (\ref{eq:DL}) governing smooth orthogonal nets. Like in
the case of conjugate nets, we will demonstrate that this analogy
can be given a qualitative content, so that smooth O-nets can be
approximated by discrete ones. However, there is a substantial
obstruction in performing this. We think of smooth rotation
coefficients as being approximated by discrete ones. But since the
discrete rotation coefficients $\beta_{kj}$ only have $i\neq k,j$
as evolution directions (that is, they are {\em plaquette}
variables attached to elementary squares parallel to $\PP_{jk}$),
there is seemingly no chance to get an approximation of such
smooth quantities as $\partial_i\beta_{ij}$ involved in the smooth
orthogonality constraint (\ref{eq:DL}). In order to be able to
achieve such an approximation, we need some discrete analogs of
the smooth rotation coefficients which would live on {\em edges}.
Such analogs will be introduced with the help of the notion of a
frame of a discrete O-net. Technical means for defining and
constructing frames are given by the apparatus of M\"obius
differential geometry.
\medskip

{\bf Transformations of discrete O-nets.} A natural generalization
of Definition \ref{dfn:rp} would be the following one.
\begin{dfn}\label{dfn:dos rp}
  A pair of $m$-dimensional discrete O-nets $f,f^+:\bbZ^m\to\bbR^N$
  is called a {\em Ribaucour pair}, if four points $f$, $\tau_i f$,
  $f^+$ and $\tau_i f^+$ are concircular at any $u\in\bbZ^m$
  and for any $1\leq i\leq m$. The net $f^+$ is called
  a {\em Ribaucour transform} of the net $f$.
\end{dfn}
But this simply means that, if we set $F(u,0)=f(u)$,
$F(u,1)=f^+(u)$, then $F:\bbZ^m\times\{0,1\}\to\bbR^N$ is an
$M$-dimensional discrete O-net, where $M=m+1$. So, in the discrete
case there is no difference between orthogonal nets and their
Ribaucour transformations. To specify a Ribaucour transform $f^+$
of a given $m$-dimensional discrete O-net $f$, one clearly needs
the following data:
\newcounter{dR}
\begin{list}{(R$_\arabic{dR}^\Delta$)}
{\usecounter{dR}\setlength{\leftmargin}{1cm}}
\item value of $f^+(0)$;
\item values of $\beta_{Mi}$ (say) on ``vertical'' elementary squares
attached to all edges of the coordinate axes $\PP_i$ for $1\le i\le m$.
\end{list}

{\bf M\"obius-geometric description of discrete O-nets.} Putting
discrete O-nets in the M\"obius-geometric model $\bbQ_{\,0}^N$ of
the Euclidean space $\bbR^N$, we observe first of all that Theorem
\ref{prop: os in K} admits an almost literal generalization to the
discrete case.

\begin{thm}\label{Th: dos in K}
A Q-net $f:\bbZ^M\to\bbR^N$ is a discrete O-net, if and only if
$|f|^2$ satisfies the same equation (\ref{eq:dcn property}) as $f$
does, in other words, if the corresponding
$\hat{f}:\bbZ^M\to\bbQ_{\,0}^N$ is a Q-net in $\bbR^{N+1,1}$.
\end{thm}
(The first claim here is due to \cite{KSch}.) Thus, in the
M\"obius-geometric picture discrete O-nets are a quadratic
reduction (to $\bbQ_{\,0}^N$) of Q-nets.

Let $\hat{v}_i$ be the (Lorentz) unit vectors  parallel to
$\delta_i\hat{f}$.
For the discrete metric coefficients $h_i=|\delta_i f|$ there
holds also $h_i=|\delta_i\hat{f}|$. One readily verifies that from
$\langle\hat{f}, \hat{f}\rangle=0$ there follows
$\langle\tau_i\hat{f}+\hat{f},\hat{v}_i\rangle=0$ and
$h_i=-2\langle\hat{f},\hat{v}_i\rangle$.  Therefore,
$\hat{v}_i$ are interpreted in the M\"obius-geometric picture as
reflections taking $\hat{f}$ to $\tau_i\hat{f}$:
\begin{equation}  \label{eq:xinthemirror}
  \tau_i\hat{f}=\hat{f}+h_i\hat{v}_i=\hat{v}_i\hat{f}\hat{v}_i.
\end{equation}
Due to Theorem \ref{Th: dos in K}, vectors $\hat{v}_i$ satisfy the
same linear relations (\ref{dos: evolve v}) as $v_i$, with the
same discrete rotations coefficients $\beta_{ji}$.
\begin{thm}[Spinor frame of a discrete O-net]   \label{thm:frame for dos}
For a discrete O-net $f:\bbZ^M\to\bbR^N$ (and the corresponding
Q-net $\hat{f}:\bbZ^M\to\bbQ_{\,0}^N$), there exists a function
$\psi:\bbZ^M\to\cH_\infty$ (called {\em a frame of} $\hat{f}$)
such that
\begin{equation}\label{eq:pinpointthem}
  \hat{f}=\psi^{-1}\ee_0\psi,
\end{equation}
satisfying the system of difference equations:
\begin{equation}\label{eq: frame dos}
  \tau_i\psi = -\ee_i\psi\hat{v}_i\;,\quad 1\le i\leq M.
\end{equation}
\end{thm}
\noindent
 {\bf Proof.} This theorem is due to \cite{BH}. Circularity of an
elementary quadrilateral of $f$ implies that the angle between
$\delta_i f$ and $\tau_i\delta_j f$ and the angle between
$\delta_j f$ and $\tau_j\delta_i f$ sum up to $\pi$:
\begin{equation}\label{circularity}
  \langle\tau_j v_i,v_j\rangle+
  \langle v_i,\tau_i v_j\rangle=0.
\end{equation}
The same relation holds for the vectors $\hat{v}_i$. Moreover,
since for these vectors there holds also eq. (\ref{dos: evolve v})
with $\nu_{ij}=\nu_{ji}$, there follows that eq.
(\ref{circularity}) for $\hat{v}_i$ is equivalent to
\begin{equation}\label{circ}
  \hat{v}_j(\tau_j\hat{v}_i)+\hat{v}_i(\tau_i\hat{v}_j)=0.
\end{equation}
But this is exactly the compatibility condition of the system
(\ref{eq: frame dos}) of linear difference equations for $\psi$.
Indeed, the solvability condition $\tau_j\tau_i\psi=\tau_i\tau_j\psi$ is
written as
\[
\ee_i\ee_j\psi\hat{v}_j(\tau_j\hat{v}_i)=
\ee_j\ee_i\psi\hat{v}_i(\tau_i\hat{v}_j),
\]
which is equivalent to (\ref{circ}). Thus, a solution to (\ref{eq:
frame dos}) exists and is uniquely defined by the choice of
$\psi(0)\in\cH_\infty$. Choose $\psi(0)$ so that eq.
(\ref{eq:pinpointthem}) holds at $u=0$. Then eqs. (\ref{eq:
frame dos}), (\ref{eq:xinthemirror}) imply that eq.
(\ref{eq:pinpointthem}) holds everywhere on $\bbZ^M$. $\Box$
\medskip

Now introduce vectors $V_i=\psi\hat{v}_i\psi^{-1}$, so that the
frame equations (\ref{eq: frame dos}) take the form
$(\tau_i\psi)\psi^{-1}=-\ee_iV_i$. Expanding these vectors with
respect to the basis vectors $\ee_k$, we have a formula analogous
to (\ref{eq: os S}):
\begin{equation} \label{eq: dos V}
  V_i=\psi\hat{v}_i\psi^{-1}=
  \sigma_i\ee_i-\frac{1}{2}\sum_{k\neq i}\rho_{ki}\ee_k+h_i\ee_\infty.
\end{equation}
The fact that the $\ee_\infty$-component here is equal to $h_i$,
is easily demonstrated. Indeed, from eq. (\ref{eq: frame dos})
there follows that
$\tau_i\hat{f}-\hat{f}=h_i\hat{v}_i=h_i(\tau_i\psi)^{-1}\ee_i\psi$.
Now eq. (\ref{eq:pinpointthem}) allows us to rewrite this
equivalently as $[\ee_0,(\tau_i\psi)\psi^{-1}]=h_i\ee_i$, which
proves the claim above. Observe also the normalization condition
\begin{equation} \label{eq: N for dos}
  \sigma_i^2=1-\frac{1}{4}\sum_{k\neq i}\rho_{ki}^2\,.
\end{equation}

Coefficients $\rho_{ki}$ are {\em edge variables} analogous to
smooth rotation coefficients. Indeed, vectors $V_i$ are defined on
edges of $\bbZ^M$ parallel to the coordinate axis $\PP_i$, but
they do not immediately reflect the local geometry near these
edges. Rather, they are obtained by integration of the frame
equations (\ref{eq: frame dos}), and thus are of a {\em non-local}
nature. Thus, in the discrete case we have two different analogs
of the rotation coefficients: local plaquette variables
$\beta_{ij}$ for $1\le i\neq j\le M$, defined on elementary
squares of $\bbZ^M$ parallel to $\PP_{ij}$, and non-local edge
variables $\rho_{ki}$ for $1\le i\le M$, $1\le k\le N$, $k\neq i$,
defined on edges of $\bbZ^M$ parallel to $\PP_i$.

Evolution equations for $V_i$ are obtained
from (\ref{dos: evolve v}) and the frame equations (\ref{eq: frame dos}):
\begin{equation*} \label{eq:master}
  \tau_iV_j=\nu_{ji}^{-1}\ee_i(V_j+\beta_{ij}V_i)\ee_i.
\end{equation*}
In the derivation one uses the identity
$\hat{v}_i(\hat{v}_j+\beta_{ji}\hat{v}_i)\hat{v}_i=
\hat{v}_j+\beta_{ij}\hat{v}_i\,$, which follows easily from
(\ref{eq: dos beta}). The resulting evolution equations for the
edge variables $\rho_{kj}$ read:
\begin{eqnarray}
    \tau_i\rho_{kj} & = & \nu_{ji}^{-1}
    (\rho_{kj}+\rho_{ki}\beta_{ij}),     \label{eq:dB}\\
    \tau_i\rho_{ij} & = & \nu_{ji}^{-1}
    (-\rho_{ij}+2\sigma_i\beta_{ij}).          \label{eq:dL}
\end{eqnarray}
Here $1\le i\neq j\leq M$, $\,1\le k\leq N$, and $k\neq i,j$.
The circularity constraint (\ref{eq: dos beta}) can be now written as
\begin{equation}\label{eq:dQ}
  \beta_{ij}+\beta_{ji}=\sigma_i\rho_{ij}+\sigma_j\rho_{ji}
   -\frac{1}{2}\sum_{k\neq i,j}\rho_{ki}\rho_{kj},
\end{equation}
and gives a relation between local plaquette variables
$\beta_{ij}$ and non-local edge variables $\rho_{kj}$. The system
consisting of (\ref{eq:dB}), (\ref{eq:dL}) and (\ref{eq:dQ}) can
be regarded as the {\em discrete Lam\'e system}.
\medskip

{\bf Continuous limit.} In Sect. \ref{subsect: os approx} we will
show how to choose these data on the lattice with mesh size $\eps$
in $m$ directions and mesh size 1 in the $M$-th direction, in
order to achieve a simultaneous approximation of a smooth
orthogonal net and its Ribaucour transformation. Performing a
smooth limit so that two or three of the directions remain
discrete, one can derive the permutability properties of Ribaucour
transformations formulated in Theorem \ref{thm:ortho permut}. This
way, permutability of Ribaucour transformations becomes a simple
consequence of properties of discrete orthogonal nets, combined
with the convergence result mentioned above.

\section{Discrete Moutard nets}
 \label{Subsect: discr Mou}

Discrete Moutard nets were introduced in \cite{NSch}.
\begin{dfn}\label{def: dmn}
A map $f:\bbZ^2\to\bbR^N$ is called a two-dimensional {\em
discrete M-net}, if it satisfies the {\em discrete Moutard
equation}
\begin{equation}\label{eq:dMou 2d}
\tau_1\tau_2 f+f=a_{12} (\tau_1 f+\tau_2 f)
\end{equation}
with some $a_{12}:\bbZ^2\to\bbR$.
\end{dfn}
Initial data that can be used to determine a discrete M-net are:
\newcounter{dM}
\begin{list}{(M$_\arabic{dM}^\Delta$)}
{\usecounter{dM}\setlength{\leftmargin}{1cm}}
\item values of $f$ on the coordinate axes $\PP_1$, $\PP_2$;
\item function $a_{12}$ defined on elementary squares of $\bbZ^2$.
\end{list}
Like in the smooth case, discrete M-nets are not Q-nets, and there
are no $M$-dimensional discrete M-nets with $M\ge 3$. There is a
construction, analogous to the smooth case, relating discrete
M-nets to a special class of Q-nets. Let $\nu:\bbZ^2\to\bbR$ be
any solution of the same discrete Moutard equation (\ref{eq:dMou
2d}) (for instance, any component of the vector $f$), then
$y=\nu^{-1}f:\bbZ^2\to\bbR^N$ is a Q-net:
\[
y_{12}-y=\frac{\nu_1(\nu_{12}+\nu)}{\nu_{12}(\nu_1+\nu_2)}\,(y_1-y)+
         \frac{\nu_2(\nu_{12}+\nu)}{\nu_{12}(\nu_1+\nu_2)}\,(y_2-y).
\]
However, there exists a relation to discrete conjugate nets of a
quite different flavor, which is of a purely discrete nature and
has no smooth analogs. To describe it, perform in (\ref{eq:dMou})
the change of variables
\begin{equation}\label{eq:fbar}
f(u)\mapsto(-1)^{u_2}\,f(u),\qquad
u=(u_1,u_2)\in\bbZ^2.
\end{equation}
Then the new function $f$ satisfies the following modified form of the
discrete Moutard equation:
\[
\tau_1\tau_2 f-f=a_{12} (\tau_2 f-\tau_1 f).
\]
This equation admits a multidimensional generalization:
\begin{dfn} \label{dfn:tn}
  A map $f:\bbZ^M\to\bbR^N$ is called an $M$-dimensional
  {\em T-net (trapezoidal net)},
  if for any $u\in\bbZ^M$ and for any pair of indices $i\neq j$ there
  holds the {\em discrete Moutard equation}
\begin{equation}\label{eq:dMou}
 \tau_i\tau_j f-f=a_{ij} (\tau_j f-\tau_i f),
\end{equation}
 with some $a_{ij}:\bbZ^M\to\bbR$, in other words, if all the elementary
 quadrilaterals $(f,\tau_i f,\tau_i\tau_j f,\tau_j f)$ are planar and have
 parallel diagonals.
\end{dfn}
Of course, coefficients $a_{ij}$ have to be skew-symmetric,
$a_{ij}=-a_{ji}$. As usual, we will consider these functions as
attached to elementary squares of $\bbZ^M$.

T-nets, unlike discrete M-nets, form a subclass of Q-nets. The
condition of parallel diagonals is expressed as
$c_{ij}+c_{ji}+2=0$ for the coefficients $c_{ij}$ of a Q-net, the
skew-symmetric coefficients of the T-net being $a_{ij}=c_{ij}+1$.
\smallskip

{\itbf M=2:} {\bf T-surfaces.} To define a two-dimensional T-net
$f:\bbZ^2\to\bbR^N$, one can prescribe two coordinate curves,
$f\rest_{\PP_1}$ and $f\rest_{\PP_2}\,$, and a real-valued
function $a_{12}$ on elementary squares of $\bbZ^2$.

{\itbf M=3:} {\bf basic 3D system.} We show that three-dimensional
T-nets are described by a well-defined three-dimensional system.
An inductive construction step of the net $f$ is as follows.
Suppose that $f$, $f_i$ and $f_{ij}$ are given for all $1\leq
i\neq j\leq 3$, satisfying eq. (\ref{eq:dMou}). Three equations
(\ref{eq:dMou}) for the facets of an elementary cube on
Fig.~\ref{Fig:cube eq} adjacent to $f_{123}$, lead to consistent
results for $f_{123}$ for arbitrary initial data, if and only if
the following conditions are satisfied:
\[
(\tau_ia_{jk})a_{ij}+(\tau_ja_{ki})(a_{jk}+a_{ij})=-1,
\]
where $(i,j,k)$ is an arbitrary permutation of $(1,2,3)$.
These conditions constitute a system of 6 (linear) equations for 3
unknown variables $\tau_i a_{jk}$ in terms of the known ones
$a_{jk}$. It turns out that this system is not overdetermined but
admits a unique solution:
\begin{equation}\label{eq:star-triang}
\frac{\tau_ia_{jk}}{a_{jk}}=
-\frac{1}{a_{ij}a_{jk}+a_{jk}a_{ki}+a_{ki}a_{ij}},
\end{equation}
With $\tau_ia_{jk}$ so defined, eqs. (\ref{eq:dMou}) are fulfilled
on all three quadrilaterals adjacent to $f_{123}$.

Eqs. (\ref{eq:star-triang}) represent a well-defined birational
map $\{a_{jk}\}\mapsto\{\tau_ia_{jk}\}$, which can be considered
as the fundamental 3D system related to T-nets. It is sometimes
called the ``star-triangle map''. Moreover, this system is a
reduction of the system describing Q-nets:
\begin{thm}
If three coordinate surfaces $f\rest_{\PP_{12}}$,
$f\rest_{\PP_{23}}$, $f\rest_{\PP_{13}}$ of a Q-net
$f:\bbZ^3\to\bbR^N$ are T-nets, then $f$ is a T-net.
\end{thm}
{\bf Proof.}  Let three quadrilaterals $(f,f_i,f_{ij},f_j)$ be
planar and have parallel diagonals. The planarity of the three
quadrilaterals $(f_i,f_{ij},f_{ijk},f_{ik})$ defines the point
$f_{123}$ as the intersection point of three planes planes
$\tau_i\Pi_{jk}$. Then these three quadrilaterals automatically
have parallel diagonals. Indeed, by the above argument, there
exists a point $f_{123}$ with this property, and it has to
coincide with the one defined by the planarity condition. $\Box$

{\itbf M \boldmath$\ge$ 4:} {\bf consistency.} The 4D consistency
of T-nets is a consequence of the analogous property of Q-nets,
since T-constraint is compatible with the Q-property. On the level
of formulas we have for T-nets with $M\ge 4$ the system
(\ref{eq:dMou}), (\ref{eq:star-triang}). All indices $i,j,k$ vary
now between 1 and $M$, and for any triple of pairwise different
indices $(i,j,k)$, equations involving these indices solely, form
a closed subset. Initial data which allow for a unique solution of
this system consist of values of $f$ on the coordinate axes
$\PP_i$ for $1\le i\le M$, and values of $M(M-1)/2$ functions
$a_{ij}$ on all elementary squares of the coordinate planes
$\PP_{ij}$ for $1\le i<j\le M$.
\medskip

{\bf Transformations of discrete M- and T-nets.} Because of the
multi-dimen\-sional consistency, transformations of T-nets do not
differ from the nets themselves, and are described by the discrete
Moutard equations. We would like, however, to consider
transformations of discrete M-nets, which have geometric
applications. For this, we consider a three-dimensional T-net
$F:\bbZ^2\times\{0,1\}$, perform the change of variables
\[
F(u)\mapsto (-1)^{u_2}F(u),\qquad u=(u_1,u_2,u_3),
\]
and then set $f=F(\cdot,0)\,$, $f^+=F(\cdot,1)$. Eqs.
(\ref{eq:dMou}) for $(i,j)=(1,2)$ turn into (\ref{eq:dMou 2d}), so
that $f$, $f^+$ are discrete M-nets, while for $(i,j)=(1,3)$ and
$(2,3)$ eqs. (\ref{eq:dMou}) turn into
\begin{equation}\label{eq:dmp}
\tau_1 f^+-f=b_1(f^+-\tau_1 f),\qquad \tau_2 f^++f=b_2(f^++\tau_2 f),
\end{equation}
where $b_1=a_{13}\,$, $b_2=a_{32}$. The quantities $b_i$, defined
on the ``vertical'' plaquettes of $\bbZ^2\times\{0,1\}$, parallel
to $\PP_{i3}$, can be also associated to the edges of $\bbZ^2$
parallel to the coordinate axes $\PP_i$. Eqs.
(\ref{eq:star-triang}) with one of the indices equal to $3$
express compatibility of eqs. (\ref{eq:dmp}) with eq.
(\ref{eq:dMou}):
\begin{equation}\label{eq:dMou 1}
\frac{\tau_2 b_1}{b_1}=\frac{\tau_1 b_2}{b_2}=
\frac{a_{12}^+}{a_{12}}=\frac{1}{(b_1+b_2)a_{12}-b_1b_2}.
\end{equation}
Eqs. (\ref{eq:dmp}), (\ref{eq:dMou 1}) define a {\em discrete
Moutard transformation} of $f$. To specify a Moutard transform
$f^+$ of a given discrete M-net $f$, one can prescribe the
following data:
\newcounter{dT}
\begin{list}{(MT$_\arabic{dT}^\Delta$)}
{\usecounter{dT}\setlength{\leftmargin}{1cm}}
\item value of $f^+(0)$;
\item values of $b_i$ on ``vertical'' elementary squares attached
to all edges of the coordinate axes $\PP_i$ for $i=1,2$.
\end{list}
The permutability properties of Moutard transformations are
governed by the discrete Moutard equations (\ref{eq:dMou}) with
$i,j\ge 3$.
\smallskip

Due to the first equation in (\ref{eq:dMou 1}), there exists a
function $\theta:\bbZ^2\to\bbR$ (associated to the points of
$\bbZ^2$) such that
\begin{equation}\label{eq:dMou b}
b_i=\frac{\theta}{\tau_i\theta},\qquad  i=1,2.
\end{equation}
The last equation in (\ref{eq:dMou 1}) implies that this function
is a scalar solution of the discrete Moutard equation
(\ref{eq:dMou 2d}). This solution is specified by its values on
the coordinate axes $\PP_i$, which are immediately obtained, via
eq. (\ref{eq:dMou b}), from the data (MT$_2^\Delta$). Recall that
$f^+$ is a solution of the discrete Moutard equation (\ref{eq:dMou
2d}) with the transformed potential $a_{12}^+$. The second
equation in (\ref{eq:dMou 1}) gives a representation of $a_{12}^+$
in terms of $\theta$:
\begin{equation}
a_{12}^+=a_{12}\,\frac{(\tau_1\theta)(\tau_2\theta)}
{\theta(\tau_1\tau_2\theta)}=
\frac{\tau_1\tau_2\theta^++\theta^+}
{\tau_1\theta^++\tau_2\theta^+}, \qquad \theta^+=\frac{1}{\theta}.
\end{equation}
This is a discrete analog of the classical formulation of the
Moutard transformation.

\section{Discrete asymptotic nets}
 \label{Subsect: discr asymp}

The following definition is due to Sauer \cite{Sauer} for $M=2$, when
it describes a discrete analog of surfaces parametrized along asymptotic
lines (A-surfaces). For $M\ge 3$, see \cite{D4, DNS, N}.
\begin{dfn} \label{dfn:dan}
  A map $f:\bbZ^M\to\bbR^3$ is called an $M$-dimensional {\em
  \,discrete A-net (discrete asymptotic net)} in $\bbR^3$,
  if for any $u\in\bbZ^M$ all the points $f(u\pm e_i)$ lie in some
  plane $\cP(u)$ through $f(u)$.
\end{dfn}
Note that this definition belongs to the projective geometry.
For $M=3$, the geometry of an elementary hexahedron
of a discrete A-net is exactly that of a {\em M\"obius pair of
tetrahedra}, i.e., a pair of tetrahedra which are inscribed in
each other. In other words, each vertex of each tetrahedron lies
in the plane of the corresponding facet of the other one (eight
conditions). In our case the two tetrahedra are those with the
vertices $(f,f_{12},f_{23},f_{31})$ and with the vertices
$(f_1,f_2,f_3,f_{123})$. Such pairs were introduced by M\"obius
\cite{M}, who demonstrated that eight conditions mentioned above
are not independent: any one of them follows from the remaining
seven. M\"obius pairs are remarkable and well-studied objects of
the projective geometry, cf. \cite{B}.

For a non-degenerate discrete A-net, all quadrilaterals $(f,\tau_i
f,\tau_i\tau_j f,\tau_j f)$ are non-planar. It would be in
principle possible to consider discrete A-nets in $\bbR^N$ with
$N>3$, however it would not lead to an essential generalization.
Indeed, for any fixed $u\in\bbZ^M$ and for any pair of indices
$i\neq j$ from $\{1,\ldots,M\}$, consider the three-dimensional
affine subspace of $\bbR^N$ through $f=f(u)$ which contains
$\tau_i f$, $\tau_j f$ and $\tau_i\tau_j f$. A simple induction
shows that the whole net $f$ lies in this subspace.

For a discrete A-net $f:\bbZ^M\to\bbR^3$ we have a field of
tangent planes $\cP:\bbZ^M\to{\rm Gr}_2(3)$, and therefore a
well-defined {\em normal direction} at every point of $\bbZ^M$.  A
remarkable way to fix a certain normal field is given by the {\em
discrete Lelieuvre representation} \cite{KP} which states:
\begin{thm}
For a non-degenerate discrete A-net $f$, there exists a normal
field $n:\bbZ^M\to\bbR^3$ such that
\begin{equation}\label{eq:dLel}
\delta_i f=\tau_i n\times n,\qquad i=1,\ldots,M,
\end{equation}
called a {\em Lelieuvre normal field}. It is uniquely defined by a
value at one point $u_0\in\bbZ^M$. All other Lelieuvre normal
fields are obtained by $n(u)\mapsto\alpha n(u)$ for
$|u|=u_1+\ldots+u_M$ even, and $n(u)\mapsto\alpha^{-1} n(u)$
for $|u|$ odd, with some $\alpha\in\bbR$ (black-white rescaling).
\end{thm}
It follows from eq. (\ref{eq:dLel}) immediately that
$(\tau_i\tau_j n-n)\times (\tau_i n-\tau_j n)=0$, that is, the
Lelieuvre normal field satisfies the discrete Moutard equations
\begin{equation}\label{eq:dMou as}
 \tau_i\tau_j n-n=a_{ij} (\tau_j n-\tau_i n)
\end{equation}
with some $a_{ij}:\bbZ^M\to\bbR$. Conversely, given a T-net
$n:\bbZ^M\to\bbR^3$, formula (\ref{eq:dLel}) produces a discrete
A-net $f:\bbZ^M\to\bbR^3$.
\begin{thm}
Discrete A-nets in $\bbR^3$ are in a one-to-one correspondence,
via the discrete Lelieuvre representation (\ref{eq:dLel}), with
T-nets in $\bbR^3$.
\end{thm}
In particular, the initial data which determine a discrete A-net
are analogous to the data (M$_{1,2}^\Delta$) for the Lelieuvre
normal field:
\newcounter{dA}
\begin{list}{(A$_\arabic{dA}^\Delta$)}
{\usecounter{dA}\setlength{\leftmargin}{1cm}}
\item values of $n$ on the coordinate axes $\PP_i$ for $1\le i\le
M$;
\item values of $M(M-1)/2$ functions $a_{ij}$ on all elementary
squares of the coordinate planes $\PP_{ij}$ for $1\le i<j\le M$.
\end{list}
\medskip

{\bf Transformations of discrete A-nets.} A natural generalization
of Definition \ref{dfn:rp} would be the following one.
\begin{dfn}\label{dfn:dwp}
A pair of discrete A-nets $f,f^+:\bbZ^m\to\bbR^3$ is called a {\em
Weingarten pair}, if, for any $u\in\bbZ^m$, the line segment
$[f(u),f^+(u)]$ lies in both tangent planes to $f$ and $f^+$ at
the points $f(u)$ and $f^+(u)$, respectively. The net $f^+$ is
called a {\em Weingarten transform} of the net $f$.
\end{dfn}
But, as usual, this definition means simply that the net
$F:\bbZ^m\times\{0,1\}\to\bbR^3$ with $F(u,0)=f(u)$ and
$F(u,1)=f^+(u)$ is an $M$-dimensional discrete A-net, where
$M=m+1$. So, once again, transformations of discrete A-nets do not
differ from the nets themselves. The Lelieuvre representation of
the net $F$ can be written now as a relation between the Lelieuvre
representations of the nets $f$, $f^+$:
\begin{equation}\label{eq:dWei}
f^+-f=n^+\times n.
\end{equation}

Definition \ref{dfn:dwp} makes sense for an arbitrary $M$;
however, we will be mainly interested in the case $M=2$, which is
an immediate discretization of the smooth A-surfaces. The change
of variables (\ref{eq:fbar}) for the Lelieuvre normal field,
\begin{equation}\label{eq:nbar}
n(u)\mapsto(-1)^{u_2}\,n(u),\qquad u=(u_1,u_2)\in\bbZ^2,
\end{equation}
leads to a replacement of the general Lelieuvre formulas (\ref{eq:dLel})
by
\begin{equation}\label{eq:dLel 2d}
\delta_1 f=\tau_1 n\times n=\delta_1 n\times n,\qquad
\delta_2 f=n\times\tau_2 n=n\times\delta_2 n.
\end{equation}
For a Moutard transformation of the Lelieuvre normal field there hold
formulas
(\ref{eq:dmp}):
\begin{equation}\label{eq:dwp}
\tau_1 n^+-n=b_1(n^+-\tau_1 n),\qquad \tau_2 n^++n=b_2(n^++\tau_2 n).
\end{equation}
Thus, a Weingarten transform $f^+$ of a given discrete A-net $f$
is determined by a Moutard transform $n^+$ of the Lelieuvre normal
field $n$, and in order to specify the latter, one can prescribe
the following data:
\newcounter{dW}
\begin{list}{(W$_\arabic{dW}^\Delta$)}
{\usecounter{dW}\setlength{\leftmargin}{1cm}}
\item value of $n^+(0)$;
\item values of $b_i$ on ``vertical'' elementary squares attached
to all edges of the coordinate axes $\PP_i$ for $i=1,2$.
\end{list}
Permutability of discrete Weingarten transformations is a direct
consequence of the 4D consistency of T-nets. Permutability of
smooth Weingarten transformations (Theorem \ref{thm:as permut})
will follow, if one combines the discrete permutability with the
convergence results of Sect. \ref{Subsect: an conv}.

\section{Discrete K-nets}
 \label{Subsect: discr K}

In discretizing the K-surfaces and their transformations, we take
as a starting point the characterization of Theorem \ref{Th: ks}.
\begin{dfn} \label{dfn:dks}
  A discrete A-net $f:\bbZ^M\to\bbR^3$ is called an M-dimensional
  {\em discrete K-net}, if for any elementary quadrilateral
  $(f,\tau_i f,\tau_i\tau_j f,\tau_j f)$
  there holds:
  \[
  |\tau_i\tau_j f-\tau_j f|=|\tau_i f-f|\quad and \quad
  |\tau_i\tau_j f-\tau_i f|=|\tau_j f-f|,
  \]
  in other words, if for any $i=1,\ldots,M$, the function
  $\alpha_i=|\delta_i f|$ (defined on the edges parallel to the
  coordinate axes $\PP_i$) depends on $u_i$ only.
\end{dfn}
This notion is due to Sauer \cite{Sauer} in the case $M=2$ and to
Wunderlich \cite{W} in the case $M=3$. A study of discrete
K-surfaces within the framework of the theory of integrable
systems was performed in \cite{BP1}.

A characterization of the Lelieuvre normal field of a discrete
K-net is analogous to the smooth case.
\begin{thm}\label{Th: dks}
 The Lelieuvre normal field $n:\bbZ^M\to\bbR^3$ of a discrete
 K-net $f:\bbZ^M\to\bbR^3$ takes values, possibly upon a black-white
 rescaling, in some sphere $S^2\subset\bbR^3$. In case when the radius
 of this sphere is equal to 1 (which we associate to the discrete
 Gaussian curvature $K=-1$), the Lelieuvre normal field coincides with
 the Gauss map. Conversely, any T-net $n$ in the unit sphere $\bbS^2$
 is the Gauss map and the Lelieuvre normal field of a discrete K-net
 $f$ with $K=-1$. The functions $\beta_i=|\delta_i n|$ depend on $u_i$
 only, and are related to the functions $\alpha_i=|\delta_i f|$ by
 $\alpha_i=\beta_i(1-\beta_i^2/4)^{1/2}$.
\end{thm}

Thus, discrete K-nets are in a one-to-one correspondence with
T-nets in $\bbS^2$. We proceed to the study of the latter object.
\smallskip

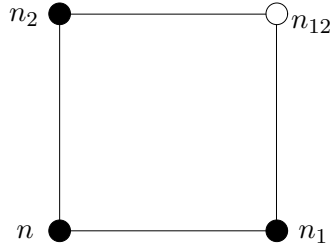
\begin{figure}[htbp]
\begin{center}
\setlength{\unitlength}{0.05em}
\begin{picture}(160,160)(0,0)
 \put(0,0){\circle*{15}}    \put(150,0){\circle*{15}}
 \put(0,150){\circle*{15}}  \put(150,150){\circle{15}}
 \path(0,0)(150,0)       \path(0,0)(0,150)
 \path(150,0)(150,142.5)   \path(0,150)(142.5,150)
 \put(-30,-5){$n$}
 \put(-35,145){$n_2$}  \put(165,-5){$n_1$}
 \put(160,140){$n_{12}$}
\end{picture}
\caption{2D system on an elementary
quadrilateral}\label{Fig:square eq}
\end{center}
\end{figure}

{\itbf M=2:} {\bf basic 2D system.} If two coordinate curves of the Gauss
map $n$ of a discrete K-surface $f$ are given, that is, $n\rest_{\PP_1}$
and $n\rest_{\PP_2}$, then one can extend $n$ to the whole of $\bbZ^2$.
An inductive step in extending $n$ to $\bbZ_+^2$ consists of computing
$\tau_1\tau_2 n=n+a_{12}(\tau_2 n-\tau_1 n)$, where the coefficient
$a_{12}$ (attached to every elementary square of the discrete surface) is
determined by the condition that $\tau_1\tau_2 n\in\bbS^2$, so that
\[
  a_{12}=\frac{\langle n,\tau_1 n-\tau_2 n\rangle}
         {1-\langle \tau_1 n,\tau_2 n\rangle}\,.
\]
This elementary construction step, i.e., finding the fourth vertex
of an elementary square out of the known three vertices, is
symbolically represented on Fig.~\ref{Fig:square eq}.

{\itbf M \boldmath$\ge$ 3:} {\bf consistency.} Turning to the case
$M\ge 3$, we see that one can prescribe all coordinate lines of a
T-net $n$ in $\bbS^2$, i.e., $n\rest_{\PP_i}$ for all $1\leq i\leq
M$. Indeed, these data are independent, and one can, by induction,
construct the whole net from them. The inductive step is
essentially two-dimensional and consists of determining
$\tau_i\tau_j n$, provided $n$, $\tau_i n$ and $\tau_j n$ are
known. In order for this inductive process to work without
contradictions, equations
\begin{equation}\label{eq:klat}
\tau_i\tau_j n-n=a_{ij}(\tau_j n-\tau_i n),\qquad
a_{ij}=\frac{\langle n,\tau_i n-\tau_j n\rangle}
       {1-\langle \tau_i n,\tau_j n\rangle}\,
\end{equation}
must have a very special property. To see this, consider in detail
the case of $M=3$; higher dimensions do not add anything new. From
$n$ and $n_i$ one determines all $n_{ij}$ uniquely. After that,
one has, in principle, three different ways to determine
$n_{123}$, from three squares adjacent to this point; see
Fig.~\ref{cube}. These three values for $n_{123}$ have to
coincide, independently of initial conditions.
\begin{dfn}
A 2D system is called 3D consistent, if it can be imposed on all
two-dimensional faces of an elementary cube of $\bbZ^3$.
\end{dfn}
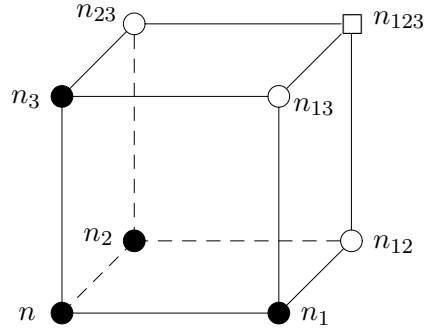
\begin{figure}[htbp]
\begin{center}
\setlength{\unitlength}{0.05em}
\begin{picture}(200,210)(0,0)
 \put(0,0){\circle*{15}}    \put(150,0){\circle*{15}}
 \put(0,150){\circle*{15}}  \put(150,150){\circle{15}}
 \put(50,50){\circle*{15}} \put(50,200){\circle{15}}
 \put(200,50){\circle{15}}
 \put(193,193){$\Box$}
 \path(0,0)(150,0)       \path(0,0)(0,150)
 \path(150,0)(150,142.5)   \path(0,150)(142.5,150)
 \path(0,150)(44.7,194.7)    \path(155.3,155.3)(193,193)
 \path(57.5,200)(193,200)
 \path(200,193)(200,57.5) \path(150,0)(194.7,44.7)
 \dashline[+30]{10}(0,0)(50,50)
 \dashline[+30]{10}(50,50)(50,192.5)
 \dashline[+30]{10}(50,50)(192.5,50)
 \put(-30,-5){$n$}
 \put(-35,145){$n_3$} \put(215,45){$n_{12}$}
 \put(165,-5){$n_1$} \put(160,140){$n_{13}$}
 \put(15,50){$n_2$}  \put(10,205){$n_{23}$}
 \put(215,200){$n_{123}$}
\end{picture}
\caption{3D consistency of 2D systems}\label{cube}
\end{center}
\end{figure}

\begin{thm}\label{Th: M-lat in quadric}
The 2D system (\ref{eq:klat}) governing T-nets in $\bbS^2$ is 3D
consistent.
\end{thm}
{\bf Proof.} This can be checked by a tiresome computation, which
can be however avoided by the following conceptual argument. Note
that T-nets in $\bbS^2$ are a result of imposing of two admissible
reductions on Q-nets, namely the T-reduction and the restriction
to a quadric $\bbS^2$. This reduces the effective dimension of the
system by 1 (allows to determine the fourth vertex of an
elementary quadrilateral from the three known ones), and transfers
the original 3D equation into the 3D consistency of the reduced 2D
equation. Indeed, after finding $n_{12}$, $n_{23}$ and $n_{13}$,
one can construct $n_{123}$ according to the Q-condition (as
intersection of three planes). Then both the T-condition and the
$\bbS^2$-condition are fulfilled for all three quadrilaterals
adjacent to $n_{123}$. Therefore, these quadrilaterals satisfy our
2D system. Note that this argument holds also for T-nets in an
arbitrary quadric (not necessarily in $\bbS^2$). $\Box$
\medskip

Thus, we have for $M\ge 3$ a system consisting of interrelated
two-dimensional building blocks (\ref{eq:klat}), for all pairs of
indices $i,j$ between 1 and $M$. The $M$-dimensional consistency
of this system means that all two-dimensional building blocks can be
imposed without contradictions. A set of initial data which determines
a solution of the system (\ref{eq:klat}) consists of
\begin{list}{(K$^\Delta$)}
{\setlength{\leftmargin}{1cm}}
\item values of $n$ on the coordinate axes $\PP_i$ for
  $1\le i\le M$.
\end{list}
\medskip

{\bf Transformations of discrete K-nets.} A natural generalization
of Definition \ref{dfn:bt} would be the following one.
\begin{dfn}\label{dfn:dbt}
A Weingarten pair of discrete K-nets $f,f^+:\bbZ^m\to\bbR^3$ is
called a {\em B\"acklund pair}, if the distance $|f^+-f|$ is
constant, i.e., does not depend on $u\in\bbZ^m$.
\end{dfn}
But, comparing this with Definition \ref{dfn:dks}, we see that the
net $F:\bbZ^m\times\{0,1\}\to\bbR^3$ with $F(u,0)=f(u)$ and
$F(u,1)=f^+(u)$ is an $M$-dimensional discrete K-net, where
$M=m+1$. Once more, transformations of discrete K-nets do not
differ from the nets themselves. To specify a B\"acklund transform
$f^+$ of a given $m$-dimensional discrete K-net $f$, or a Moutard
transform $n^+$ of the Gauss map $n$, one can prescribe the
following data:
\begin{list}{(B$^\Delta$)}
{\setlength{\leftmargin}{1cm}}
\item value of $n^+(0)$.
\end{list}
We will be mainly interested in transformations of discrete K-surfaces.
After the usual change of variables (\ref{eq:nbar}) we will have for the
Gauss map the equation
\begin{equation}\label{eq:klat 2d}
\tau_1\tau_2 n+n=a_{12}(\tau_1 n+\tau_2 n),\qquad
a_{12}=\frac{\langle n,\tau_1 n+\tau_2 n\rangle}
       {1+\langle \tau_1 n,\tau_2 n\rangle}\,,
\end{equation}
while for its transformation we will have the equations
\begin{eqnarray}
\tau_1 n^+-n=b_1(n^+-\tau_1 n), & \quad &
b_1=\frac{\langle n,\tau_1 n-n^+\rangle}{1-\langle\tau_1 n,n^+\rangle},
\label{eq: dbt 2d 1}\\
\tau_2 n^++n=b_2(n^++\tau_2 n), & \quad &
b_2=\frac{\langle n,\tau_2 n+n^+\rangle}{1+\langle\tau_2 n,n^+\rangle}.
\label{eq: dbt 2d 2}
\end{eqnarray}
Of course, this has to be supplied with the Lelieuvre formulas
(\ref{eq:dLel 2d}) and (\ref{eq:dwp}). Permutability of B\"acklund
transformations for discrete K-surfaces is a direct consequence of
the 3D consistency of T-nets in $\bbS^2$, and the smooth result
(Theorem \ref{thm:ks permut}) will be derived in Sect.
\ref{Subsect: ks conv} by the continuous limit.

\section{Discrete isothermic nets}
 \label{Subsect: discr iso}

Discrete isothermic surfaces were introduced in \cite{BP2}. See also
\cite{HHP, HJ1, HJ2, Sch1} for further developments.
\begin{dfn} \label{dfn:dis}
  A discrete O-surface $f:\bbZ^2\to\bbR^N$ is called a
  {\em discrete I-surface (discrete isothermic surface)}, if
  the cross-ratios of its elementary quadrilaterals satisfy
\begin{equation}\label{eq:dis property}
  q(f,f_1,f_{12},f_2)=-\frac{\alpha_1}{\alpha_2}\,,
\end{equation}
  where the functions $\alpha_i$ depend on $u_i$ only ($i=1,2$), and the
  usual notations $f_i=\tau_i f$, $f_{12}=\tau_1\tau_2 f$ are used.
\end{dfn}
It is natural to consider the functions $\alpha_i$ as defined on
edges of $\bbZ^2$ parallel to the coordinate axes $\PP_i$, with
the property that any two opposite edges of any elementary square
carry the same value of the corresponding $\alpha_i$. Such
functions $\alpha_i$  will be called a {\em labelling} of the
edges of $\bbZ^2$. Actually, it will be more convenient to
associate $-\alpha_2$ to the edges parallel to $\PP_2$. To
determine a labelling, one has to prescribe it on the coordinate
axes $\PP_i$. Thus, discrete I-surfaces are governed by a 2D
equation (\ref{eq:dis property}), and the following initial data
determine such a surface completely:
\newcounter{dI}
\begin{list}{(I$_\arabic{dI}^\Delta$)}
{\usecounter{dI}\setlength{\leftmargin}{1cm}}
\item values of $f$ in the coordinate axes $\PP_i\,$ $(i=1,2)$, i.e.,
two discrete curves $f\rest_{\PP_i}$ with a common intersection point
$f(0)$;
\item two functions $\alpha_i:\PP_i\to\bbR$ on the edges of the
coordinate axes $\PP_i$ for $i=1,2$, which yield an edge
labelling, $\alpha_i=\alpha_i(u_i)$.
\end{list}
A discrete analog of the function $s$ from eq. (\ref{eq:is prop})
is defined on the vertices of $\bbZ^2$.
\begin{thm}\label{Th: dis metric}
For a discrete isothermic surface $f$, there exists a function
$s:\bbZ^2\to\bbR$ (positive if both $\alpha_i>0$), such that
\begin{equation}\label{eq:dis metric}
|f_i-f|^2=\alpha_i ss_i\qquad (i=1,2).
\end{equation}
It is defined uniquely, up to a black-white rescaling (which is
fixed by prescribing $s$ at one point, say at $u=0$). Conversely,
existence of such a function $s$ for a discrete O-surface $f$
implies that $f$ is isothermic.
\end{thm}
The following property actually characterizes discrete isothermic
surfaces.
\begin{thm}[Dual discrete I-surface] \label{Th: dis dual}
Let $f:\bbZ^2\to\bbR^N$ be a discrete isothermic surface. Then the
$\bbR^N$-valued discrete one-form $\delta f^*$ defined by
\begin{equation}\label{eq: dis dual}
\delta_1 f^*=\alpha_1\frac{\delta_1 f}{|\delta_1 f|^2}=
\frac{\delta_2 f}{ss_1}\,,  \qquad
\delta_2 f^*=-\alpha_2\frac{\delta_2 f}{|\delta_2 f|^2}=
-\frac{\delta_2 f}{ss_2}\,,
\end{equation}
is closed. Its integration defines (up to a translation)
a surface $f^*:\bbZ^2\to\bbR^N$, called {\em dual} to the surface $f$,
or {\em Christoffel transform} of the surface $f$. The surface $f^*$
is discrete isothermic, with
\begin{equation}\label{eq:dis prop dual}
 q(f^*,f_1^*,f_{12}^*,f_2^*)=-\frac{\alpha_1}{\alpha_2}\,.
\end{equation}
\end{thm}
{\bf Proof.} We rewrite eq. (\ref{eq:dis property}) in several
equivalent forms. Recall that in this equation we regard $f$ as
belonging to the Clifford algebra ${\cal C}\ell(\bbR^N)$, and that
the cross-ratio is defined by eq. (\ref{eq:q}). One checks
straightforwardly that eq. (\ref{eq:dis property}) is equivalent
to:
\begin{equation}\label{eq:dis 3leg I}
(\alpha_1+\alpha_2)(f_{12}-f)^{-1}=
\alpha_1(f_1-f)^{-1}+\alpha_2(f_2-f)^{-1},
\end{equation}
or, interchanging the roles of $f$ and $f_{12}$, to:
\begin{equation}\label{eq:dis 3leg II}
(\alpha_1+\alpha_2)(f_{12}-f)^{-1}=\alpha_1(f_{12}-f_2)^{-1}+
\alpha_2(f_{12}-f_1)^{-1}.
\end{equation}
Taking into account that $(f_i-f)^{-1}=-(f_i-f)/|f_i-f|^2$, we find:
\begin{equation}\label{eq:dis dual 1}
\alpha_1\frac{f_1-f}{|f_1-f|^2}+\alpha_2\frac{f_2-f}{|f_2-f|^2}=
\alpha_1\frac{f_{12}-f_2}{|f_{12}-f_2|^2}+
\alpha_2\frac{f_{12}-f_1}{|f_{12}-f_1|^2}\,,
\end{equation}
or, due to (\ref{eq:dis metric}),
\begin{equation}\label{eq:dis dual 2}
\frac{f_1-f}{ss_1}+\frac{f_2-f}{ss_2}=
\frac{f_{12}-f_2}{s_2s_{12}}+\frac{f_{12}-f_1}{s_1s_{12}}\,,
\end{equation}
This is equivalent to the first statement of the theorem. The second one
follows readily from the Clifford algebra expressions
$f_1^*-f^*=\alpha_1(f_1-f)^{-1}$ and $f_2^*-f^*=-\alpha_2(f_2-f)^{-1}$.
Note that
$|f_i^*-f^*|^2=\alpha_i^2|f_i-f|^{-2}=\alpha_i(ss_i)^{-1},$
so that $s^*=s^{-1}$ can be taken as the analog of the function $s$ for
the dual surface $f^*$. $\Box$

\begin{dfn}[Discrete Darboux transformation] \label{dfn:ddt}
  \quad A Ribaucour \linebreak transform
  $f^+:\bbZ^2\to\bbR^N$ of a given discrete isothermic surface
  $f:\bbZ^2\to\bbR^N$ is called a {\em Darboux transform}, if the
  cross-ratios of its elementary quadrilaterals can be likewise
  factorized:
\begin{equation}\label{eq:ddt property}
  q(f^+,f_1^+,f_{12}^+,f_2^+)=-\frac{\alpha_1}{\alpha_2}\,.
\end{equation}
\end{dfn}
It can be demonstrated that, apart from the trivial case when $f^+$ is a
M\"obius transform of $f$, the Darboux transform is given by the following
formulas:
\begin{equation}\label{eq:ddt}
  q(f,f_1,f_1^+,f^+)=\frac{\alpha_1}{c}, \qquad
  q(f,f_2,f_2^+,f^+)=-\frac{\alpha_2}{c}.
\end{equation}
where $c\in\bbR$ is its {\em parameter}.
Clearly, following data determine a Darboux transform $f^+$ of a given
discrete isothermic surface $f$ uniquely:
\newcounter{dD}
\begin{list}{(D$_\arabic{dD}^\Delta$)}
{\usecounter{dD}\setlength{\leftmargin}{1cm}}
\item a point $f^+(0)$;
\item a real number $c$, the parameter of the transformation.
\end{list}
Comparing formulas (\ref{eq:dis property}), (\ref{eq:ddt property}), and
(\ref{eq:ddt}), we see that iterating a Darboux transformation is nothing
but adding a third dimension to a two-dimensional discrete isothermic net,
$f^+=f_3$. The parameter $c$ plays the role of the function $\alpha_3$,
attached to all edges parallel to the third lattice direction. Formula
(\ref{eq:dis metric}) for the third direction reads:
\begin{equation}
 |f^+-f|^2=css^+
\end{equation}
(which literally coincides with the corresponding formula
(\ref{eq: Darboux int}) for the smooth case). All this is possible
due to the following fundamental statement about multidimensional
discrete isothermic nets:
\begin{thm}\label{Thm 3D cr}
The 2D cross-ratio equation
\begin{equation}\label{eq:i/j}
q(f,f_i,f_{ij},f_j)=\frac{\alpha_i}{\alpha_j}
\end{equation}
is 3D consistent for any labelling $\alpha_i$ of the edges.
\end{thm}
This statement is due to \cite{HHP}; its generalization, where
$\cC\ell(\bbR^N)$ is replaced by an arbitrary associative algebra,
was proven in \cite{BS2}, and will be given in Theorem \ref{Thm 3d
cr noncommut}.

{\bf Remark}. In view of the fundamental importance of the
consistency property which holds if the right-hand side of the
cross-ratio system is factorized as in (\ref{eq:i/j}), one might
wonder why we introduced the minus sign in eq. (\ref{eq:dis
property}) for discrete I-surfaces, which propagated also into eq.
(\ref{eq:ddt}) for the Darboux transformations. This choice is
motivated by the convenience of passing to the continuous limit
only. The fundamental factorization form (\ref{eq:i/j}) is
restored by the change $\alpha_2\to -\alpha_2$.
\medskip

{\bf M\"obius-geometric characterization of discrete isothermic
surfaces.} Cross-ratio of four concircular points is a M\"obius
invariant quantity. Therefore, the notions of discrete isothermic
surfaces and of their Darboux transformations are also M\"obius
invariant. To give their characterization within the
M\"obius-geometric formalism, we note first of all that, according
to Theorem \ref{Th: dos in K}, eq. (\ref{eq:dis dual 2}) holds
also for the image $\hat{f}$ of the net $f$ in the quadric
$\bbQ_{\,0}^N$:
\[
\frac{\hat{f}_1-\hat{f}}{ss_1}+\frac{\hat{f}_2-\hat{f}}{ss_2}=
\frac{\hat{f}_{12}-\hat{f}_2}{s_2s_{12}}+
\frac{\hat{f}_{12}-\hat{f}_1}{s_1s_{12}}\,.
\]
This yields a discrete analog of Theorem \ref{Th: iso Mob}:
\begin{thm}\label{Th: dis Mob}
The lift $\hat{s}=s^{-1}\hat{f}:\bbZ^2\to\cn$ of $f$ to the light
cone of $\bbR^{N+1,1}$ satisfies the {\em discrete Moutard
equation}
\begin{equation}\label{eq: dis Mou}
\tau_1\tau_2\hat{s}+\hat{s}=a_{12}(\tau_1\hat{s}+\tau_2\hat{s}),
\end{equation}
with $a_{12}=(\tau_1\tau_2 s^{-1}+s^{-1})/(\tau_1 s^{-1}+\tau_2 s^{-1})$.

Conversely, given a discrete M-net $\hat{s}:\bbZ^2\to\cn$ in the
light cone, let the functions  $s:\bbZ^2\to\bbR_+$ and
$f:\bbZ^2\to\bbR^N$ be defined by
\[
\hat{s}=s^{-1}(f+\ee_0+|f|^2\ee_\infty)
\]
(so that $s^{-1}$ is the $\ee_0$-component, and $s^{-1}f$ is the
$\bbR^N$-part of $\hat{s}$ in the basis
$\ee_1,\ldots,\ee_N,\ee_0,\ee_\infty$). Then $f$ is a discrete
isothermic surface.
\end{thm}

Eq. (\ref{eq: dis Mou}) is a 2D equation in $\cn$ which determines
$\tau_1\tau_2\hat{s}$ from $\hat{s}$, $\tau_1\hat{s}$, $\tau_2\hat{s}$.
This yields:
\begin{equation}\label{eq: dis Mou a}
a_{12}=\frac{\langle \hat{s},\tau_1\hat{s}+\tau_2\hat{s}\rangle}
 {\langle \tau_1\hat{s},\tau_2\hat{s}\rangle}\,.
\end{equation}
>From this one readily sees that for a discrete M-net in $\cn$ the
quantities $\langle \hat{s},\tau_i\hat{s}\rangle$ depend on $u_i$
only. They are related to the labelling $\alpha_i$ of Definition
\ref{dfn:dis} by $\langle
\hat{s},\tau_i\hat{s}\rangle=-\alpha_i/2$. Thus, the labelling is
already encoded in the lift $\hat{s}$ of a discrete I-surface $f$.

For the consistent cross-ratio system (\ref{eq:i/j}) one can
introduce the scalar function $s:\bbZ^{M}\to\bbR$ (not necessarily
positive) according to eq. (\ref{eq:dis metric}) with $1\le i\le
M$, which yields a multidimensional version of the lift
$\hat{s}=s^{-1}\hat{f}:\bbZ^{M}\to\cn$. It satisfies a discrete
Moutard equation with minus signs:
\begin{equation} \label{eq:minusMou}
\tau_i\tau_j\hat{s}-\hat{s}=a_{ij}(\tau_j\hat{s}-\tau_i\hat{s}).
\end{equation}
This is an instance of 3D consistent T-nets in a quadric (see the
proof of Theorem \ref{Th: M-lat in quadric}), and describes
Darboux transformations of discrete I-surfaces and their
permutability properties. To perform a smooth limit a usual change
of signs (see Sect.\ref{Subsect: discr Mou}) is required. This
results in eq. (\ref{eq: dis Mou}) for the surface and in
\begin{equation}
\tau_1\hat{s}^+-\hat{s}=b_1(\hat{s}^+-\tau_1\hat{s}),\qquad
\tau_2\hat{s}^++\hat{s}=b_2(\hat{s}^++\tau_2\hat{s})
\end{equation}
for its Darboux transformations.

\section{Discrete surfaces made of spheres}
 \label{Subsect S-iso}

It is common to think of discrete surfaces as of maps
$f:\bbZ^2\to\bbR^N$. However, it is not the only possibility. One
of the further ideas is to ``blow up'' the points $f$, making them
to (small) hyperspheres in $\bbR^N$. We discuss here two
realizations of this idea, which generalize the notions of
discrete O-nets and of discrete I-surfaces, maintaining important
features of these notions, including transformations with
permutability properties.
\medskip

{\bf S-orthogonal nets.} We first generalize $M$-dimensional
discrete O-nets $f$ in $\bbR^N$. Blowing up the points to spheres
is done in the following way:
\begin{dfn}\label{dfn:S-circ}
An\, {\em S-circular net} is a map
\begin{equation}\label{eq:S-map}
S:\,\bbZ^M\to\big\{\mbox{\rm
non-oriented\;\;hyperspheres\;\;in\;\;}\bbR^N\big\}
\end{equation}
such that, for any elementary square of $\,\bbZ^M$, the
corresponding four hyperspheres $S$, $\tau_i S$, $\tau_j S$ and
$\tau_i\tau_j S$ have a common orthogonal circle.
\end{dfn}
Clearly, this reduces to Definition \ref{dfn:dos}, if the radii of
all hyperspheres become infinitely small. Observe also that the
plane of the common orthogonal circle contains the Euclidean
centers of all four hyperspheres.

Recall (see Appendix to Chapter \ref{Sect: smooth diffe geo}) that
hyperspheres in $\bbR^N$ can be represented as elements of
$\bbP(\bbR^{N+1,1}_{\,\rm out})$, where
\begin{equation}\label{eq: Rout}
 \bbR^{N+1,1}_{\,\rm out}=\big\{\hat{s}\in\bbR^{N+1,1}:\,\langle
 \hat{s},\hat{s}\rangle >0\big\}
\end{equation}
is the space-like part of $\bbR^{N+1,1}$. Thus, blowing up points
to spheres corresponds to the step away from $\cn$ to
$\bbR^{N+1,1}_{\,\rm out}$. To establish a characterization of
S-circular nets similar to that of Theorem \ref{Th: dos in K}, the
following statement is needed.
\begin{thm}
Four points of $\,\bbR^{N+1,1}_{\,\rm out}$ are linearly dependent
if and only if the corresponding four hyperspheres in $\bbR^N$
\begin{itemize}
\item[(i)] have a common orthogonal circle, or
\item[(ii)] intersect along an $(N-3)$-sphere, or else
\item[(iii)] intersect at exactly one point.
\end{itemize}
Case (iii) can be regarded as a degenerate case of both (i) and
(ii).
\end{thm}
The distinction between these three cases depends on the signature
of the quadratic form obtained as the restriction of the Minkowski
scalar product to the $(N-1)$-dimensional (Lorentz-)orthogonal
complement of the three-dimensional vector subspace of
$\bbR^{N+1,1}$ spanned by our four points. This form is
positive-definite in the case (i), it has the signature $(N-2,1)$
in the case (ii), and it is degenerate in the case (iii). It is
understood that for $N=3$ the 0-sphere in the case (ii) is a point
pair.

Since it is difficult to control the signature of this quadratic
form, it makes sense to introduce a notion which is somewhat wider
than that of a S-circular net.
\begin{dfn}\label{dfn:S-O}
An\, {\em S-orthogonal net} is a map (\ref{eq:S-map}) such that,
for any elementary square of $\,\bbZ^M$, the points $\hat{s}$,
$\tau_i\hat{s}$, $\tau_j\hat{s}$ and $\tau_i\tau_j\hat{s}$,
corresponding to four spheres $S$, $\tau_i S$, $\tau_j S$ and
$\tau_i\tau_j S$, are linearly dependent.
\end{dfn}

For $N=3$, S-orthogonal nets with all intersections of type (ii)
are natural discrete analogs of sphere congruences parametrized
along principal lines, because four infinitesimally neighboring spheres
of such a congruence intersect this way, the pairs of intersection points
comprising two enveloping surfaces of the congruence (see, e.g.,
\cite{E2}).

\begin{thm}
Six circles (corresponding to the faces) of an elementary cube of
a S-circular net lie on a 2-sphere. For $N=3$ these spheres
comprise an S-orthogonal net.
\end{thm}

A convenient choice of representatives $\hat{s}$ of hyperspheres
$S$ in a fixed affine hyperplane of $\bbR^{N+1,1}_{\,\rm out}$ is
given by eq. (\ref{eq:sph0}). Therefore:
\begin{thm}
A map (\ref{eq:S-map}) is an S-orthogonal net, if and only if the
corresponding map
\[
\hat{s}:\,\bbZ^M\to\bbR^{N+1,1}_{\,\rm out}\cap\{\xi_0=1\},\qquad
\hat{s}=c+\ee_0+\big(|c|^2-r^2\big)\ee_\infty,
\]
is a discrete Q-net in $\bbR^{N+1,1}$, i.e., if the centers
$c:\bbZ^M\to\bbR^N$ of the spheres $S$ form a discrete Q-net in
$\bbR^N$, and their radii $r:\bbZ^M\to\bbR_+$  are such that the
function $|c|^2-r^2$ satisfies the same equation (\ref{eq:dcn
property}) as the centers $c$.
\end{thm}

A smooth version of this theorem for arbitrary sphere congruences
in $\bbR^3$ parametrized along principal lines can be found, e.g., in
\cite{E2}, \S 99.

The multidimensional consistency of discrete quadrilateral nets
can be immediately transferred into the corresponding property of
S-orthogonal nets, with the following reservation: given seven
points $\hat{s}$, $\tau_i\hat{s}$, $\tau_i\tau_j\hat{s}$ in
$\bbR^{N+1,1}_{\,\rm out}\cap\{\xi_0=1\}$, the Q-property
(planarity condition) uniquely defines the eighth point
$\tau_1\tau_2\tau_3\hat{s}$ in $\bbR^{N+1,1}\cap\{\xi_0=1\}$,
which, however, might get outside of $\bbR^{N+1,1}_{\,\rm out}$,
and therefore might not represent a real hypersphere. Thus, the
corresponding discrete 3D system is well-defined only on a proper
open subset of the definition domain. As long as it is defined, it
can be used to produce ``S-Ribaucour'' transformations of
S-orthogonal nets, with usual permutability properties.
\medskip

{\bf S-isothermic surfaces.}  Recall that, according to Theorem
\ref{Th: dis Mob}, discrete I-surfaces are characterized as
discrete M-surfaces (or T-surfaces) in $\cn$. In this context, we
blow up points to spheres by replacing $\cn$ by $\cn_{\,\kappa}$
(see eqs. (\ref{eq:sph can}), (\ref{eq:Lkappa})):
\begin{dfn}\label{dfn:S-iso}
A map $S:\bbZ^2\to\{{\mbox {\rm oriented hyperspheres in\;\,}}
\bbR^N\}$ is called an\, {\em S-isothermic surface}, if the
corresponding map $\hat{s}:\bbZ^2\to\cn_{\,\kappa}$ is a discrete
M-surface (T-surface).
\end{dfn}
S-isothermic surfaces, along with their dual surfaces were
originally introduced in \cite{BP3} for the special case of
touching spheres. The general class of Definition \ref{dfn:S-iso},
together with Darboux transformations and dual surfaces, is due to
\cite{Hof}.

S-isothermic surfaces are governed by the equation (\ref{eq: dis
Mou}) with
\begin{equation}\label{eq: S-iso Mou a}
a_{12}=\frac{\langle \hat{s},\tau_1\hat{s}+\tau_2\hat{s}\rangle}
{\kappa^2+\langle \tau_1\hat{s},\tau_2\hat{s}\rangle}\,,
\end{equation}
and have the property that the quantities $\alpha_i=\langle
\hat{s},\tau_i\hat{s}\rangle$ depend on $u_i$ only.
If radii of all hyperspheres become uniformly small, $r(u)\sim\kappa s(u)$,
$\kappa\to 0$, then in the limit we recover discrete isothermic surfaces,
as characterized in Theorem \ref{Th: dis Mob}.

Consistency of discrete T-nets in $\cn_{\,\kappa}$ implies, in
particular, Darboux transformations for S-isothermic surfaces,
governed by eq. (\ref{eq:minusMou}). A Darboux transform
$\hat{s}^+:\bbZ^2\to\cn_{\,\kappa}$ of a given S-isothermic
surface $\hat{s}$ is uniquely specified by a choice of one of its
spheres $\hat{s}^+(0)$.

Geometric properties of S-isothermic surfaces are the following.
First of all, S-isothermic surfaces form a subclass of
S-orthogonal nets. Further, the quantities
$\langle\hat{s},\tau_i\hat{s}\rangle$ depend on $u_i$ only, and
have the meaning of cosines of the intersection angles of the
neighboring spheres (resp., of their so called inversive distances
if they do not intersect). An important characterization is the
following generalization of Theorem \ref{Th: dis dual}.
\begin{thm}[Dual S-isothermic surface] \label{Th: S-iso dual} \quad

Let $S:\bbZ^2\to\{{\mbox {\rm oriented hyperspheres
in\;\;}}\bbR^N\}$ be an S-isothermic surface. Denote the Euclidean
centers and radii of $S$ by $c:\bbZ^2\to\bbR^N$ and
$r:\bbZ^2\to\bbR$, respectively. Then the $\bbR^N$-valued discrete
one-form $\delta c^*$ defined by
\begin{equation}\label{eq: S-iso dual}
\delta_1 c^*=\frac{\delta_1 c}{rr_1}\,,\qquad
\delta_2 c^*=-\frac{\delta_2 c}{rr_2}
\end{equation}
is closed, so that its integration defines (up to a translation) a
function $c^*:\bbZ^2\to\bbR^N$. Define also $r^*:\bbZ^2\to\bbR$ by
$r^*=r^{-1}$. Then the hyperspheres $S^*$ with the centers $c^*$
and radii $r^*$ form an S-isothermic surface, called {\em dual} to
$S$.
\end{thm}
{\bf Proof.} Consider eq. (\ref{eq: dis Mou}), in terms of
\[
\hat{s}=r^{-1}\big(c+\ee_0+(|c|^2-r^2)\ee_\infty\big).
\]
Its $\ee_0$-part yields:
$a_{12}=(r_{12}^{-1}+r^{-1})/(r_1^{-1}+r_2^{-1})$. This allows us
to rewrite eq. (\ref{eq: dis Mou}) as
\begin{equation}\label{eq: S-iso dual eq}
(r_1^{-1}+r_2^{-1})(\hat{s}_{12}+\hat{s})=
(r_{12}^{-1}+r^{-1})(\hat{s}_1+\hat{s}_2).
\end{equation}
A direct computation shows that the $\bbR^N$-part of
this equation can be rewritten as
\begin{equation}
\frac{c_1-c}{rr_1}+\frac{c_2-c}{rr_2}=
\frac{c_{12}-c_2}{r_2r_{12}}+\frac{c_{12}-c_1}{r_1r_{12}},
\end{equation}
which is equivalent to closeness of the form $\delta c^*$ defined
by (\ref{eq: S-iso dual}). In the same way, the $\ee_\infty$-part
of eq. (\ref{eq: S-iso dual eq}) is equivalent to closeness of the
discrete form $\delta w$ defined by
\begin{equation*}\label{eq: S-iso dual w}
\delta_1 w=\frac{\delta_1(|c|^2-r^2)}{rr_1}\,,\qquad
\delta_2 w=-\frac{\delta_2(|c|^2-r^2)}{rr_2}\,.
\end{equation*}
For similar reasons, the second claim of the theorem is equivalent
to the closeness of the form
\begin{equation*}\label{eq: S-iso dual z}
\delta_1 z=\frac{\delta_1\big(|c^*|^2-(r^*)^2\big)}{r^*r_1^*}\,,\qquad
\delta_2 z=-\frac{\delta_2\big(|c^*|^2-(r^*)^2\big)}{r^*r_2^*}\,.
\end{equation*}
But one easily checks that
\[
\delta_1 w=\langle c_1^*-c^*,c_1+c\rangle-\frac{r_1}{r}+\frac{r}{r_1}\,,
\qquad
\delta_2 w=\langle c_2^*-c^*,c_2+c\rangle+\frac{r_2}{r}-\frac{r}{r_2}\,,
\]
\[
\delta_1 z=\langle c_1-c,c_1^*+c^*\rangle-\frac{r}{r_1}+\frac{r_1}{r}\,,
\qquad
\delta_2 z=\langle c_2-c,c_2^*+c^*\rangle+\frac{r}{r_2}-\frac{r_2}{r}\,.
\]
The sum of these two one-forms is closed:
\[
\delta_1(w+z)=2\langle c_1^*,c_1\rangle-2\langle c^*,c\rangle,\quad
\delta_2(w+z)=2\langle c_2^*,c_2\rangle-2\langle c^*,c\rangle,
\]
therefore they are closed simultaneously. $\Box$

\chapter{Approximation} \label{Sect: approx}

\section{Discrete hyperbolic systems}
 \label{Sect: subsect discr hyp}

To formulate the most general scheme covering all the situations
encountered so far, we have to put our hyperbolic systems into the
first order form. It should be stressed that this is necessary
only for general theoretical considerations, and will never be
done for concrete examples.
\begin{dfn} A {\em hyperbolic system of first order partial difference
equations} is a system of the form
\begin{equation}\label{eq:prehyperbolic}
  \delta_i x_k=g_{k,i}(x), \quad i\in\cE_k.
\end{equation}
for functions $x_k:\bbZ^M\to\cX_k$ with values in Banach spaces $\cX_k$.
For each $x_k$, eqs. (\ref{eq:prehyperbolic}) are given for
$i\in\cE_k\subset\{1,\ldots,M\}$, the {\em evolution directions} of $x_k$.
The complement $\cS_k=\{1,\ldots,M\}\setminus \cE_k$ consists of
{\em static directions} of $x_k$.
\end{dfn}
We think of the variable $x_k(u)$ as attached to the elementary cell
$\cC_k$ of dimension $\#\cS_k$ adjacent to the point
$u\in\bbZ^M$ and parallel to $\PP_{\cS_k}$:
\[
\cC_k=\Big\{u+\sum_{i\in\cS_k}\mu_ie_i:\;\mu_i\in[0,1]\Big\}.
\]
Here, recall,
\[
  \PP_\cS=\{u\in\bbZ^M: \;u_i=0 \quad\mbox{if}\quad i\notin\cS\},
\]
for an index set $\cS\subset\{1,\ldots,M\}$.
\begin{dfn}\quad

 1) A {\em local Goursat problem} for the hyperbolic system
 (\ref{eq:prehyperbolic}) consists of finding a solution $x_k$
 for all $k$ and for all cells $\,\cC_k$ within the elementary
 cube of $\bbZ^M$ at the origin from the prescribed values $x_k(0)$.
 A {\em global Goursat problem} consists of finding a solution
 of (\ref{eq:prehyperbolic}) on $\bbZ^M$ subject to the following
 initial data:
\begin{equation}\label{eq:preGoursat}
  x_k\rest_{\PP_{\cS_k}}=X_k,
\end{equation}
 where $X_k:\PP_{\cS_k}\to\cX_k$ are given functions.

 2) The system (\ref{eq:prehyperbolic}) is called
 {\em consistent}, if the local Goursat problem for it is uniquely
 solvable for arbitrary initial data $x(0)$.
\end{dfn}
The following rather obvious but extremely important statement
holds:
\begin{thm}  A Goursat problem for a consistent hyperbolic system
(\ref{eq:prehyperbolic}) has a unique solution $x$ on all of $\bbZ^M$.
\end{thm}

Consistency conditions read:
$\delta_j\delta_i x_k=\delta_i\delta_j x_k$ for all $i\neq j$.
Substituting eqs. (\ref{eq:prehyperbolic}), one gets the
following equations:
\begin{equation}  \label{eq:comp2}
  \delta_j g_{k,i}(x)=\delta_i g_{k,j}(x),\quad i\neq j,
\end{equation}
or $ g_{k,i}\big(x+g_j(x)\big)-g_{k,i}(x)=
g_{k,j}\big(x+g_i(x)\big)-g_{k,j}(x)$, where $g_i(x)$ is a vector function
whose $\ell$-th component is equal to $g_{\ell,i}(x)$, if $i\in\cE_\ell$,
and is undefined otherwise.
\begin{lem}\label{subsyst}
For a consistent system of hyperbolic equations
(\ref{eq:prehyperbolic}), the function $g_{k,i}$ depends on those
components $x_\ell$ only for which $\cS_\ell\subset\cS_k\cup\{i\}$.
\end{lem}
{\bf Proof.}
Eqs. (\ref{eq:comp2}) have to hold identically in $x$. This
implies that the function $g_{k,i}$ can depend on those components
$x_{\ell}$ only, for which $\delta_j x_\ell$ is defined,
i.e., for which $j\in\cE_\ell$. As (\ref{eq:comp2}) has to be
satisfied for all $j\in\cE_k$, $j\neq i$, one obtains that for
these $\ell$ there holds $\cE_k\setminus\{i\}\subset\cE_\ell$.
$\Box$ \medskip

It follows from Lemma \ref{subsyst} that for any subset
$\cS\subset\{1,\ldots,M\}$, equations of (\ref{eq:prehyperbolic})
for $k$ with $\cS_k\subset\cS$ and for $i\in\cS$ form a closed
subsystem, in the sense that $g_{k,i}$ depend on $x_\ell$ with
$\cS_\ell\subset\cS$ only.

\begin{dfn}
The {\em essential dimension} $d$ of the system
(\ref{eq:prehyperbolic}) is given by
\begin{equation}\label{eq: dim}
d=1+\max_k\big(\#\cS_k\big).
\end{equation}
\end{dfn}
If $d=M$, system (\ref{eq:prehyperbolic}) has no lower-dimensional
hyperbolic subsystems. If $d<M$, then $d$-dimensional subsystems
corresponding to $\cS$ with $\#\cS=d$ are hyperbolic. In this
case, {\itbf consistency} of system (\ref{eq:prehyperbolic}) is a
manifestation of a very special property of its $d$-dimensional
subsystems, which we suggest to treat as the {\itbf discrete
integrability} (at least under some further conditions, excluding
certain non-interesting situations, like trivial evolution in some
of the directions). Sect. \ref{Sect: consist} will be devoted to
giving a solid background for this suggestion in the case $d=2$.

\section{Discrete approximation in hyperbolic systems}

To handle with approximation results for discrete geometric models, we
need to introduce small parameters into hyperbolic systems of partial
difference equations. The definition domain of our functions becomes
\[
\PP^\veps=\eps_1\bbZ\times\cdots\times\eps_M\bbZ.
\]
If $\eps_i=0$ for some
index $i$, the respective component in $\PP^{\veps}$ is replaced by $\bbR$.
For instance, if $\veps=(0,\ldots,0)$, then $\PP^\veps=\bbR^M$.
So the domains $\PP^\veps$ posses continuous and discrete directions,
with mesh sizes depending on the parameters $\eps_i$.
Definitions of translations and difference quotients are modified
for functions on $\PP^\veps$ in an obvious way:
\[
  (\tau_i f)(u)=f(u+\eps_ie_i),\quad
  (\delta_i f)(u)=\frac{1}{\eps_i}
  \big(f(u+\eps_ie_i)-f(u)\big).
\]
If $\eps_i=0$, then $\delta_i$ is naturally replaced by the partial derivative
$\partial_i$. For a multi-index $\alpha=(\alpha_1,\ldots,\alpha_M)$, we set
$\delta^\alpha=\delta_1^{\alpha_1}\ldots\delta_M^{\alpha_M}$.

The definition of elementary cells $\cC_k$, carrying the variables $x_k$,
is modified as follows:
\[
\cC_k=\Big\{u+\sum_{i\in\cS_k}\mu_ie_i:\;\mu_i\in[0,\eps_i]\Big\}
\]
(so that the cell size shrinks to zero in the directions with $\eps_i=0$).
We see how the discreteness helps to organize the ideas: in the continuous
case, when all $\eps_i=0$, all the functions $x_k$ live at {\em points},
independently on the dimensions $\#\cS_k$ of their static spaces. In the
discrete case, when all $\eps_i>0$, one can clearly distinguish between
functions living on vertices (those without static directions), on edges
(those with exactly one static direction), on elementary squares (those
with exactly two static directions), etc.

Having in mind the limit $\eps\ra 0$,
we will treat only the case when the first $m\le M$ parameters go to
zero, $\eps_1=\cdots=\eps_m=\eps$, while the other $M-m$ ones remain
constant, $\eps_{m+1}=\cdots=\eps_M=1$. Thus, in this case
$\PP^\veps=(\eps\bbZ)^m\times\bbZ^{M-m}$, and we set
$\PP=\PP^0=\bbR^m\times\bbZ^{M-m}$.  Assuming that the functions
$g_{k,i}=g_{k,i}^\eps$ on the right-hand sides of (\ref{eq:prehyperbolic})
depend on $\eps$ smoothly and have limits as $\eps\to 0$, we will study the
convergence of solutions $x^\eps$ of the difference hyperbolic system
(\ref{eq:prehyperbolic}) towards solutions $x^0$ of the limiting
differential(-difference) hyperbolic system
\begin{eqnarray}
    \partial_i x_k & = & g^0_{k,i}(x), \quad 1\le i\leq m,
    \label{eq:prehyperbolic2}\\
    \delta_i x_k & = & g^0_{k,i}(x), \quad m+1\le i\le M.
    \label{eq:prehyperbolic3}
\end{eqnarray}
Naturally, (\ref{eq:prehyperbolic2}), (\ref{eq:prehyperbolic3}) describe
the respective $m$-dimensional smooth geometry with $M-m$ permutable
transformations.

Throughout this section, a {\em smooth} function
$g:\mathcal{D}\ra\cX$ is one that is infinitely often
differentiable on its domain, $g\in C^\infty(\mathcal{D})$. For a
compact set $\cK\subset\mathcal{D}$, we say that a sequence of
smooth functions $g^\eps$ converges toward a smooth function
$g$ with the order $\OH(\eps)$ in $C^\infty(\cK)$, if
\[
  \|g^\eps-g\|_{C^\ell(\cK)}\;\leq c_\ell\eps
\]
with suitable constants $c_\ell$ for any $\ell\in\bbN$.
Convergence in $C^\infty(\mathcal{D})$ means convergence in
$C^\infty(\cK)$ for all compact sets $\cK\subset\mathcal{D}$.

Convergence of discrete functions (defined on lattices $\PP^\eps$
with different $\eps$) is understood as follows. We say that a
family of discrete functions $x^\eps: \PP^\eps\to\cX$ converges to
a function $x:\PP\to\cX$ with the order $\OH(\eps)$ in
$C^\infty(\PP)$, if for any multi-index
$\alpha=(\alpha_1,\ldots,\alpha_M)$ there holds
\[
 \sup_{u\in\PP^\eps}|\delta^\alpha(x^\eps-x)(u)|\leq c_\alpha\eps
\]
with some constants $c_\alpha$. The symbol $x$ on the left-hand side is
understood as a restriction of $x$ to $\PP^\eps\subset\PP$.

Finally, we mention that we are mainly concerned with {\em local}
problems, so that actually we work with bounded domains of the
lattices $\PP^\eps$,
\[
 \PP^\eps(r)=\{u\in\PP^\eps:\;u_i\in[0,r],\;\;u_j\in\{0,1\},\;\;
 1\leq i\leq m<j\leq M\}.
\]
Each $\PP^\eps(r)$ contains only finitely many points (though
their number grows infinitely with $\eps\to 0$, if $r$ remains
fixed). For functions $u^\eps$ defined on a bounded lattice domain
$\PP^\eps(r)$ only, the notion of convergence is modified in an
obvious way: the supremum is taken only over those lattice sites
$u$, where the respective difference quotient $\delta^\alpha
x^\eps(u)$ exists.

\begin{thm}  \label{thm:approx}
Consider a Goursat problem (\ref{eq:prehyperbolic}), (\ref{eq:preGoursat})
for a hyperbolic system  of difference  equations. Suppose that:
\begin{itemize}
  \item[i)] the discrete system (\ref{eq:prehyperbolic}) is consistent for
    all $\eps>0$;
  \item[ii)] functions $g^\eps_{k,i}$ on the right-hand sides of
    eqs. (\ref{eq:prehyperbolic}) converge with the order $\OH(\eps)$
    in $C^\infty(\cX_k)$ to smooth functions $g^0_{k,i}$;
  \item[iii)] the Goursat data $X^\eps_k$ converge with the order
    $\OH(\eps)$ in $C^\infty(\PP_{\cS(k)})$ to smooth functions $X^0_k$.
  \end{itemize}
Then, for $\eps>0$ small enough, solution $x^\eps$ of the Goursat
problem exists and is unique on $\PP^\eps(r)$ for a suitable
($\eps$-independent) $r>0$; moreover, solutions $x^\eps$ converge
to a smooth function $x^0$ with the order $\OH(\eps)$ in
$C^\infty(\PP(r))$; this function $x^0$ is a unique solution on
$\PP(r)$ of the Goursat problem for (\ref{eq:prehyperbolic2}),
(\ref{eq:prehyperbolic3}) with the Goursat data $X_k^0$.
\end{thm}
{\bf Proof} of this theorem is technical, based on the discrete
Gronwall-type estimates. Details can be found in \cite{BMS1}.
$\Box$
\smallskip

In condition {\em ii)}, convergence $g_{k,i}^\eps\to g_{k,i}^0$ in
$C^\infty(\cX_k)$, i.e., on every compact subset of $\cX_k$, is
assumed for simplicity of presentation only. In applications,
functions $g^\eps_{k,i}$ are often defined on certain subdomains
$\mathcal{D}_k^\eps\subset\cX_k$, with the property that
$\mathcal{D}_k^0$ is open and dense in $\cX_k$. In such a case,
one requires in {\em ii)} the convergence $g^\eps_{k,i}\to
g^0_{k,i}$ in $C^\infty(\mathcal{D}^0)$. Then conclusions of
Theorem \ref{thm:approx} hold for generic initial data.

As for condition {\em iii)}, smooth data
$X^0_k:\PP_{\cS_k}\to\cX_k$ are usually given a priori, and
discrete data $X^\eps_k$ are obtained by restriction to the
lattice: $X_k^\eps=X_k^0\rest_{\PP_{\cS_k}^\eps}$. In such a
situation, condition {\em iii)} is fulfilled automatically.

\section{Conjugate nets}

{\bf Discretization of a conjugate net.} Recall that a smooth
conjugate net $f:\bbR^m\to\bbR^N$ is determined by the initial
data (Q$_{1,2}$) (see Sect. \ref{Sect: conj}), while a discrete
Q-net $f^\eps:(\eps\bbZ)^m\to\bbR^N$ is determined by the initial
data (Q$_{1,2}^\Delta$) (see Sect. \ref{Sect: subsect discr
conj}). We now demonstrate how to produce from the data
(Q$_{1,2}$) certain discrete data (Q$_{1,2}^\Delta$), which will
assure the convergence of the corresponding discrete Q-nets to the
smooth conjugate net.

Define the discrete curves $f^\eps\rest_{\PP_i^\eps}$ by restricting
the curves $f\rest_{\PP_i}$ to the lattice points:
\[
f^\eps(u)=f(u),\qquad u\in\PP^\eps_i,\qquad 1\le i\le m.
\]
Similarly, define the plaquette functions
$c_{ij}^\eps\rest_{\PP_{ij}^\eps}$ by restricting
$c_{ij}\rest_{\PP_{ij}}$ to the lattice points:
\[
c_{ij}^\eps(u)=c_{ij}(u),\qquad u\in\PP_{ij}^\eps,\qquad
1\le i\neq j\le m.
\]
An even better option is to read off the values of
$c_{ij}\rest_{\PP_{ij}}$ at the middlepoints of the corresponding
plaquettes of $\PP_{ij}^\eps$:
\[
c_{ij}^\eps(u)=
c_{ij}(u+\et e_i+\et e_j),\qquad u\in\PP_{ij}^\eps,\qquad
1\le i\neq j\le m.
\]
Either choice gives the data (Q$_{1,2}^\Delta$) which define an
$\eps$-dependent family of discrete Q-nets
$f^\eps\!:(\eps\bbZ)^m\to\bbR^N$, called {\em canonical discrete
Q-nets} corresponding to the initial data (Q$_{1,2}$).
\begin{thm} \label{thm:cnapprox}
For some $r>0$, the canonical discrete Q-nets
$f^\eps\!:\PP^\eps(r)\to\bbR^N$ converge, as $\eps\to 0$, to the
unique conjugate net $f:\PP(r)\to\bbR^N$ with the initial data
{\em (Q$_{1,2}$)}. Convergence is with the order $\OH(\eps)$ in
$C^\infty(\PP(r))$.
\end{thm}
This follows directly from Theorem \ref{thm:approx}, since eqs.
(\ref{eq:dcn property}), (\ref{eq:dcn c}) and (\ref{eq:cn property}),
(\ref{eq:cn c}) are manifestly hyperbolic (and can be easily re-written
in the first order form).
\vspace{0pt}

{\bf Discretization of a Jonas pair.} Recall that a Jonas
transform of a given conjugate net is determined by the initial
data (J$_{1,2}$) (see Sect. \ref{Sect: conj}). We now produce out
of these the initial data (J$_{1,2}^\Delta$) (see Sect. \ref{Sect:
subsect discr conj}) for an $\eps$-dependent family of Jonas
transforms of canonical discrete Q-nets corresponding to the
initial data (Q$_{1,2}$).

Take the point $f^+(0)$ from (J$_1$). Define the edge functions
$a_i^\eps\rest_{\PP_i^\eps}$, $b_i^\eps\rest_{\PP_i^\eps}$ by
restricting the functions $a_i\rest_{\PP_i}$, $b_i\rest_{\PP_i}$
to the lattice points, or, better, to the middlepoints of the
corresponding edges of $\PP_i^\eps$. This gives the data set
(J$_{1,2}^\Delta$); along with the data (Q$_{1,2}^\Delta$)
produced above this yields in a canonical way an $\eps$-dependent
family of discrete Q-nets
$F^\eps:(\eps\bbZ)^m\times\{0,1\}\to\bbR^N$, which will be called
the canonical ones for the initial data (Q$_{1,2}$), (J$_{1,2}$).

\begin{thm}  \label{thm:jonapprox}
 The canonical Q-nets $(f^\eps)^+=F^\eps(\cdot,1):
 \PP^\eps(r)\to\bbR^N$ converge to the unique Jonas transform
 $f^+:\PP(r)\to\bbR^N$ of $f$ with the initial data {\em (J$_{1,2}$)}.
 Convergence is with the order $\OH(\eps)$ in $C^\infty(\PP(r))$.
\end{thm}
Again, this follows directly from Theorem \ref{thm:approx},
applied to the hyperbolic systems consisting of eqs. (\ref{eq: djp
a})--(\ref{eq: djp c}) in the discrete case and of eqs.
(\ref{eq:jp dA})--(\ref{eq:jp C+}) in the smooth case. Note that
the discrete equations are implicit, and their solvability for
$\eps$ small enough is guaranteed on the set $\{a_j\neq 0: 1\le
j\le m\}$, which is open dense in the phase space.

\section{Orthogonal nets}
 \label{subsect: os approx}

We start with approximation of a single orthogonal net
$f:\bbR^m\to\bbR^N$. For the approximating discrete O-nets, we
have $M=m$ and all $\eps_i=\eps$. In all formulas of Sect.
\ref{Subsect: discr ortho} one has to replace the lattice
functions $h_i$, $\beta_{ij}$, $\rho_{kj}$ by $\eps h_i$,
$\eps\beta_{ij}$, $\eps\rho_{kj}$, respectively. Observe that
formulas (\ref{eq: dos n}), (\ref{eq: N for dos}) become
\begin{eqnarray*}
 & \nu_{ji}=\nu_{ij}=(1-\eps^2\beta_{ij}\beta_{ji})^{1/2}=1+\OH(\eps^2),&
 \\
 & \sigma_i=\Big(1-\displaystyle\frac{\eps^2}{4}
 \sum_{k\neq i}\rho_{ki}^2\Big)^{1/2}=1+\OH(\eps^2). &
\end{eqnarray*}
Under this re-scaling, eqs. (\ref{dos: evolve x}), (\ref{dos:
evolve v}), (\ref{dos: evolve h}) and (\ref{dos: evolve beta}) can
be put into the standard form (\ref{eq:prehyperbolic}) with the
functions on the right-hand sides approximating, as $\eps\to 0$,
the corresponding functions in eqs. (\ref{eq:DX})--(\ref{eq:DB})
with the order $\OH(\eps)$.

Nevertheless, formally speaking, Theorem \ref{thm:approx} cannot
be applied to orthogonal nets. The reason for this is that the
full system of differential equations describing orthogonal nets,
consisting of eqs. (\ref{eq:DX})--(\ref{eq:DB}) {\em and} the
constraint (\ref{eq:ddL}), is non-hyperbolic. Its
non-hyperbolicity rests on the fact that the constraint
(\ref{eq:ddL}) is not resolved with respect to the derivatives
$\partial_i\beta_{ij}$. Note, however, that constraint
(\ref{eq:ddL}) does not take part in the evolution of solutions
starting with the data given in the coordinate planes $\PP_{ij}$:
it is satisfied automatically, provided it is fulfilled for the
coordinate surfaces $f\rest_{\PP_{ij}}$. Therefore, we will obtain
a convergence result for orthogonal nets as soon as it will be
established for coordinate surfaces.
\medskip

{\bf Discretization of an O-surface.} Initial data for a smooth O-surface
$f:\PP_{12}\to\bbR^N$ are:
\begin{enumerate}
\item[(i)] two smooth curves $f\rest_{\PP_i}$ $(i=1,2)$,
intersecting orthogonally at $f(0)$;
\item[(ii)] a smooth function $\gamma_{12}:\PP_{12}\to\bbR$, whose
designated meaning is
$\gamma_{12}=\frac{1}{2}(\partial_1\beta_{12}-\partial_2\beta_{21})$.
\end{enumerate}
Let $\hat{f}\rest_{\PP_i}$ be the images of the curves $f\rest_{\PP_i}$
in the M\"obius-geometric model $\bbQ_{\,0}^N$. Let
$h_i=|\partial_i\hat{f}|$ and $\hat{v}_i=h_i^{-1}\partial_i\hat{f}$
be the metric coefficients and unit tangent vectors
of the coordinate curves.
Choose an initial frame $\psi(0)\in\cH_\infty$ such that
\[
\hat{f}(0)=\psi^{-1}(0)\ee_0\psi(0),\quad
\hat{v}_i(0)=\psi^{-1}(0)\ee_i\psi(0)\quad (i=1,2).
\]
Define the frames $\psi:\PP_i\to\cH_\infty$ of the curves
$\hat{f}\rest_{\PP_i}$ as the solutions of eqs. (\ref{eq: frame
os}) for $i=1,2$ (considered as ordinary differential equations)
with the initial value $\psi(0)$. {\em Rotation coefficients of
the curves $\hat{f}\rest_{\PP_i}$} are the functions
$\beta_{ki}:\PP_i\to\bbR$ defined by the formula (\ref{eq: os S})
for $i=1,2$.

Define the discrete coordinate curves
$\hat{f}^\eps\rest_{\PP_i^\eps}$ by restricting the functions
$\hat{f}\rest_{\PP_i}$ to the lattice points. Let
$h_i^\eps=|\delta_i\hat{f}^\eps|$ and
$\hat{v}_i^\eps=(h_i^\eps)^{-1}\delta_i\hat{f}^\eps$ be the
discrete metric coefficients and unit vectors along the discrete
curves. Define the frame $\psi^\eps:\PP_i^\eps\to\cH_\infty$ by
iterating the difference equation (\ref{eq: frame dos}) for
$i=1,2$ with the initial condition $\psi^\eps(0)=\psi(0)$. Then
{\em canonical rotation coefficients of the discrete curves}
$\hat{f}^\eps\rest_{\PP_i^\eps}$ are the coefficients
$\rho^\eps_{ki}:\PP_i^\eps\to\bbR$ in the expansions
\[
 V_i^\eps=\psi^\eps\hat{v}_i^\eps(\psi^\eps)^{-1}
 =\sigma_i^\eps\ee_i-\frac{\eps}{2}\,
 \sum_{k\neq i}\rho_{ki}^\eps\ee_k+\eps h_i^\eps\ee_\infty.
\]
Finally, let the plaquette function
$\gamma_{12}^\eps:\PP_{12}^\eps\to\bbR$ be obtained by restricting
$\gamma_{12}$ to the lattice points (or to the middlepoints of the
corresponding plaquettes of $\PP_{12}^\eps$).

Thus, we get valid Goursat data for a hyperbolic system of
first-order difference equations for the variables $\hat{f}^\eps$,
$\hat{v}_i^\eps$, $h_i^\eps$, $\rho_{ki}^\eps$, consisting of eqs.
(\ref{dos: evolve x}), (\ref{dos: evolve v}), (\ref{dos: evolve h}),
(\ref{eq:dB}) and (\ref{eq:dL}) with distinct $i,j\in\{1,2\}$ and
$1\le k\le N$, where the following expressions should be inserted:
\[
\beta_{12}^\eps=\sigma_1^\eps\rho_{12}^\eps-\frac{\eps}{2}\Big(\frac{1}{2}
\sum_{k>2}\rho_{k1}^\eps\rho_{k2}^\eps-\gamma_{12}^\eps\Big),\quad
\beta_{21}^\eps=\sigma_2^\eps\rho_{21}^\eps-\frac{\eps}{2}\Big(\frac{1}{2}
\sum_{k>2}\rho_{k1}^\eps\rho_{k2}^\eps+\gamma_{12}^\eps\Big).
\]
The discrete nets $\hat{f}^\eps:\PP_{12}^\eps\to\bbQ_{\,0}^N$
defined as solutions of the Goursat problem just described are
discrete O-surfaces, since they fulfill the circularity constraint
(\ref{eq:dQ}). They will be called {\em canonical discrete
O-surfaces} constructed from the above initial data.

\begin{thm}  \label{thm: approx C-surf}
  For some fixed $r>0$, the canonical discrete O-surfaces
  $\hat{f}^\eps:\PP_{12}^\eps(r)\to\bbQ_{\,0}^N$ converge, with
  the order $\OH(\eps)$ in $C^\infty(\PP_{12}(r))$,
  to the unique O-surface $\hat{f}:\PP_{12}(r)\to\bbQ_{\,0}^N$ with
  the initial data $f\rest_{\PP_i}\,$ ($i=1,2$) and
  $\frac{1}{2}(\partial_1\beta_{12}-\partial_2\beta_{21})=\gamma_{12}$.
  Edge rotation coefficients $\rho_{ki}^\eps$ and plaquette rotation
  coefficients $\beta_{12}^\eps$, $\beta_{21}^\eps$ of the discrete
  O-surfaces $\hat{f}^\eps$ converge to the corresponding
  rotation coefficients $\beta_{ki}$ of the O-surface $\hat{f}$.
\end{thm}
{\bf Proof.} First, we show the convergence of the frames,
$\psi^\eps\to\psi$, and of the rotation coefficients,
$\rho_{ki}^\eps\to\beta_{ki}\,$, along the discrete curves
$\hat{f}^\eps\rest_{\PP_i^\eps}$. This follows from two
observations. First,
$\hat{v}_i^\eps(0)=\hat{v}_i(0)+\frac{\eps}{2}(\partial_i\hat{v}_i)(0)+
\OH(\eps^2)$, so that there holds:
\[
 (\tau_i-1)\psi^\eps(0)=
 -\frac{\eps}{2}\,\ee_i\psi(0)(\partial_i\hat{v}_i)(0)+\OH(\eps^2).
\]
Second, combining frame equations on two neighboring edges of
$\PP_i^\eps$, one finds that everywhere on $\PP_i^\eps$ there
holds:
\[
 (\tau_i-\tau_i^{-1})\psi^\eps=
 -\ee_i\psi^\eps(1-\tau_i^{-1})\hat{v}_i^\eps=
 -\eps\ee_i\psi^\eps(\partial_i\hat{v}_i)+\OH(\eps^2).
\]
The claim follows by standard methods of the ODE theory.

Now an application of Theorem \ref{thm:approx} shows that
functions $\hat{f}^\eps:\PP_{12}^\eps\to\bbQ_{\,0}^N$ converge to
the functions $\hat{f}:\PP_{12}\to\bbQ_{\,0}^N$ which solve the
Goursat problem for the hyperbolic system of first order
differential equations, consisting of eqs.
(\ref{eq:DX})--(\ref{eq:DB}) with distinct $i,j\in\{1,2\}$ and
$1\le k\le N$, and
\[
\partial_1\beta_{12}=-\frac{1}{2}\sum_{k>2}\beta_{k1}\beta_{k2}+
\gamma_{12},\quad
\partial_2\beta_{21}=-\frac{1}{2}\sum_{k>2}\beta_{k1}\beta_{k2}-
\gamma_{12}.
\]
The solutions $\beta_{ki}$ satisfy the orthogonality constraint
(\ref{eq:DL}) and the relation
$\frac{1}{2}(\partial_1\beta_{12}-\partial_2\beta_{21})=\gamma_{12}$. $\Box$
\medskip

{\bf Discretization of an $m$-dimensional orthogonal net.}
Given the initial data (O$_{1,2}$) for an $m$-dimensional orthogonal net
(see Sect. \ref{Sect: subsect os}), we can apply the procedure described in
the previous paragraph, with an initial frame $\psi(0)\in\cH_\infty$ such that
\[
\hat{f}(0)=\psi^{-1}(0)\ee_0\psi(0),\quad
\hat{v}_i(0)=\psi^{-1}(0)\ee_i\psi(0)\quad (1\le i\le m),
\]
to produce, in a canonical way, the discrete O-surfaces
$\hat{f}^\eps\rest_{\PP_{ij}^\eps}$ and their plaquette rotation
coefficients $\beta_{ij}^\eps$. Thus, we get data
(O$_{1,2}^\Delta$) (see Sect. \ref{Subsect: discr ortho}) for an
$\eps$-dependent family of discrete O-nets
$\hat{f}^\eps:(\eps\bbZ)^m\to\bbQ_{\,0}^N$. These nets will be
called the {\em canonical discrete O-nets} corresponding to the
initial data (O$_{1,2}$).

\begin{thm}  \label{thm:cDupin}
  The canonical discrete O-nets $\hat{f}^\eps:\PP^\eps(r)\to\bbR^N$
  converge, as $\eps\to 0$, to the unique orthogonal net
  $\hat{f}:\PP(r)\to\bbR^N$ with the initial data {\em (O$_{1,2}$)}.
  Convergence is with the order $\OH(\eps)$ in $C^\infty(\PP(r))$.
\end{thm}
{\bf Proof.} The data (O$_{1,2}^\Delta$) yield a well-posed Goursat problem
for the hyperbolic system of first-order difference equations for the variables
$\hat{f}^\eps$, $\hat{v}_i^\eps$, $h_i^\eps$, $\beta_{ij}^\eps$, consisting
of eqs. (\ref{dos: evolve x}), (\ref{dos: evolve v}), (\ref{dos: evolve h}),
(\ref{dos: evolve beta}). The convergence of these Goursat data is assured
by Theorem \ref{thm: approx C-surf}. Now the claim of the theorem follows
directly from Theorem \ref{thm:approx}. $\Box$
\medskip

{\bf Discretization of a Ribaucour pair.} Given the initial data
(R$_{1,2}$) for a Ribaucour transform of an orthogonal net (see
Sect. \ref{Sect: subsect os}), define the plaquette rotation
coefficients $\beta_{Mi}^\eps$ on the ``vertical'' plaquettes
along the edges of the coordinate axes $\PP_i^\eps$ by restricting
the corresponding functions $\eps\theta_i$ to lattice points or,
alternatively, to middlepoints of the corresponding edges of
$\PP_i^\eps$:
\[
\beta_{Mi}^\eps(u)=\eps\theta_i(u)\;\;{\rm or}\;\;
\eps\theta_i(u+\et),\qquad u\in\PP_i^\eps,
\qquad 1\le i\le m.
\]
Thus, we get the data (R$_{1,2}^\Delta$) (see Sect. \ref{Subsect:
discr ortho}), which, together with (O$_{1,2}^\Delta$), allow us
to construct in a canonical way discrete O-nets
$F^\eps:(\eps\bbZ)^m\times\{0,1\}\to\bbR^N$. They will be called
the canonical ones corresponding to the initial data (O$_{1,2}$),
(R$_{1,2}$).
\begin{thm}  \label{thm:ribapprox}
  The canonical discrete O-nets $(f^\eps)^+=F^\eps(\cdot,1):
 \PP^\eps(r)\to\bbR^N$ converge to the unique Ribaucour transform
 $f^+:\PP(r)\to\bbR^N$ of $f$ with the initial data {\em(R$_{1,2}$)}.
 Convergence is with the order $\OH(\eps)$ in $C^\infty(\PP(r))$.
\end{thm}
{\bf Proof.} Define $v_M^\eps(0)$ as the unit vector parallel to
$\delta f(0)=f^+(0)-f(0)$, and set $h_M^\eps(0)=|\delta f(0)|$.
These data along with $\beta_{Mi}^\eps$ on the coordinate axes,
added to the previously found ones $f^\eps(0)$, $v_i^\eps$,
$h_i^\eps$, $\beta_{ij}^\eps$ for $1\le i,j\le m$, form valid
Goursat data for the system (\ref{dos: evolve x}), (\ref{dos:
evolve v}), (\ref{dos: evolve h}), (\ref{dos: evolve beta}).

The circularity constraint (\ref{eq: dos beta}) implies that
$\beta_{iM}^\eps=-2\langle v_i^\eps,v_M^\eps\rangle-\eps\theta_i$
on all edges of $\PP_i^\eps$.
Perform the substitution
\[
v_M^\eps=y+\OH(\eps),\quad \!h_M^\eps=\ell+\OH(\eps),\quad
\beta_{Mi}^\eps=\eps\theta_i+\OH(\eps^2),\quad \!
\beta_{iM}^\eps=-2\langle v_i,y\rangle+\OH(\eps)
\]
in eqs. (\ref{dos: evolve v}), (\ref{dos: evolve h}),
(\ref{dos: evolve beta}) with one of the indices equal to $M$.
Taking into account that in this limit one has
\[
\nu_{iM}^{-1}=\nu_{Mi}^{-1}=1-\eps\langle v_i,y\rangle\theta_i+\OH(\eps^2),
\]
one sees that the limiting equations coincide with eqs. (\ref{eq:rp v}),
(\ref{eq:rp h}), (\ref{eq:rp beta}).
A reference to Theorem \ref{thm:approx} finishes the proof. $\Box$

\section{Moutard nets}
 \label{Subsect: mn conv}

{\bf Discretization of an M-net.} Given the initial data
(M$_{1,2}$) for an M-net (see Sect. \ref{Subsect: mn}), we produce
initial data (M$_{1,2}^\Delta$) for an $\eps$-dependent family of
discrete M-nets with $\eps_1=\eps_2=\eps$. Discrete curves
$f^\eps\rest_{\PP_i^\eps}$ are obtained from the smooth curves
$f\rest_{\PP_i}$ by restricting to the lattice points:
\[
f^\eps(u)=f(u),\qquad u\in\PP_i^\eps,\quad i=1,2.
\]
The plaquette function $a_{12}^\eps:(\eps\bbZ)^2\to\bbR$
is obtained from the function $q_{12}$ restricted to the lattice points:
\[
a_{12}^\eps(u)=1+{\textstyle\frac{1}{2}}\eps^2 q_{12}(u),
\qquad u\in(\eps\bbZ)^2
\]
(one could also restrict $q_{12}$ to middlepoints of the
corresponding plaquettes). Now discrete M-nets
$f^\eps:(\eps\bbZ)^2\to\bbR^N$ are defined as solutions of the
difference equation (\ref{eq:dMou 2d}) with the above data
(M$_{1,2}^\Delta$).

\begin{thm} \label{thm:mnapprox}
Canonical discrete M-nets $f^\eps\!:\PP^\eps(r)\to\bbR^N$
converge, as $\eps\to 0$, to the unique M-net $f:\PP(r)\to\bbR^N$
with the initial data {\em (M$_{1,2}$)}. Convergence is with the
order $\OH(\eps)$ in $C^\infty(\PP(r))$.
\end{thm}
{\bf Proof.} Eq. (\ref{eq:dMou 2d}) is manifestly hyperbolic (and can be
easily put in the first order form). It approximates eq. (\ref{eq:Mou}),
because it can be re-written as
\[
\delta_1\delta_2 f={\textstyle\frac{1}{2}}q_{12}(\tau_1 f+\tau_2 f)=
q_{12}\big(f+\et\delta_1 f+\et\delta_2 f\big).
\]
Now Theorem \ref{thm:approx} can be applied. $\Box$
\smallskip

{\bf Discretization of a Moutard pair.}
Let the initial data (MT$_{1,2}$) for a Moutard transformation be given.
Define the edge variables $b_i^\eps\rest_{\PP_i^\eps}$  from the
functions $p_i\rest_{\PP_i}$ restricted to the lattice points:
\[
b_1^\eps(u_1,0)=1+\eps p_1(u_1,0),\qquad
b_2^\eps(0,u_2)=1+\eps p_2(0,u_2),\qquad u_i\in\eps\bbZ
\]
(one could restrict $p_i\rest_{\PP_i}$ to the middlepoints of the
corresponding edges, as well). This gives us the data
(MT$_{1,2}^\Delta$), which canonically generate discrete M-nets
$(f^\eps)^+:(\eps\bbZ)^2\to\bbR^N$.
\begin{thm}  \label{thm:mtapprox}
 Canonical discrete M-nets $(f^\eps)^+:\PP^\eps(r)\to\bbR^N$
 converge to the unique Moutard transform
 $f^+:\PP(r)\to\bbR^N$ of $f$ with the initial data {\em (MT$_{1,2}$)}.
 Convergence is with the order $\OH(\eps)$ in $C^\infty(\PP(r))$.
\end{thm}
{\bf Proof.} The system consisting of eqs. (\ref{eq:dmp}),
(\ref{eq:dMou 1}) is hyperbolic. Upon substituting
$b_i=1+\eps p_i$ and $a_{12}=1+\frac{1}{2}\eps^2q_{12}\,$,
these equations can be re-written as
\[
\delta_1 f^++\delta_1 f=p_1(f^+-\tau_1 f),\qquad
\delta_2 f^+-\delta_2 f=p_2(f^++\tau_2 f),
\]
and
\[
\frac{1+\eps\tau_2p_1}{1+\eps p_1}=\frac{1+\eps\tau_1p_2}{1+\eps p_2}=
\frac{1+(\eps^2/2)\,q_{12}^+}{1+(\eps^2/2)\,q_{12}}=
\frac{1}{1+\eps^2(q_{12}-p_1p_2)+\OH(\eps^3)}.
\]
Clearly, they aproximate, as $\eps\to 0$, eqs.
(\ref{eq:mp 1})--(\ref{eq:mp 2}) and (\ref{eq:Mou 1})--(\ref{eq:Mou 2}),
respectively. It remains to apply Theorem \ref{thm:approx}. $\Box$

\section{A-surfaces}
 \label{Subsect: an conv}

{\bf Discretization of an A-surface.} Initial data (A$_{1,2}$)
for an A-surface (see Sect. \ref{Subsect: an}), are nothing but initial
data (M$_{1,2}$) for the Lelieuvre normal field $n:\bbR^2\to\bbR^3$.
In Theorem \ref{thm:mnapprox}, we described the canonical construction
of the initial data (A$_{1,2}^\Delta$), which give a converging family
of the discrete Lelieuvre normal fields $n^\eps:(\eps\bbZ)^2\to\bbR^3$.
Eqs. (\ref{eq:dLel 2d}) define the discrete A-surfaces
$f^\eps:(\eps\bbZ)^2\to\bbR^3$, called the {\em canonical discrete A-surfaces}
corresponding to the initial data (A$_{1,2}$).

\begin{thm} \label{thm:anapprox}
Canonical discrete A-surfaces $f^\eps\!:\PP^\eps(r)\to\bbR^3$ converge,
as $\eps\to 0$, to the unique A-surface $f:\PP(r)\to\bbR^3$ with the initial
data {\em (A$_{1,2}$)}. Convergence is with the order $\OH(\eps)$ in
$C^\infty(\PP(r))$.
\end{thm}
{\bf Proof.} Eqs. (\ref{eq:dLel 2d}) are hyperbolic and approximate eqs.
(\ref{eq:Lel}). Theorem \ref{thm:approx} can be applied to prove the
convergence of $f^\eps$, after the convergence of $n^\eps$ has been already
proved. $\Box$
\smallskip

{\bf Discretization of a Weingarten pair.}
Initial data (W$_{1,2}$) for a Weingarten transformation are nothing
but initial data (MT$_{1,2}$) for a Moutard transformation of the Lelieuvre
normal field. We already described, in Theorem \ref{thm:mtapprox}, the
canonical construction of the initial data (W$_{1,2}^\Delta$) for
a converging family of the discrete Lelieuvre normal fields
$(n^\eps)^+:(\eps\bbZ)^2\to\bbR^3$. Now the transformed A-surfaces
$(f^\eps)^+:(\eps\bbZ)^2\to\bbR^3$ are obtained by eq. (\ref{eq:dWei}).
\begin{thm}  \label{thm:weiapprox}
  Canonical discrete A-nets $(f^\eps)^+\!:\PP^\eps(r)\to\bbR^3$ converge
  to the unique Weingarten transform $f^+:\PP(r)\to\bbR^3$ of $f$ with
  the initial data {\em (W$_{1,2}$)}. Convergence is with the order
  $\OH(\eps)$ in $C^\infty(\PP(r))$.
\end{thm}
{\bf Proof} follows by comparing the (identical) formulas
(\ref{eq:Wei}) and (\ref{eq:dWei}), after the convergence of
$n^\eps$ and $(n^\eps)^+$ has been proved. $\Box$

\section{K-surfaces}
 \label{Subsect: ks conv}

{\bf Discretization of a K-surface.} Given the initial data (K)
for a K-surface (see Sect. \ref{Subsect: Ksurf}), we define
initial data (K$^\Delta$) (see Sect. \ref{Subsect: discr K}) for
an $\eps$-dependent family of discrete K-surfaces with
$\eps_1=\eps_2=\eps$ by restricting $n\rest_{\PP_i}$ to the
lattice points, as for general A-surfaces. Define discrete M-nets
$n^\eps:(\eps\bbZ)^2\to\bbS^2$ as solutions of the difference
equations (\ref{eq:klat 2d}) with the initial data (K$^\Delta$).
Finally, define the discrete K-surfaces
$f^\eps:(\eps\bbZ)^2\to\bbR^3$ with the help of the discrete
Lelievre representation (\ref{eq:dLel 2d}). These will be called
the {\em canonical discrete K-surfaces} corresponding to the
initial data (K).

\begin{thm} \label{thm:ksapprox}
Canonical discrete K-surfaces $f^\eps\!:\PP^\eps(r)\to\bbR^3$ converge,
as $\eps\to 0$, to the unique K-surface $f:\PP(r)\to\bbR^3$ with the initial
data {\em (K)}. Convergence is with the order $\OH(\eps)$ in
$C^\infty(\PP(r))$.
\end{thm}
{\bf Proof.} We have for $n=n^\eps$:
\[
a_{12}^\eps=\frac{\langle n,\tau_1 n+\tau_2 n\rangle}
  {1+\langle \tau_1 n,\tau_2 n\rangle}
 =\frac{2+\eps\langle n,\delta_1 n+\delta_2 n\rangle}
  {2+\eps\langle n,\delta_1 n+\delta_2 n\rangle+\eps^2
  \langle \delta_1 n,\delta_2 n\rangle}.
\]
Since $\langle n,\delta_i n\rangle=\OH(\eps)$, we find that
\[
 a_{12}^\eps=1-\textstyle\frac{1}{2}\eps^2
 \langle \delta_1 n,\delta_2 n\rangle +\OH(\eps^4).
\]
Comparing this with eq. (\ref{eq: Mou n q}), we see that Theorem
\ref{thm:approx} can be applied in order to prove approximation of
$n:\bbR^2\to\bbS^2$ by $n^\eps:(\eps\bbZ)^2\to\bbS^2$. Finally,
approximation of $f$ by $f^\eps$ follows exactly as for general A-surfaces.
$\Box$
\smallskip

{\bf Discretization of a B\"acklund pair.} Let the initial data
(B) for a B\"acklund transformation of a given K-surface $f$, i.e.
the point $n^+(0)$, be given. Take it as the initial data
(B$^\Delta$) for the discrete B\"acklund transformations
$(f^\eps)^+:(\eps\bbZ)^2\to\bbR^3$ of the family $f^\eps$ of
discrete K-surfaces constructed in Theorem \ref{thm:ksapprox}.
\begin{thm}  \label{thm:backapprox}
  Canonical discrete K-surfaces $(f^\eps)^+\!:\PP^\eps(r)\to\bbR^3$
  converge to the unique B\"acklund transform $f^+:\PP(r)\to\bbR^3$ of
  $f$ with the initial data {\em (B)}. Convergence is with the order
  $\OH(\eps)$ in $C^\infty(\PP(r))$.
\end{thm}
{\bf Proof.} In eqs. (\ref{eq: dbt 2d 1}), (\ref{eq: dbt 2d 2}) we have:
\begin{eqnarray*}
b_1=\frac{\langle n,\tau_1 n-n^+\rangle}
{1-\langle\tau_1 n,n^+\rangle} & = &
1+\eps\,\frac{\langle\delta_1 n,n^+\rangle}
{1-\langle n,n^+\rangle}+\OH(\eps^2),\\
b_2=\frac{\langle n,\tau_2 n+n^+\rangle}
{1+\langle\tau_2 n,n^+\rangle} & = &
1-\eps\,\frac{\langle\delta_2 n,n^+\rangle}
{1+\langle n,n^+\rangle}+\OH(\eps^2).
\end{eqnarray*}
Comparing this with eqs. (\ref{eq:wp 1})--(\ref{eq:wp 2}) and
applying Theorem \ref{thm:approx}, we prove approximation of the
Gauss maps. $\Box$

\chapter{Consistency as integrability} \label{Sect: consist}

\section[From 3D consistency to BT and ZCR]
{From 3D consistency of discrete 2D systems to B\"acklund
transformations and zero curvature representations}
\label{Sect: 3d consist}

Equations of a discrete 2D hyperbolic system are associated to
elementary squares of $\bbZ^2$. Such a system may contain fields
on vertices and/or on edges of $\bbZ^2$. In our
differential-geometric considerations we encountered two types of
such systems: the cross-ratio system, with fields on the vertices
and the parameters on the edges of $\bbZ^2$, and discrete
Lorentz-harmonic maps into quadrics, with fields on the vertices
of $\bbZ^2$. One can imagine further types of discrete hyperbolic
2D systems, for instance, those with fields on edges only. We will
mainly consider systems of the cross-ratio type, and will discuss
other types briefly.
\smallskip

\begin{figure}[htbp]
\begin{center}
\setlength{\unitlength}{0.06em}
\begin{picture}(200,140)(-50,-20)
  \put(100,  0){\circle{10}} \put(0  ,100){\circle{10}}
  \put(  0,  0){\circle*{10}}  \put(100,100){\circle*{10}}
  \thicklines
  \path( 5,  0)(95,0)
  \path( 5,100)(95,100)
  \path(  0, 5)(0,95)
  \path(100, 5)(100,95)
  \put(-10,-20){$f$}
  \put(99,-20){$f_1$}
  \put(99,113){$f_{12}$}
  \put(-10,113){$f_2$}
  \put(43,-13){$\alpha_1$}
  \put(43,105){$\alpha_1$}
  \put(-20,47){$\alpha_2$}
  \put(105,47){$\alpha_2$}
\end{picture}
\caption{Elementary quadrilateral; fields on vertices}
\label{Fig:quadrilateral}
\end{center}
\end{figure}
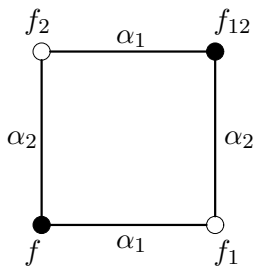

{\bf 2D systems with fields on vertices and with labelled edges.}
A typical representative of this class of equations is the
cross-ratio system:
\begin{equation}\label{eq:2d cr}
q(f,f_1,f_{12},f_2)=\frac{\alpha_1}{\alpha_2}\,.
\end{equation}
The variables $f:\bbZ^2\to\cC\ell(\bbR^N)$ are associated to the
vertices of $\bbZ^2$, while the variables $\alpha_i:\bbZ^2\to\bbR$
satisfy the labelling property
\begin{equation}\label{eq:2d lab}
\delta_2\alpha_1=0,\qquad \delta_1\alpha_2=0.
\end{equation}
Therefore $\alpha_i$ play the role of parameters naturally
associated to edges of $\bbZ^2$ parallel to $\PP_i$ and constant
along the strips in the complementary direction. (As a disclaimer,
we mention that further edge variables have to be introduced in
order to put this equation in the first-order form.) One more
example of such a system with vertex variables and parameters
sitting on edges and having the labelling property is given by the
{\em Hirota equation},
\begin{equation}\label{eq:2d H}
\frac{f_{12}}{f}=\frac{\alpha_1f_1-\alpha_2f_2}{\alpha_2f_1-\alpha_1f_2},
\end{equation}
see \cite{BP1, BMS1} for the geometric interpretation of this
system in terms of discrete K-surfaces. A general system of this
class consists of equations
\begin{equation}\label{eq:2d Q}
Q(f,f_1,f_{12},f_2;\alpha_1,\alpha_2)=0,
\end{equation}
see Fig.~\ref{Fig:quadrilateral}. In the simplest situation one
has complex fields $f:\bbZ^2\to\bbC$, and complex parameters
$\alpha_i$ on the edges of $\bbZ^2$ parallel to $\PP_i$,
satisfying the labelling condition (\ref{eq:2d lab}). We require
that eq. (\ref{eq:2d Q}) be uniquely solvable for any one of its
arguments $f,f_1,f_2,f_{12}\in\widehat{\mathbb{C}}$. Therefore,
the solutions have to be fractional-linear in each of their
arguments. This naturally leads to the following assumption.
\smallskip

{\it Linearity.} The function $Q(x,u,y,v;\alpha,\beta)$ is a polynomial
of degree 1 in each argument (affine linear):
\begin{equation*}\label{eq:Qlin}
  Q(x,u,y,v;\alpha,\beta)= a_1(\alpha,\beta)xuyv+\dots
  +a_{16}(\alpha,\beta).
\end{equation*}

Note that for the cross-ratio equation (\ref{eq:2d cr}) with
complex-valued arguments, $Q=\beta(x-u)(y-v)-\alpha(u-y)(v-x)$,
while for the Hirota equation (\ref{eq:2d H}),
$Q=\alpha(xu+yv)-\beta(xv+yu)$.

Actually, this setup admits an important generalization:
elementary quadrilaterals carrying eqs. (\ref{eq:2d Q}) can be
attached one to another with the combinatorics more complicated
then that of $\bbZ^2$.
\begin{dfn}
A {\,\rm quad-graph} $\cD$ is a strongly regular polytopal cell
decomposition of a surface with all quadrilateral faces.
\end{dfn}
We denote by $V(\cD)$, $E(\cD)$, $F(\cD)$ the sets of vertices, of
edges and of faces of $\cD$, respectively. We consider eq.
(\ref{eq:2d Q}) for fields $f:V(\cD)\to\bbC$, with
$\alpha:E(\cD)\to\bbC$ being a labelling of edges of $\cD$, i.e.,
a function taking equal values on any pair of opposite edges of
any quadrilateral from $F(\cD)$. In the context of equations on
quad-graphs, there are no distinguished coordinate directions,
nevertheless it will be convenient to continue to use notations of
eq. (\ref{eq:2d Q}), with the understanding that indices are used
locally (within one quadrilateral), and do not stand for shifts
into the globally defined coordinate directions. So,
$f,f_1,f_{12},f_2$ can be {\em any} cyclic enumeration of the
vertices of an elementary quadrilateral. Eq. (\ref{eq:2d Q})
should not depend on the enumeration of vertices, therefore the
following assumption is natural when considering equations on
general quad-graphs.
\smallskip

{\it Symmetry.} The function $Q$ has the symmetry properties
\begin{equation*}\label{eq:Qsym}
 Q(x,u,y,v;\alpha,\beta)=\eps Q(x,v,y,u;\beta,\alpha)=
 \sigma Q(u,y,v,x;\beta,\alpha),\quad \eps,\sigma=\pm 1.
\end{equation*}
(For the complex cross-ratio equation and for the Hirota equation
the functions $Q$ given above possess this property with
$\eps=\sigma=-1$).
\smallskip

\begin{figure}[htbp]
\begin{center}
\setlength{\unitlength}{0.05em}
\begin{picture}(200,220)(0,0)
 \put(0,0){\circle*{10}}    \put(150,0){\circle{10}}
 \put(0,150){\circle{10}}   \put(150,150){\circle*{10}}
 \put(50,200){\circle*{10}} \put(200,200){\circle{10}}
 \put(50,50){\circle{10}}   \put(200,50){\circle*{10}}
 \path(0,0)(145,0)       \path(0,0)(0,145)
 \path(150,5)(150,150)   \path(5,150)(150,150)
 \path(3.53,153.53)(50,200)    \path(150,150)(196.47,196.47)
 \path(50,200)(195,200)
 \path(200,195)(200,50) \path(200,50)(153.53,3.53)
 \dashline[+30]{10}(0,0)(46.47,46.47)
 \dashline[+30]{10}(50,55)(50,200)
 \dashline[+30]{10}(55,50)(200,50)
 \put(-30,-5){$f$}
 \put(-35,145){$f_3$} \put(215,45){$f_{12}$}
 \put(165,-5){$f_1$} \put(160,140){$f_{13}$}
 \put(15,50){$f_2$}  \put(10,205){$f_{23}$}
 \put(215,205){$f_{123}$}
 \put(80,7){$\alpha_1$} \put(33,15){$\alpha_2$}
 \put(110,57){$\alpha_1$} \put(180,15){$\alpha_2$}
 \put(0,177){$\alpha_2$} \put(75,135){$\alpha_1$}
 \put(140,170){$\alpha_2$} \put(110,205){$\alpha_1$}
 \put(-28,70){$\alpha_3$}\put(155,75){$\alpha_3$}
 \put(22,110){$\alpha_3$}\put(205,115){$\alpha_3$}
\end{picture}
\caption{3D consistency of 2D systems}\label{cube again}
\end{center}
\end{figure}
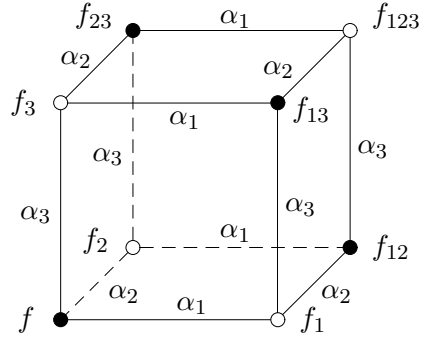

Assume now that eq. (\ref{eq:2d Q}) possesses the property of
{\itbf 3D consistency}. Recall that this means that this equation
can be consistently imposed on all 2D faces of a combinatorial
cube on Fig.~\ref{cube again}, with the parallel edges of the
$i$-th direction all carrying the parameter $\alpha_i$. More
precisely, given the initial data $f$, $f_i$, one determines
$f_{ij}$ in virtue of eq. (\ref{eq:2d Q}), and then three
equations for $f_{123}$ (for three faces adjacent to this vertex)
lead to identical results.
This is the case for the
complex cross-ratio equation, with
\begin{equation}\label{eq:cr 123}
  f_{123}=\frac{(\alpha_1-\alpha_2)f_1f_2+
  (\alpha_2-\alpha_3)f_2f_3+(\alpha_3-\alpha_1)f_3f_1}
  {(\alpha_2-\alpha_1)f_3+(\alpha_3-\alpha_2)f_1+(\alpha_1-\alpha_3)f_2},
\end{equation}
as well as for the Hirota equation, with
\begin{equation}\label{eq:H 123}
  f_{123}=\frac{\alpha_3(\alpha_1^2-\alpha_2^2)f_1f_2+
  \alpha_1(\alpha_2^2-\alpha_3^2)f_2f_3+
  \alpha_2(\alpha_3^2-\alpha_1^2)f_3f_1}
  {\alpha_3(\alpha_2^2-\alpha_1^2)f_3+
  \alpha_1(\alpha_3^2-\alpha_2^2)f_1+
  \alpha_2(\alpha_1^2-\alpha_3^2)f_2}.
\end{equation}
We will demonstrate that this condition automatically leads to two
basic structures associated in the soliton theory with {\itbf
integrability}: B\"acklund transformations and zero curvature
representation.

\begin{figure}[htbp]
\begin{center}
\setlength{\unitlength}{0.05em}
\begin{picture}(200,240)(0,0)
 \put(0,0){\circle*{10}}    \put(150,0){\circle{10}}
 \put(0,150){\circle{10}}   \put(150,150){\circle*{10}}
 \put(50,200){\circle*{10}} \put(200,200){\circle{10}}
 \put(50,50){\circle{10}}   \put(200,50){\circle*{10}}
 \path(0,0)(145,0)       \path(0,0)(0,145)
 \path(150,5)(150,150)   \path(5,150)(150,150)
 \path(3.53,153.53)(50,200)    \path(150,150)(196.47,196.47)
 \path(50,200)(195,200)
 \path(200,195)(200,50) \path(200,50)(153.53,3.53)
 \dashline[+30]{10}(0,0)(46.47,46.47)
 \dashline[+30]{10}(50,55)(50,200)
 \dashline[+30]{10}(55,50)(200,50)
 \put(-25,-5){$f$} \put(-30,145){$f^+$}
 \put(162,-6){$f_1$} \put(160,140){$f_1^+$}
 \put(210,45){$f_{12}$} \put(210,200){$f_{12}^+$}
 \put(20,50){$f_2$}  \put(25,205){$f_2^+$}
 \put(80,7){$\alpha_1$} \put(33,15){$\alpha_2$}
 \put(110,57){$\alpha_1$} \put(180,15){$\alpha_2$}
 \put(0,177){$\alpha_2$} \put(75,135){$\alpha_1$}
 \put(140,170){$\alpha_2$} \put(110,205){$\alpha_1$}
 \put(-15,75){$\lambda$}\put(155,75){$\lambda$}
 \put(35,115){$\lambda$}\put(205,115){$\lambda$}
\end{picture}
\caption{Elementary cube of $\bf D$}\label{back cube}
\end{center}
\end{figure}
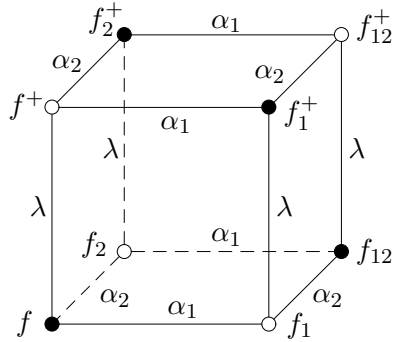

\begin{thm}
For any solution $f:V(\cD)\to\bbC$ of a 3D consistent equation
(\ref{eq:2d Q}) on the quad-graph $\cD$, there is a two-parameter
family of solutions $f^+:V(\cD)\to\bbC$ of the same equation,
satisfying
\begin{equation}\label{eq:2d back}
Q(f,f_i,f_i^+,f^+;\alpha_i,\lambda)=0
\end{equation}
for all edges $(f,f_i)\in E(\cD)$. Such a solution $f^+$ is called a
{\itbf B\"acklund transform} of $f$, and is determined by its value
at one vertex of $\cD$ and by the parameter $\lambda$.
\end{thm}
{\bf Proof.} Extend the planar quad-graph $\cD$ into the third
dimension. Formally speaking, we consider the second copy $\cD^+$
of $\cD$ and add edges connecting each vertex $f\in V(\cD)$ with
its copy $f^+\in V(\cD^+)$. (We slightly abuse the notations here,
by using the same letter $f$ for vertices of the quad-graph and
for the fields assigned to these vertices.) On this way we obtain
a ``3D quad--graph'' ${\bf D}$, with the set of vertices
\[
V({\bf D})=V(\cD)\cup V(\cD^+),
\]
with the set of edges
\[
E({\bf D})=E(\cD)\sqcup E(\cD^+)\sqcup\{(f,f^+):f\in V(\cD)\},
\]
and with the set of faces
\[
F({\bf D})=F(\cD)\sqcup F(\cD^+)\sqcup
\{(f,f_1,f_1^+,f^+):f,f_1\in V(\cD)\}.
\]
Elementary building blocks of $\bf D$ are  combinatorial cubes
$(f,f_1,f_{12},f_2, \linebreak f^+,f_1^+,f_{12}^+,f_2^+)$, as
shown on Fig.~\ref{back cube}. The labelling on $E(\bf D)$ is
defined in the natural way: each edge $(f^+,f_i^+)\in E(\cD^+)$
carries the same label $\alpha_i$ as its counterpart $(f,f_i)\in
E(\cD)$, while all ``vertical'' edges $(f,f^+)$ carry one and the
same label $\lambda$. Clearly, the content of Fig.~\ref{back cube}
is the same as of Fig.~\ref{cube again}, up to notations. Now, a
solution $f^+:V(\cD^+)\to\bbC$ on the first flour of $\bf D$ is
well defined due to the 3D consistency, and is determined by its
value at one vertex of $\cD^+$ and by $\lambda$. We can assume
that $f^+$ is defined on $V(\cD)$ rather than on $V(\cD^+)$, since
these two sets are in a one-to-one correspondence. $\Box$

\begin{thm}\label{Th zcr}
A 3D consistent system (\ref{eq:2d Q}) on the quad-graph $\cD$
admits a {\itbf zero curvature representation} with spectral
parameter dependent $2\times 2$ matrices: there exist matrices
attached to directed edges of $\cD$,
\begin{equation}\label{eq:2d L}
L({\mathfrak e},\alpha({\mathfrak e});\lambda):\vec{E}(\cD)\to
{\rm GL}_2(\bbC)[\lambda],
\end{equation}
such that for any quadrilateral face
$(f,f_1,f_{12},f_2)\in F(\cD)$ there holds
\begin{equation}\label{eq:2d zc}
L(f_{12},f_1,\alpha_2;\lambda)L(f_1,f,\alpha_1;\lambda)=
L(f_{12},f_2,\alpha_1;\lambda)L(f_2,f,\alpha_2;\lambda),
\end{equation}
identically in $\lambda$.
\end{thm}
{\bf Proof.} Due to the linearity assumption, equations
(\ref{eq:2d back}) can be solved for $f_i^+$ in terms of a
M\"obius (fractional-linear) transformation of $f^+$ with
coefficients depending on $f$, $f_i$:
\begin{equation}\label{eq:2d zc Mob}
f_i^+=L(f_i,f,\alpha_i;\lambda)\,[f^+],
\end{equation}
Here we use the standard matrix notation for the action of M\"obius
transformations:
\begin{equation}\label{eq:Mob act}
L[z]=(az+b)(cz+d)^{-1},\quad{\rm where}\quad
L=\left(\begin{array}{cc} a & b \\ c & d \end{array}\right).
\end{equation}
Now 3D consistency for $f_{12}^+$ yields that for any $f^+$ there
holds:
\begin{equation}\label{eq:2d zc weak}
L(f_{12},f_1,\alpha_2;\lambda)L(f_1,f,\alpha_1;\lambda)\,[f^+]=
L(f_{12},f_2,\alpha_1;\lambda)L(f_2,f,\alpha_2;\lambda)\,[f^+].
\end{equation}
Therefore, eq. (\ref{eq:2d zc}) holds at least projectively, i.e.,
up to a scalar factor. A normalization of determinants of $L$ (or
any other suitable normalization) allows one to achieve that eq.
(\ref{eq:2d zc}) holds in the usual sense. $\Box$
\medskip

As an example, we derive a zero curvature representation for the
complex cross-ratio equation (\ref{eq:2d cr}). It will be
convenient to re-define the spectral parameter in this case by
$\lambda\mapsto\lambda^{-1}$, so that equations on the vertical
faces of Fig.~\ref{back cube} read:
\[
\frac{(f_i^+-f^+)(f-f_i)}{(f^+-f)(f_i-f_i^+)}=\lambda\alpha_i.
\]
This gives the M\"obius transformation (\ref{eq:2d zc Mob}) with
\begin{equation}\label{eq:2d cr L1}
L(f_i,f,\alpha_i,\lambda)=I+\frac{\lambda\alpha_i}{f-f_i}
\left(\begin{array}{cc} f_i & -ff_i \\ 1 & -f\end{array}\right)
\end{equation}
The determinant of this matrix is constant (equal to
$1-\lambda\alpha_i$), therefore no further normalization is
required. A more usual form of the transition matrices of the
zero-curvature representation for the complex cross-ratio equation
is obtained by the gauge transformation
\[
L(f_i,f,\alpha_i;\lambda)\mapsto
A^{-1}(f_i)L(f_i,f,\alpha_i;\lambda)A(f), \qquad
A(f)=\left(\begin{array}{cc} 1 & f \\ 0 & 1 \end{array}\right),
\]
which leads to the matrices
\begin{equation}\label{eq:2d cr L2}
L(f_i,f,\alpha_i;\lambda)=\left(\begin{array}{cc} 1 & f-f_i \\
\lambda\alpha_i(f-f_i)^{-1} & 1 \end{array}\right).
\end{equation}
These matrices (\ref{eq:2d cr L2}) are interpreted as matrices of the
M\"obius transformations of the shifted quantities:
\[
f_i^+-f_i=L(f_i,f,\alpha_i;\lambda)\,[f^+-f].
\]

Thus, we have seen that 3D consistency of a 2D equation on a
quad-graph with complex fields at vertices and with labelled edges
implies existence of B\"acklund transformations and of the zero
curvature representation. This is not a pure existence statement
but rather a construction: both attributes can be derived in a
systematic way starting with no more information than the equation
itself, they are in a sense encoded in the equation provided it is
3D consistent.
\smallskip

{\bf Non-commutative equations.} The truth contained in the last
paragraph is by no means restricted to the situation for which it
was demonstrated. For instance, in \cite{BS2} it was extended to
equations with fields on vertices taking values in some
associative but {\em non-commutative} algebra $\cA$ with unit over
a field $\cK$, and with edge labels with values in $\cK$. The
transition matrices of the zero curvature representation are in
this case $2\times 2$ matrices with entries from $\cA$. They act
on $\cA$ according to eq. (\ref{eq:Mob act}), where now the order
of the factors is essential. The examples worked out in \cite{BS2}
include the non-commutative analogs of the Hirota equation
(\ref{eq:2d H}) and of the cross-ratio equation (\ref{eq:2d cr})
(recall that the latter equation with fields in
$\cA=\cC\ell(\bbR^N)$ and with parameters $\alpha_i$ from
$\cK=\bbR$ governs discrete isothermic surfaces in $\bbR^N$). In
particular:
\begin{thm}\label{Thm 3d cr noncommut}
The cross-ratio equation in an associative algebra $\cA$ is 3D
consistent. It possesses a zero curvature representation with the
transition matrices (\ref{eq:2d cr L2}), where the inversion is
treated in $\cA$.
\end{thm}
{\bf Proof.} The transition matrices are obtained from the
equations themselves, in a full analogy with the scalar case.
Verification of the zero curvature representation thus obtained is
a matter of a direct computation. The non-commutative 3D
consistency is a consequence of the zero curvature representation,
see \cite{BS2}. $\Box$

We consider here one more equation of this kind:
\begin{equation}\label{eq: Schief}
(f_{12}-f)(f_2-f_1)=\alpha_2-\alpha_1,
\end{equation}
with the vertex variables $f$ taking values in $\cA$ and with the
edge labels $\alpha$ from $\cK$. Clearly, in the case when
$\cA=\cC\ell(\bbR^N)$, solutions of this equation are special
T-nets in $\bbR^N$:
 \begin{equation}\label{eq: Schief vec}
f_{12}-f=\frac{\alpha_1-\alpha_2}{|f_2-f_1|^2}(f_2-f_1).
\end{equation}
In the latter vector form this equation was introduced in
\cite{Sch1} under the name of {\em discrete Calapso equation}.
\begin{thm}
Equation (\ref{eq: Schief}) in an associative algebra $\cA$ is 3D
consistent. It possesses a zero curvature representation with the
transition matrices
\begin{equation}\label{eq: Schief L}
 L(f_i,f,\alpha_i;\lambda)=\left(\begin{array}{cc}
 f & \lambda-\alpha_i-ff_i \\ 1 & -f_i \end{array}\right).
\end{equation}
\end{thm}
{\bf Proof.} Equations (\ref{eq: Schief}) on the vertical faces of
Fig.~\ref{back cube} read:
\[
f_i^+=f+(\lambda-\alpha_i)(f^+-f_i)^{-1}=L(f_i,f,\alpha_i;\lambda)\,[f^+].
\]
This gives the transition matrices, which can then be used to
prove the 3D consistency, cf. \cite{BS2}. $\Box$
\smallskip

{\bf 2D systems with fields on vertices.} Regarding T-nets in quadrics,
we encounter an equation with vertex variables
$f:\bbZ^2\to\bbQ=\{f\in\bbR^N:\langle f,f\rangle=\kappa^2\}$,
and with no edge variables:
\begin{equation}\label{eq:2d Lh}
f_{12}-f=a(f_2-f_1),\qquad
a=\frac{\langle f,f_1-f_2\rangle}{\kappa^2-\langle f_1,f_2\rangle}=
\frac{2\langle f,f_1-f_2\rangle}{|f_1-f_2|^2}.
\end{equation}
For the quantities $\alpha_i=2\langle f,f_i\rangle$ the labelling
property (\ref{eq:2d lab}) is fulfilled, but now they are
functions of the vertex variables $f$ rather than parameters of
equation. Comparing eq. (\ref{eq:2d Lh}) with eq. (\ref{eq: Schief
vec}), we see that the former can be regarded as a particular
instance of the latter.
\begin{thm}
Equation (\ref{eq:2d Lh}), describing T-nets in quadrics, is 3D
consistent. It possesses a zero curvature representation with the
transition matrices with entries from $\cC\ell(\bbR^N)$:
\begin{equation}\label{eq:2d Lh L}
 L(f_i,f;\lambda)=\left(\begin{array}{cc}
 f & \lambda+f_i f \\ 1 & -f_i \end{array}\right).
\end{equation}
\end{thm}
{\bf Proof.} 3D consistency has been proven geometrically in
Theorem \ref{Th: M-lat in quadric}. As for the transition
matrices, we can take those from eq. (\ref{eq: Schief L}) with
\[
\alpha_i=2\langle f,f_i\rangle=-ff_i-f_if.
\]
Note the geometrical meaning of the spectral parameter:
$\lambda=2\langle f,f^+\rangle$ for the B\"acklund transformation
$f^+$ from which the transition matrices are constructed. $\Box$
\smallskip

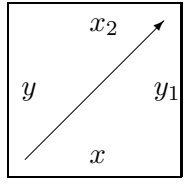
\begin{figure}[htbp]
\setlength{\unitlength}{0.06em}
\begin{center}
\begin{picture}(100,120)
  \put( 0,  0){\line(1,0){100}}
  \put( 0,100){\line(1,0){100}}
  \put(  0, 0){\line(0,1){100}}
  \put(100, 0){\line(0,1){100}}
  \put(10,10){\vector(1,1){80}}
  \put(47,8){$x$}
  \put(47,84){$x_2$}
  \put(8,47){$y$}
  \put(84,47){$y_1$}
\end{picture}
\end{center}
\caption{Elementary quadrilateral; fields on edges}
\label{Fig:quadrilateral edges}
\end{figure}

\begin{figure}[htbp]
\setlength{\unitlength}{0.07em}
\begin{center}
\begin{picture}(450,170)(-30,-20)
  \path(0,0)(100,0)(150,30)(150,130)(50,130)(0,100)(0,0)
  \path(0,0)(50,30)(50,130)
  \path(50,30)(150,30)
  \put(40,35){\vector(-1,2){30}}
  \put(140,40){\vector(-1,1){80}}
  \put(95,5){\vector(-2,1){40}}
  \put(105,85){$R_{13}$}
  \put(23,80){$R_{23}$}
  \put(40,5){$R_{12}$}
  \put(40,-11){$x$}
  \put(130,6){$y$}
  \put(155,75){$z$}
  \put(90,135){$x_{23}$}
  \put(5,120){$y_{13}$}
  \put(-25,50){$z_{12}$}
  \put(90,40){$x_2$}
  \put(20,25){$y_1$}
  \put(55,75){$z_1$}
  \put(190,65){$=$}
  \path(250,0)(350,0)(400,30)(400,130)(300,130)(250,100)(250,0)
  \path(250,100)(350,100)(350,0)
  \path(350,100)(400,130)
  \put(390,35){\vector(-1,2){30}}
  \put(340,10){\vector(-1,1){80}}
  \put(345,105){\vector(-2,1){40}}
  \put(280,30){$R_{13}$}
  \put(355,45){$R_{23}$}
  \put(290,105){$R_{12}$}
  \put(290,-11){$x$}
  \put(380,6){$y$}
  \put(405,75){$z$}
  \put(340,135){$x_{23}$}
  \put(255,120){$y_{13}$}
  \put(225,50){$z_{12}$}
  \put(300,87){$x_3$}
  \put(368,100){$y_3$}
  \put(332,50){$z_2$}
\end{picture}
\end{center}
\caption{Yang--Baxter relation} \label{Fig:YB}
\end{figure}
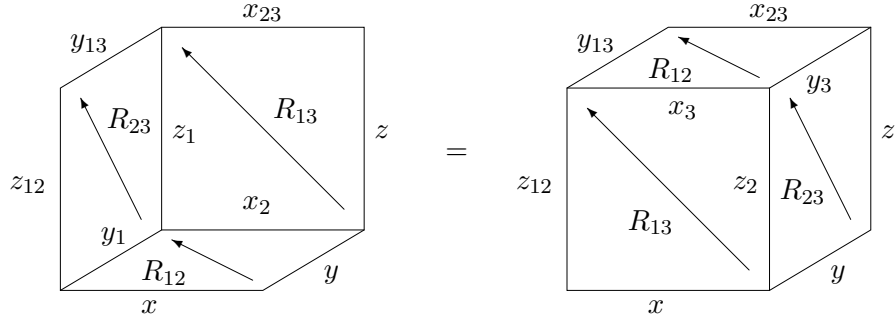

{\bf 2D systems with fields on edges.} Another large class of 2D
systems on quad-graphs build those with fields assigned to the
{\it edges}, see Fig.~\ref{Fig:quadrilateral edges}.
In this situation it is natural to assume that each elementary
quadrilateral carries a map $R:\cX\times\cX\mapsto\cX\times\cX$,
with $\cX$ being the set where the fields $x,y$ take values, so
that $(x_2,y_1)=R(x,y)$. The 3D consistency of such maps can be
encoded in the formula
\begin{equation}\label{YB map}
R_{23}\circ R_{13}\circ R_{12}=R_{12}\circ R_{13}\circ R_{23},
\end{equation}
where each $R_{ij}:\cX^3\mapsto\cX^3$ acts as the map $R$ on the
factors $i,j$ of the cartesian product $\cX^3$, and acts
identically on the third factor. This equation should be
understood as follows. The fields $x,y$ are supposed to be
attached to the edges parallel to the 1st and the 2nd coordinate
axes, respectively. Additionally, consider the fields $z$ attached
to the edges parallel to the 3rd coordinate axis.
Fig.~\ref{Fig:YB} illustrates eq. (\ref{YB map}), its left-hand
side corresponding to the chain of maps along the three rear faces
of the cube:
\[
(x,y)\To{R_{12}}(x_2,y_1),\quad (x_2,z)\To{R_{13}}(x_{23},z_1),\quad
(y_1,z_1)\To{R_{23}}(y_{13},z_{12}),
\]
and the right-hand side corresponding to the chain of maps along
the three front faces of the cube:
\[
(y,z)\To{R_{23}}(y_3,z_2),\quad (x,z_2)\To{R_{13}}(x_3,z_{12}),\quad
(x_3,y_3)\To{R_{12}}(x_{23},y_{13}).
\]
So, eq. (\ref{YB map}) assures that two ways of obtaining
$(x_{23},y_{13},z_{12})$ from the initial data $(x,y,z)$ lead to
the same results.
Maps with the property (\ref{YB map}) were introduced by Drinfeld
\cite{D} under the name of {\em set-theoretical solutions of the
Yang-Baxter equation}, an alternative name {\em Yang--Baxter maps}
was proposed by Veselov in the recent study \cite{V}, see also
references therein. The notion of the zero curvature
representation makes perfect sense for Yang--Baxter maps, and is
expressed as
\[
L(y_1;\lambda)L(x_2;\lambda)=L(x;\lambda)L(y;\lambda).
\]
There is a construction of zero curvature representations for
Yang-Baxter maps with no more input information than the maps
themselves \cite{SV}, close in spirit to Theorem \ref{Th zcr}.
Consider {\it parameter-dependent Yang-Baxter maps}
$R(\alpha,\beta)$, with the parameters $\alpha,\beta\in\bbC$
assigned to the same edges of the quadrilateral on Fig.
\ref{Fig:quadrilateral edges} as the fields $x,y$, opposite edges
carrying the same parameters. Although this can be considered as a
particular case of the general notion, by introducing $\tilde\cX
=\cX\times\bbC$ and $\tilde R(x,\alpha; y,
\beta)=R(\alpha,\beta)(x,y)$, it is convenient for us to keep the
parameter separately. Thus, on Fig. \ref{Fig:YB} all edges
parallel to the $x$ (resp. $y,z$) axis, carry the parameter
$\alpha$ (resp. $\beta,\gamma$), and the corresponding version of
the Yang-Baxter relation reads:
\begin{equation}\label{eq:sYB}
R_{23}(\beta,\gamma) R_{13}(\alpha,\gamma) R_{12}(\alpha,\beta) =
R_{12}(\alpha,\beta) R_{13}(\alpha,\gamma) R_{23}(\beta,\gamma).
\end{equation}

\begin{thm}
Suppose that there is an effective action of the linear group
$G={\rm GL}_N(\bbC)$ on the set $\cX$ (i.e., $A\in G$ acts
identically on $\cX$ only if $A=I$), and that the Yang-Baxter map
$R(\alpha,\beta)$ has the following special form:
\begin{equation}\label{map}
x_2=B(y,\beta,\alpha)[x], \quad  y_1=A(x,\alpha,\beta)[y].
\end{equation}
Here $A,B:\cX\times\bbC\times\bbC\to G$ are some matrix-valued
functions on $\cX$ depending on parameters $\alpha$ and $\beta$
and $A[x]$ denotes the action of $A\in G$ on $x\in\cX$. Then,
whenever $(x_2,y_1)=R(\alpha,\beta)(x,y)$, there holds:
\begin{eqnarray}
 A(y_1,\beta,\lambda)A(x_2,\alpha,\lambda) & = &
 A(x,\alpha,\lambda)A(y,\beta,\lambda),
\label{eq:YBLax 1}\\
 B(x_2,\alpha,\lambda)B(y_1,\beta,\lambda) & = &
 B(y,\beta,\lambda)B(x,\alpha,\lambda). \label{eq:YBLax 2}
\end{eqnarray}
In other words, both $A(x,\alpha,\lambda)$ and
$B^{-1}(x,\alpha,\lambda)$ form zero curvature representations for
$R$.
\end{thm}

\section{Classification of 3D consistent systems}

The notion of 3D consistency proves useful also in various
classification problems of the integrable systems theory. We give
here a presentation of results of \cite{ABS1} on the
classification of integrable quad-graph equations (\ref{eq:2d Q})
with complex fields on vertices and a complex-valued labelling of
edges, and those of \cite{ABS2} on the classification of
integrable equations with complex fields on edges (Yang-Baxter
maps).

\subsection{2D systems with fields on vertices and labelled
edges}\label{Subsect:classif vert}

The classification of integrable equations of the type (\ref{eq:2d
Q}) is performed in \cite{ABS1} under the linearity and symmetry
assumptions of the previous subsection, and one additional
assumption. This latter one is less natural but it is fulfilled
for the vast majority of interesting examples, including the
complex cross-ratio and Hirota equations. An attentive look at
eqs. (\ref{eq:cr 123}), (\ref{eq:H 123}) tells that for the
cross-ratio and for the Hirota equations there holds:
\smallskip

{\it Tetrahedron property.} The value $f_{123}$, existing due to 3D
consistency, depends on $f_1$, $f_2$ and $f_3$, but not on $f$.
\smallskip

Thus, the fields $f_1$, $f_2$, $f_3$ and $f_{123}$ (vertices of a
white tetrahedron on Fig.~\ref{back cube}) are related by a
well-defined equation. Of course, due to symmetry the same holds
for the fields $f$, $f_{12}$, $f_{23}$ and $f_{13}$ (vertices of
the black tetrahedron).

\begin{thm}\label{th:list}
The 3D consistent quad-graph equations (\ref{eq:2d Q}) with the
linearity, symmetry, and tetrahedron properties are exhausted, up
to common M\"obius transformations of the variables $f$ and point
transformations of the parameters $\alpha$, by the following three
lists Q, H, A (we use the abbreviations $x=f$, $u=f_1$, $v=f_2$,
$y=f_{12}$, $\alpha=\alpha_1$, $\beta=\alpha_2$).
\smallskip

List $Q$:
\begin{itemize}
  \item[{\rm(Q1)}] \quad $\alpha(x-v)(u-y)-\beta(x-u)(v-y)+
  \delta^2\alpha\beta(\alpha-\beta)=0$,
  \item[{\rm(Q2)}] \quad $\begin{array}{r} \\
           \alpha(x-v)(u-y)-\beta(x-u)(v-y)+\alpha\beta
             (\alpha-\beta)(x+y+u+v)\\
           -\alpha\beta(\alpha-\beta)(\alpha^2-\alpha\beta+\beta^2)=0,
                         \end{array}$
  \item[{\rm(Q3)}] \quad $\begin{array}{r} \\
                   \sin(\alpha)(xu+vy)-\sin(\beta)(xv+uy)
           -\sin(\alpha-\beta)(xy+uv) \\
                   +\delta^2\sin(\alpha-\beta)\sin(\alpha)\sin(\beta)=0,
                         \end{array}$
  \item[{\rm(Q4)}] \quad $\begin{array}{r} \\
                   {\rm sn}(\alpha)(xu+vy)-{\rm sn}(\beta)(xv+uy)
           -{\rm sn}(\alpha-\beta)(xy+uv) \\
                   +{\rm sn}(\alpha-\beta){\rm sn}(\alpha){\rm sn}(\beta)
           (1+k^2xyuv)=0,
                         \end{array}$
  \end{itemize}
where ${\rm sn}(\alpha)={\rm sn}(\alpha;k)$ is the Jacobi elliptic
function with the modulus $k$.
\smallskip

List $H$:
\begin{itemize}
  \item[{\rm(H1)}] \quad $(x-y)(u-v)+\beta-\alpha=0$,
  \item[{\rm(H2)}] \quad $(x-y)(u-v)+(\beta-\alpha)(x+y+u+v)+
                   \beta^2-\alpha^2=0$,
  \item[{\rm(H3)}] \quad $\alpha(xu+vy)-\beta(xv+uy)+
                   \delta(\alpha^2-\beta^2)=0$.
\end{itemize}

List $A$:
\begin{itemize}
  \item[{\rm(A1)}] \quad $\alpha(x+v)(u+y)-\beta(x+u)(v+y)-
                   \delta^2\alpha\beta(\alpha-\beta)=0$,
  \item[{\rm(A2)}] \quad $\sin(\alpha)(xv+uy)-\sin(\beta)(xu+vy)
                   -\sin(\alpha-\beta)(1+xyuv)=0.$
\end{itemize}
\end{thm}

\noindent
{\bf Remarks.}

1) Parameter $\delta$ in eqs. (Q1), (Q3), (H3), (A1) can be scaled away,
so that one can assume without loss of generality that $\delta=0$
or $\delta=1$.

2) If one extends the transformation group of equations by
allowing M\"obius transformations to act on the variables $x,y$
differently than on $u,v$ (white and black subgraphs of a
bipartite quad-graph), then eq. (A1) is related to (Q1) by the
change $u\to -u$, $v\to -v$, and eq. (A2) is related to
(Q3)$_{\delta=0}$ by the change $u\to 1/u$, $v\to 1/v$. So, really
independent equations are given by the lists Q and H.

3) Note that the above lists contain the complex versions of the
cross-ratio equation (Q1)$_{\delta=0}$, of the Hirota equation
(H3)$_{\delta=0}$, and of the discrete Calapso equation (H1). The
latter two equations are, probably, the oldest ones in our lists,
they can be found in the work of Hirota \cite{H}. Eqs. (Q1) and
${\rm (Q3)}_{\delta=0}$ go back to \cite{QNCL}. Eq. (Q4) was found
in \cite{A} (in the Weierstrass normalization of an elliptic
curve; the observation that the formulas become much nicer in the
Jacobi normalization is due to J.~Hietarinta). Eqs. (Q2),
(Q3)$_{\delta=1}$, (H2) and (H3)$_{\delta=1}$ seem to have
appeared explicitly for the first time in \cite{ABS1}.
\smallskip

The classification in Theorem \ref{th:list} is performed modulo
simultaneous M\"obius transformations of all the variables.
Applying this theorem, one would like to identify a 2D equation at
hand with one of the equations of the lists Q, H, A. This task is
greatly simplified by using some further information. For all 3D
consistent equations with the linearity, symmetry and the
tetrahedron properties, the following holds. The symmetric
biquadratic polynomial
\begin{equation}\label{gxu}
   g(x,u;\alpha,\beta)=QQ_{yv}-Q_yQ_v
\end{equation}
admits a representation
\begin{equation}\label{ghk}
   g(x,u;\alpha,\beta)=k(\alpha,\beta)h(x,u;\alpha),
\end{equation}
where the factor $k$ is antisymmetric,
$k(\beta,\alpha)=-k(\alpha,\beta)$, and the coefficients of the
polynomial $h(x,u;\alpha)$ depend on a single parameter $\alpha$
in such a way that its discriminant
\begin{equation}\label{r}
   r(x)=h^2_u-2hh_{uu}
\end{equation}
does not depend on $\alpha$ at all. Thus, we can characterize 3D
consistent equations $Q=0$ of Theorem \ref{th:list} by much less
complicated objects, namely the biquadratic symmetric polynomials
$h(x,u;\alpha)$ of two variables and the quartic polynomials
$r(x)$ of one variable. The action of the simultaneous M\"obius
transformations on the variables $x\mapsto (ax+b)/(cx+d)$
transforms the polynomials $h$, $r$ as follows:
\begin{eqnarray*}
 & h(x,u;\alpha)\mapsto
(cx+d)^2(cu+d)^2h\!\left(
\dfrac{ax+b}{cx+d},\dfrac{au+b}{cu+d};\,\alpha\!\right),
& \\ & r(x)\mapsto (cx+d)^4\,r\!\left(\dfrac{ax+b}{cx+d}\right). &
\end{eqnarray*}
Using such transformations one can put the polynomial $r(x)$ into
one of the following canonical forms, depending on the
distribution of its zeroes: either $r(x)=0$, or $r(x)=1$ (one
quadruple zero), or $r(x)=x$ (one simple zero and one triple
zero), or $r(x)=x^2$ (two double zeroes), or $r(x)=1-x^2$ (two
simple zeroes and one double zero), or, finally,
$r(x)=(1-x^2)(1-k^2x^2)$ with $k^2\neq 1$ (four simple zeroes).
So, the first step in identifying a 3D consistent equation is
computing the polynomial $r(x)$ and putting it by a M\"obius
transformation into one of the canonical forms above. After that,
one has to identify the polynomial $h(x,u;\alpha)$ with one of
those corresponding to equations of our lists. This might require
a further M\"obius transformation in the cases $r(x)=0$, $r(x)=1$,
and $r(x)=x^2$. In other three cases, $r(x)=x$, $r(x)=1-x^2$, and
$r(x)=(1-x^2)(1-k^2x^2)$ (corresponding to equations (Q2), (Q3)
and (Q4)), the polynomials $h(x,u;\alpha)$ are uniquely determined
by $r(x)$, as a one-parameter family of symmetric biquadratic
polynomials with the discriminant $r(x)$. In the non-degenerate
case $r(x)=(1-x^2)(1-k^2x^2)$ this family is given by
\begin{equation}
h(x,u;\alpha)=\frac{1}{2\,{\rm sn}(\alpha)}\Big(
x^2+u^2-2\,{\rm cn}(\alpha){\rm dn}(\alpha)xu-
{\rm sn}^2(\alpha)(1+k^2x^2u^2)\Big).
\end{equation}
One can recognize this polynomial as the addition theorem for the
elliptic function ${\rm sn}(x;k)$; more precisely,
$h(x,u;\alpha)=0$ if and only if $x={\rm sn}(\xi;k)$ and $u={\rm
sn}(\eta;k)$ with $\xi-\eta=\pm\alpha$. This is the origin of the
elliptic parametrization of the equation (Q4). Eqs. (Q1)--(Q3) are
obtained from (Q4) by degenerations of an elliptic curve into
rational ones. Similarly, eqs. (H1)--(H2) are limiting cases of
(H3). One could be tempted to spoil down the lists Q, H to one
item each. However, the limit procedures necessary for that are
outside of our group of admissible (M\"obius) transformations,
and, on the other hand, in many situations the ``degenerate''
equations (Q1)--(Q3) and (H1)--(H2) are of interest for
themselves. This resembles the situation with the list of six
Painlev\'e equations and the coalescences connecting them.

\subsection{2D systems with fields on edges}

Consider Yang-Baxter maps $R:\cX\times\cX\to\cX\times\cX$,
$(x,y)\mapsto(u,v)$ in the following special framework. Suppose
that $\cX$ is an irreducible algebraic variety, and that $R$ is a
birational automorphism of $\cX\times\cX$. Thus, the birational
map $R^{-1}:\cX\times\cX\to\cX\times\cX$, $(u,v)\mapsto(x,y)$ is
defined. This is depicted on the left square on
Fig.~\ref{fig:YBmaps}. Further, a non-degeneracy condition is
imposed on $R$: rational maps $u(\cdot,y):\cX\to\cX\,$ and
$v(x,\cdot):\cX\to\cX\,$ should be well defined for generic $x$,
resp. $y$. In other words, birational maps $\bar
R:\cX\times\cX\to\cX\times\cX$, $(x,v)\mapsto(u,y)$ and $\bar
R^{-1}:\cX\times\cX\to\cX\times\cX$, $(u,y)\mapsto(x,v)$, called
companion maps to $R$, should be defined. This requirement is
visualized on the right square on Fig.~\ref{fig:YBmaps}.
Birational maps $R$ satisfying this condition are called {\em
quadrirational} \cite{ABS2}.

\def\tmpg{\path(0,0)(0,100)(100,100)(100,0)(0,0)
 \put(47,-13){$x$}\put(-13,47){$v$}\put(47,105){$u$}\put(105,47){$y$}}
\begin{figure}[htbp]
\begin{center}
\setlength{\unitlength}{0.06em}
\begin{picture}(300,140)(0,-15)
  \tmpg \put(13,90){\vector(1,-1){77}}  \put(30,30){$R$}
        \put(87,10){\vector(-1,1){77}}  \put(60,55){$R^{-1}$}
 \put(200,0){
  \tmpg \put(10,13){\vector(1,1){77}}   \put(35,55){$\bar R$}
        \put(90,87){\vector(-1,-1){77}} \put(55,30){$\bar R^{-1}$} }
\end{picture}
\caption{A map $R$ on $\cX\times\cX$, its inverse and its companions}
\label{fig:YBmaps}
\end{center}
\end{figure}
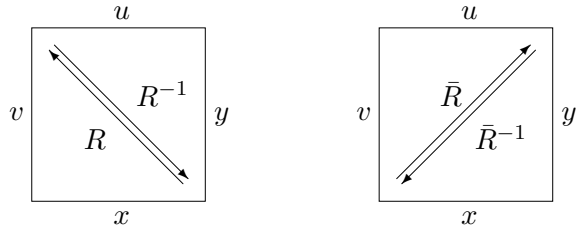

It turns out to be possible to classify all quadrirational maps in
the case $\cX=\bbC\bbP^1$; we give a short presentation of the
corresponding results \cite{ABS2}. Birational isomorphisms of
$\bbC\bbP^1\times\bbC\bbP^1$ are with necessity of the form
\begin{equation}\label{eq:biMoeb}
  R:\;\; u=\frac{a(y)x+b(y)}{c(y)x+d(y)}\,,\quad
         v=\frac{A(x)y+B(x)}{C(x)y+D(x)}\,,
\end{equation}
where $a(y),\ldots,d(y)$ are polynomials in $y$, and
$A(x),\ldots,D(x)$ are polynomials in $x$. For quadrirational
maps, the degrees of all these polynomials are $\le 2$. Dependent
on the highest degree of the coefficients of each fraction in
(\ref{eq:biMoeb}), we say that the map is [1:1], [1:2], [2:1], or
[2:2]. The most rich and interesting subclass is [2:2]. A
necessary condition for a map of this subclass to be
quadrirational is that quartic polynomials
$\delta(y)=a(y)d(y)-b(y)c(y)$ and $\Delta(x)=A(x)D(x)-B(x)C(x)$
are {\em simultaneously} of one of the following five types: they
have either (I) four simple roots, or (II) two simple and two
double roots, or (III) two double roots, or (IV) one simple and
one triple root, or, finally, (V) one quadruple root.
\begin{thm}\label{Th: quadrirat}
Any quadrirational map [2:2] on $\bbC\bbP^1\times\bbC\bbP^1$ is
equivalent, under some change of variables acting by M\"obius
transformations independently on each field $x,y,u,v$, to exactly
one of the following five maps: \bigskip

\noindent
 {\rm(R$_\I$):}\quad  $u= \alpha yP,\quad  v= \beta xP,\quad
  P= \dfrac{(1-\beta)x+\beta-\alpha+(\alpha-1)y}
           {\beta(1-\alpha)x+(\alpha-\beta)yx+\alpha(\beta-1)y}$,
\medskip

\noindent
 {\rm(R$_\II$):}\quad $u= \dfrac{y}{\alpha}\,P,\quad
  v= \dfrac{x}{\beta}\,P,\quad
  P= \dfrac{\alpha x-\beta y+\beta-\alpha}{x-y}$,
\medskip

\noindent
 {\rm(R$_\III$):}\quad $u= \dfrac{y}{\alpha}\,P,\quad
  v= \dfrac{x}{\beta}\,P,\quad
  P= \dfrac{\alpha x-\beta y}{x-y}$,
\medskip

\noindent
 {\rm(R$_\IV$):}\quad $u= yP,\quad v= xP,\quad
  P= 1+\dfrac{\beta-\alpha}{x-y}$,
\medskip

\noindent
{\rm(R$_\V$):}\quad $u= y+P,\quad  v= x+P,\quad
  P= \dfrac{\alpha-\beta}{x-y}$,
\medskip

\noindent
with some suitable constants $\alpha,\beta$.
\end{thm}

Each one of these maps is an involution and coincides with its
companion maps, so that all four arrows on Fig.~\ref{fig:YBmaps}
are described by the same formulas. Note also that these maps come
with the intrisically built-in parameters $\alpha,\beta$. Neither
their existence nor a concrete dependence on parameters is
pre-supposed in Theorem \ref{Th: quadrirat}. A geometric
interpretation of these parameters can be given in terms of
singularities of the map; it turns out that the parameter $\alpha$
is naturally assigned to the edges $x,u$, while $\beta$ is
naturally assigned to the edges $y,v$.

The most remarkable fact about the maps (R$_\I$)--(R$_\V$) is
their 3D consistency. For $\cT=\I,\II,\III,\IV$ or $\V$, denote
the corresponding map $R_\cT$ of Theorem \ref{Th: quadrirat} by
$R_{\cT}(\alpha,\beta)$, indicating the parameters explicitly.
Moreover, for any $\alpha_1,\alpha_2,\alpha_3\in\bbC$, denote by
$R_{ij}=R_{\cT}(\alpha_i,\alpha_j)$ the corresponding maps acting
nontrivially on the $i$-th and the $j$-th factors of
$(\bbC\bbP^1)^3$.

\begin{thm}\label{Th: YB}
For any $\cT=\I,\II,\III,\IV$ or $\V$, the maps $R_{ij}$
satisfy the Yang--Baxter equation (\ref{YB map}).
\end{thm}
Actually, the 3D consistency of quadrirational maps on
$\bbC\bbP^1\times \bbC\bbP^1$ takes place not only for the normal
forms $R_{\cT}(\alpha_i, \alpha_j)$ but under much more general
circumstances. The only condition for quadrirational maps [2:2]
consists in matching singularities along all edges of the cube
(see details in \cite{ABS2}). Similar statements hold also for
quadrirational maps [1:1] and [1:2], so that in the case
$\cX=\bbC\bbP^1$ the properties  of being quadrirational and of
satisfying the Yang-Baxter equation (\ref{YB map}) are related
very closely.

\begin{figure}[htbp]
\begin{center}\includegraphics[width=8cm]{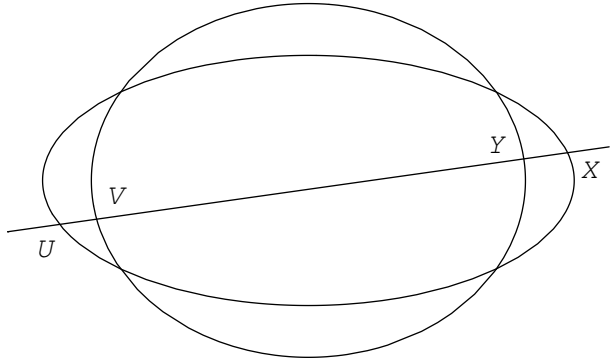}\end{center}
\caption{A quadrirational map on a pair of conics}\label{fig:greeks}
\end{figure}

The maps $R_\cT$ of Theorem \ref{Th: quadrirat} admit a very nice
geometric interpretation. Consider a pair of nondegenerate conics
$Q_1$, $Q_2$ on the plane $\bbC\bbP^2$, so that both $Q_i$ are
irreducible algebraic varieties isomorphic to $\bbC\bbP^1$.  Take
$X\in Q_1$, $Y\in Q_2$, and let $\ell=\overline{XY}$ be the line
through $X,Y$ (well-defined if $X\neq Y$). Generically, the line
$\ell$ intersects $Q_1$ at one further point $U\neq X$, and
intersects $Q_2$ at one further point $V\neq Y$. This defines a
map $\cF:(X,Y)\mapsto(U,V)$, see Fig.~\ref{fig:greeks} for the
$\bbR^2$ picture. The map $\cF:Q_1\times Q_2\mapsto Q_1\times Q_2$
is quadrirational, is an involution and moreover coincides with
its both companions. This follows immediately from the fact that,
knowing one root of a quadratic equation, the second one is a
rational function of the input data. Intersection points $X\in
Q_1\cap Q_2$ correspond to the singular points $(X,X)$ of the map
$\cF$.

Generically, two conics intersect in four points, however,
degeneracies can happen. Five possible types $\I-\V$ of
intersection of two conics are described in detail in \cite{Be}:
\begin{itemize}\setlength{\itemsep}{-0.1pt}
\item[\I:] four simple intersection points;
\item[\II:] two simple intersection points and one point of tangency;
\item[\III:] two points of tangency;
\item[\IV:] one simple intersection point and one point of the
second order tangency;
\item[\V:] one point of the third order tangency.
\end{itemize}
Using rational parametrizations of the conics:
\[
\bbC\bbP^1\ni x\mapsto X(x)\in Q_1\subset\bbC\bbP^2,
 \quad {\rm resp.} \quad
\bbC\bbP^1\ni y\mapsto Y(y)\in Q_2\subset\bbC\bbP^2,
\]
it is easy to see that $\cF$ pulls back to the map
$F:(x,y)\mapsto(x_2,y_1)$ which is quadrirational on
$\bbC\bbP^1\times\bbC\bbP^1$. One shows by a direct computation
that the maps $F$ for the above five situations are exactly the
five maps listed in Theorem \ref{Th: quadrirat}. Now, we obtain
the following geometric interpretation of the statement of Theorem
\ref{Th: YB}.
\begin{figure}[htbp]\label{fig:greeks3}
\begin{center}\includegraphics[width=10cm]{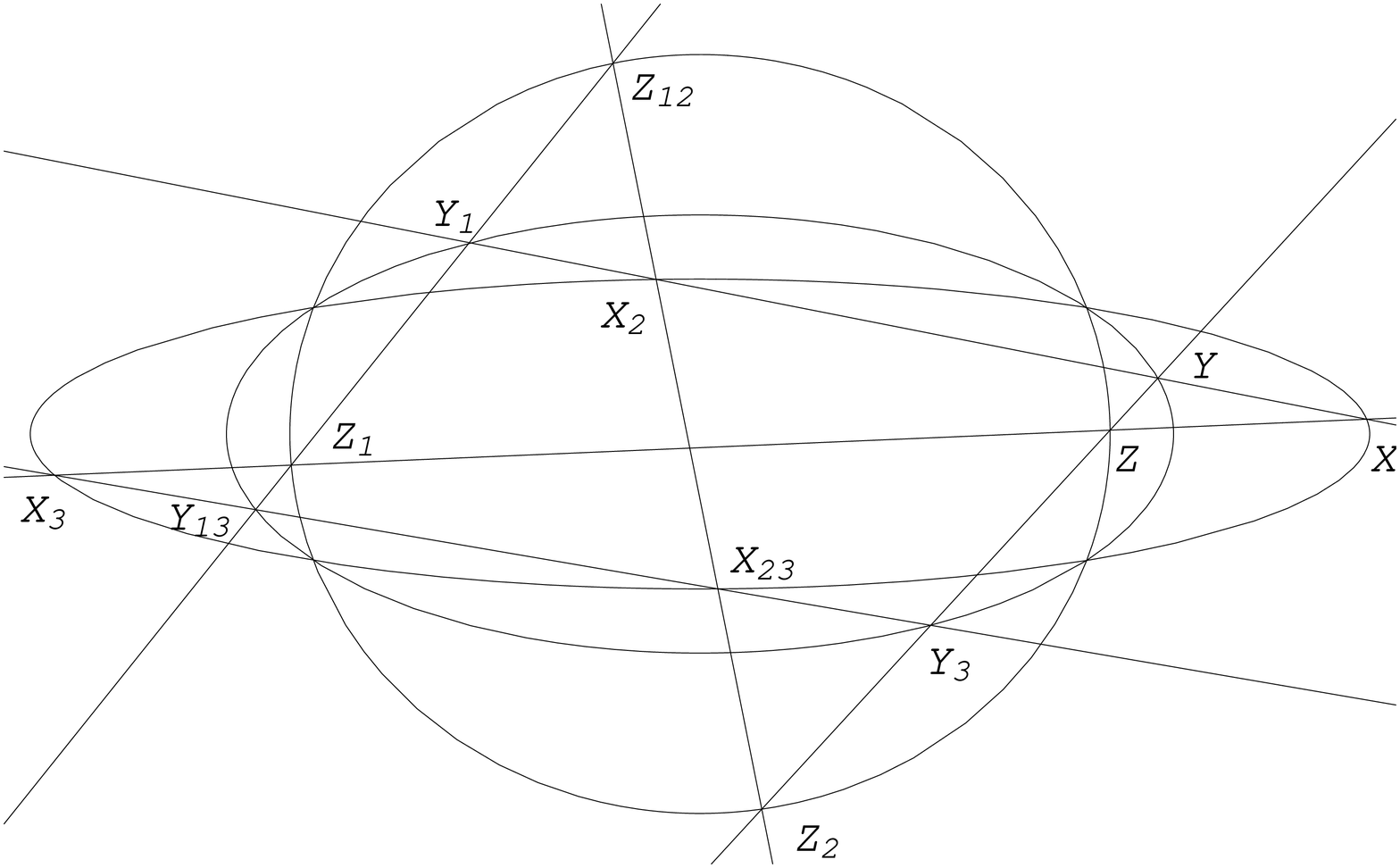}\end{center}
\caption{3D consistency at a linear pencil of conics}
\end{figure}
\begin{thm}\label{th:greeks}
Let $Q_i$, $i=1,2,3$ be three non-degenerate members of a linear
pencil of conics. Let $X\in Q_1$, $Y\in Q_2$ and $Z\in Q_3$ be
arbitrary points on these conics. Define the maps $\cF_{ij}$ as
above, corresponding to the pair of conics $(Q_i,Q_j)$. Set
$(X_2,Y_1)=\cF_{12}(X,Y)$, $(X_3,Z_1)=\cF_{13}(X,Z)$, and
$(Y_3,Z_2)=\cF_{23}(Y,Z)$. Then
\begin{gather*}
\label{geom 1}
 X_{23}=\overline{X_3Y_3}\cap\overline{X_2Z_2}\in Q_1,\quad
 Y_{13}=\overline{X_3Y_3}\cap\overline{Y_1Z_1}\in Q_2, \\
\label{geom 2}
 Z_{12}=\overline{Y_1Z_1}\cap\overline{X_2Z_2}\in Q_3.
\end{gather*}
\end{thm}

\section[Discrete Laplace type equations]
{Integrability of discrete Laplace type equations}
 \label{Sect: nonlin Laplace}

There exist discrete equations on graphs which are not covered by
the theory we developed up to now.
\begin{dfn}\label{Dfn:dLaplace}
Let $\cG$ be a graph, with the set of vertices $V(\cG)$ and the
set of edges $E(\cG)$. {\itbf Discrete Laplace type equations} on
the graph $\cG$ for a function $f:V(\cG)\to\bbC$ read:
\begin{equation}\label{eq:nlin Laplace}
\sum_{x\sim x_0}\phi\big(f(x_0),f(x);\nu(x_0,x)\big)=0.
\end{equation}
There is one equation for every vertex $x_0\in V(\cG)$; the
summation is extended over the set of vertices of $\cG$ connected
to $x_0$ by an edge; the function $\phi$ depends on some
parameters $\nu:E(\cG)\to\bbC$.
\end{dfn}
The classical {\em discrete Laplace equations} on $\cG$ are a
particular case of this definition:
\begin{equation}\label{eq: Laplace}
\sum_{x\sim x_0}\nu(x_0,x)\big(f(x)-f(x_0)\big)=0,
\end{equation}
with some weights $\nu:E(\cG)\to\bbR_+$ attached to (undirected)
edges of $\cG$. We use the notation ${\rm star}(x_0)={\rm
star}(x_0;\cG)$ for the set of edges of $\cG$ incident to $x_0$.
In view of this, discrete Laplace type equations live on stars of
$\cG$.

The notion of {\em integrability} of discrete Laplace type
equations is not well established yet. We propose here a
construction which works under additional assumptions about the
graph $\cG$ and is based on the relation to 3D consistent
quad-graph equations (with fields on vertices). The additional
structural assumption on $\cG$ is that it comes from a strongly
regular polytopal cell decomposition of an oriented surface. To
any such $\cG$ there corresponds canonically a dual cell
decomposition $\cG^*$ (it is only defined up to isotopy, but can
be fixed uniquely with the help of the Voronoi/Delaunay
construction). If one assigns a direction to an edge $e\in
E(\cG)$, then it will be assumed that the dual edge $e^*\in
E(\cG^*)$ is also directed, in a way consistent with the
orientation of the underlying surface, namely so that the pair
$(e,e^*)$ is positively oriented at its crossing point. This
orientation convention implies that $e^{**}=-e$. The {\it double}
$\cD$ is a quad-graph, constructed from $\cG$ and its dual $\cG^*$
as follows. The set of vertices of the double $\cD$ is
$V(\cD)=V(\cG)\sqcup V(\cG^*)$. Each pair of dual edges, say
$e=(x_0,x_1)\in E(\cG)$ and $e^*=(y_0,y_1)\in E(\cG^*)$, defines a
quadrilateral $(x_0,y_0,x_1,y_1)$. These quadrilaterals constitute
the faces of the cell decomposition (quad-graph) $\cD$. Thus, a
star of a vertex in $x_0\in V(\cG)$ generates a flower of adjacent
quadrilaterals from $F(\cD)$ around $x_0$, see
Figs.~\ref{Fig:star}, \ref{Fig:flower}.

This construction can be reversed. Start with a bipartite
quad-graph $\cD$, whose vertices $V(\cD)$ are decomposed into two
complementary halves (``black'' and ``white'' vertices), such that
the ends of each edge from $E(\cD)$ are of different colours. For
instance, any quad-graph embedded in $\bbC$, or in an open disc,
is automatically bipartite. Any bipartite quad-graph produces two
dual polytopal (not necessarily quadrilateral) cell decompositions
$\cG$ and $\cG^*$, with $V(\cG)$ being the ``black'' vertices of
$\cD$ and $V(\cG^*)$ being the ``white'' ones, and edges of $\cG$
(resp. of $\cG^*$) connecting ``black'' (resp. ``white'') vertices
along the diagonals of each face of $\cD$.

\begin{figure}[htbp]
    \setlength{\unitlength}{16pt}
\begin{minipage}[t]{165pt}
\begin{picture}(9,9)(-1,0)
\put(4,4){\circle*{0.3}} \put(4.2,8){\circle*{0.3}}
\put(7,5.9){\circle*{0.3}} \put(7.5,2.6){\circle*{0.3}}
\put(4,0.2){\circle*{0.3}} \put(0.4,3.9){\circle*{0.3}}
\path(4,4)(4.2,8) \path(4,4)(7,5.9) \path(4,4)(7.5,2.6)
\path(4,4)(4,0.2) \path(4,4)(0.4,3.9)
\put(3.1,4.4){$x_0$} \put(7.4,5.9){$x_1$} \put(4.1,8.4){$x_2$}
\put(0.2,4.3){$x_3$} \put(3.9,-0.4){$x_4$} \put(7.7,2.2){$x_5$}
\end{picture}
    \caption{The star of the vertex $x_0$ in the graph $\cG$.}
    \label{Fig:star}
\end{minipage}\hfill
\begin{minipage}[t]{160pt}
\begin{picture}(8,9)(-1,0)
\put(4,4){\circle*{0.3}} \put(4.2,8){\circle*{0.3}}
\put(5,6){\circle{0.3}} \put(3,5.7){\circle{0.3}}
\put(7,5.9){\circle*{0.3}} \put(6,4.2){\circle{0.3}}
\put(7.5,2.6){\circle*{0.3}} \put(5.1,2.5){\circle{0.3}}
\put(4,0.2){\circle*{0.3}} \put(3,2.1){\circle{0.3}}
\put(0.4,3.9){\circle*{0.3}}
\path(4.93,5.86)(4,4)\path(4.95,6.125)(4.2,8)\path(5.15,5.9875)(7,5.9)
\path(3.09,5.55)(4,4)\path(3.07,5.84)(4.2,8)\path(2.87,5.61)(0.4,3.9)
\path(3.07,2.24)(4,4)\path(2.87,2.19)(0.4,3.9)\path(3.07,1.96)(4,0.2)
\path(5.01,2.635)(4,4)\path(5.01,2.36)(4,0.2)\path(5.25,2.5)(7.5,2.6)
\path(5.85,4.185)(4,4)\path(6.09,4.1)(7.5,2.6)\path(6.09,4.35)(7,5.9)
\put(3.1,4){$x_0$} \put(6.3,4.2){$y_0$} \put(7.3,5.9){$x_1$}
\put(5.1,6.3){$y_1$} \put(4.1,8.3){$x_2$} \put(2.1,6){$y_2$}
\put(-0.5,4.3){$x_3$} \put(2.4,1.6){$y_3$} \put(3.9,-0.3){$x_4$}
\put(5.1,2.0){$y_4$} \put(7.5,2.1){$x_5$}
\end{picture}
   \caption{Faces of $\cD$ around the vertex $x_0$.}
   \label{Fig:flower}
\end{minipage}
    \end{figure}

In order to extract from this geometric construction some
consequences for discrete Laplace type equations, we need the
following deep and somewhat mysterious property, which quite often
accompanies the 3D consistency of quad-graph equations with fields
on vertices (\ref{eq:2d Q}). We write here this equation in
slightly modified notations as
\begin{equation}\label{eq:2d for 3leg}
Q(x_0,y_0,x_1,y_1;\alpha_0,\alpha_1)=0.
\end{equation}
For notational simplicity, vertices $x$ stand here for the
corresponding fields $f(x)$; the edges $(x_0,y_0)$, $(x_0,y_1)$
carry the labels $\alpha_0$, $\alpha_1$, respectively.

\begin{dfn}\label{Def:3leg}
An equation (\ref{eq:2d for 3leg}) possesses a {\itbf three-leg
form} centered at the vertex $x_0$, if it is equivalent to the
equation
\begin{equation}\label{eq:3leg add}
\psi(x_0,y_0;\alpha_0)-\psi(x_0,y_1;\alpha_1)=
\phi(x_0,x_1;\alpha_0,\alpha_1)
\end{equation}
with some functions $\psi, \varphi$. The terms on the left-hand
side correspond to ``short'' legs $(x_0,y_0), (x_0,y_1)\in
E(\cD)$, while the right-hand side corresponds to the ``long'' leg
$(x_0,x_1)\in E(\cG)$.
\end{dfn}

A summation of quad-graph equations for the flower of
quadrilaterals adjacent to the ``black'' vertex $x_0\in V(\cG)$
immediately leads, due to the telescoping effect, to the following
statement.
\begin{thm}\label{Th: Laplace for 3legs}
a) Suppose that eq. (\ref{eq:2d for 3leg}) on a bipartite
quad-graph $\cD$ possesses a three-leg form. Then a restriction of
any solution $f:V(\cD)\to\bbC$ to the ``black'' vertices $V(\cG)$
satisfies the discrete equations of the Laplace type,
\begin{equation}\label{eq:Laplace for 3legs}
\sum_{k=1}^n\phi(x_0,x_k;\alpha_{k-1},\alpha_k)=0,
\end{equation}
where $n$ is the valence of the vertex $x_0$ in $\cG$.

b) Conversely, given a solution $f:V(\cG)\to\bbC$ of these Laplace
type equations on a graph $\cG$ coming from a cell decomposition
of a simply-connected surface, there exists a one-parameter family
of its extensions $f:V(\cD)\to\bbC$ satisfying eq. (\ref{eq:2d Q})
on the double $\cD$. Such an extension is uniquely defined by a
value at one arbitrary vertex from $V(\cG^*)$.
\end{thm}
Sometimes it is more convenient to write the three-leg equation
(\ref{eq:3leg add}) in the multiplicative form:
\begin{equation}\label{eq:3leg mult}
\Psi(x_0,y_0;\alpha_0)/\Psi(x_0,y_1;\alpha_1)=
\Phi(x_0,x_1;\alpha_0,\alpha_1)
\end{equation}
with some functions $\Psi$, $\Phi$, so that the Laplace type
equations (\ref{eq:Laplace for 3legs}) also become multiplicative:
\begin{equation}\label{eq:Laplace for 3legs mult}
\prod_{k=1}^n\Phi(x_0,x_k;\alpha_{k-1},\alpha_k)=1.
\end{equation}

This relation between integrable quad-graph equations and Laplace
type equations was discovered in \cite{BS1}. The three-leg forms
for all quad-graph equations of Theorem \ref{th:list} were found
in \cite{ABS2}. The next theorem provides three-leg forms for all
equations of the lists Q and H (the results for the list A follow
from these ones). The functions $\phi(x,y;\alpha,\beta)$, resp.
$\Phi(x,y;\alpha,\beta)$, corresponding to the ``long'' legs,
yield integrable additive (resp. multiplicative) Laplace type
equations on arbitrary planar graphs.

\begin{thm}\label{Th: 3leg forms}
Three-leg forms exist for all equations from Theorem
\ref{th:list}:
\medskip

${\rm (Q1)_{\delta=0}}$: Additive three-leg form with
$\phi(x,y;\alpha,\beta)=\psi(x,y;\alpha-\beta)$,
\begin{equation}\label{Q10 3leg}
\psi(x,u;\alpha)=\frac{\alpha}{x-u}\,.
\end{equation}

${\rm (Q1)_{\delta=1}}$: Multiplicative three-leg form with
$\Phi(x,y;\alpha,\beta)=\Psi(x,y;\alpha-\beta)$,
\begin{equation}\label{Q11 3leg}
\Psi(x,u;\alpha)=\frac{x+\alpha-u}{x-\alpha-u}\,.
\end{equation}

{\rm (Q2)}: Multiplicative three-leg form with
$\Phi(x,y;\alpha,\beta)=\Psi(x,y;\alpha-\beta)$,
\begin{equation}\label{Q2 3leg Psi}
\Psi(x,u;\alpha)=\frac{(X+\alpha)^2-u}{(X-\alpha)^2-u}\,,
\end{equation}
where $x=X^2$.

${\rm (Q3)_{\delta=0}}$: Multiplicative three-leg form with
$\Phi(x,y;\alpha,\beta)=\Psi(x,y;\alpha/\beta)$,
\begin{equation}\label{Q30 3leg Psi}
\Psi(x,u;\alpha)=\frac{\alpha x-u}{x-\alpha u}\,.
\end{equation}

${\rm (Q3)_{\delta=1}}$: Multiplicative three-leg form with
$\Phi(x,y;\alpha,\beta)=\Psi(x,y;\alpha-\beta)$,
\begin{equation}\label{Q31 3leg Psi}
\Psi(x,u;\alpha)=\frac{\sin(X+\alpha)-u}{\sin(X-\alpha)-u}\,,
\end{equation}
where $x=\sin(X)$.

{\rm (Q4)}: Multiplicative three-leg form with
$\Phi(x,y;\alpha,\beta)=\Psi(x,y;\alpha-\beta)$,
\begin{equation}\label{Q4 3leg Psi}
\Psi(x,u;\alpha)=
  \frac{\Theta(X+\alpha)}{\Theta(X-\alpha)}\cdot
  \frac{{\rm sn}(X+\alpha)-u}{{\rm sn}(X-\alpha)-u}\,,
\end{equation}
where $x={\rm sn}(X)$, and $\Theta(X)$ is the Jacobi
theta-function.

{\rm (H1)}: Additive three-leg form with
\begin{equation}\label{H1 3leg phi psi}
\phi(x,y;\alpha,\beta)=\frac{\alpha-\beta}{x-y}\,, \quad
\psi(x,u;\alpha)=x+u.
\end{equation}

{\rm (H2)}: Multiplicative three-leg form with
\begin{equation}\label{H2 3leg phi psi}
\Phi(x,y;\alpha,\beta)=\frac{x-y+\alpha-\beta}{x-y-\alpha+\beta}\,,
\quad  \Psi(x,u;\alpha)=x+u+\alpha.
\end{equation}

{\rm (H3)}: Multiplicative three-leg form with
\begin{equation}\label{H3 3leg phi psi}
\Phi(x,y;\alpha,\beta)=\frac{\beta x-\alpha y}{\alpha x-\beta
y}\,,  \quad  \Psi(x,u;\alpha)=xu+\delta\alpha.
\end{equation}
\end{thm}

One sees that there are only six ``long'' legs functions
$\phi(x,y;\alpha,\beta)$, resp. $\Phi(x,y;\alpha,\beta)$, leading
to integrable Laplace type equations on arbitrary planar graphs.
Three of them are rational in $x, y$. Each of the corresponding
Laplace type equations admits two extensions to a 3D consistent
quad-graph equation: one from the list Q, where the ``short'' legs
essentially coincide with the ``long'' ones -- (Q1)$_{\delta=0}$,
(Q1)$_{\delta=1}$, and (Q3)$_{\delta=0}$, and another one from the
list H, with different ``short'' legs -- (H1), (H2), and (H3).
Other three functions $\Phi$ are rational in $y$ only, and require
for a uniformizing change of the variable $x$. Corresponding
Laplace type equations admit only one extension to 3D consistent
quad-graph equations -- (Q2), (Q3)$_{\delta=1}$, and (Q4). Thus,
restriction to a subgraph allows us to derive from 3D consistent
quad-graph equations much more general systems, which clearly
inherit the integrability. See a detailed realization of this idea
for the (Q4) equation in \cite{AS}. \smallskip

{\bf Remark.} It should be mentioned that existence of a three-leg
form allows us to derive (and, in some sense, to explain) the
tetrahedron property of Sect. \ref{Subsect:classif vert}. Indeed,
consider three faces adjacent to the vertex $f_{123}$ on
Fig.~\ref{cube again}, namely the quadrilaterals
$(f_1,f_{12},f_{123},f_{13})$, $(f_2,f_{23},f_{123},f_{12})$, and
$(f_3,f_{13},f_{123},f_{23})$. A summation of the three-leg forms
(centered at $f_{123}$) of equations corresponding to these three
faces leads to the equation
\begin{equation}
\phi(f_{123},f_1;\alpha_2,\alpha_3)+\phi(f_{123},f_2;\alpha_3,\alpha_1)
+\phi(f_{123},f_3;\alpha_1,\alpha_2)=0.
\end{equation}
This equation relates the fields at the vertices of the ``white''
tetrahedron on Fig.~\ref{cube again}. It can be interpreted as a
discrete Laplace type equation coming from a spatial flower with
three petals.


\section[Geometry of 3D consistent equations]
{Geometry of boundary value problems for 3D consistent equations}
\label{Sect: discr geom}

We discuss here several aspects of the problem of embedding of a
quad-graph into a regular multi-dimensional square lattice. As a
combinatorial problem, it was studied in a more general setting of
arbitrary cubic complexes in \cite{DSh1, DSh2, ShSh}. We are
interested here in its relation to integrable equations.

We start with the question about correct initial value problems
for discrete 2D equations on quad-graphs. Let $P$ be a {\em path}
in the quad-graph $\cD$, i.e., a sequence of edges ${\mathfrak
e}_j=(v_j,v_{j+1})\in E(\cD)$. We denote by $E(P)=\{{\mathfrak
e}_j\}$ and $V(P)=\{v_j\}$ the set of edges and the set of
vertices of the path $P$, respectively. One says that the Cauchy
problem for the path $P$ is well-posed, if for any set of data
$f_P:V(P)\to\bbC$ there exists a unique solution $f:V(\cD)\to\bbC$
such that $f\rest_{V(P)}=f_P$.

The answer to this question was given in \cite{AV} with the help
of the notion of a strip in $\cD$.
\begin{dfn}\label{dfn:strip}
A {\,\em strip} in $\cD$ is a sequence of quadrilateral faces
$q_j\in F(\cD)$ such that any pair $q_{j-1}$, $q_j$ is adjacent
along the edge $\mathfrak e_j= q_{j-1}\cap q_j$, and $\mathfrak
e_j$, $\mathfrak e_{j+1}$ are opposite edges of $q_j$. The edges
$\mathfrak e_j$ are called {\,\em traverse edges} of the strip.
\end{dfn}
So, any strip in $\cD$ can be associated to a label $\alpha$
sitting on all its traverse edges $\mathfrak e_j$.
\begin{thm}\label{Th: AdlerVeselov}
Consider one of the 3D consistent systems
(\ref{eq:2d Q}) of Theorem \ref{th:list}.
Let $\cD$ be a finite simply connected quad-graph
without self-crossing strips, and let $P$ be a path without
self-crossings in $\cD$. Then:
\begin{enumerate}
\item[i)] If each strip in $\cD$ intersects $P$ exactly once, then
the Cauchy problem for $P$ is well-posed.
\item[ii)] If some strip in $\cD$ intersects $P$ more than once,
then the Cauchy problem for $P$ is overdetermined (has in general
no solutions).
\item[iii)] If some strip in $\cD$ does not intersect $P$, then the
Cauchy problem for $P$ is underdetermined (has in general more
than one solution).
\end{enumerate}
\end{thm}
It should be mentioned that this theorem is not valid for
equations without the 3D consistency property. One of the proofs
of the statement i) in \cite{AV} is based on the embedding of
$\cD$ into the unit cube of $\bbZ^N$, where $N$ is the number of
edges in $P$ (the number of distinct strips in $\cD$). In this
embedding the path $P$ turns into the path $(0,0,0,\ldots,0)$,
$(1,0,0,\ldots,0)$, $(1,1,0,\ldots,0)$, $\ldots$ ,
$(1,1,1,\ldots,1)$. It is clear that for a 3D consistent equation
the Cauchy problem for this path is well-posed.

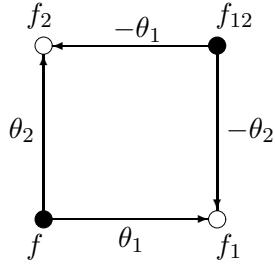
\begin{figure}[htbp]
\begin{center}
\setlength{\unitlength}{0.06em}
\begin{picture}(200,140)(-50,-20)
  \put(100,  0){\circle{10}} \put(0  ,100){\circle{10}}
  \put(  0,  0){\circle*{10}}  \put(100,100){\circle*{10}}
  \put(5,0){\vector(1,0){90}}
  \put(0,5){\vector(0,1){90}}
  \put(95,100){\vector(-1,0){90}}
  \put(100,95){\vector(0,-1){90}}
  \put(-10,-20){$f$}
  \put(99,-20){$f_1$}
  \put(99,113){$f_{12}$}
  \put(-10,113){$f_2$}
  \put(43,-16){$\theta_1$}
  \put(40,105){$-\theta_1$}
  \put(-20,47){$\theta_2$}
  \put(105,47){$-\theta_2$}
\end{picture}
\caption{Labelling of directed edges} \label{Fig:labelling dir}
\end{center}
\end{figure}

A different aspect of embeddability of the quad-graph $\cD$ into a
regular multi-dimensional cubic lattice was studied in \cite{BMS},
based on the following result \cite{KS}.
\begin{thm}\label{theorem KS}
A quad-graph $\cD$ admits an embedding in $\bbC$ with all rhombic
faces if and only if the following two conditions are satisfied:
\begin{itemize}
\item[i)] No strip crosses itself or is periodic.
\item[ii)] Two distinct strips cross each other at most once.
\end{itemize}
\end{thm}
Given a rhombic embedding $p:V(\cD)\to\bbC$, one defines the following
function on the directed edges of $\cD$ with values in
$\bbS^1=\{\theta\in\bbC:|\theta|=1\}$:
\begin{equation}\label{label canonical}
\theta(x,y)=p(y)-p(x),\quad \forall (x,y)\in\vec{E}(\cD).
\end{equation}
This function can be called a {\em labelling of directed edges},
since it satisfies $\theta(-\mathfrak e)=-\theta(\mathfrak e)$ for
any $\mathfrak e\in\vec{E}(\cD)$, and the values of $\theta$ on
two opposite and equally directed edges of any quadrilateral from
$F(\cD)$ are equal. See. Fig.~\ref{Fig:labelling dir}. For any
labelling $\theta:\vec{E}(\cD)\to\bbS^1$ of directed edges, the
function $\alpha=\theta^2:E(\cD)\to\bbS^1$ is a {\it labelling} of
(undirected) edges of $\cD$ in our usual sense.

\begin{dfn}
A rhombic embedding $p:V(\cD)\to\bbC$ of a quad-graph $\cD$
is called {\itbf quasicrystallic}, if the set of values of the
function $\theta:\vec{E}(\cD)\to\bbS^1$ defined by
(\ref{label canonical}), is finite, say
$\Theta=\{\pm\theta_1,\ldots,\pm\theta_d\}$.
\end{dfn}
It is of a central importance that any quasicrystallic
rhombic embedding $p$ can be seen as a sort of a projection
of a certain two-dimensional subcomplex (combinatorial surface)
$\Omega_\cD$ of a multi-dimensional regular square lattice
$\bbZ^d$. The vertices of $\Omega_\cD$ are given by a map
$P:V(\cD)\to\bbZ^d$ constructed as follows. Fix some $x_0\in
V(\cD)$, and set $P(x_0)=\bf 0$. The images in $\bbZ^d$ of all other
vertices of $\cD$ are defined recurrently by the property:
\smallskip

for any two neighbors $x,y\in V(\cD)$, if
$p(y)-p(x)=\pm\theta_i\in\Theta$, then

$P(y)-P(x)=\pm {\boldsymbol e}_i\,$, where ${\boldsymbol e}_i$ is
the $i$-th coordinate vector of $\bbZ^d$. \smallskip

\noindent
Edges and faces of $\Omega_\cD$ correspond to edges and faces of
$\cD$, so that the combinatorics of $\Omega_\cD$ is that of $\cD$.

\begin{figure}[htbp]
\setlength{\unitlength}{0.08em}
\begin{center}
\begin{picture}(350,120)(-30,-10)
  \path(0,0)(50,0)(90,30)(90,80)(40,80)(0,50)(0,0)
  \path(0,0)(40,30)(40,80)
  \path(40,30)(90,30)
  \put(45,38){$f$}
  \put(-5,-15){$f_1$}
  \put(95,30){$f_2$}
  \put(40,87){$f_3$}
  \put(45,-15){$f_{12}$}
  \put(95,87){$f_{23}$}
  \put(-10,60){$f_{13}$}
  \put(130,45){\vector(1,0){20}}
  \path(200,0)(250,0)(290,30)(290,80)(240,80)(200,50)(200,0)
  \path(200,50)(250,50)(250,0)
  \path(250,50)(290,80)
  \put(225,35){$f_{123}$}
  \put(195,-15){$f_1$}
  \put(295,30){$f_2$}
  \put(240,87){$f_3$}
  \put(245,-15){$f_{12}$}
  \put(295,87){$f_{23}$}
  \put(190,60){$f_{13}$}
\end{picture}
\end{center}
\caption{Elementary flip}
\label{Fig:flip}
\end{figure}
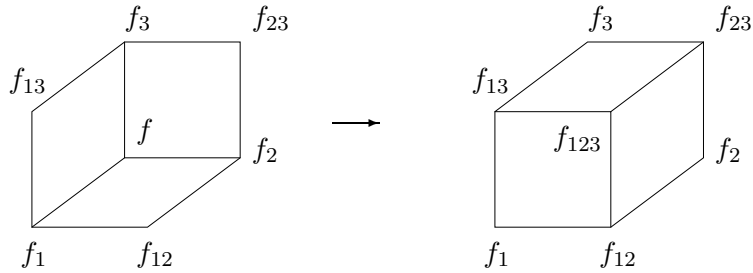

To exploit possibilities provided by the 3D consistency,
extend the labeling $\theta:\vec{E}(\cD)\to\bbS^1$ to all edges of
$\bbZ^d$, assuming that all edges parallel to (and directed as)
$\boldsymbol e_k$ carry the label $\theta_k$. This gives, of course,
also the labelling $\alpha=\theta^2$ of undirected edges of $\bbZ^d$.
Now, any 3D consistent equation can be imposed not only on $\Omega_\cD$,
but on the whole of $\bbZ^d$:
\begin{equation}\label{eq: square Q}
  Q(f,f_j,f_{jk},f_k;\alpha_j,\alpha_k)=0,\qquad 1\le j\neq k\le d.
\end{equation}
Here indices stand for the shifts into the coordinate directions.
Obviously, for any solution $f:\bbZ^d\to\bbC$ of eq. (\ref{eq:
square Q}), its restriction to $V(\Omega_\cD)\sim V(\cD)$ gives a
solution of the corresponding equation on the quad-graph $\cD$. As
for the reverse procedure, i.e., for the extension of an arbitrary
solution of eq. (\ref{eq:2d Q}) from $\cD$ to $\bbZ^d$, more
thorough considerations are necessary. An elementary step of such
an extension consists of finding $f$ at the eighth vertex of an
elementary 3D cube from the known values at seven vertices, see
Fig.~\ref{Fig:flip}. This can be alternatively viewed as a flip
(elementary transformation) on the set of rhombically embedded
quad-graphs $\cD$, or on the set of the corresponding surfaces
$\Omega_\cD$ in $\bbZ^d$. The 3D consistency assures that any
quad-graph $\cD$ (or any corresponding surface $\Omega_\cD$)
obtainable from the original one by such flips, carries a unique
solution of eq. (\ref{eq: square Q}) which is an extension of the
original one.

\begin{dfn}
For a given set $V\subset\bbZ^d$, its {\itbf hull} $\cH(V)$  is
the minimal set $\cH\subset\bbZ^d$ containing $V$ and satisfying
the condition: if three vertices of an elementary square belong to
$\cH$, then so does the fourth vertex.
\end{dfn}
One shows by induction that for an arbitrary connected subcomplex
of $\bbZ^d$ with the set of vertices $V$, its hull is a {\it
brick}
\begin{equation}\label{brick}
 \Pi_{\ba,\bb}=\big\{\bn=(n_1,\ldots,n_d)\in\bbZ^d:\;
 a_k\le n_k\le b_k,\;\, k=1,\ldots,d\,\big\},
\end{equation}
where
\begin{equation}
 a_k=a_k(V)=\min_{\bn\in V} n_k,\quad
 b_k=b_k(V)=\max_{\bn\in V} n_k,\quad
 k=1,\ldots,d,
\end{equation}
and in the case that $n_k$ are unbounded from below or from above on
$V$, we set $a_k(V)=-\infty$, resp. $b_k(V)=\infty$.

Combinatorially, all points of the hull $\cH(V(\Omega_\cD))$ can
be reached by flips of Fig.~\ref{Fig:flip}, starting from
$\Omega_\cD$. However, there might be obstructions, having nothing
to do with 3D consistency, for extending solutions of eq.
(\ref{eq:2d Q}) from a combinatorial surface (two-dimensional
subcomplex of $\bbZ^d$) to its hull. For instance, the surface
$\Omega$ shown on Fig.~\ref{fig nonconvex} supports solutions of
eq. (\ref{eq:2d Q}) which {\em cannot} be extended to solutions of
eq. (\ref{eq: square Q}) on the whole of $\cH(V(\Omega))$: the
recursive extension will lead to contradictions. The reason for
this is a non-monotonicity of $\Omega$: it contains pairs of
points which cannot be connected by a monotone path, i.e., by a
path in $\Omega$ with all directed edges lying in one octant of
$\bbZ^d$. However, such surfaces $\Omega$ {\em do not} come from
rhombic embeddings, and in the case of $\Omega_\cD$ there will be
no contradictions.
\begin{thm}\label{th extend}
Let the combinatorial surface $\Omega_\cD$ in $\bbZ^d$ come from a
rhombic embedding of a quad-graph $\cD$, and let its hull be
$\cH(V(\Omega_\cD))=\Pi_{\ba,\bb}$. An arbitrary solution of eq.
(\ref{eq:2d Q}) on $\Omega_\cD$ can be uniquely extended to a
solution of eq. (\ref{eq: square Q}) on $\Pi_{\ba,\bb}$.
\end{thm}

\begin{figure}[htbp]
\begin{center}
\setlength{\unitlength}{0.04em}
\begin{picture}(500,250)(-80,-50)
 \path(10,10)(-70,-50)(70,-50)(70,100)(280,100)(280,0)(350,0)(430,60)
      (360,60)(360,160)(150,160)(70,100)
 \path(110,130)(320,130)(320,30)(390,30)
 \path(140,100)(220,160)
 \path(210,100)(290,160)
 \path(280,100)(360,160)
 \path(280, 50)(360,110)
 \path(280,  0)(360, 60)
 \dashline[+30]{10}(150,160)(150,100)
 \path(150,100)(150,10)
 \path(70,-50)(150,10)
 \path(10,10)(70,10)
 \dashline[+30]{10}(70,10)(150,10)
 \path(-30,-20)(70,-20)
 \dashline[+35]{10}(70,-20)(110,-20)
 \path(110,-20)(110,100)
 \dashline[+30]{10}(110,100)(110,130)
 \path(0,-50)(70,2.5)
 \path(70,50)(136.6,100)
 \path(70, 0)(150, 60)
\end{picture}
\caption{A non-monotone surface in $\bbZ^3$}\label{fig nonconvex}
\end{center}
\end{figure}

Note that intersections of $\Omega_\cD$ with bricks correspond to
{\em combinatorially convex subsets} of $\cD$, as defined in
\cite{M2}.

\chapter[Discrete linear complex analysis]
{Discrete complex analysis. Linear theory}
\label{Chap:_discr_analysis-linear}

\section[Basics of discrete complex analysis]
{Basic notions of the discrete linear complex analysis}
\label{Sect: discr analysis Intro}

Many constructions in the discrete complex analysis
are parallel to the discrete differential geometry
in the space of the real dimension 2.

Recall that a harmonic function $u:\bbR^2\simeq\bbC\to\bbR$ is
characterized by the relation:
\[
\Delta u=\frac{\partial^2 u}{\partial x^2}+
 \frac{\partial^2 u}{\partial y^2}=0.
\]
A conjugate harmonic function
$v:\bbR^2\simeq\bbC\to\bbR$ is defined by the Cauchy-Riemann
equations:
\[
\frac{\partial v}{\partial y}=\frac{\partial u}{\partial x}\,,
\qquad
 \frac{\partial v}{\partial x}=-\frac{\partial u}{\partial y}\,.
\]
Equivalently, $f=u+iv:\bbR^2\simeq\bbC\to\bbC$ is holomorphic, i.e.,
satisfies the Cauchy-Riemann equation:
\[
\frac{\partial f}{\partial y}=i\,\frac{\partial f}{\partial x}\,.
\]
The real and the imaginary parts of a holomorphic function
are harmonic, and any real-valued harmonic function can be
considered as a real part of a holomorphic function.

A standard classical way to discretize these notions, going back
to Ferrand \cite{F} and Duffin \cite{Du1}, is the following. A
function $u:\bbZ^2\to\bbR$ is called {\em discrete harmonic}, if
\[
(\Delta u)_{m,n}=
 u_{m+1,n}+u_{m-1,n}+u_{m,n+1}+u_{m,n-1}-4u_{m,n}=0.
\]
A natural definition domain of a conjugate discrete harmonic
function $v:(\bbZ^2)^*\to\bbR$ is the {\em dual lattice:}

\begin{figure}[htbp]
    \setlength{\unitlength}{14pt}
\begin{center}
\begin{picture}(8,9)
\put(0,0){\circle*{0.5}}\put(0,4){\circle*{0.5}}
\put(0,8){\circle*{0.5}}
\put(4,0){\circle*{0.5}}\put(4,4){\circle*{0.5}}
\put(4,8){\circle*{0.5}}
\put(8,0){\circle*{0.5}}\put(8,4){\circle*{0.5}}
\put(8,8){\circle*{0.5}}
\put(2,2){\circle{0.5}}\put(2,6){\circle{0.5}}
\put(6,2){\circle{0.5}}\put(6,6){\circle{0.5}}
\path(0,0)(0,8)(8,8)(8,0)(0,0) \path(4,0)(4,8) \path(0,4)(8,4)
\dashline[+30]{0.4}(2.25,2)(5.75,2)
\dashline[+30]{0.4}(2.25,6)(5.75,6)
\dashline[+30]{0.4}(2,2.25)(2,5.75)
\dashline[+30]{0.4}(6,2.25)(6,5.75)
\put(8.5,4){$u$}\put(6.5,6){$v$}
\end{picture}
\end{center}
    \end{figure}

\noindent The defining {\em discrete Cauchy--Riemann equations}
read:
\begin{figure}[htbp]
    \setlength{\unitlength}{15pt}
    \begin{minipage}[t]{170pt}
\begin{picture}(6,7)(-5,-1)
\put(0,2){\circle*{0.5}} \put(4,2){\circle*{0.5}}
\put(-1.2,2){$u_0$}      \put(4.5,2){$u_1$} \path(0,2)(4,2)
\put(2,0){\circle{0.5}} \put(2,4){\circle{0.5}}
\put(2.7,-0.2){$v_0$}   \put(2.7,4){$v_1$}
\dashline[+30]{0.4}(2,0.25)(2,3.75)
\end{picture}
\[
\quad\qquad v_1-v_0=u_1-u_0
\]
    \end{minipage}\hfill
\begin{minipage}[t]{170pt}
    \setlength{\unitlength}{15pt}
\begin{picture}(6,7)(-2,-1)
\put(2,0){\circle*{0.5}} \put(2,4){\circle*{0.5}}
\put(2.7,-0.2){$u_0$}    \put(2.7,4){$u_1$} \path(2,0)(2,4)
\put(0,2){\circle{0.5}} \put(4,2){\circle{0.5}}
\put(-1.2,2){$v_0$} \put(4.5,2){$v_1$}
\dashline[+30]{0.4}(0.25,2)(3.75,2)
\end{picture}
\[
 v_1-v_0=-(u_1-u_0)\qquad\qquad
\]
    \end{minipage}
\end{figure}
\medskip

\noindent The corresponding {\em discrete holomorphic} function
$f:\bbZ^2\cup(\bbZ^2)^*\to\bbC$ is defined on the superposition of
the original square lattice $\bbZ^2$ and the dual one
$(\bbZ^2)^*$, by the formula $f=\left\{\begin{array}{cc} u, &
\bullet\,, \\ iv, & \circ\,, \end{array}\right.$ which comes to
replace the smooth one $f=u+iv$. Remarkably, the discrete
Cauchy--Riemann equation for $f$ is one and the same for the both
pictures above:
\begin{figure}[htbp]
    \setlength{\unitlength}{15pt}
    \begin{minipage}[t]{170pt}
\begin{picture}(6,5)(-5,0)
\put(0,2){\circle*{0.5}} \put(4,2){\circle*{0.5}}
\put(-1.2,2){$f_1$}      \put(4.5,2){$f_3$}
\put(2,0){\circle{0.5}} \put(2,4){\circle{0.5}}
\put(2.7,-0.2){$f_2$}   \put(2.7,4){$f_4$} \path(0,2)(1.824,0.176)
\path(2.176,0.176)(4,2)(2.176,3.824) \path(1.824,3.824)(0,2)
\end{picture}
    \end{minipage}\hfill
\begin{minipage}[t]{170pt}
    \setlength{\unitlength}{15pt}
\begin{picture}(6,5)(-3,0)
\put(2,0){\circle*{0.5}} \put(2,4){\circle*{0.5}}
\put(2.7,-0.2){$f_2$}    \put(2.7,4){$f_4$}
\put(0,2){\circle{0.5}} \put(4,2){\circle{0.5}}
\put(-1.2,2){$f_1$} \put(4.5,2){$f_3$}
\path(0.176,1.824)(2,0)(3.824,1.824)
\path(3.824,2.176)(2,4)(0.176,2.176)
\end{picture}
    \end{minipage}
\vspace{0.4cm}
\[
 f_4-f_2=i(f_3-f_1).
\]
\end{figure}
\newpage
\noindent This discretization of the Cauchy-Riemann equations
apparently preserves the most number of important structural
features, and the corresponding theory has been developed in
\cite{F, D1}. A pioneering step in the direction of a further
generalization of the notions of discrete harmonic and discrete
holomorphic functions was undertaken by Duffin \cite{D2}, where
the combinatorics of $\bbZ^2$ was given up in favor of arbitrary
planar graphs with rhombic faces. A far reaching generalization of
these ideas was given by Mercat \cite{M1}, who extended the theory
to discrete Riemann surfaces.

\medskip
Discrete harmonic functions can be defined for an arbitrary graph
$\cG$ with the set of vertices $V(\cG)$ and the set of edges
$E(\cG)$.
\begin{dfn}\label{Dfn discr harmonic}
For a given weight function $\nu:E(\cG)\to\bbR_+$ on edges of
$\cG$, the {\itbf Laplacian} is the operator acting on functions
$f:V(\cG)\to\bbC$ by
\begin{equation}\label{eq: Laplacian}
(\Delta f)(x_0)=\sum_{x\sim x_0}\nu(x_0,x)(f(x)-f(x_0)),
\end{equation}
where the summation is extended over the set of vertices $x$
connected to $x_0$ by an edge. A function $f:V(\cG)\to\bbC$ is
called {\itbf discrete harmonic} (with respect to the weights
$\nu$), if $\Delta f=0$.
\end{dfn}

\noindent
The positivity of weights $\nu$ in this definition is important
from the analytic point of view, as it guarantees, e.g., the maximum
principle for the discrete Laplacian under suitable boundary conditions
(so that discrete harmonic functions come as minimizers of a convex
functional). However, from the pure algebraic point of view, one might
consider at times also arbitrary real (or even complex) weights.

If $\cG$ comes from a cellular decomposition of an oriented
surface, let $\cG^*$ be its dual graph, and let the quad-graph
$\cD$ be its double, see Sect. \ref{Sect: nonlin Laplace}. Extend
the weight function to the edges of $\cG^*$ according to the rule
\begin{equation}\label{eq: nu*}
\nu(e^*)=1/\nu(e).
\end{equation}

\begin{figure}[htbp]
\setlength{\unitlength}{0.04em}
\begin{picture}(200,220)(-300,-110)
 \put(0,0){\circle*{13}}\put(200,0){\circle*{13}}
 \put(100,-75){\circle{13}} \put(100,75){\circle{13}}
 \put(4,-3){\vector(4,-3){91}}
 \put(95,71.25){\vector(-4,-3){91}}
 \put(196,3){\vector(-4,3){91}}
 \put(105,-71.25){\vector(4,3){91}}
 \put(-37,-7){$x_0$} \put(215,-7){$x_1$}
 \put(94,-100){$y_0$} \put(94,97){$y_1$}
 \put(20,55){$-\theta_1$} \put(150,-60){$\theta_1$}
 \put(145,55){$-\theta_0$} \put(25,-60){$\theta_0$}
\end{picture}
\caption{Positively oriented quadrilateral, with a labelling of
directed edges}\label{Fig: diamond dir}
\end{figure}
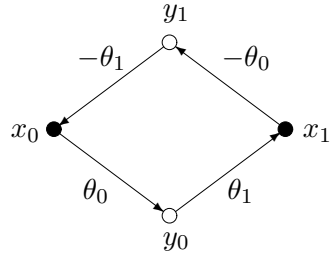

\begin{dfn}\label{Dfn discr holo}
A function $f:V(\cD)\to\bbC$ is called {\itbf discrete
holomorphic} (with respect to the weights $\nu$), if for any
positively oriented quadrilateral $(x_0,y_0,x_1,y_1)\in F(\cD)$
(see Fig.~\ref{Fig: diamond dir}) there holds:
\begin{equation}\label{eq: discr CR}
  \frac{f(y_1)-f(y_0)}{f(x_1)-f(x_0)}=i\nu(x_0,x_1)=
  -\frac{1}{i\nu(y_0,y_1)}.
\end{equation}
These equations are called {\itbf discrete Cauchy-Riemann
equations}.
\end{dfn}

The relation between discrete harmonic and discrete holomorphic functions
is the same as in the smooth case. It is given by the following
statement which is a particular case of Theorem \ref{Th: Laplace
for 3legs}.

\begin{thm}\label{Th: dholo vs dharm}
a) If a function $f:V(\cD)\to\bbC$ is discrete holomorphic, then
its restrictions to $V(\cG)$ and to $V(\cG^*)$ are discrete
harmonic.

b) Conversely, any discrete harmonic function $f:V(\cG)\to\bbC$
admits a family of discrete holomorphic extensions to $V(\cD)$,
differing by an additive constant on $V(\cG^*)$. Such an extension
is uniquely defined by a value at one arbitrary vertex $y\in
V(\cG^*)$.
\end{thm}

\section[M-transformation for CR equations]
{Moutard transformation for Cauchy-Riemann equations on quad-graphs}

Observe that discrete Cauchy-Riemann equations (\ref{eq: discr
CR}) is formally not different from the Moutard equations
(\ref{eq:dMou}) for T-nets. One only has to fix the
orientation of all quadrilateral faces $(x_0,y_0,x_1,y_1)\in F(\cD)$.
We assume that it is inherited from the orientation of the underlying
surface.

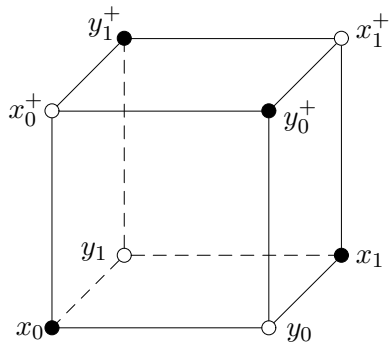
\begin{figure}[htbp]
\begin{center}
\setlength{\unitlength}{0.05em}
\begin{picture}(200,240)(0,0)
 \put(0,0){\circle*{10}}    \put(150,0){\circle{10}}
 \put(0,150){\circle{10}}   \put(150,150){\circle*{10}}
 \put(50,200){\circle*{10}} \put(200,200){\circle{10}}
 \put(50,50){\circle{10}}   \put(200,50){\circle*{10}}
 \path(0,0)(145,0)       \path(0,0)(0,145)
 \path(150,5)(150,150)   \path(5,150)(150,150)
 \path(3.53,153.53)(50,200)    \path(150,150)(196.47,196.47)
 \path(50,200)(195,200)
 \path(200,195)(200,50) \path(200,50)(153.53,3.53)
 \dashline[+30]{10}(0,0)(46.47,46.47)
 \dashline[+30]{10}(50,55)(50,200)
 \dashline[+30]{10}(55,50)(200,50)
 \put(-25,-5){$x_0$} \put(-30,145){$x_0^+$}
 \put(162,-6){$y_0$} \put(160,140){$y_0^+$}
 \put(210,45){$x_1$} \put(210,200){$x_1^+$}
 \put(20,50){$y_1$}  \put(25,205){$y_1^+$}
\end{picture}
\caption{Elementary cube of $\bf D$}\label{Fig: Mou cube}
\end{center}
\end{figure}

One can apply now the Moutard transformation of Sect.
\ref{Subsect: discr Mou} to discrete holomorphic functions.
To this aim, one has to choose the orientation of {\em all}
elementary quadrilaterals on Fig.~\ref{Fig: Mou cube}.
This can be done, for example, as follows: for the quadrilaterals
$(x_0^+,y_0^+,x_1^+,y_1^+)\in F(\cD^+)$, choose the orientation
to coincide with that of the corresponding $(x_0,y_0,x_1,y_1)\in F(\cD)$.
For a ``vertical'' quadrilateral over an edge $(x,y)\in E(\cD)$,
assume that $x\in V(\cG)$, $y\in V(\cG^*)$, and choose the positive
orientation corresponding to the cyclic order $(x,y,y^+,x^+)$ of its
vertices. Observe that under this convention, two opposite ``vertical''
quadrilaterals are always oriented differently.

In the case of arbitrary quad-graphs, one has to generalize
one more ingredient of the Moutard transformation, namely the
data (MT$_2^\Delta$).

\begin{thm}
On an arbitrary bipartite quad-graph $\cD$, valid initial conditions
for a Moutard transformation of the discrete Cauchy-Riemann equations
consist of
\newcounter{dCR}
\begin{list}{{\rm(MCR$_\arabic{dCR}^\Delta$)}}
{\usecounter{dCR}\setlength{\leftmargin}{0.5cm}}
\item value of $f$ at one point $x^+(0)\in V(\cD^+)$;
\item values of weights on ``vertical'' quadrilaterals $(x,y,y^+,x^+)$
attached to all edges $(x,y)$ of a Cauchy path in $\cD$.
\end{list}
\end{thm}
See Theorem \ref{Th: AdlerVeselov} for necessary and
sufficient conditions for a path to be a Cauchy path, i.e., to
support initial data for a well-posed Cauchy problem. It is natural
to assign the weights on the ``vertical'' quadrilaterals to the
underlying edges of $\cD$.

Weights $\nu$ on the faces of $\cD$ together with the data
(MCR$_2^\Delta$) yield the transformed weights $\nu^+$ on the faces
of $\cD^+$, as well as the weights over all edges of $E(\cD)$.
This can be considered as a Moutard transformation for
{\em Cauchy-Riemann equations} on $\cD$. Finding a {\em solution}
$f:V(\cD^+)\to\bbC$ of the transformed equations
requires additionally the datum (MCR$_1^\Delta$).

Note that the system of weights $\nu$ is highly redundant, due to
eq. (\ref{eq: nu*}). To fix the ideas in writing the equations, we
stick to the weights attached to the ``black'' diagonals of the
quadrilateral faces of the complex $\bf D$. On the ground floor,
these are the edges of the ``black'' graph $\cG$; on the first
floor, these are the  edges of the ``black'' graph  which is a
copy of $\cG^*$; and for the ``vertical'' faces, these are the
edges $(x,y^+)$, with $x\in  V(\cG)$ and $y\in V(\cG^*)$. Needless
to say that the latter weights can be assigned to the quad-graph
edges $(x,y)\in E(\cD)$. So, we write the discrete Cauchy-Riemann
equations as follows:
\begin{eqnarray}
f(y_1)-f(y_0) & = & i\nu(x_0,x_1)\big(f(x_1)-f(x_0)\big),\\
f(x_0^+)-f(x_1^+) & = &
   i\nu(y_0^+,y_1^+)\big(f(y_1^+)-f(y_0^+)\big),\\
f(x^+)-f(y) & = & i\nu(x,y^+)\big(f(y^+)-f(x)\big).
\end{eqnarray}
Denote, for the sake of brevity,
\[
\nu=\nu(x_0,x_1),\quad \nu^+=\nu(y_0^+,y_1^+),\quad
 \mu_{jk}=\nu(x_j,y_k^+).
\]
Regarding the weights $\nu$, $\mu_{00}$, and $\mu_{01}$ as the input
of the Moutard transformation on an elementary hexahedron of $\bf D$,
its output consists of the weights $\nu^+$, $\mu_{10}$, and $\mu_{11}$,
given by (cf. eq. (\ref{eq:dMou 1}))
\begin{equation}\label{eq: dCR Darboux}
\nu^+\nu=-\mu_{11}\mu_{00}=-\mu_{10}\mu_{01}
=\frac{\nu\mu_{00}\mu_{01}}{\mu_{00}-\mu_{01}-\nu}\,.
\end{equation}
This transformation is well defined for real weights $\nu$, $\mu_{jk}$,
but it {\em does not} preserve, in general, positivity of the weights
$\nu$.

To give a different form of this transformation, observe that
relation $\mu_{11}\mu_{00}=\mu_{10}\mu_{01}$ for each elementary
quadrilateral $(x_0,y_0,x_1,y_1)$ of $\cD$ yields the existence of
the function $\theta:V(\cD)\to\bbC$, defined up to a
constant factor, such that $i\mu_{jk}=\theta(y_k)/\theta(x_j)$.
Moreover, choosing $\theta(x_0)$ real at some point $x_0\in V(\cG)$,
one sees that $\theta$ takes real values on $V(\cG)$ and
imaginary values on $V(\cG^*)$.  An easy computation shows
that the last equation in (\ref{eq: dCR Darboux}) is equivalent to
\[
\theta(y_1)-\theta(y_0)=i\nu(x_0,x_1)\big(\theta(x_1)-\theta(x_0)\big),
\]
so that the function $\theta$ is discrete holomorphic with respect
to the weights $\nu$. For the transformed weights $\nu^+$ one finds:
\begin{equation}\label{eq:dCR Darboux theta}
\nu^+\nu=\frac{\theta(y_0)\theta(y_1)}{\theta(x_0)\theta(x_1)}.
\end{equation}
Conversely, an arbitrary discrete holomorphic function
$\theta:V(\cD)\to\bbC$ defines, via eq. (\ref{eq:dCR Darboux theta}),
a Moutard transformation of the discrete Cauchy-Riemann equations.
It should be mentioned that the data
(MCR$_2^\Delta$) can be re-formulated in terms of the function
$\theta$:
\medskip

(MCR$_2^\Delta$)\quad values of $\theta$ at all vertices along
a Cauchy path in $\cD$.
\medskip

\noindent
{\bf Remark.} Moutard transformation for discrete Cauchy-Riemann
equations yields,
by restriction to the ``black'' graphs, a sort of Darboux transformation
of arbitrary discrete Laplacians on $\cG$ into discrete Laplacians on
$\cG^*$. A further discussion of this transformation for discrete Laplace
operators can be found in \cite{NSD, DGNS}.

\section{Integrable discrete Cauchy-Riemann equations}

We now turn to a fruitful question about ``stationary points'' of the
Moutard transformation discussed in the previous section. More precisely,
this is the question about conditions on the weights $\nu:E(\cG)\to\bbR_+$
such that there exists a Moutard transformation for which the opposite
faces of any elementary hexahedron of $\bf D$ (see Fig.~\ref{Fig: Mou cube})
carry {\em identical} equations.

\begin{thm}\label{Th 3D CR}
A system of discrete Cauchy-Riemann equations with the function
$\nu:E(\cG)\sqcup E(\cG^*)\to\bbR_+$ satisfying (\ref{eq: nu*}) admits
a Moutard transformation into itself, if and only if for all
$x_0\in V(\cG)$ and for all $y_0\in V(\cG^*)$ the following
conditions are fulfilled:
\begin{equation}\label{eq:prod=1}
\prod_{e\in{\rm star}(x_0;\cG)}
 \frac{1+i\nu(e)}{1-i\nu(e)}=1,\quad
 \prod_{e^*\in{\rm star}(y_0;\cG^*)}
 \frac{1+i\nu(e^*)}{1-i\nu(e^*)}=1.
\end{equation}
\end{thm}
{\bf Proof.} Opposite faces of $\cD$ and $\cD^+$ carry identical
equations, if $\nu^+\nu=1$ in eq. (\ref{eq: dCR Darboux}). Clearly,
this yields also $\mu_{11}\mu_{00}=\mu_{10}\mu_{01}=-1$, which means
that the opposite ``vertical'' faces also support identical equations
(recall that opposite ``vertical'' faces carry different orientation).
Moreover, given $\nu=\nu(x_0,x_1)$ for an elementary quadrilateral
$(x_0,y_0,x_1,y_1)$ of $\cD$, we find that the input data $\mu_{00}$,
$\mu_{01}$ of the Moutard transformation should be related as follows:
\[
\frac{\nu\mu_{00}\mu_{01}}{\mu_{00}-\mu_{01}-\nu}=1
\quad\Leftrightarrow\quad
\mu_{01}=\frac{\mu_{00}-\nu}{\mu_{00}\nu+1}=
\begin{pmatrix} 1 & -\nu \\ \nu & 1 \end{pmatrix}[\mu_{00}],
\]
where the standard notation for the action of $PGL_2(\bbC)$ on
$\bbC$ by M\"obius transformations is used. This means that all
the weights on the vertical faces of a ``stationary'' Moutard
transformation are completely defined by just one of them, so that
such transformations form a one-parameter family. To derive a
condition for $\nu$ for the existence of a ``stationary'' Moutard
transformation, consider a flower of quadrilaterals
$(x_0,y_{k-1},x_k,y_k)$ around $x_0\in V(\cG)$ (see
Fig.~\ref{Fig:flower}). In the natural notations, we find:
\[
\mu_{0,k}=\frac{\mu_{0,k-1}-\nu_k}{\mu_{0,k-1}\nu_k+1}=
\begin{pmatrix} 1 & -\nu_k \\ \nu_k & 1\end{pmatrix}[\mu_{0,k-1}].
\]
Running around $x_0$ should for any $\mu_{00}$ return its value,
which means that the matrix product
$\begin{pmatrix} A & -B \\
B & A \end{pmatrix}=
\displaystyle{\prod_k^{\curvearrowleft}}\begin{pmatrix} 1 & -\nu_k \\
\nu_k & 1 \end{pmatrix}$
should be proportional to the identity matrix. This matrix product
is easily computed:
\[
A=\frac{1}{2}\Big(\prod_k(1+i\nu_k)+\prod_k(1-i\nu_k)\Big),\quad
B=\frac{1}{2i}\Big(\prod_k(1+i\nu_k)-\prod_k(1-i\nu_k)\Big),
\]
and the condition $B=0$ is
equivalent to the first equality in (\ref{eq:prod=1}). The second
condition in (\ref{eq:prod=1}) is proved similarly, by considering
a flower of quadrilaterals around $y_0\in V(\cG^*)$.  $\Box$
\medskip

Thus, existence of a ``stationary'' Moutard transformation singles
out a special class of discrete Cauchy-Riemann equations, which
have to be considered as 2D systems with the 3D consistency
property, see Sect. \ref{Sect: 3d consist}. In other words, such
Cauchy-Riemann equations should be termed {\em integrable}. The
main difference as compared with the examples in Sect. \ref{Sect:
3d consist}, is that discrete Cauchy-Riemann equations naturally
depend on the {\em orientation} of the elementary quadrilaterals,
and that their parameters $\nu$ are apparently assigned not to the
edges of the quad-graph, but rather to diagonals of its faces.

The integrability condition (\ref{eq:prod=1}) admits a nice geometric
interpretation. It is convenient (especially for positive real-valued
$\nu$) to use the notation
\begin{equation}\label{eq:phi}
\nu(e)=\tan\frac{\phi(e)}{2},\quad \phi(e)\in(0,\pi).
\end{equation}
The condition $\nu(e^*)=1/\nu(e)$ is translated into
\begin{equation}\label{eq:phi*}
\phi(e^*)=\pi-\phi(e),
\end{equation}
while the condition (\ref{eq:prod=1}) says that for
all $x_0\in V(\cG)$ and for all $y_0\in V(\cG^*)$ there holds:
\begin{equation}\label{eq:prodexp=1}
\prod_{e\in{\rm star}(x_0;\cG)}\exp(i\phi(e))=1, \quad
\prod_{e^*\in{\rm star}(y_0;\cG^*)}\exp(i\phi(e^*))=1.
\end{equation}
These conditions should be compared with conditions \cite{KS}
characterizing the angles $\phi:E(\cG)\sqcup E(\cG^*)\to(0,\pi)$
of a rhombic embedding of a quad-graph $\cD$, which consist of
(\ref{eq:phi*}) and
\begin{equation}\label{eq:sum=2pi}
\sum_{e\in{\rm star}(x_0;\cG)}\phi(e)=2\pi, \quad
\sum_{e^*\in{\rm star}(y_0;\cG^*)}\phi(e^*)=2\pi,
\end{equation}
for all $x_0\in V(\cG)$ and for all $y_0\in V(\cG^*)$. Thus, the
integrability condition (\ref{eq:prodexp=1}) says that the system
of angles $\phi:E(\cG)\sqcup E(\cG^*)\to(0,\pi)$ comes from a
realization of the quad-graph $\cD$ as a {\it rhombic ramified
embedding} in $\bbC$. Flowers of such an embedding can wind around
its vertices more than once.

Another formulation of the integrability conditions
is given in terms of the {\em edges} of the rhombic realizations.

\begin{thm}
Integrability condition (\ref{eq:prod=1}) for the weight function
$\nu:E(\cG)\sqcup E(\cG^*)\to\bbR_+$ is equivalent to the
following one: there exists a labeling
$\theta:\vec{E}(\cD)\to\bbS^1$ of directed edges of $\cD$ such
that, in notations of Fig.~\ref{Fig: diamond dir},
\begin{equation}\label{nu paral}
\nu(x_0,x_1)=\frac{1}{\nu(y_0,y_1)}=
i\,\frac{\theta_0-\theta_1}{\theta_0+\theta_1}.
\end{equation}
Under this condition, the 3D consistency of the discrete
Cauchy-Riemann equations is assured by the following values of the
weights $\nu$ on the diagonals of the vertical faces of $\bf D$:
\begin{equation}\label{nu vert}
\nu(x,y^+)=i\,\frac{\theta-\lambda}{\theta+\lambda},
\end{equation}
where $\theta=\theta(x,y)$, and $\lambda\in\bbC$ is an arbitrary
number having the interpretation of the label carried by all
vertical edges of $\bf D$: $\lambda=\theta(x,x^+)=\theta(y,y^+)$.
\end{thm}
So, integrable discrete Cauchy-Riemann equations can be given a form
with parameters attached to {\em directed} edges of $\cD$:
\begin{equation}\label{CR paral}
\frac{f(y_1)-f(y_0)}{f(x_1)-f(x_0)}=
\frac{\theta_1-\theta_0}{\theta_1+\theta_0},
\end{equation}
where
\[
\theta_0=p(y_0)-p(x_0)=p(x_1)-p(y_1), \quad
\theta_1=p(y_1)-p(x_0)=p(x_1)-p(y_0),
\]
and $p:V(\cG)\to\bbC$ is a rhombic realization of the quad-graph $\cD$.
Since
\[
\frac{\theta_1-\theta_0}{\theta_1+\theta_0}=
\frac{p(y_1)-p(y_0)}{p(x_1)-p(x_0)},
\]
we see that for a discrete holomorphic function $f:V(\cG)\to\bbC$,
the quotient of diagonals of the $f$-image of any quadrilateral
$(x_0,y_0,x_1,y_1)\in F(\cD)$ is equal to the quotient of
diagonals of the corresponding rhombus.

A standard construction of zero curvature representation for 3D
consistent equations, given in Theorem \ref{Th zcr}, leads in the
present case to the following result.
\begin{thm}\label{Thm: CR zcr}
The discrete Cauchy-Riemann equations (\ref{CR paral}) admit a
zero curvature representation with spectral parameter dependent
$2\times 2$ matrices along $(x,y)\in\vec{E}(\cD)$ given by
\begin{equation}\label{eq:L CR}
\renewcommand{\arraystretch}{1.4}
L(y,x,\alpha;\lambda)=\begin{pmatrix} \lambda+\theta &
-2\theta(f(x)+f(y)) \\ 0 & \lambda-\theta \end{pmatrix},
\end{equation}
where $\theta=p(y)-p(x)$.
\end{thm}
Linearity of the discrete Cauchy-Riemann equations is reflected in the
triangular structure of the transition matrices.

Also, all constructions of Sect. \ref{Sect: discr geom} can be
applied to integrable discrete Cauchy-Riemann equations. In
particular, for weights coming from a quasicrystallic rhombic
embedding of the quad-graph $\cD$, with the labels $\Theta=\{\pm
\theta_1,\ldots,\pm\theta_d\}$, discrete holomorphic functions can
be extended from the corresponding surface
$\Omega_\cD\subset\bbZ^d$ to its hull, preserving discrete
holomorphy. Here we have in mind the following natural definition:
\begin{dfn}\label{def multidim holo}
A function $f:\bbZ^d\to\bbC$ is called {\itbf discrete
holomorphic}, if it satisfies, on each elementary square of
$\,\bbZ^d$, the equation
\begin{equation}\label{CR square}
  \frac{f({\boldsymbol n}+{\boldsymbol e}_j+{\boldsymbol e}_k)-
  f({\boldsymbol n})}{f({\boldsymbol n}+{\boldsymbol e}_j)-
  f({\boldsymbol n}+{\boldsymbol e}_k)}=
  \frac{\theta_j+\theta_k}{\theta_j-\theta_k}\,.
\end{equation}
\end{dfn}
For discrete holomorphic functions in $\bbZ^d$, the transition matrices
along the edges $(\bn,\bn+\be_k)$ of $\bbZ^d$ are given by
\begin{equation}\label{L CR d}
\renewcommand{\arraystretch}{1.4}
L_k({\boldsymbol n};\lambda)=\begin{pmatrix} \lambda+\theta_k &
-2\theta_k(f({\boldsymbol n}+{\boldsymbol e}_k)+f({\boldsymbol n})) \\
0 & \lambda-\theta_k \end{pmatrix}.
\end{equation}

All observations of this section hold also in the case of generic
complex weights $\nu$, which leads to $\theta\in\bbC$ and to parallelogram
realizations of $\cD$.

\section{Discrete exponential functions}
\label{Sect exp}

A discrete exponential function on quad-graphs $\cD$ was defined
and studied in \cite{M1,K}. Its definition, for a given rhombic
embedding $p:V(\cD)\to\bbC$, is as follows: fix a point $x_0\in
V(\cD)$. For any other point $x\in V(\cD)$, choose some path
$\{\mathfrak e_j\}_{j=1}^n\subset \vec{E}(\cD)$ connecting $x_0$ to $x$,
so that $\mathfrak e_j=(x_{j-1},x_j)$ and $x_n=x$. Let the slope of the
$j$th edge be $\theta_j=p(x_j)-p(x_{j-1})\in\bbS^1$. Then
\[
e(x;z)=\prod_{j=1}^n \frac{z+\theta_j}{z-\theta_j}.
\]
Clearly, this definition depends on the choice of the point
$x_0\in V(\cD)$, but not on the path connecting $x_0$ to $x$. A
question posed in \cite{K} is whether discrete exponential
functions form a basis in the space of discrete holomorphic functions
on $\cD$.

An extension of the discrete exponential function from
$\Omega_\cD$ to the whole of $\bbZ^d$ is the {\it discrete
exponential function}, given by the following simple formula:
\begin{equation}\label{discr exp}
  e(\boldsymbol n;z)=\prod_{k=1}^d
  \Big(\frac{z+\theta_k}{z-\theta_k}\Big)^{n_k}.
\end{equation}
For $d=2$, this function was considered in \cite{F,D1}. The
discrete Cauchy-Riemann equations for the discrete exponential
function are easily checked: they are equivalent to a simple
identity
\[
\Big(\frac{z+\theta_j}{z-\theta_j}\cdot
\frac{z+\theta_k}{z-\theta_k}-1\Big)\Big/
\Big(\frac{z+\theta_j}{z-\theta_j}-
\dfrac{z+\theta_k}{z-\theta_k}\Big)=
\frac{\theta_j+\theta_k}{\theta_j-\theta_k}\,.
\]
At a given $\boldsymbol n\in\bbZ^d$, the discrete exponential
function is rational with respect to the parameter $z$, with poles
at the points $\epsilon_1\theta_1,\ldots,\epsilon_d\theta_d$,
where $\epsilon_k={\rm sign}\, n_k$.

Equivalently, one can identify the discrete exponential function
by its initial values on the axes:
\begin{equation}\label{discr exp ini}
  e(n\boldsymbol e_k;z)=\Big(\frac{z+\theta_k}{z-\theta_k}\Big)^n.
\end{equation}
A still another characterization says that $e(\cdot;z)$ is the
B\"acklund transformation of the zero solution of discrete
Cauchy-Riemann equations on $\bbZ^d$, with the ``vertical''
parameter $z$.

We now show that the discrete exponential functions form a basis
in some natural class of functions (growing not faster than
exponentially), thus answering in affirmative the above mentioned
question by Kenyon.
\begin{thm}\label{Th: exp dense}
Let $f$ be a discrete holomorphic function on $V(\cD)\sim
V(\Omega_\cD)$, satisfying
\begin{equation}\label{exp growth}
  |f(\boldsymbol n)|\le \exp(C(|n_1|+\ldots+|n_d|)),\qquad \forall \bn\in
  V(\Omega_\cD),
\end{equation}
with some $C\in\bbR$. Extend it to a discrete holomorphic function
on $\cH(V(\Omega_\cD))$. There exists a function $g$ defined
on the disjoint union of small neighborhoods around the points
$\pm\theta_k\in\bbC$ and holomorphic on each one of these
neighborhoods, such that
\begin{equation}\label{dense}
f(\boldsymbol n)-f(\boldsymbol 0)=\frac{1}{2\pi i}\int_{\Gamma}
g(\lambda)e(\boldsymbol n;\lambda)d\lambda, \qquad \forall
\boldsymbol n\in\cH(V(\Omega_\cD)),
\end{equation}
where $\Gamma$ is a collection of $2d$ small loops, each one
running counterclockwise around one of the points $\;\pm\theta_k$.
\end{thm}
\noindent
{\bf Proof} is constructive and consists of three steps.
\begin{itemize}
\item[(i)] Extend $f$ from $V(\Omega_\cD)$ to
$\cH(V(\Omega_\cD))$; inequality (\ref{exp
growth}) propagates in the extension process, if the constant $C$
is chosen large enough.
\item[(ii)] Introduce the restrictions $f_n^{(k)}$ of
$f:\cH(V(\Omega_\cD))\to\bbC$ to the coordinate axes:
\[
f_n^{(k)}=f(n\be_k), \qquad a_k(\Omega_\cD)\le n\le
b_k(\Omega_\cD).
\]
\item[(iii)] Set
$g(\lambda)=\sum_{k=1}^d (g_k(\lambda)+g_{-k}(\lambda))$, where
the functions $g_{\pm k}(\lambda)$ vanish everywhere except in
small neighborhoods of the points $\pm\theta_k$, respectively, and
are given there by convergent series
\begin{equation}\label{dense g_i}
g_k(\lambda)=\frac{1}{2\lambda}\left(f_1^{(k)}-f(\boldsymbol 0)
+\sum_{n=1}^\infty\Big(\frac{\lambda-\theta_k}{\lambda+\theta_k}\Big)^n
\big(f_{n+1}^{(k)}-f_{n-1}^{(k)}\big)\right),
\end{equation}
and a similar formula for $g_{-k}(\lambda)$. Formula (\ref{dense})
is then easily verified by computing the residues at
$\lambda=\pm\theta_k$.
\end{itemize}
It is important to observe that the data $f_n^{(k)}$, necessary
for the construction of $g(\lambda)$, are {\em not among} the
values of $f$ on $V(\cD)\sim V(\Omega_\cD)$ known initially, but
are encoded in the extension process.

\section{Discrete logarithmic function}
\label{Sect discr log}

The discrete logarithmic function on $\cD$ can be defined as
follows \cite{K}. Fix some point $x_0\in V(\cD)$, and set
\begin{equation}\label{dLog}
 \ell(x)=\frac{1}{2\pi i}\int_{\Gamma}
 \frac{\log(\lambda)}{2\lambda}\,e(x;\lambda)d\lambda, \qquad
 \forall x\in V(\cD).
\end{equation}
Here the integration path $\Gamma$ is the same as in Theorem
\ref{Th: exp dense}, and fixing $x_0$ is necessary for the definition
of the discrete exponential function on $\cD$. To make
(\ref{dLog}) a valid definition, one has to specify a branch
of $\log(\lambda)$ in a neighborhood of each point $\pm\theta_k$. This
choice depends on $x$, and is done as follows.

Assume, without loss of generality, that the circular order of the
points $\pm\theta_k$ on the positively oriented unit circle
$\bbS^1$ is the following: $\theta_1$, $\ldots$, $\theta_d$,
$-\theta_1$, $\ldots$, $-\theta_d$. We set
$\theta_{k+d}=-\theta_k$ for $k=1,\ldots,d$, and then define
$\theta_r$ for all $r\in\bbZ$ by $2d$-periodicity. For each
$r\in\bbZ$, assign to $\theta_r=\exp(i\gamma_r)\in\bbS^1$ a
certain value of argument $\gamma_r\in\bbR$: choose a value
$\gamma_1$ of the argument of $\theta_1$ arbitrarily, and then
extend it according to the rule
\[
\gamma_{r+1}-\gamma_r\in (0,\pi),\qquad \forall r\in\bbZ.
\]
Clearly, there holds $\gamma_{r+d}=\gamma_{r}+\pi$, and therefore
also $\gamma_{r+2d}=\gamma_{r}+2\pi$. It will be convenient to
consider the points $\theta_r$, supplied with the arguments
$\gamma_r$, as belonging to the Riemann surface $\Lambda$ of the
logarithmic function (a branched covering of the complex
$\lambda$-plane).

For each $m\in\bbZ$, define the ``sector'' $U_m$ on the embedding
plane $\bbC$ of the quad-graph $\cD$ as the set of all points of
$V(\cD)$ which can be reached from $x_0$ along paths with all
edges from $\{\theta_m,\ldots,\theta_{m+d-1}\}$. Two sectors
$U_{m_1}$ and $U_{m_2}$ have a non-empty intersection, if and only
if $|m_1-m_2|<d$. The union $U=\bigcup_{m=-\infty}^\infty U_m$ is
a branched covering of the quad-graph $\cD$, and serves as the
definition domain of the discrete logarithmic function.

The definition (\ref{dLog}) of the latter should be read as
follows: for $x\in U_m$, the poles of $e(x;\lambda)$ are exactly
the points $\theta_m,\ldots,\theta_{m+d-1}\in\Lambda$. The
integration path $\Gamma$ consists of $d$ small loops on $\Lambda$
around these points, and $\arg(\lambda)=\Im\log(\lambda)$ takes
values in a small open neighborhood (in $\bbR$) of the interval
\begin{equation}\label{ineq S}
[\gamma_m,\gamma_{m+d-1}]
\end{equation}
of length less than $\pi$. If $m$ increases by $2d$, the interval
(\ref{ineq S}) is shifted by $2\pi$. As a consequence, the
function $\ell$ is discrete holomorphic, and its restriction to
the ``black'' points $V(\cG)$ is discrete harmonic everywhere on
$U$ except at the point $x_0$:
\begin{equation}\label{Green 1}
\Delta \ell(x)=\delta_{x_0x}.
\end{equation}

Thus, the functions $g_k$ in the integral representation
(\ref{dense}) of an arbitrary discrete holomorphic function,
defined originally in disjoint neighborhoods of the points
$\alpha_r$, in the case of the discrete logarithmic function are
actually restrictions of a single analytic function
$\log(\lambda)/(2\lambda)$ to these neighborhoods. This allows one
to deform the integration path $\Gamma$ into a connected contour
lying on a single leaf of the Riemann surface of the logarithm,
and then to use standard methods of the complex analysis in order
to obtain asymptotic expressions for the discrete logarithmic
function. In particular, one can show \cite{K} that at the
``black'' points $V(\cG)$ there holds:
\begin{equation}\label{Green 2}
\ell(x)\sim \log|x-x_0|,\qquad x\to\infty.
\end{equation}
Properties (\ref{Green 1}), (\ref{Green 2}) characterize
the {\itbf discrete Green's function} on $\cG$. Thus:
\begin{thm}
The discrete logarithmic function on $\cD$, restricted to
the set of vertices $V(\cG)$ of the ``black'' graph $\cG$, coincides
with discrete Green's function on $\cG$.
\end{thm}

Now we extend the discrete logarithmic function to $\bbZ^d$, which
will allow us to gain significant additional information about it
\cite{BMS}. Introduce, in addition to the unit vectors
$\be_k\in\bbZ^d$ (corresponding to $\theta_k\in\bbS^1$), their
opposites $\be_{k+d}=-\be_k$, $k\in[1,d]$ (corresponding to
$\theta_{k+d}=-\theta_k$), and then define $\be_r$ for all
$r\in\bbZ$ by $2d$-periodicity. Then
\begin{equation}\label{sector}
  S_m=\bigoplus_{r=m}^{m+d-1}\bbZ\be_{r}\subset\bbZ^d
\end{equation}
is a $d$-dimensional octant containing exactly the part of
$\Omega_\cD$ which is the $P$-image of the sector $U_m\subset\cD$.
Clearly, only $2d$ different octants appear among $S_m$ (out of
$2^d$ possible $d$-dimensional octants). Define $\widetilde S_m$
as the octant $S_m$ equipped with the interval (\ref{ineq S}) of
values for $\Im\log(\theta_r)$. By definition, $\widetilde
S_{m_1}$ and $\widetilde S_{m_2}$ intersect, if the underlying
octants $S_{m_1}$ and $S_{m_2}$ have a non-empty intersection
spanned by the common coordinate semi-axes $\bbZ\be_r$, and
$\Im\log(\theta_r)$ for these common semi-axes match. It is easy
to see that $\widetilde S_{m_1}$ and $\widetilde S_{m_2}$
intersect, if and only if $|m_1-m_2|<d$. The union $\widetilde
S=\bigcup_{m=-\infty}^\infty \widetilde S_m$ is a branched
covering of the set $\bigcup_{m=1}^{2d}S_m\subset\bbZ^d$.
\begin{dfn}
The {\itbf discrete logarithmic function} on $\widetilde S$ is
given by the formula
\begin{equation}\label{Green f prelim}
\ell(\bn)=\frac{1}{2\pi i}
\int_\Gamma\frac{\log(\lambda)}{2\lambda}\,e(\bn;\lambda)d\lambda,\qquad
\forall \bn\in\widetilde S,
\end{equation}
where for $\bn\in\widetilde S_m$ the integration path $\Gamma$ consists
of $d$ loops around $\theta_m$, $\ldots$, $\theta_{m+d-1}$ on $\Lambda$,
and $\Im\log(\lambda)$ on $\Gamma$ is chosen in a small open neighborhood
of the interval (\ref{ineq S}).
\end{dfn}
The discrete logarithmic function on $\cD$ can be described as the
restriction of the discrete logarithmic function on $\widetilde S$
to a branched covering of $\Omega_\cD\sim\cD$. This holds for an
{\em arbitrary} quasicrystallic quad-graph $\cD$ whose set of edge
slopes coincides with $\Theta=\{\pm\theta_1,\ldots,\pm\theta_d\}$.

Now we are in a position to give an alternative definition of the
discrete logarithmic function. Clearly, it is completely
characterized by its values $\ell(n\be_r)$, $r\in[m,m+d-1]$, on the
coordinate semi-axes of an arbitrary octant $\widetilde S_m$. Let us
stress once more that the points $n\be_r$ do not lie, in general, on
the original quad-surface $\Omega_\cD$.
\begin{thm}\label{prop initial values}
The values $\ell_n^{(r)}=\ell(n\be_r)$, $r\in[m,m+d-1]$, of the discrete
logarithmic function on $\widetilde S_m\subset\widetilde S$ are given by:
\begin{equation}\label{Green axes}
\ell_{n}^{(r)}=\left\{\begin{array}{ll}
2\Big(1+\frac{1}{3}+\ldots+\frac{1}{n-1}\Big), & n\;\;{\rm even},\\
\log(\theta_r)=i\gamma_r, & n\;\;{\rm odd}.\end{array}\right.
\end{equation}
Here the values $\log(\theta_r)=i\gamma_r$ are chosen in the interval
(\ref{ineq S}).
\end{thm}
{\bf Proof.} Comparing formula (\ref{Green f prelim}) with (\ref{dense g_i}),
we see that the values $\ell_n^{(r)}$ can be obtained from the expansion of
$\log(\lambda)$ in a neighbourhood of $\lambda=\theta_r$ into the power
series with respect to the powers of $(\lambda-\theta_r)/(\lambda+\theta_r)$.
This expansion reads:
\[
\log(\lambda)=\log(\theta_r)+\sum_{n=1}^\infty\frac{1-(-1)^n}{n}\Big(
\frac{\lambda-\theta_r}{\lambda+\theta_r}\Big)^n.
\]
Thus, we come to a simple difference equation
\begin{equation}\label{CR axis recur}
 n(\ell_{n+1}^{(r)}-\ell_{n-1}^{(r)})=1-(-1)^n,
\end{equation}
with the initial conditions
\begin{equation}\label{init eps}
 \ell_0^{(r)}=\ell(\boldsymbol 0)=0,\qquad
 \ell_1^{(r)}=\ell(\be_r)=\log(\theta_r),
\end{equation}
which yields eq. (\ref{Green axes}). $\Box$
\smallskip

Observe that values (\ref{Green axes}) at even (resp. odd) points
imitate the behaviour of the real (resp. imaginary) part of the
function $\log(\lambda)$ along the semi-lines ${\rm
arg}(\lambda)={\rm arg}(\theta_r)$. This can be easily
extended to the whole of $\widetilde S$. Restricted to black
points $\bn\in\widetilde S$ (those with $n_1+\ldots+n_d$ even),
the discrete logarithmic function models the real part of the
logarithm. In particular, it is real-valued
and does not branch: its values on $\widetilde S_m$ depend on
$m\pmod{2d}$ only. In other words, it is a well defined function
on $S_m$. On the contrary, the discrete logarithmic function
restricted to white points $\bn\in\widetilde S$ (those with
$n_1+\ldots+n_d$ odd) takes purely imaginary values, and increases
by $2\pi i$, as $m$ increases by $2d$. Hence, this restricted
function models the imaginary part of the logarithm.

It turns out that recurrent relations (\ref{CR axis recur}) are
characteristic for an important class of solutions of the discrete
Cauchy-Riemann equations, namely for the {\em isomonodromic} ones.
In order to introduce this class, recall that discrete holomorphic
functions in $\bbZ^d$ possess a zero curvature representation with
the transition matrices (\ref{L CR d}). The {\em moving frame}
$\Psi(\cdot,\lambda):\bbZ^d\to GL_2(\bbC)[\lambda]$ is defined by
prescribing some $\Psi(\boldsymbol 0;\lambda)$, and by extending
it recurrently according to the formula
\begin{equation}\label{Psi recur}
  \Psi(\bn+\be_k;\lambda)=L_k(\bn;\lambda)\Psi(\bn;\lambda).
\end{equation}
Finally, define the matrices $A(\cdot;\lambda):\bbZ^d\to
GL_2(\bbC)[\lambda]$ by
\begin{equation}\label{A}
  A(\bn;\lambda)=\frac{d\Psi(\bn;\lambda)}{d\lambda}\Psi^{-1}(\bn;\lambda).
\end{equation}
These matrices satisfy a recurrent relation, which results by
differentiating (\ref{Psi recur}),
\begin{equation}\label{A recur}
  A(\bn+\be_k;\lambda)=
  \frac{dL_k(\bn;\lambda)}{d\lambda}L_k^{-1}(\bn;\lambda)+
  L_k(\bn;\lambda)A(\bn;\lambda)L_k^{-1}(\bn;\lambda),
\end{equation}
and therefore they are defined uniquely after fixing some
$A(\boldsymbol 0;\lambda)$.
\begin{dfn}\label{Dfn: isomonodrom}
A discrete holomorphic function $f:\bbZ^d\to\bbC$ is called {\itbf
isomonodromic}, if, for some choice of $A(\boldsymbol 0;\lambda)$,
the matrices $A(\bn;\lambda)$ are meromorphic in $\lambda$, with poles
whose positions and orders do not depend on $\bn\in\bbZ^d$.
\end{dfn}
This term originates in the theory of integrable nonlinear
differential equations, where it is used for solutions with a
similar analytic characterization \cite{IN}.

It is clear how to extend Definition \ref{Dfn: isomonodrom} to functions
on the covering $\widetilde S$. In the following statement, we restrict
ourselves to the octant $S_1=(\bbZ_+)^d$ for notational simplicity.
\begin{thm}\label{prop CR isomonodromic}
The discrete logarithmic function is isomonodromic: for a proper choice
of $A(\boldsymbol 0;\lambda)$, the matrices $A(\bn;\lambda)$ at any point
$\bn\in(\bbZ_+)^d$ have simple poles only:
\begin{equation}\label{A CR simple poles}
  A(\bn;\lambda)=\frac{A^{(0)}(\bn)}{\lambda}+
  \sum_{l=1}^d\Big(\frac{B^{(l)}(\bn)}{\lambda+\theta_l}+
  \frac{C^{(l)}(\bn)}{\lambda-\theta_l}\Big),
\end{equation}
with
\begin{eqnarray}
  A^{(0)}(\bn) & = &
 \renewcommand{\arraystretch}{1.4}
 \begin{pmatrix}
 0 & (-1)^{n_1+\ldots+n_d} \\ 0 & 0 \end{pmatrix},
\label{A CR} \\ \nonumber\\
 B^{(l)}(\bn) & = &
 \renewcommand{\arraystretch}{1.4}
  n_l\begin{pmatrix} 1 & -(\ell(\bn)+\ell(\bn-\be_l)) \\
  0 & 0 \end{pmatrix},\label{B CR} \\ \nonumber\\
  C^{(l)}(\bn) & = &
  \renewcommand{\arraystretch}{1.4}
  n_l\begin{pmatrix} 0 & \ell(\bn+\be_l)+\ell(\bn) \\
  0 & 1 \end{pmatrix}.  \label{C CR}
\end{eqnarray}
At any point $\bn\in\widetilde S$, the following constraint holds:
\begin{equation}\label{CR constr}
\sum_{l=1}^d n_l\Big(\ell(\bn+\be_l)-\ell(\bn-\be_l)\Big)=
1-(-1)^{n_1+\ldots+n_d}.
\end{equation}
\end{thm}
{\bf Proof.} The proper choice of $A(\boldsymbol 0;\lambda)$
mentioned in the Theorem, can be read off formula (\ref{A CR}):
$A(\boldsymbol 0;\lambda)=\frac{1}{\lambda}\begin{pmatrix} 0 & 1
\\ 0 & 0 \end{pmatrix}$. The proof consists of two parts.
\begin{itemize}
 \item[(i)] First, one proves the claim for the points of
the coordinate semi-axes. For any $r=1,\ldots,d$, construct the
matrices $A(n\be_r;\lambda)$ along the $r$-th coordinate semi-axis
via formula (\ref{A recur}) with the transition matrices (\ref{L
CR d}). This formula shows that singularities of
$A(n\be_r;\lambda)$ are poles at $\lambda=0$ and at
$\lambda=\pm\theta_r$, and that the pole $\lambda=0$ remains
simple for all $n>0$. By a direct computation and using
mathematical induction, one shows that it is exactly the recurrent
relation (\ref{CR axis recur}) for $f_n^{(r)}=f(n\be_r)$ that
assures the poles $\lambda=\pm\theta_r$ to remain simple for all
$n>0$. Thus, eq. (\ref{A CR simple poles}) holds on the $r$-th
coordinate semi-axis, with $B^{(l)}(n\be_r)=C^{(l)}(n\be_r)=0$ for
$l\neq r$.
 \item[(ii)] The second part of the proof is conceptual, and is based
upon the multidimensional consistency only. Proceed by induction,
whose scheme follows filling out the hull of the coordinate
semi-axes: each new point is of the form $\bn+\be_j+\be_k$, $j\neq
k$, with three points $\bn$, $\bn+\be_j$ and $\bn+\be_k$ known
from the previous steps, where the statements of the proposition
are assumed to hold. Suppose that (\ref{A CR simple poles}) holds
at $\bn+\be_j$, $\bn+\be_k$. The new matrix
$A(\bn+\be_j+\be_k;\lambda)$ is obtained by two alternative
formulas,
\begin{eqnarray}\label{prop CR jk}
\lefteqn{A(\bn+\be_j+\be_k;\lambda)
 = \frac{dL_k(\bn+\be_j;\lambda)}{d\lambda}\,
 L_k^{-1}(\bn+\be_j;\lambda)}\nonumber\\
 && +L_k(\bn+\be_j;\lambda)A(\bn+\be_j;\lambda)
     L_k^{-1}(\bn+\be_j;\lambda),
\end{eqnarray}
and the one with interchanged roles of $k$ and $j$. Eq. (\ref{prop
CR jk}) shows that all poles of $A(\bn+\be_j+\be_k;\lambda)$
remain simple, with the possible exception of
$\lambda=\pm\theta_k$, whose orders might increase by 1. The same
statement holds with $k$ replaced by $j$. Therefore, {\em all}
poles remain simple, and (\ref{A CR simple poles}) holds at
$\bn+\be_j+\be_k$. Formulas (\ref{A CR})--(\ref{C CR}) and
constraint (\ref{CR constr}) follow by direct computations based
on eq. (\ref{prop CR jk}). $\Box$
\end{itemize}

\chapter[Integrable circle patterns]{Discrete complex analysis.
Integrable circle patterns} \label{Chap:_discr_analysis-nonlinear}

\section{Circle patterns}
\label{Sect: patterns intro}

The idea that circle packings and, more generally, circle patterns
serve as a discrete counterpart of analytic functions is by now
well established, see a popular review \cite{St1} and a monograph
\cite{St2}. The origin of this idea is connected with Thurston's
approach to the Riemann mapping theorem via circle packings, see
\cite{T1}. Since then the theory bifurcated to several areas.

One of them is mainly dealing with approximation problems, and in
this context it is advantageous to stick from the beginning with
some fixed regular combinatorics. The most popular are hexagonal
packings, for which the $C^{\infty}$ convergence to the Riemann
mapping has been established \cite{RS, HS}. Similar results are
available also for circle patterns with the combinatorics of the
square grid introduced by Schramm \cite{S}.

Another area concentrates around the uniformization theorem of
Koebe-Andreev-Thurston, and is dealing with circle packing
realizations of cell complexes of prescribed combinatorics,
rigidity properties, constructing hyperbolic 3-manifolds, etc.
\cite{T2, MR, BeaS, He}. A new variational approach to this area
is given in \cite{BSp}, where also a review of related results can
be found. An important differential-geometric application of this
approach is the construction of discrete minimal surfaces through
circle patterns \cite{BHSp}.

The Schramm circle patterns with the combinatorics of the square
grid constitute also a prominent example in the area we are mainly
interested here, which deals with the interrelations of circle
patterns with integrable systems. We give here a presentation of
several results in this area, based on \cite{BP3, BHo, AB1}.

\begin{dfn}\label{Dfn patterns}
Let $\cG$ be an arbitrary cell decomposition of an open or closed
disk in $\bbC$. A map $z:V(\cG)\mapsto\bbC$ defines a {\itbf
circle pattern with the combinatorics of} $\cG$, if the following
condition is satisfied. Let $y\in F(\cG)\sim V(\cG^*)$ be an
arbitrary face of $\cG$, and let $x_1,x_2,\ldots,x_n$ be its
consecutive vertices. Then the points
$z(x_1),z(x_2),\ldots,z(x_n)\in\bbC$ lie on a circle, and their
circular order is just the listed one. We denote this circle by
$C(y)$, thus putting it into a correspondence with the face $y$,
or, equivalently, with the respective vertex of the dual cell
decomposition $\cG^*$.
\end{dfn}
As a consequence of this condition, if two faces $y_0,y_1\in
F(\cG)$ have a common edge $(x_0,x_1)$, then the circles $C(y_0)$
and $C(y_1)$ intersect in the points $z(x_1),z(x_2)$. In other
words, the edges from $E(\cG)$ correspond to pairs of neighboring
(intersecting) circles of the pattern. Similarly, if several faces
$y_1,y_2,\ldots,y_m\in F(\cG)$ meet in one point $x_0\in V(\cG)$,
then the corresponding circles $C(y_1),C(y_2),\ldots,C(y_m)$ also
have a common intersection point $z(x_0)$.

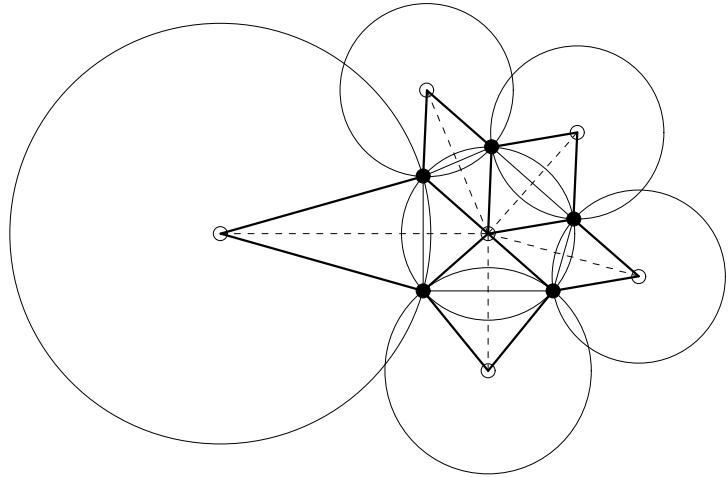
\begin{figure}[htbp]
\begin{center}
\setlength{\unitlength}{0.025em}
\begin{picture}(1000,800)(-700,-400)
 \put(0,0){\circle{20}}
 \put(0,0){\circle{240}}
 \put(5,119.8958){\circle*{20}}
 \put(118.3216,20){\circle*{20}}
 \put(90,-79.3725){\circle*{20}}
 \put(-90,-79.3725){\circle*{20}}
 \put(-90,79.3725){\circle*{20}}
 \put(0,-190){\circle{20}}
 \put(0,-190){\circle{285.2258}}
 \put(208.3216,-59.3725){\circle{20}}
 \put(208.3216,-59.3725){\circle{240}}
 \put(-370,0){\circle{20}}
 \put(-370,0){\circle{582.06527}}
 \put(-85,198.2683){\circle{20}}
 \put(-85,198.2683){\circle{240}}
 \put(123.3216,139.8958){\circle{20}}
 \put(123.3216,139.8958){\circle{240}}
 \path(5,119.8958)(118.3216,20)
 \path(118.3216,20)(90,-79.3725)
 \path(90,-79.3725)(-90,-79.3725)
 \path(-90,-79.3725)(-90,79.3725)
 \path(-90,79.3725)(5,119.8958)
 \dashline[+30]{10}(0,0)(0,-190)
 \dashline[+30]{10}(0,0)(208.3216,-59.3725)
 \dashline[+30]{10}(0,0)(-370,0)
 \dashline[+30]{10}(0,0)(-85,198.2683)
 \dashline[+30]{10}(0,0)(123.3216,139.8958)
 \thicklines
 \path(-370,0)(-90,-79.3725)  \path(-370,0)(-90,79.3725)
 \path(0,0)(-90,-79.3725)  \path(0,0)(-90,79.3725)
 \path(-85,198.2683)(-90,79.3725)  \path(-85,198.2683)(5,119.8958)
 \path(0,0)(5,119.8958)
 \path(123.3216,139.8958)(5,119.8958)  \path(123.3216,139.8958)(118.3216,20)
 \path(0,0)(118.3216,20)
 \path(208.3216,-59.3725)(118.3216,20) \path(208.3216,-59.3725)(90,-79.3725)
 \path(0,0)(90,-79.3725)
 \path(0,-190)(90,-79.3725)  \path(0,-190)(-90,-79.3725)
    \end{picture}
\caption{Circle pattern}\label{circle pattern}
\end{center}
\end{figure}

Given a circle pattern with the combinatorics of $\cG$, we can extend the
function $z$ to the vertices of the dual graph, setting
\[
z(y)={\rm \;center\;\;of\;\;the\;\;circle}\;\;C(y),
\quad y\in F(\cG)\simeq V(\cG^*).
\]
After this extension, the map $z$ is defined on all of
$V(\cD)=V(\cG)\sqcup V(\cG^*)$, where $\cD$ is the double of
$\cG$. Consider a face of the double. Its $z$-image is a
quadrilateral of the {\em kite} form, whose vertices correspond to
the intersection points and the centers of two neighboring circles
$C_0,C_1$ of the pattern. Denote the radii of $C_0,C_1$ by
$r_0,r_1$, respectively. Let $x_0,x_1$ correspond to the
intersection points, and let $y_0,y_1$ correspond to the centers
of the circles. Give the circles $C_0,C_1$ a positive orientation
(induced by the orientation of the underlying $\bbC$), and let
$\phi\in(0,\pi)$ stand for the intersection angle of these
oriented circles. This angle $\phi$ is equal to the kite angles at
the ``black'' vertices $z(x_0),z(x_1)$, see Fig.~\ref{two
circles}, where the complementary angle $\phi^*=\pi-\phi$ is also
shown. It will be convenient to assign the intersection angle
$\phi=\phi(e)$ to the ``black'' edge $e=(x_0,x_1)\in E(\cG)$, and
to assign the complementary angle $\phi^*=\phi(e^*)$ to the dual
``white'' edge $e^*=(y_0,y_1)\in E(\cG^*)$. Thus, the function
$\phi:E(\cG)\sqcup E(\cG^*)\to(0,\pi)$ satisfies eq.
(\ref{eq:phi*}).
\begin{figure}[htbp]
\begin{center}
\setlength{\unitlength}{0.04em}
\begin{picture}(500,340)(-300,-160)
 \put(-150,0){\circle{10}}
 \put(-170,-30){$z(y_0)$}
 \put(50,0){\circle{10}}
 \put(45,-25){$z(y_1)$}
 \put(0,50){\circle*{10}}
 \put(-67,55){$z(x_0)$}
 \put(0,-50){\circle*{10}}
 \put(-8,-90){$z(x_1)$}
 \put(50,0){\circle{141.42136}}
 \put(130,0){$C_1$}
 \put(-150,0){\circle{316.22777}}
 \put(-295,0){$C_0$}
 \path(0,50)(0,-50)
 \dashline[+30]{9}(-150,0)(-110,0)
 \dashline[+30]{10}(-85,0)(50,0)
 \dottedline{6}(0,50)(50,100)
 \dottedline{6}(0,50)(-25,125)
 \thicklines
 \path(-150,0)(0,50)  \path(0,50)(50,0)
 \path(50,0)(0,-50)  \path(0,-50)(-150,0)
  \put(1,75){$\phi^*$}
  \put(-105,-5){$\psi_{01}$}
    \end{picture}
\caption{Two intersecting circles}\label{two circles}
\end{center}
\end{figure}
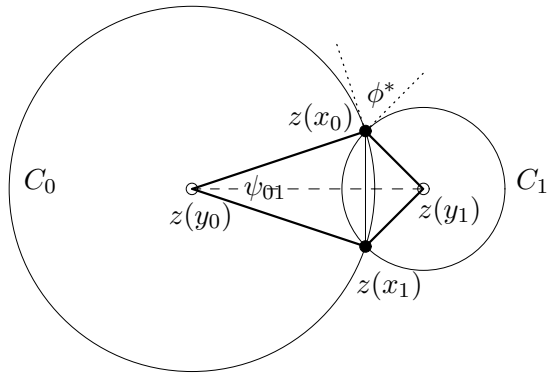

The geometry of Fig.~\ref{two circles} yields following relations.
First of all, the cross-ratio of the four points corresponding to
the vertices of a quadrilateral face of $\cD$ is expressed through
the intersection angle of the circles $C_0,C_1$:
\begin{equation}\label{pat cr}
q(z(x_0),z(y_0),z(x_1),z(y_1))=\exp(2i\phi^*).
\end{equation}
Further, running around a ``black'' vertex of $\cD$ (a common
intersection point of several circles of the pattern), we see that
the sum of the consecutive kite angles vanishes $\pmod {2\pi}$,
hence:
\begin{equation}\label{pat angles int points}
  \prod_{e\in{\rm star}(x_0;\cG)}\exp(i\phi(e))=1,
  \quad \forall x_0\in V(\cG).
\end{equation}
Finally, let $\psi_{01}$ be the angle of the kite
$(z(x_0),z(y_0),z(x_1),z(y_1))$ at the ``white'' vertex $z(y_0)$,
i.e., the angle between the half-lines from the center $z(y_0)$ of
the circle $C_0$ to the intersection points $z(x_0),z(x_1)$ with
its circle $C_1$. It is not difficult to calculate this angle:
\begin{equation}\label{pat psi}
\exp(i\psi_{01})=\frac{r_0+r_1\exp(i\phi^*)}{r_0+r_1\exp(-i\phi^*)}.
\end{equation}
Running around the ``white'' vertex of $\cD$, we come to the
relation
\begin{equation}\label{pat prodpsi=1}
\prod_{j=1}^m\frac{r_0+r_j\exp(i\phi_j^*)}
{r_0+r_j\exp(-i\phi_j^*)}=1,\quad \forall y_0\in V(\cG^*),
\end{equation}
where the product is extended over all edges
$e_j^*=(y_0,y_j)\in{\rm star}(y_0;\cG^*)$, and
$\phi_j^*=\phi(e_j^*)$, while $r_j$ are the radii of the circles
$C_j=C(y_j)$. This formula has been used in \cite{BSp} as the
basis of a variational proof of the existence and uniqueness of
Delaunay circle patterns with prescribed combinatorics and
intersection angles.

\section{Integrable cross-ratio and Hirota systems}
\label{Sect: q int}

Our main interest is in the circle patterns with prescribed
combinatorics and with prescribed intersection angles for all
pairs of neighboring angles. (This is not the only class of
patterns deserving a study from the point of view of
integrability, see, e.g., \cite{BHS} for a different integrable
class.) According to eq. (\ref{pat cr}), prescribing all
intersection angles amounts to prescribing cross-ratios for all
quadrilateral faces of the quad-graph $\cD$. Thus, we come to the
study of cross-ratio equations on arbitrary quad-graphs.

Let there be given a function $Q:E(\cG)\sqcup E(\cG^*)\to\bbC$
satisfying the condition
\begin{equation}\label{eq: Q*}
Q(e^*)=1/Q(e),\qquad\forall e\in E(\cG).
\end{equation}
\begin{dfn}
The {\itbf cross-ratio system} on $\cD$ corresponding to the
function $Q$, consists of the following equations for a function
$z:V(\cD)\to\bbC$, one for any quadrilateral face
$(x_0,y_0,x_1,y_1)$ of $\cD$:
\begin{equation}\label{cross-rat eq}
q(z(x_0),z(y_0),z(x_1),z(y_1))=Q(x_0,x_1)=1/Q(y_0,y_1).
\end{equation}
\end{dfn}
An important distinction from the discrete Cauchy-Riemann
equations is that the cross-ratio equations actually do not feel
the orientation of quadrilaterals.

We have already encountered 3D consistent cross-ratio systems on
$\bbZ^d$ in Sect. \ref{Sect: 3d consist} (see eq. (\ref{eq:2d
cr})), in the version with labelled edges. A natural
generalization to the case of arbitrary quad-graphs is this:
\begin{figure}[htbp]
\setlength{\unitlength}{0.04em}
\begin{picture}(200,220)(-300,-110)
 \put(0,0){\circle*{10}}\put(200,0){\circle*{10}}
 \put(100,-75){\circle{10}} \put(100,75){\circle{10}}
 \path(4,3)(96,72)
 \path(4,-3)(96,-72)
 \path(196,3)(104,72)
 \path(196,-3)(104,-72)
 \put(-37,-7){$x_0$} \put(215,-7){$x_1$}
 \put(94,-100){$y_0$} \put(94,97){$y_1$}
 \put(25,55){$\alpha_1$} \put(25,-60){$\alpha_0$}
 \put(150,55){$\alpha_0$} \put(150,-60){$\alpha_1$}
\end{picture}
\caption{Quadrilateral, with a labelling of undirected
edges}\label{Fig: diamond undir}
\end{figure}
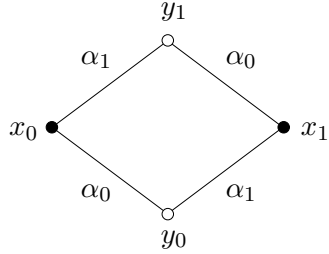
\begin{dfn}\label{Dfn integrable cr}
A cross-ratio system is called {\itbf integrable}, if there exists
a labelling $\alpha:E(\cD)\to\bbC$ of {\em undirected} edges of
$\cD$ such that the function $Q$ admits the following
factorization (in notations of Fig.\,\ref{Fig: diamond undir}):
\begin{equation}\label{eq: Q int}
Q(x_0,x_1)=\frac{1}{Q(y_0,y_1)}=\frac{\alpha_0}{\alpha_1}.
\end{equation}
\end{dfn}
Clearly, integrable cross-ratio systems are 3D consistent (see
Theorem \ref{Thm 3D cr}), admit {\em B\"acklund transformations},
and possess {\em zero curvature representation} with the
transition matrices (\ref{eq:2d cr L2}). It is not difficult to
give an equivalent re-formulation of the integrability condition
(\ref{eq: Q int}).
\begin{thm}\label{Thm 3D cr quad}
A cross-ratio system with the function $Q:E(\cG)\sqcup
E(\cG^*)\to\bbC$ is integrable, if and only if for all $x_0\in
V(\cG)$ and for all $y_0\in V(\cG^*)$ the following conditions are
fulfilled:
\begin{equation}\label{eq:cr prod=1}
 \prod_{e\in{\rm star}(x_0;\cG)} Q(e)=1,\quad
 \prod_{e^*\in{\rm star}(y_0;\cG^*)} Q(e^*)=1.
\end{equation}
\end{thm}
For a labelling of undirected edges $\alpha:E(\cD)\to\bbC$, there
can be found a labelling $\theta:\vec{E}(\cD)\to\bbC$ of directed
edges such that $\alpha=\theta^2$. The function $p:V(\cD)\to\bbC$
defined by $p(y)-p(x)=\theta(x,y)$ gives, according to eq.
(\ref{eq:cr prod=1}), a parallelogram realization (ramified
embedding) of the quad-graph $\cD$. The cross-ratio equations are
written as
\begin{equation}\label{cr paral}
q(z(x_0),z(y_0),z(x_1),z(y_1))= \frac{\theta_0^2}{\theta_1^2}=
q(p(x_0),p(y_0),p(x_1),p(y_1)),
\end{equation}
in other words, the cross-ratio of the vertices of the $f$-image
of any quadrilateral $(x_0,y_0,x_1,y_1)\in F(\cD)$ is equal to the
cross-ratio of the vertices of the corresponding parallelogram. In
particular, one always has the {\em trivial} solution $z(x)\equiv
p(x)$ for all $x\in V(\cD)$.

A very useful transformation of the cross-ratio system is given by
the following construction.
\begin{dfn}\label{Dfn Hirota oriented}
The {\itbf Hirota system} for a given labelling of directed edges
$\theta:\vec{E}(\cD)\to\bbC$ consists of the following equations
for the function $w:V(\cD)\to\bbC$, one for every quadrilateral
face $(x_0,y_0,x_1,y_1)\in F(\cD)$:
\begin{equation}\label{Hirota}
\theta_0w(x_0)w(y_0)+\theta_1w(y_0)w(x_1)-\theta_0w(x_1)w(y_1)-
\theta_1w(y_1)w(x_0)=0.
\end{equation}
\end{dfn}
Note that this Hirota system is different from eq. (\ref{eq:2d H})
of Sect. \ref{Sect: 3d consist}, in that the present version has
parameters defined on {\em directed} edges. In terms of the
parallelogram realization $p:V(\cD)\to\bbC$ of the quad-graph
$\cD$ corresponding to the labelling $\theta$, eq. (\ref{Hirota})
reads:
\begin{eqnarray}\label{Hirota in p}
\lefteqn{w(x_0)w(y_0)\big(p(y_0)-p(x_0)\big)+
         w(y_0)w(x_1)\big(p(x_1)-p(y_0)\big)+}\nonumber\\
      && w(x_1)w(y_1)\big(p(y_1)-p(x_1)\big)+
         w(y_1)w(x_0)\big(p(x_0)-p(y_1)\big)=0.\qquad
\end{eqnarray}
Obviously, a transformation $w\to cw$ on $V(\cG)$ and $w\to
c^{-1}w$ on $V(\cG^*)$ with a constant $c\in\bbC$, called
hereafter a black-white scaling, maps solutions of the Hirota
system into solutions. A relation between the cross-ratio and the
Hirota system is based on the following observation:
\begin{thm}
Let $w:V(\cD)\to\bbC$ be a solution of the Hirota system. Then the
relation
\begin{equation}\label{w vs z}
z(y)-z(x)=\theta(x,y) w(x)w(y)=w(x)w(y)\big(p(y)-p(x)\big)
\end{equation}
for all directed edges $(x,y)\in\vec{E}(\cD)$ correctly defines a
unique (up to an additive constant) function $z:V(\cD)\to\bbC$
which is a solution of the cross-ratio system (\ref{cr paral}).
Conversely, for any solution $z$ of the cross-ratio system
(\ref{cr paral}), relation (\ref{w vs z}) defines a function
$w:V(\cD)\to\bbC$ correctly and uniquely (up to a black-white
scaling); this function $w$ solves the Hirota system
(\ref{Hirota}).
\end{thm}
In particular, the trivial solution $z(x)=p(x)$ of the cross-ratio
system corresponds to the trivial solution of the Hirota system,
$w(x)\equiv 1$ for all $x\in V(\cD)$. By a direct computation one
can establish the following fundamental property.
\begin{thm}
The Hirota system (\ref{Hirota}) is 3D consistent.
\end{thm}
As a usual consequence, the Hirota system admits B\"acklund
transformations and possesses zero curvature representation with
transition matrices along the edge $(x,y)\in\vec{E}(\cD)$ given by
\begin{equation}\label{L Hirota}
\renewcommand{\arraystretch}{1.4}
L(y,x,\theta;\lambda)=\begin{pmatrix} 1 & -\theta w(y) \\
-\lambda\theta/w(x) & w(y)/w(x) \end{pmatrix}, \quad {\rm where}
\quad\theta=p(y)-p(x).
\end{equation}

\section{Integrable circle patterns}

Returning to circle patterns, let $\{z(x):x\in V(\cG)\}$ be
intersection points of the circles of a pattern, and let
$\{z(y):y\in V(\cG^*)\}$ be their centers. Due to eq. (\ref{pat
cr}), the function $z:V(\cD)\to\bbC$ satisfies a cross-ratio
system with $Q:E(\cG)\sqcup E(\cG^*)\to\bbS^1$ defined as
$Q(e)=\exp(2i\phi(e))$. Because of eq. (\ref{pat angles int
points}), the first one of the integrability conditions
(\ref{eq:cr prod=1}) is fulfilled for an {\em arbitrary} circle
pattern. Therefore, integrability of the cross-ratio system for
circle patterns with the prescribed intersection angles
$\phi:E(\cG^*)\to(0,\pi)$ is equivalent to:
\begin{equation}\label{CP angles centers 2}
  \prod_{e^*\in{\rm star}(y_0;\cG^*)}\exp(2i\phi(e^*))=1,
  \quad \forall y_0\in V(\cG^*),
\end{equation}
This is equivalent to the existence of the edge labelling
$\alpha:E(\cD)\to\bbC$ such that, in notations of Fig.~\ref{two
circles},
\begin{equation}\label{eq: int pat fact 2}
\exp(2i\phi^*)=\frac{\alpha_0}{\alpha_1}.
\end{equation}
Moreover, one can assume that the labelling $\alpha$ takes values
in $\bbS^1$.

Our definition of integrable circle patterns will require somewhat
more than integrability of the corresponding cross-ratio system.
\begin{dfn}\label{dfn integrable pattern}
A circle pattern with the prescribed intersection angles
$\phi:E(\cG^*)\to(0,\pi)$ is called {\itbf integrable}, if
\begin{equation}\label{CP angles centers}
  \prod_{e^*\in{\rm star}(y_0;\cG^*)}\exp(i\phi(e^*))=1,
  \quad \forall y_0\in V(\cG^*),
\end{equation}
i.e., if for any circle of the pattern the sum of its intersection
angles with all neighboring circles vanishes $\pmod{2\pi}$.
\end{dfn}
This requirement is equivalent to a somewhat sharper factorization
than (\ref{eq: int pat fact 2}), namely, to the existence of a
labelling of directed edges $\theta:\vec E(\cD)\to\bbS^1$ such
that, in notations of Fig.~\ref{two circles},
\begin{equation}\label{eq: int pat fact 1}
\exp(i\phi)=\frac{\theta_1}{\theta_0} \quad\Leftrightarrow\quad
\exp(i\phi^*)=-\frac{\theta_0}{\theta_1}.
\end{equation}
(Of course, the last condition yields (\ref{eq: int pat fact 2})
with $\alpha=\theta^2$.) The parallelogram realization
$p:V(\cD)\to\bbC$ corresponding to the labelling $\theta\in\bbS^1$
is actually a {\em rhombic} one.
\begin{thm}
Combinatorial data $\cG$ and intersection angles
$\phi:E(\cG)\to(0,\pi)$ belong to an integrable circle pattern, if
and only if they admit an isoradial realization. In this case, the
dual combinatorial data $\cG^*$ and intersection angles
$\phi:E(\cG^*)\to(0,\pi)$ admit a realization as an isoradial
circle pattern, as well.
\end{thm}
{\bf Proof.} The rhombic realization $p:V(\cD)\to\bbC$ of the
quad-graph $\cD$ corresponds to a circle pattern with the same
combinatorics and the same intersection angles as the original one
and with all radii equal to 1, and, simultaneously, to an
analogous dual circle pattern. $\Box$
\medskip

Consider a rhombic realization $p:V(\cD)\to\bbC$ of $\cD$.
Solutions $z:V(\cD)\to\bbC$ of the corresponding integrable
cross-ratio system which come from integrable circle patterns are
characterized by the property that the $z$-image of any
quadrilateral $(x_0,y_0,x_1,y_1)$ from $F(\cD)$ is a {\it kite}
with the prescribed angle $\phi$ at the black vertices $z(x_0)$,
$z(x_1)$ (cf. Fig.~\ref{two circles}). It turns out that the
description of this class of kite solutions admits a more
convenient analytic characterization in terms of the corresponding
solutions $w:V(\cD)\to\bbC$ of the Hirota system defined by eq.
(\ref{w vs z}).

\begin{thm}\label{Thm Hirota pat reduction}
The solution $z$ of the cross-ratio system corresponds to a circle
pattern, if and only if the solution $w$ of the Hirota system,
corresponding to $z$ via (\ref{w vs z}), satisfies the condition
\begin{equation}\label{Hirota red}
w(x)\in\bbS^1,\quad w(y)\in\bbR_+,\quad \forall x\in V(\cG),\;y\in
V(\cG^*).
\end{equation}
The values $w(y)\in\bbR_+$ have then the interpretation of the
radii of the circles $C(y)$, while the (arguments of the) values
$w(x)\in\bbS^1$ measure the rotation of the tangents to the
circles intersecting at $z(x)$ with respect to the isoradial
realization of the pattern.
\end{thm}
{\bf Proof.} As easily seen, the kite conditions are equivalent
to:
\[
\frac{|w(x_0)|}{|w(x_1)|}=1 \quad{\rm and}\quad
\frac{w(y_0)}{w(y_1)}\in\bbR_+.
\]
This yields (\ref{Hirota red}), possibly upon a black-white
scaling. $\Box$
\medskip

The conditions (\ref{Hirota red}) form an {\em admissible
reduction} of the Hirota system with $\theta\in\bbS^1$, in the
following sense: if any three of the four points $w(x_0)$,
$w(y_0)$, $w(x_1)$, $w(y_1)$ satisfy the condition (\ref{Hirota
red}), then so does the fourth one. This is immediately seen, if
one rewrites the Hirota equation (\ref{Hirota}) in one of the
equivalent forms:
\begin{equation}\label{Hirota alt}
 \frac{w(x_1)}{w(x_0)}=
 \frac{\theta_1w(y_1)-\theta_0w(y_0)}{\theta_1w(y_0)-\theta_0w(y_1)}
 \quad\Leftrightarrow\quad
 \frac{w(y_1)}{w(y_0)}=
 \frac{\theta_0w(x_0)+\theta_1w(x_1)}{\theta_0w(x_1)+\theta_1w(x_0)}.
\end{equation}
As a consequence of this remark, we obtain B\"acklund
transformations for integrable circle patterns.
\begin{thm}
Let all $\theta\in\bbS^1$, and let $p:V(\cD)\to\bbC$ be the
corresponding rhombic realization of $\cD$. Let the solution
$w:V(\cD)\to\bbC$ of the Hirota system correspond to a circle
pattern with the combinatorics of $\,\cG$, i.e., satisfy
(\ref{Hirota red}). Consider its B\"acklund transformation
$w^+:V(\cD)\to\bbC$ with an arbitrary parameter $\lambda\in\bbS^1$
and with an arbitrary initial value $w^+(x_0)\in\bbR_+$ or
$w^+(y_0)\in\bbS^1$. Then there holds:
\begin{equation}\label{Hirota red hat}
w^+(x)\in\bbR_+,\quad w^+(y)\in\bbS^1,\quad \forall x\in
V(\cG),\;y\in V(\cG^*),
\end{equation}
so that $w^+$ corresponds to a circle pattern with the
combinatorics of $\,\cG^*$, which we call a B\"acklund transform
of the original circle pattern.
\end{thm}

We close this section by mentioning several Laplace type equations
which can be used to describe integrable circle patterns. First of
all, the restriction of the function $z$ to $V(\cG)$ (i.e., the
intersection points of the circles) satisfies the equations
\[
\sum_{k=1}^{n}\frac{\alpha_k-\alpha_{k+1}}{z(x_k)-z(x_0)}=0.
\]
Here $z(x_0)$ is any intersection point, where $n$ circles
$C(y_1),\ldots,C(y_n)$ meet, $z(x_k)$ are the second intersection
points of $C(y_k)$ with $C(y_{k+1})$, and $\alpha_k$ are the
labels on the edges $(x_0,y_k)\in E(\cD)$. Analogously, the
restriction of the function $z$ to $V(\cG^*)$ (i.e., the centers
of the circles) satisfies the equation
\[
\sum_{j=1}^m\frac{\alpha_{j-1}-\alpha_j}{z(y_j)-z(y_0)}=0.
\]
Here $z(y_0)$ is the center of any circle $C(y_0)$, which
intersects with $m$ circles $C(y_1),\ldots,C(y_{m})$ with the
centers at the points $z(y_j)$; the intersection of $C(y_0)$ with
$C(y_j)$ consists of two points $z(x_{j-1})$, $z(x_j)$, and
$\alpha_j$ are the labels on the edges $(y_0,x_j)\in E(\cD)$.
These both Laplace type equations follow from the first claim of
Theorem \ref{Th: Laplace for 3legs}, applied to the cross-ratio
system, which is nothing but the case (Q1)$_{\delta=0}$ of Theorem
\ref{Th: 3leg forms}.

A similar construction can be applied to the Hirota system,
written in the three-leg form (\ref{Hirota alt}). Again, it yields
two multiplicative Laplace type equations -- on $\cG$ and on
$\cG^*$. It is instructive to look at the equation on $\cG^*$ (for
the radii $r_j=w(y_j)$ of the circles):
\[
\prod_{j=1}^m\frac{\theta_jr_j-\theta_{j-1}r_0}
{\theta_jr_0-\theta_{j-1}r_j}=1.
\]
Due to eq. (\ref{eq: int pat fact 1}), this equation can be
written in terms of the intersection angles $\phi_j$ of $C(y_0)$
with $C(y_j)$, and takes the form of eq. (\ref{pat prodpsi=1}).
Interestingly, the latter equation holds for {\em any} circle
pattern and is not specific for integrable ones (as opposed to the
similar Laplace type equation on $\cG$).

\section{$z^a$ and $\log z$ circle patterns}
\label{Sect z^a}

Due to the 3D consistency of the cross-ratio and the Hirota
systems, we can follow the procedure of Sect. \ref{Sect: discr
geom} and to extend solutions to this systems from a
quasicrystallic quad-graph $\cD$, realized as a quad-surface
$\Omega_\cD\subset\bbZ^d$, to the whole of $\bbZ^d$ (more
precisely, to the hull of $\Omega_\cD$). Then, one can ask about
isomonodromic solutions. As shown in \cite{AB1, BHo}, this leads to
discrete analogs of the power function. Naturally, these discrete
power functions are defined on the same branched covering
$\widetilde S$ of the set $\bigcup_{m=1}^{2d}S_m\subset\bbZ^d$ as
the discrete logarithmic function of Sect. \ref{Sect discr log}.

\begin{dfn}\label{Def z isomonodr}
For $a\in(0,1)$, the {\itbf discrete} {\boldmath $z^{2a}$} is the
solution of the discrete cross-ratio system on $\widetilde S$
defined by the values on the coordinate semi-axes
$z_n^{(r)}=z(n\be_r)$, $r\in[m,m+d-1]$, which solve the recurrent
relation
\begin{equation}\label{cr axis recur}
 n\,\frac{(z_{n+1}-z_n)(z_n-z_{n-1})}
 {z_{n+1}-z_{n-1}}=a z_n.
\end{equation}
with the initial conditions:
\begin{equation}\label{cr init}
z_0^{(r)}=z(\boldsymbol 0)=0,\qquad
z_1^{(r)}=z(\be_r)=\theta_r^{2a}=\exp(2a\log\theta_r),
\end{equation}
where $\log\theta_r$ is chosen in the interval (\ref{ineq S}).
\end{dfn}
By induction, one can derive the following explicit expressions
for the solutions $z_n^{(r)}$:
\begin{equation}
z_{2n}^{(r)}=\prod_{k=1}^{n-1}\frac{k+a}{k-a}\cdot\frac{n}{n-a}\cdot
\theta_r^{2a},\qquad
z_{2n+1}^{(r)}=\prod_{k=1}^{n}\frac{k+a}{k-a}\cdot \theta_r^{2a}.
\end{equation}
Observe the asymptotic relation
\begin{equation}
z_n^{(r)}=c(a)(n\theta_r)^{2a}\big(1+O(n^{-1})\big).
\end{equation}

The functions $z^{2a}$ correspond to (intersection points and
circle centers of) circle patterns. In order to prove this, it is
convenient to study the corresponding solutions of the Hirota
equation. Therefore, we introduce the functions $w^{2a-1}$ related
to $z^{2a}$ by
\[
z(\bn+\be_j)-z(\bn)=\theta_jw(\bn)w(\bn+\be_j).
\]
\begin{dfn}\label{Def w isomonodr}
For $a\in(0,1)$, the {\itbf discrete} {\boldmath $w^{2a-1}$} is
the solution of the discrete Hirota system on $\widetilde S$
defined by the values on the coordinate semi-axes
$w_n^{(r)}=w(n\be_r)$, $r\in[m,m+d-1]$, which solve the recurrent
relation
\begin{equation}\label{Hirota axis recur}
 n\,\frac{w_{n+1}-w_{n-1}}{w_{n+1}+w_{n-1}}=
 \big(a-\tfrac{1}{2}\big)\big(1-(-1)^n\big).
\end{equation}
with the initial conditions:
\begin{equation}\label{Hirota init}
w_0^{(r)}=w(\boldsymbol 0)=0,\qquad
w_1^{(r)}=w(\be_r)=\theta_r^{2a-1}=\exp((2a-1)\log\theta_r),
\end{equation}
where $\log\theta_r$ is chosen in the interval (\ref{ineq S}).
\end{dfn}
One can easily find a closed expression for $w_n^{(r)}$:
\begin{equation}\label{Hirota axes}
w_{2n}^{(r)}=\prod_{k=1}^n\frac{k-1+a}{k-a}\,,\qquad
w_{2n+1}^{(r)}=\theta_r^{2a-1}.
\end{equation}
Observe the asymptotics at $n\to\infty$,
\begin{equation}
w_{2n}^{(r)}=c(a)n^{2a-1}\big(1+O(n^{-1})\big).
\end{equation}
The main technical advantage of the $w$ variables is seen from the
following observation.
\begin{thm}
The function $w^{2a-1}$ takes values from $\bbR_+$ at the white
points and values from $\bbS^1$ at the black points. Therefore,
the function $z^{2a}$ defines a circle pattern.
\end{thm}
{\bf Proof.} The claim for $w^{2a-1}$ on the coordinate axes is
obvious from the explicit formulas (\ref{Hirota axes}), and can be
extended to the whole of $\widetilde S$ according to the remark
after Theorem \ref{Thm Hirota pat reduction}. The statement for
$z^{2a}$ is now a consequence of Theorem \ref{Thm Hirota pat
reduction}. $\Box$

\begin{figure}[htbp]
\begin{center}
\includegraphics[width=0.40\textwidth]{sqSchramm.eps}

\includegraphics[width=0.45\textwidth]{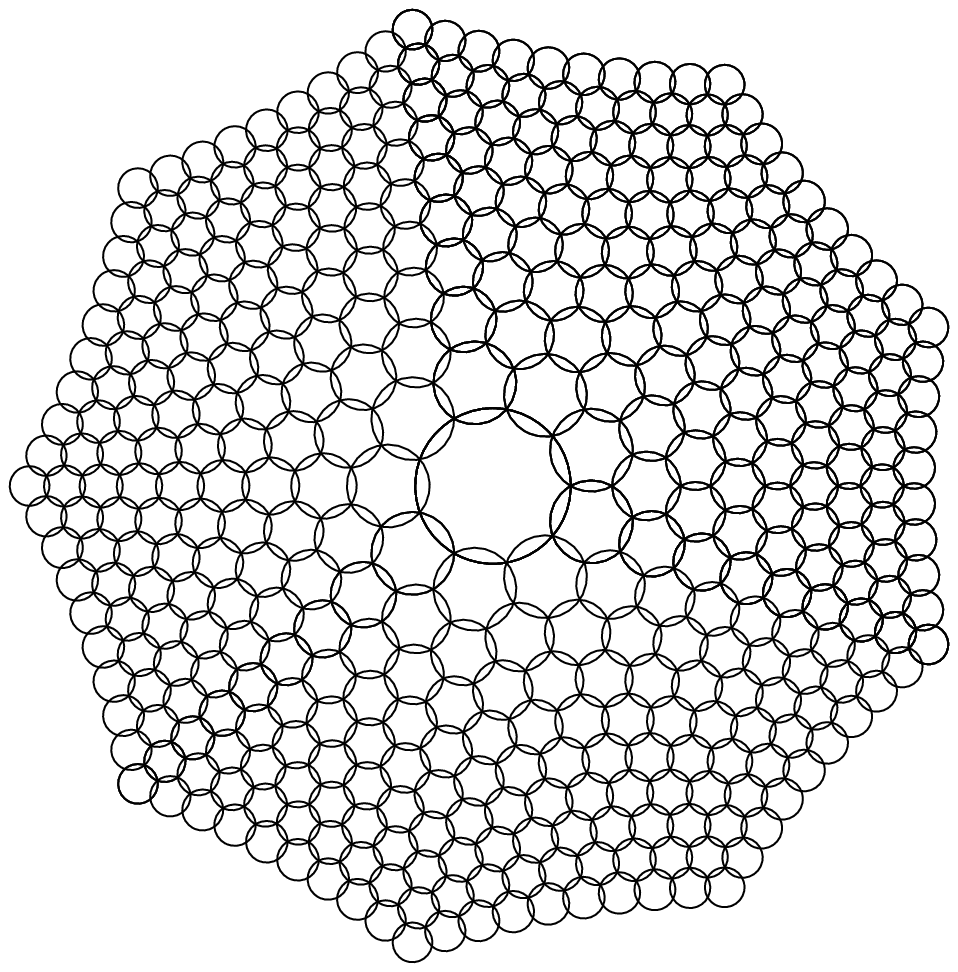}\hfill
\includegraphics[width=0.45\textwidth]{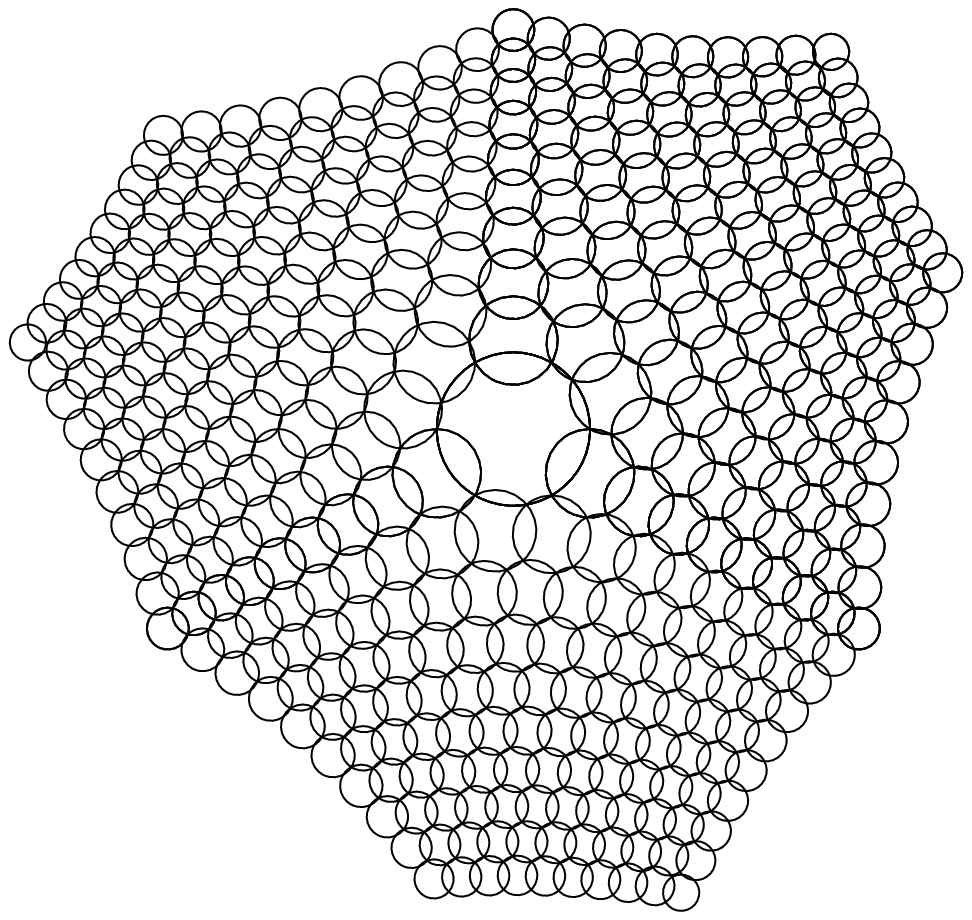}
\end{center}
\caption{Circle patterns $z^{4/5}$ with the combinatorics of the
square grid, and $z^{2/3}$ with the combinatorics of the regular
hexagonal lattice (isotropic and non-isotropic).}
 \label{fig: z^c}
\end{figure}

The restriction of $z^{2a}$ to various quad-surfaces $\Omega_\cD$
give the discrete analogs of the power function on the
corresponding quasicrystallic quad-graphs $\cD$ with the set
$\Theta=\{\pm\theta_1,\ldots,\pm\theta_d\}$ of edge slopes (see
Fig.~\ref{fig: z^c}).

These pictures lead to the conjecture \cite{BP3} that the circle
patterns $z^{2a}$ are embedded, i.e., interiors of different kites
are disjoint. This conjecture has been proven in \cite{Ag} for the
case of the square grid combinatorics. The fact that the circle
patterns $z^{2a}$ are immersed, i.e., the neighboring kites do not
overlap, was proven in \cite{AB1} for the square grid and in
\cite{AB2} for the hexagonal grid combinatorics. One possible
approach to the general case could be based on applying the well
developed techniques of the theory of isomonodromic solutions
\cite{IN} to the asymptotic study of the discrete $z^{2a}$,
because of the following statement.
\begin{thm}\label{Thm cr isomonodrom}
The discrete $z^{2a}$ is isomonodromic: for a proper choice of
$A(\boldsymbol 0;\lambda)$, the matrices $A(\bn;\lambda)$ at any
point $\bn\in(\bbZ_+)^d$ have simple poles only:
\begin{equation}\label{A cr simple poles}
  A(\bn;\lambda)=\frac{A^{(0)}(\bn)}{\lambda}+
  \sum_{l=1}^d\frac{B^{(l)}(\bn)}{\lambda-\theta_l^{-2}}\;,
\end{equation}
with
\begin{equation}\label{A cr}
 A^{(0)}(\bn) =
 \renewcommand{\arraystretch}{1.4}
 \begin{pmatrix}
 -a/2 & -az(\bn) \\ 0 & a/2 \end{pmatrix},
\end{equation}
\begin{eqnarray}\label{B cr}
 \lefteqn{B^{(l)}(\bn) =
  \dfrac{n_l}{z(\bn+\be_l)-z(\bn-\be_l)}}\nonumber\\
 && \times\renewcommand{\arraystretch}{1.4}
 \begin{pmatrix} z(\bn+\be_l)-z(\bn) &
  (z(\bn+\be_l)-z(\bn))(z(\bn)-z(\bn-\be_l))\\
  1 & z(\bn)-z(\bn-\be_l) \end{pmatrix}\!.\qquad
\end{eqnarray}
At any point $\bn\in\widetilde S$, the discrete $z^{2a}$ satisfies
the following constraint:
\begin{equation}\label{cr constr}
\sum_{j=1}^d
n_j\,\frac{(z(\bn+\be_j)-z(\bn))(z(\bn)-z(\bn-\be_j))}
{z(\bn+\be_j)-z(\bn-\be_j)}=a z(\bn).
\end{equation}
\end{thm}
{\bf Proof} follows the same scheme as the proof of Theorem
\ref{prop CR isomonodromic}: one first shows that the poles of
$A(n\be_r;\lambda)$ remain simple, due to the recurrent relations
(\ref{cr axis recur}), and then shows that the order of poles does
not increase at the points $\bn$ away from the coordinate axes,
due to the multidimensional consistency. Note that this scheme
works also for Theorem \ref{Thm Hirota isomonodrom}. $\Box$
\smallskip

The transition between $z$ and $w$ variables is a matter of
straightforward computations. Actually, the next theorem is
dealing with the same matrices as Theorem \ref{Thm cr isomonodrom}
but written in different variables.

\begin{thm}\label{Thm Hirota isomonodrom}
The discrete $w^{2a-1}$ is isomonodromic: for a proper choice of
$A(\boldsymbol 0;\lambda)$, the matrices $A(\bn;\lambda)$ at any
point $\bn\in(\bbZ_+)^d$ have simple poles only:
\begin{equation}\label{A Hirota simple poles}
  A(\bn;\lambda)=\frac{A^{(0)}(\bn)}{\lambda}+
  \sum_{l=1}^d\frac{B^{(l)}(\bn)}{\lambda-\theta_l^{-2}}\;,
\end{equation}
with
\begin{equation}\label{A Hirota}
 A^{(0)}(\bn) =
 \renewcommand{\arraystretch}{1.4}
 \begin{pmatrix}
 -a/2 & * \\ 0 & a/2 \end{pmatrix},
\end{equation}
\begin{equation}\label{B Hirota}
 B^{(l)}(\bn) =
 \renewcommand{\arraystretch}{1.4}
  \dfrac{n_l}{w(\bn+\be_l)+w(\bn-\be_l)}
  \begin{pmatrix} w(\bn+\be_l) & \theta_l w(\bn+\be_l)w(\bn-\be_l)\\
  1/\theta_l & w(\bn-\be_l) \end{pmatrix}.
\end{equation}
The upper right entry of the matrix $A^{(0)}(\bn)$, denoted by the
asterisk in (\ref{A Hirota}), is given by $
A_{12}^{(0)}(\bn)=-\sum_{l=1}^d B_{12}^{(l)}(\bn)$. At any point
$\bn\in\widetilde S$, the discrete $w^{2a-1}$ satisfies the
following constraint:
\begin{equation}\label{Hirota constr}
\sum_{l=1}^dn_l\,\frac{w(\bn+\be_l)-w(\bn-\be_l)}
 {w(\bn+\be_l)+w(\bn-\be_l)}=
 \big(a-\tfrac{1}{2}\big)\big(1-(-1)^{n_1+\ldots+n_d}\Big).
\end{equation}
\end{thm}

\begin{figure}[htbp]
\begin{center}
\includegraphics[width=0.48\textwidth]{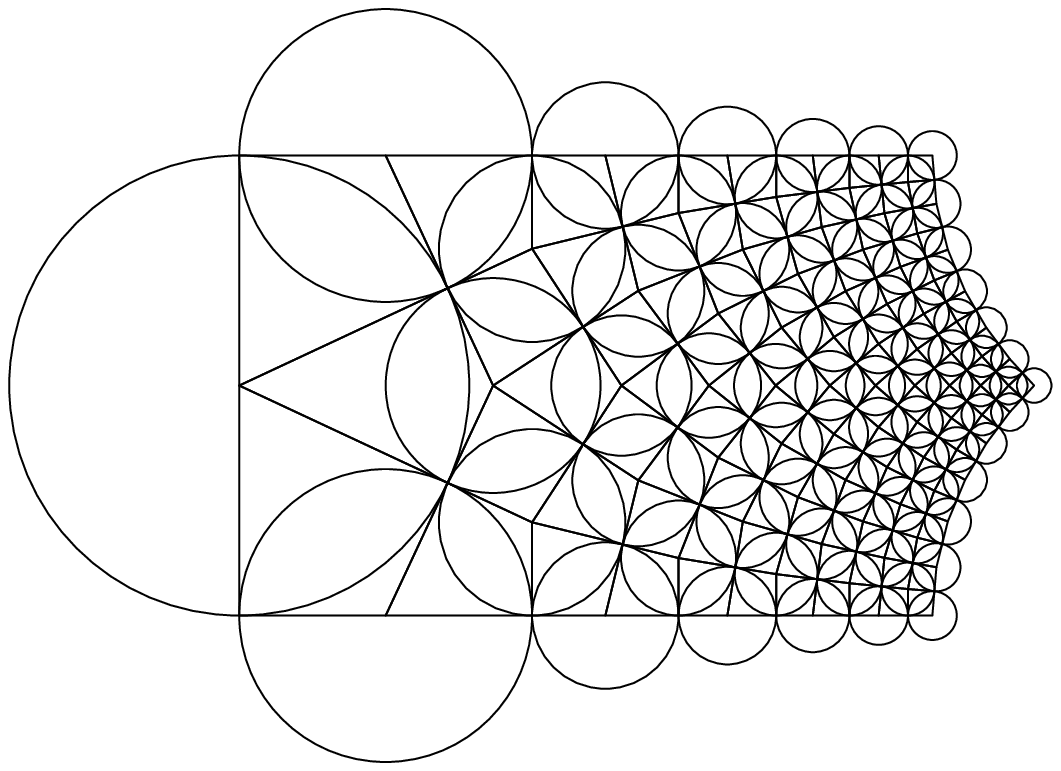}
\includegraphics[width=0.48\textwidth]{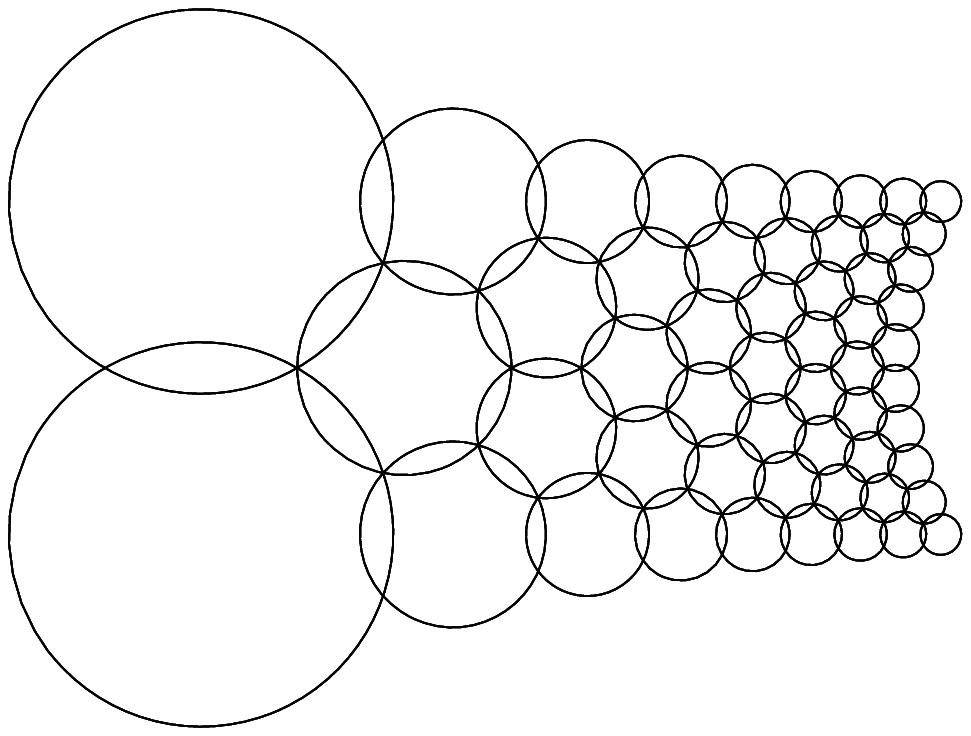}
\end{center}
\caption{Discrete logarithm circle patterns with the combinatorics
of the regular square and hexagonal lattices} \label{fig:hex Log}
\end{figure}

It is interesting do study the limiting behaviour of the function
$z^{2a}$ as $a\to 0$. It is not difficult to see that for all
$\bn\neq \boldsymbol 0$ one has $z^{2a}(\bn)\to 1$. Denote
\begin{equation}\label{eq: nonlin log}
{\rm L}(\bn)=\lim_{a\to 0}\frac{z^{2a}(\bn)-1}{2a}.
\end{equation}
This function is called the {\em discrete logarithmic function};
it should not be confused with the namesake function $\ell(\bn)$
in the linear theory (Sect. \ref{Sect discr log}). From eq.
(\ref{eq: nonlin log}) the following characterization is found:
the discrete logarithmic function ${\rm L}$ is the solution of the
discrete cross-ratio system on $\widetilde S$ defined by the
values on the coordinate semi-axes ${\rm L}_n^{(r)}={\rm
L}(n\be_r)$, $r\in[m,m+d-1]$, which solve the recurrent relation
\begin{equation}\label{dL axis recur}
 n\,\frac{({\rm L}_{n+1}-{\rm L}_n)
 ({\rm L}_n-{\rm L}_{n-1})}
 {{\rm L}_{n+1}-{\rm L}_{n-1}}=\frac{1}{2}
\end{equation}
with the initial conditions:
\begin{equation}\label{dL init}
{\rm L}_0^{(r)}={\rm L}(\boldsymbol 0)=\infty,\qquad {\rm
L}_1^{(r)}={\rm L}(\be_r)=\log \theta_r,
\end{equation}
where $\log\theta_r$ is chosen in the interval (\ref{ineq S}).
Explicit expressions:
\begin{equation}
 {\rm L}_{2n}^{(r)}=
\log\theta_r+\sum_{k=1}^{n-1}\frac{1}{k}+\frac{1}{2n}\,, \qquad
 {\rm L}_{2n+1}^{(r)}=
\log\theta_r+\sum_{k=1}^{n}\frac{1}{k}\,.
\end{equation}
\begin{thm}\label{Thm log nonlinear}
The discrete logarithm is isomonodromic and satisfies, at any
point $\bn\in\widetilde S$, the following constraint:
\begin{equation}\label{dL constr}
\sum_{j=1}^d
n_j\,\frac{({\rm L}(\bn+\be_j)-{\rm L}(\bn))
({\rm L}(\bn)-{\rm L}(\bn-\be_j))}
{{\rm L}(\bn+\be_j)-{\rm L}(\bn-\be_j)}=\frac{1}{2}.
\end{equation}
\end{thm}
By restriction to quad-surfaces $\Omega_\cD$, we come to the
discrete logarithmic function on arbitrary quasicrystallic
quad-graphs $\cD$. By construction, they all correspond to circle
patterns. A conjecture that these circle patterns are embedded
seems plausible (see Fig.~\ref{fig:hex Log}).

\section{Linearization}
\label{Sect linearization}

Let $\theta:\vec{E}(\cD)\to\bbC$ be an edge labelling, and let
$p:V(\cD)\to\bbC$ be the corresponding parallelogram realization
of $\cD$ defined by $p(y)-p(x)=\theta(x,y)$. Consider the trivial
solutions
\[
z_0(x)=p(x),\qquad w_0(x)=1,\qquad\forall x\in V(\cD)
\]
of the cross-ratio system (\ref{cr paral}) and the corresponding
Hirota system (\ref{Hirota in p}). Suppose that
$z_0:V(\cD)\to\bbC$ belongs to a differentiable one-parameter
family of solutions $z_\epsilon:V(\cD)\to\bbC$, $\epsilon\in
(-\epsilon_0,\epsilon_0)$, of the same cross-ratio system, and
denote by $w_\epsilon:V(\cD)\to\bbC$ the corresponding solutions
of the Hirota system. Denote
\begin{equation}\label{linearization}
g=\frac{dz_{\epsilon}}{d\epsilon}\bigg|_{\epsilon=0}\;,\qquad
f=\bigg(w_\epsilon^{-1}\,
\frac{dw_{\epsilon}}{d\epsilon}\bigg)_{\epsilon=0}\;.
\end{equation}
\begin{thm}\label{prop linearization}
Both functions $f,g:V(\cD)\to\bbC$ solve discrete Cauchy-Riemann
equations (\ref{CR paral}).
\end{thm}
{\bf Proof.} By differentiating (\ref{w vs z}), we obtain a
relation between the functions $f,g:V(\cD)\to\bbC$:
\begin{equation}\label{f vs g}
g(y)-g(x)=\big(f(x)+f(y)\big)\big(p(y)-p(x)\big), \qquad \forall
(x,y)\in\vec{E}(\cD).
\end{equation}
The proof of proposition is based on this relation solely. Indeed,
the closeness condition for the form on the right-hand side reads:
\begin{eqnarray*}
\lefteqn{\big(f(x_0)+f(y_0)\big)\big(p(y_0)-p(x_0)\big)+
         \big(f(y_0)+f(x_1)\big)\big(p(x_1)-p(y_0)\big)+}\nonumber\\
      && \big(f(x_1)+f(y_1)\big)\big(p(y_1)-p(x_1)\big)+
         \big(f(y_1)+f(x_0)\big)\big(p(x_0)-p(y_1)\big)=0,
\end{eqnarray*}
which is equivalent to (\ref{CR paral}) for the function $f$.
Similarly, the closeness condition for $f$, that is,
\[
\big(f(x_0)+f(y_0)\big)-\big(f(y_0)+f(x_1)\big)+
\big(f(x_1)+f(y_1)\big)-\big(f(y_1)+f(x_0)\big)=0,
\]
yields:
\[
\frac{g(y_0)-g(x_0)}{p(y_0)-p(x_0)}-\frac{g(x_1)-g(y_0)}{p(x_1)-p(y_0)}
+\frac{g(y_1)-g(x_1)}{p(y_1)-p(x_1)}-\frac{g(x_0)-g(y_1)}{p(x_0)-p(y_1)}=0.
\]
Under the condition $p(y_0)-p(x_0)=p(x_1)-p(y_1)$, this is
equivalent to (\ref{CR paral}) for $g$.   $\Box$
\medskip

{\bf Remark.} This proof shows that, given a discrete holomorphic
function $f:V(\cD)\to\bbC$, relation (\ref{f vs g}) correctly
defines a unique, up to an additive constant, function
$g:V(\cD)\to\bbC$, which is also discrete holomorphic. Conversely,
for any $g$ satisfying the discrete Cauchy-Riemann equations
(\ref{CR paral}), relation (\ref{f vs g}) defines a function $f$
correctly and uniquely (up to an additive black-white constant);
this function $f$ also solves the discrete Cauchy-Riemann
equations (\ref{CR paral}). Actually, formula (\ref{f vs g})
expresses that the discrete holomorphic function $f$ is the {\it
discrete derivative} of $g$, so that $g$ is obtained from $f$ by
{\it discrete integration}. This operation was considered in
\cite{D1, D2, M1}.
\medskip

Summarizing, we have the following statement.
\begin{thm}
a) A tangent space to the set of solutions of an integrable
cross-ratio system, at a point corresponding to a rhombic
embedding of a quad-graph, consists of discrete holomorphic
functions on this embedding. This holds in both descriptions of
the above set: in terms of variables $z$ satisfying the
cross-ratio equations, and in terms of variables $w$ satisfying
the Hirota equations. The corresponding two descriptions of the
tangent space are related by taking the discrete derivative (resp.
anti-derivative) of discrete holomorphic functions.

b) A tangent space to the set of integrable circle patterns of a
given combinatorics, at a point corresponding to an isoradial
pattern, consists of discrete holomorphic functions on the rhombic
embedding of the corresponding quad-graph, which take real values
at white vertices and purely imaginary values at black ones. This
holds in the description of circle patterns in terms of circle
radii and rotation angles at intersection points (Hirota system).
\end{thm}

A spectacular instance of this linearization property is delivered
by the isomonodromic discrete logarithm studied in Sect. \ref{Sect
discr log} and isomonodromic $z^{2a}$ circle patterns of Sect.
\ref{Sect z^a}.
\begin{thm}
The tangent vector to the space of integrable circle patterns
along the curve consisting of patterns $w^{2a-1}$, at the
isoradial point corresponding to $a=1/2$, is the discrete
logarithmic function $\ell$ of Sect. \ref{Sect discr log}.
\end{thm}
{\bf Proof.} We have to prove that the discrete logarithm $\ell$
and the discrete power function $w^{2a-1}$ are related by
\[
\ell(\bn)=\Big(\frac{1}{2}\frac{d}{da}
\,w^{2a-1}(\bn)\Big)_{a=1/2}\,.
\]
Due to Theorem \ref{prop linearization}, it is enough to prove
this for the initial data on the coordinate semi-axes. But this
follows by differentiating with respect to $a$ the initial values
(\ref{Hirota axes}) at the point $a=1/2$, where all $w=1$: the
result coincides with (\ref{Green axes}). $\Box$


\begin{thebibliography}{WWWW}
\addcontentsline{toc}{chapter}{Bibliography}

\bibitem[Ad]{A} V.E.~Adler.
  B\"acklund transformation for the Krichever-Novikov equation.
  {\em Internat. Math. Res. Notices}, 1998, No 1, 1--4.

\bibitem[AdBS1]{ABS1} V.E.~Adler, A.I.~Bobenko, Yu.B.~Suris.
    Classification of integrable equations on quad-graphs.
    The consistency approach.
    {\em Commun. Math. Phys.}, 2003, {\bf 233}, 513--543.

\bibitem[AdBS2]{ABS2} V.E.~Adler, A.I.~Bobenko, Yu.B.~Suris.
    Geometry of Yang--Baxter maps: pencils of conics and
    quadrirational mappings.
    {\em Commun. Anal. Geom}, 2004, {\bf 12}, 967-1007.

\bibitem[AdS]{AS} V.E.~Adler, Yu.B.~Suris.
   Q4: Integrable master equation related to an elliptic curves.
   {\em Internat. Math. Res. Notices}, 2004, No 47, 2523--2553.

\bibitem[AdV]{AV} V.E.~Adler, A.P.~Veselov.
     Cauchy problem for integrable discrete equations on quad-graphs.
     {\em Acta Appl. Math.} {\bf 84} (2004), 237--262.

\bibitem[Ag]{Ag} S.I.~Agafonov.
    Imbedded circle patterns with the combinatorics of the square grid
    and discrete Painlev\'e equations.
    {\em Discrete Comput. Geom.} {\bf 29} (2003), 305--319.

\bibitem[AgB1]{AB1} S.I.~Agafonov, A.I.~Bobenko.
   Discrete ${\rm Z}^\gamma$  and Painlev\'e equations.
   {\em Internat. Math. Res. Notices}, 2000, no. 4, 165--193.

\bibitem[AgB2]{AB2} S.I.~Agafonov, A.I.~Bobenko.
   Hexagonal circle patterns with constant intersection angles and
   discrete Painlev\'e and Riccati equations.
   {\em J. Math. Phys.} {\bf 44} (2003), 3455--3469.

\bibitem[AkhKV]{AKV} A.A.~Akhmetshin, I.M.~Krichever, Y.S.~Volvovski.
  Discrete analogues of the Darboux-Egoroff metrics.
  {\em Tr. Mat. Inst. Steklova} {\bf 225} (1999), 21--45 (Russian);
  English translation: {\it Proc. Steklov Inst. Math.}, {\bf 2 (225)}
  (1999), 16--39.

\bibitem[BeaS]{BeaS} A.F.~Beardon, K.~Stephenson.
   The uniformization theorem for circle packings.
   {\it Indiana Univ. Math. J.}, 1990, {\bf 39}, 1383--1425.

\bibitem[Ber]{Be} M.~Berger.
   {\em Geometry.} Berlin etc.: Springer-Verlag, 1987.

\bibitem[Bi]{Bi} L.~Bianchi.
   {\em Lezioni di geometria differenziale.}
   3rd edition. Pisa: Enrico Spoerri,
   1923 (Italian). iv+806, xi+832 pp.

\bibitem[Bl]{Bl} W.~Blaschke.
  {\em Projektive Geometrie}. 3rd edition. Basel etc.: Birkh\"auser, 1954.
  197 pp.

\bibitem[BlIII]{BlIII} W.~Blaschke.
  Vorlesungen \"uber Differentialgeometrie III:
  Differentialgeometrie der Kreise und Kugeln.
  Bearbeitet von G. Thomsen. Berlin: Springer,
  1929 (German). x+474 pp.

\bibitem[Bob]{B} A.I.~Bobenko.
   Discrete conformal maps and surfaces. -- In:
  {\em SIDE III -- Symmetries and integrability of difference equations
  (Sabaudia, 1998),  CRM Proc. Lecture Notes}, {\bf 225} (2000), pp.
   97--108.

\bibitem[BobHe]{BH} A.I.~Bobenko, U.~Hertrich-Jeromin.
   Orthogonal nets and Clifford algebras.
   {\em T\^ohoku Math. Publ.} {\bf 20} (2001), 7--22.

\bibitem[BobHo]{BHo} A.I.~Bobenko, T.~Hoffmann.
   Hexagonal circle patterns and integrable systems. Patterns with
   constant angles. {\em Duke Math. J.} {\bf 116:3} (2003) 525--566.

\bibitem[BobHoSp]{BHSp} A.I.~Bobenko, T.~Hoffmann, B.A.~Springborn.
   Discrete minimal surfaces: geometry from combinatorics.
   {\em Annals of Math.}, 2005 (to appear).

\bibitem[BobHoSu]{BHS} A.I.~Bobenko, T.~Hoffmann, Yu.B.~Suris.
   Hexagonal circle patterns and integrable systems: patterns with the
   multi-ratio property and Lax equations on the regular triangular
   lattice.
   {\em Int. Math. Res. Not.}  2002,  no. 3, 111--164.

\bibitem[BobMaS1]{BMS1} A.I.~Bobenko, D.~Matthes, Yu.B.~Suris.
  Nonlinear hyperbolic equations in surface theory: integrable
  discretizations and approximation results.
  {\it Algebra Anal.}, 2005 (to appear).

\bibitem[BobMaS2]{BMS2} A.I.~Bobenko, D.~Matthes, Yu.B.~Suris.
  Discrete and smooth orthogonal systems: $C\sp \infty$-approximation.
  {\it Int. Math. Res. Not.}, {\bf 45} (2003), 2415--2459.

\bibitem[BobMeS]{BMS} A.I.~Bobenko, Ch.~Mercat, Yu.B.~Suris.
  Linear and nonlinear theories of discrete analytic functions. Integrable
  structure and isomonodromic Green's function.
  {\em J. Reine Angew. Math.}, 2005 (to appear).

\bibitem[BobP1]{BP1} A.I.~Bobenko, U.~Pinkall.
   Discrete surfaces with constant negative Gaussian curvature
   and the Hirota equation.
   {\it J. Differential Geom.}, {\bf 43} (1996), no. 3, 527--611.

\bibitem[BobP2]{BP2} A.I.~Bobenko, U.~Pinkall.
   Discrete isothermic surfaces.
  {\it J. Reine Angew. Math.}, {\bf 475} (1996), 187--208.

\bibitem[BobP3]{BP3} A.I.~Bobenko, U.~Pinkall.
   Discretization of surfaces and integrable systems. -- In: {\em
   Discrete integrable geometry and physics}, Eds. A.I.~Bobenko and
   R.~Seiler, Oxford: Clarendon Press, 1999, pp. 3--58.

\bibitem[BobSch1]{BSch1} A.I.~Bobenko, W.K.~Schief.
   Discrete indefinite affine spheres. -- In: {\em
   Discrete integrable geometry and physics}, Eds. A.I.~Bobenko and
   R.~Seiler, Oxford: Clarendon Press, 1999, pp. 113--138.

\bibitem[BobSch2]{BSch2} A.I.~Bobenko, W.K.~Schief.
   Affine spheres: discretization via duality relations.
   {\em Experiment. Math.}, {\bf 8} (1999), 261--280.

\bibitem[BobSp]{BSp} A.I.~Bobenko,B.A.~Springborn.
   Variational principles for circle patterns and Koebe's theorem.
   {\em Trans. Amer. Math. Soc.} {\bf  356}  (2004), 659--689.

\bibitem[BobSu1]{BS1} A.I.~Bobenko, Yu.B.~Suris.
    Integrable systems on quad-graphs.
    {\em Internat. Math. Res. Notices}, 2002, Nr. 11, 573--611.

\bibitem[BobSu2]{BS2} A.I.~Bobenko, Yu.B.~Suris.
    Integrable noncommutative equations on quad-graphs. The consistency
    approach.
    {\em Lett. Math. Phys.}, 2002, {\bf 61}, 241--254.

\bibitem[BogK]{BK} L.V.~Bogdanov, B.G.~Konopelchenko.
    Lattice and $q$-difference Darboux-Zakharov-Manakov systems via
    $\overline\partial$-dressing method.
    {\em J. Phys. A}, {\bf 28} (1995), L173--L178.

\bibitem[CDS]{CDS} J.~Cieslinski, A.~Doliwa, P.M.~Santini.
   The integrable discrete analogues of orthogonal coordinate systems
   are multi-dimensional circular lattices.
   {\em Phys. Lett. A} {\bf 235} (1997), 480--488.

\bibitem[Da1]{Da1} G.~Darboux.
   {\em Le\c{c}ons sur la th\'eorie g\'en\'erale des surfaces et les
   applications g\'eom\'etriques du calcul infinit\'esimal.}
   T.I--IV. 3rd edition.
   Paris: Gauthier-Villars, 1914--1927. (French)

\bibitem[Da2]{Da2} G.~Darboux.
  {\em Le\c{c}ons sur les syst\'emes orthogonaux et les coordonn\'ees
  curvilignes.}
  2nd edition. Paris: Gauthier-Villars, 1910. vii+567 pp.

\bibitem[DoSh1]{DSh1} N.P.~Dolbilin, M.A.~Shtan'ko, M.I.~Shtogrin.
  Cubic subcomplexes in regular lattices. (Russian)
  {\em Dokl. Akad. Nauk SSSR} {\bf 291} (1986), 277--279.

\bibitem[DoSh2]{DSh2} N.P.~Dolbilin, M.A.~Shtan'ko, M.I.~Shtogrin.
   Cubic manifolds in lattices. (Russian)
   {\em Izv. Ross. Akad. Nauk Ser. Mat.}, {\bf 58} (1994), 93--107;
   english translation in: {\em Russian Acad. Sci. Izv. Math.}
   {\bf 44} (1995), 301--313.

\bibitem[Do1]{D1} A.~Doliwa.
  Geometric discretization of the Toda system.
  {\em Phys. Lett. A.} {\bf 234} (1997), 187--192.

\bibitem[Do2]{D2} A.~Doliwa.
  Quadratic reductions of quadrilateral lattices.
  {\em J. Geom. Phys.} {\bf  30}  (1999), 169--186.

\bibitem[Do3]{D3} A.~Doliwa.
  Lattice geometry of the Hirota equation. In:
  {\em SIDE III -- Symmetries and integrability of difference equations
  (Sabaudia, 1998),  CRM Proc. Lecture Notes}, {\bf 225} (2000), pp.
  93--100.

\bibitem[Do4]{D4} A.~Doliwa.
  Discrete asymptotic nets and $W$-congruences
  in Pl{\"u}cker line geometry.
  {\it J. Geom. Phys.} {\bf 39} (2001), no. 1, 9--29.

\bibitem[Do5]{D5} A.~Doliwa.
  Geometric discretization of the Koenigs nets.
  {\em J. Math. Phys.} {\bf 44} (2003), 2234--2249.

\bibitem[DoMS]{DMS} A.~Doliwa,  S.V.~Manakov, P.M.~Santini.
  $\bar\partial$-reductions of the multidimensional quadrilateral
  lattice. The multidimensional circular lattice.
  {\it Commun. Math. Phys.} {\bf 196} (1998), 1--18.

\bibitem[DoNS1]{DNS} A.~Doliwa, M.~Nieszporski,  P.M.~Santini.
   Asymptotic lattices and their integrable reductions. I.
   The Bianchi-Ernst and the Fubini-Ragazzi lattices.
   {\em J. Phys. A} {\bf  34}  (2001),  10423--10439.

\bibitem[DoNS2]{DNS2} A.~Doliwa, M.~Nieszporski, P.M.~Santini.
   Geometric discretization of the Bianchi system.
   {\em J. Geom. Phys.} {\bf 52} (2004), 217--240.

\bibitem[DoS1]{DS1} A.~Doliwa, P.M.~Santini.
  Multidimensional quadrilateral lattices are integrable.
  {\it Phys. Lett. A}, {\bf 233} (1997), 265--372.

\bibitem[DoS2]{DS2} A.~Doliwa, P.M.~Santini.
  The symmetric, $d$-invariant and Egorov reductions of the
  quadrilateral lattice.
  {\em J. Geom. Phys.} {\bf  36}  (2000),  60--102.

\bibitem[DoSM]{DSM} A.~Doliwa, P.M.~Santini, M.~Ma\~{n}as.
  Transformations of quadrilateral lattices.
  {\em J. Math. Phys.} {\bf  41}  (2000),  944--990.

\bibitem[DoGNS]{DGNS}  A.~Doliwa, P.~Grinevich, M.~Nieszporski, P.M.~Santini.
   Integrable lattices and their sub-lattices: from the discrete Moutard
   (discrete Cauchy-Riemann) 4-point equation to the self-adjoint 5-point
   scheme.
   {\em Preprint} {\tt http://arxiv.org/abs/nlin.SI/0410046}.

\bibitem[Dr]{D} V.G.~Drinfeld.
  On some unsolved problems in quantum group theory.
  {\it Lecture Notes Math.}, 1992, {\bf 1510}, 1--8.

\bibitem[Du1]{Du1} R.J.~Duffin.
  Basic properties of discrete analytic functions.
  {\em Duke Math. J.}, 1956, {\bf 23}, 335--363.

\bibitem[Du2]{Du2} R.J.~Duffin.
  Potential theory on a rhombic lattice.
  {\em J. Combinatorial Theory}, 1968, {\bf 5}, 258--272.

\bibitem[DyN]{DN} I.A.~Dynnikov, S.P.~Novikov.
  Geometry of the triangle equation on two-manifolds.
  {\em Mosc. Math. J.} {\bf 3} (2003), 419--438.

\bibitem[E1]{E1} L.P.~Eisenhart.
  {\em A treatise on the differential geometry of curves and surfaces.}
  Boston: Ginn and Co., 1909. ix+379 pp.

\bibitem[E2]{E2} L.P.~Eisenhart.
  {\em Transformations of surfaces.}
  Princeton University Press, 1923. ix+379 pp.

\bibitem[F]{F} J.~Ferrand.
  Fonctions preharmoniques et functions preholomorphes.
  {\em Bull. Sci. Math.}, 2nd ser., 1944, {\bf 68}, 152--180.

\bibitem[GT]{GTs} E.I.~Ganzha, S.P.~Tsarev.
    An algebraic superposition formula and the completeness of B\"acklund
    transformations of $(2+1)$-dimensional integrable systems.
    {\em Uspekhi Mat. Nauk} {\bf 51} (1996), no. 6, 197--198 (Russian);
    English translation in: {\em Russian Math. Surveys} {\bf 51} (1996),
    1200--1202.

\bibitem[He]{He} Z.-X.~He.
   Rigidity of infinite disk patterns.
   {\it Ann. of Math.}, 1999, {\bf 149}, p. 1--33.

\bibitem[HeS]{HS} Z.-X.~He, O.~Schramm.
   The $C^{\infty}$ convergence of hexagonal disc packings to Riemann map.
   {\it Acta Math.}, 1998, {\bf 180}, 219--245.

\bibitem[He1]{HJ1} U.~Hertrich-Jeromin.
   Transformations of discrete isothermic nets and discrete cmc-1
   surfaces in hyperbolic space.
   {\em Manuscr. Math.} {\bf 102} (2000), 465--486.

\bibitem[He2]{HJ2} U.~Hertrich-Jeromin.
   {\em Introduction to M\"obius differential geometry.}
   Cambridge University Press, 2003. xii+413 pp.

\bibitem[HeHP]{HHP} U.~Hertrich-Jeromin, T.~Hoffmann, U.~Pinkall.
   A discrete version of the Darboux transform for isothermic
   surfaces.  -- In: {\em
   Discrete integrable geometry and physics}, Eds. A.I.~Bobenko and
   R.~Seiler, Oxford: Clarendon Press, 1999, pp. 59--81.

\bibitem[HeMNP]{HMNP} U.~Hertrich-Jeromin, I.~McIntosh, P.~Norman, F.~Pedit.
   Periodic discrete conformal maps.
   {\em J. Reine Angew. Math.} {\bf 534} (2001), 129--153.

\bibitem[Hi]{H} R.~Hirota.
  Nonlinear partial difference equations. I. A
  difference analog of the Korteweg--de Vries equation.
  III. Discrete sine-Gordon equation.
  {\em J. Phys. Soc. Japan}, 1977, {\bf 43}, 1423--1433, 2079--2086.

\bibitem[Ho]{Hof} T.~Hoffmann.
  Darboux transformation for S-isothermic surfaces. {\em In
  preparation.}

\bibitem[IN]{IN} A.~Its, V.~Novokshenov.
  The isomonodromic deformation  method in the theory of Painlev\'e
  equations.
  {\em Lecture Notes Math.}, {\bf 1191}, 313 pp. Berlin:  Springer, 1986.

\bibitem[J]{Jonas} H.~Jonas.
   \"Uber die Transformation der konjugierten systeme und \"uber den
   gemeinsamen Ursprung der Bianchischen Permutabilit\"atstheoreme.
   {\em Sitzungsber. Berl. Math. Ges.}, 1915, {\bf 14}, 96--118.

\bibitem[Ke]{K} R.~Kenyon.
   The Laplacian and Dirac operators on critical planar graphs.
   {\it Invent. Math.}, 2002, {\bf 150}, 409--439.

\bibitem[KeS]{KS} R.~Kenyon, J.-M.~Schlenker.
   Rhombic embeddings of planar graphs with faces of degree 4.
   {\em Trans. AMS}, 2004 (to appear).

\bibitem[KoP]{KP} B.~Konopelchenko, U.~Pinkall.
   Projective generalizations of Lelieuvre's formula.
   {\em Geom. Dedicata} {\bf 79} (2000), 81--99.

\bibitem[KoSch]{KSch} B.G.~Konopelchenko, W.K.~Schief.
   Three-dimensional integrable lattices in Euclidean spaces: conjugacy and
   orthogonality.
   {\em R. Soc. Lond. Proc. Ser. A}, {\bf  454} (1998), 3075--3104.

\bibitem[Kor]{Ko} V.E.~Korepin.
   Completely integrable models in quasicrystals.
   {\em Commun. Math. Phys.} {\bf 110} (1987), 157--171.

\bibitem[Kr]{Kr} I.M.~Krichever.
  Algebraic-geometric n-orthogonal curvilinear coordinate systems
  and the solution of associativity equations.
  {\em Funct. Anal. Appl.}, {\bf 31} (1997), 25--39.

\bibitem[MDS]{MDS} M.~Ma\~{n}as, A.~Doliwa, P.M.~Santini.
  Darboux transformations for multidimensional quadrilateral lattices.
  {\em  Phys. Lett. A}, {\bf 232} (1997), 99--105.

\bibitem[MR]{MR} A.~Marden, B.~Rodin.
   On Thurston's formulation and proof of Andreev's theorem.
   {\it Lect. Notes Math.}, 1990, {\bf 1435}, 103--115.

\bibitem[MPS]{MPS} R.R.~Martin, J.~de Pont, T.J.~Sharrock.
  Cyclide surfaces in computer aided design. In:
  {\em The mathematics of surfaces}, Ed. J.A.~Gregory, Oxford:
  Clarendon Press, 1986, pp.253--268.

\bibitem[Me1]{M1} Ch.~Mercat.
   Discrete Riemann surfaces and the Ising model.
  {\it Commun. Math. Phys.}, 2001, {\bf 218}, 177--216.

\bibitem[Me2]{M2} Ch.~Mercat.
   Exponentials form a basis of discrete holomorphic functions.
   {\em Bull. Soc. Math. France} {\bf 132} (2004) 305--326.

\bibitem[Mou]{Mou} Th.F.~Moutard.
   Sur la construction des \'equations de la forme $\frac{1}{z}\,
   \frac{\partial^2 z}{\partial x \partial y}=\lambda(x,y)$ qui
   admettenent une int\'egrale g\'en\'erale explicite.
   {\em J. \'Ec. Pol.} {\bf 45} (1878), 1--11.

\bibitem[M\"o]{M} A.F.~M\"obius.
  Kann von zwei dreiseitigen Pyramiden eine jede in Bezug auf
  andere um- und eingeschrieben zugleich heissen?
  {\em J. Reine und Angew. Math.}, {\bf 3} (1828), 273--278.

\bibitem[Nie]{N} M.~Nieszporski.
   On a discretization of asymptotic nets.
   {\em J. Geom. Phys.} {\bf  40} (2002), 259--276.

\bibitem[NieSD]{NSD} M.~Nieszporski, P.M.~Santini, A.~Doliwa.
  Darboux transformations for 5-point and 7-point self-adjoint
  schemes and an integrable discretization of the 2D Schrödinger operator.
  {\em Phys. Lett. A}, {\bf 323} (2004), 241--250.

\bibitem[Nij]{Nij} F.W.~Nijhoff.
   Lax pair for the Adler (lattice Krichever-Novikov) sastem.
   {\em Phys. Lett. A} {\bf 297} (2002), 49--58.

\bibitem[NijC]{NC} F.W.~Nijhoff, H.W.~Capel.
    The discrete Korteweg-de Vries equation.
   {\em Acta Appl. Math.} { \bf 39} (1995) 133--158.

\bibitem[NiSch]{NSch} J.J.C.~Nimmo, W.K.~Schief.
   Superposition principles associated with the Moutard transformation:
   an integrable discretization of a $(2+1)$-dimensional sine-Gordon
   system.
   {\em Proc. Roy. Soc. London Ser. A} {\bf 453} (1997), 255--279.

\bibitem[No1]{No1} S.P.~Novikov.
   Schr\"odinger operators on graphs and symplectic geometry.
   In: {\em The Arnoldfest}, Fields Inst. Commun., 24, Providence: AMS,
   1999, 397--413.

\bibitem[No2]{No2} S.P.~Novikov.
  The discrete Schr\"odinger operator.
  {\em Tr. Mat. Inst. Steklova} {\bf 224} (1999) 275--290 (Russian);
  English translation: {\em Proc. Steklov Inst. Math.} {\bf 224} (1999)
  250--265.

\bibitem[NoD]{ND} S.P.~Novikov, I.A.~Dynnikov.
   Discrete spectral symmetries of small-dimensional differential
   operators and difference operators on regular lattices and
   two-dimensional manifolds.
   {\em Uspekhi Mat. Nauk} {\bf 52} (1997), no. 5(317), 175--234 (Russian);
   English translation: {\em Russian Math. Surveys} {\bf 52} (1997),
   1057--1116.

\bibitem[Nu]{Nu} A.W.~Nutbourne.
  The solution of frame matching equation.
  In: {\em The mathematics of surfaces}, Ed. J.A.~Gregory,
  Oxford: Clarendon Press, 1986, pp.233--252.


\bibitem[QNCL]{QNCL} G.R.W.~Quispel, F.W.~Nijhoff, H.W.~Capel, J. van
   der Linden.
   Linear integral equations and nonlinear difference-difference
   equations.
   {\em Physica A}, 1984, {\bf 125}, 344--380.

\bibitem[RSch]{RSch} C.~Rogers, W.K.~Schief.
   {\em B\"acklund and Darboux transformations. Geometry and modern
   applications in soliton theory.}
   Cambridge University Press, 2002, xviii+413 pp.

\bibitem[RSu]{RS} B.~Rodin, D.~Sullivan.
   The convergence of circle packings to Riemann mapping.
   {\it J. Diff. Geom.}, 1987, {\bf 26}, 349--360.

\bibitem[Sa]{Sauer} R.~Sauer.
   {\em Differenzengeometrie}.
   Berlin etc.: Springer, 1970, 234 pp.

\bibitem[Sch1]{Sch1} W.K.~Schief.
   Isothermic surfaces in spaces of arbitrary dimension: integrability,
   discretization and B\"acklund transformations. A discrete Calapso
   equation.
   {\em Stud. Appl. Math.}, {\bf 106} (2001), 85--137.

\bibitem[Sch2]{Sch2} W.K.~Schief.
   On the unification of classical and novel integrable surfaces. II.
   Difference geometry.
   {\em R. Soc. Lond. Proc. Ser. A} {\bf 459} (2003), 373--391.

\bibitem[Schr]{S} O.~Schramm.
   Circle patterns with the combinatorics of the square grid.
   {\it Duke Math. J.}, 1997, {\bf 86}, p. 347--389.

\bibitem[ShSh]{ShSh} M.A.~Shtan'ko, M.I.~Shtogrin.
   Embedding cubic manifolds and complexes into a cubic lattice.
   {\em Uspekhi Mat. Nauk} {\bf  47} (1992),  no. 1 (283),
   219--220 (Russian); English translation: {\em Russian Math. Surveys}
   {\bf  47} (1992), 267--268.

\bibitem[St1]{St1} K.~Stephenson.
   Circle packing: a mathematical tale.
   {\em Notices Amer. Math. Soc.} {\bf  50}  (2003),  1376--1388.

\bibitem[St2]{St2} K.~Stephenson.
   {\em Introduction to the theory of circle packing: discrete analytic
   functions.}
   Cambridge Univ. Press, 2005. 400 pp.

\bibitem[SuV]{SV} Yu.B.~Suris, A.P.~Veselov.
   Lax matrices for Yang-Baxter maps.
   {\em J. Nonlinear Math. Phys.}, {\bf 10, suppl. 2}  (2003), 223--230.

\bibitem[Th1]{T1} W.P.~Thurston.
  The finite Riemann mapping theorem.
  {\em Invited talk at the international symposium on the occasion of
  the proof of the Bieberbach conjecture}, Purdue University, 1985.

\bibitem[Th2]{T2} W.P.~Thurston.
   {\em Three-dimensional geometry and topology. Vol.1}, Princeton Univ.
   Press, 1997. x+311 pp.

\bibitem[Tzi]{Tzi} G.~Tzitz\'eica.
   {\em G\'eom\'etrie diff\'erentielle projective des r\'eseaux.}
   Bucharest: Cultura Nationala, 1924.

\bibitem[V]{V} A.P.~Veselov.
   Yang-Baxter maps and integrable dynamics.
   {\em Phys. Lett. A}, {\bf 314} (2003), 214--221.

\bibitem[W]{W} W.~Wunderlich.
  Zur Differenzengeometrie der Fl{\"a}chen konstanter negativer
  Kr{\"u}mmung.
  {\it {\"O}sterreich. Akad. Wiss. Math.-Nat. Kl.}, {\bf 160} (1951),
  39--77.

\bibitem[Z]{Z} V.E.~Zakharov.
  {\em What is integrability?}
  Berlin: Springer, 1991.
\end{thebibliography}
\end{document}